\documentclass[leqno,a4paper,12pt]{article}
\usepackage{amsmath,amsthm,bbm}
\usepackage[sans]{dsfont}

\newtheorem{Thm}{\indent Theorem}[section]
\newtheorem{Prop}[Thm]{\indent Proposition}
\newtheorem{Lem}[Thm]{\indent Lemma}
\newtheorem{Cor}[Thm]{\indent Corollary}
\newtheorem{Comp}[Thm]{\indent Complement}

\theoremstyle{definition}
\newtheorem{Def}[Thm]{\indent Definition}
\newtheorem{Rem}[Thm]{\indent Remark}
\newtheorem{Ex}[Thm]{\indent Example}
\newtheorem{Exs}[Thm]{\indent Examples}
\newtheorem{Cons}[Thm]{\indent Construction}

\def\qed{{\hskip0pt\unskip\unskip\nobreak\hfil\penalty50
          \hskip1em\hbox{}\nobreak\hfil
          {\bf q.e.d.}%
          \parfillskip=0pt\finalhyphendemerits=0
          \par}\medskip}

\newenvironment{Proof}
               {{\it Proof.}\quad}
               {\qed}

\newenvironment{Proofof}[1]
               {{\it Proof of #1.}\quad}
               {\qed}

\newcommand{\Prime}{\kern3\fontdimen1\font$'$\kern-7\fontdimen1\font}

\long\def\forget#1{}

\long\def\beginSIDEREMARK#1\endSIDEREMARK
    {{\par\bigskip\advance\leftskip by 2cm
                  \advance\rightskip by -2cm\noindent
      {\bf Our own side remark:} #1
      \par\bigskip\noindent}}

\long\def\beginFORGET#1\endFORGET{#1}
\long\def\beginFORGET#1\endFORGET{}

\def\?{\ ???\ \immediate\write16{}%
\immediate\write16{Warning: There was still a question mark . . . }%
\immediate\write16{}}

\usepackage{amsmath}

\usepackage{exscale}
\usepackage{amssymb}

\input cyracc.def

\newcommand{\BA}{{\mathbb{A}}}

\newcommand{\BC}{{\mathbb{C}}}

\newcommand{\BG}{{\mathbb{G}}}

\newcommand{\BQ}{{\mathbb{Q}}}
\newcommand{\BR}{{\mathbb{R}}}
\newcommand{\BS}{{\mathbb{S}}}

\newcommand{\BZ}{{\mathbb{Z}}}

\newcommand{\Fc}{{\mathfrak{c}}}

\newcommand{\Fn}{{\mathfrak{n}}}

\newcommand{\FH}{{\mathfrak{H}}}

\newcommand{\FU}{{\mathfrak{U}}}
\newcommand{\FV}{{\mathfrak{V}}}
\newcommand{\FW}{{\mathfrak{W}}}
\newcommand{\FX}{{\mathfrak{X}}}
\newcommand{\FY}{{\mathfrak{Y}}}
\newcommand{\FZ}{{\mathfrak{Z}}}

\newcommand{\CD}{{\cal D}}

\newcommand{\CL}{{\cal L}}

\newcommand{\CO}{{\cal O}}

\newcommand{\Ad}{\mathop{\rm Ad}\nolimits}

\newcommand{\adm}{\mathop{\rm adm}\nolimits}

\newcommand{\Aut}{\mathop{\rm Aut}\nolimits}

\newcommand{\Int}{\mathop{\rm Int}\nolimits}
\newcommand{\Lie}{\mathop{\rm Lie}\nolimits}

\newcommand{\Cent}{\mathop{\rm Cent}\nolimits}

\newcommand{\ext}{\mathop{\rm exd}\nolimits}
\newcommand{\GL}{\mathop{\rm GL}\nolimits}

\newcommand{\GQm}{\mathop{\BG_{m,\BQ}}\nolimits}
\newcommand{\GRm}{\mathop{\BG_{m,\BR}}\nolimits}

\newcommand{\Hom}{\mathop{\rm Hom}\nolimits}
\newcommand{\im}{\mathop{\rm Im}\nolimits}

\newcommand{\re}{\mathop{\rm Re}\nolimits}
\newcommand{\Rad}{\mathop{\rm Rad}\nolimits}

\newcommand{\Stab}{\mathop{\rm Stab}\nolimits}

\newcommand{\loccit}{[loc.$\;$cit.]}

\def\tei{\, | \,}
\def\halb{\frac{1}{2}}

\def\id{{\rm id}}

\newbox\mybox
\def\arrover#1{\mathrel{
       \setbox\mybox=\hbox spread 1.4em{\hfil$\scriptstyle#1$\hfil}
       \vbox{\offinterlineskip\copy\mybox
             \hbox to\wd\mybox{\rightarrowfill}}}}
\def\larrover#1{\mathrel{
       \setbox\mybox=\hbox spread 1.4em{\hfil$\scriptstyle#1$\hfil}
       \vbox{\offinterlineskip\copy\mybox
             \hbox to\wd\mybox{\leftarrowfill}}}}

\def\ontoover#1{\mathrel{
       \setbox\mybox=\hbox spread 1.4em{\hfil$\scriptstyle#1$\hfil}
       \vbox{\offinterlineskip\copy\mybox
             \hbox to\wd\mybox{\rightarrowfill\hskip-2.8mm
                               $\rightarrow$}}}}
\def\leftontoover#1{\mathrel{
       \setbox\mybox=\hbox spread 1.4em{\hfil$\scriptstyle#1$\hfil}
       \vbox{\offinterlineskip\copy\mybox
             \hbox to\wd\mybox{$\leftarrow$\hskip-2.8mm
                               \leftarrowfill}}}}
\def\longto{\longrightarrow}
\def\into{\hookrightarrow}
\def\onto{\ontoover{\ }}
\def\longonto{\ontoover{\ }}
\def\isoto{\arrover{\sim}}

\def\longinto{\lhook\joinrel\longrightarrow}

\usepackage[curve,matrix,arrow,cmtip]{xy}
\NoComputerModernTips

\def\myxymessage{\def\messagetext
   {Here an xy-pic diagram was omitted to speed up compilation . . . }
   \immediate\write16{\messagetext}
   \hbox{\bf \messagetext}}
\def\filxymatrix#1{\myxymessage}
\def\filxyarray#1{\myxymessage}

\newdir^{ (}{{}*!/-3pt/\dir^{(}}
\newdir_{ (}{{}*!/-3pt/\dir_{(}}
\newdir^{ )}{{}*!/+3pt/\dir^{)}}
\newdir_{ )}{{}*!/+3pt/\dir_{)}}

\def\rscript#1{\hbox to 0pt{$\scriptstyle#1$\hss}}

\let\oldbullet\bullet
\def\bullet{{\mathchoice{\oldbullet}%
                        {\oldbullet}%
                        {\scriptscriptstyle\oldbullet}%
                        {\oldbullet}}}

\newcommand{\argdot}{{\;\bullet\;}}

\begin{document}

%

\hfuzz=3pt
\overfullrule=10pt                   


\setlength{\abovedisplayskip}{6.0pt plus 3.0pt}
\setlength{\belowdisplayskip}{6.0pt plus 3.0pt}
\setlength{\abovedisplayshortskip}{6.0pt plus 3.0pt}
\setlength{\belowdisplayshortskip}{6.0pt plus 3.0pt}

\setlength{\baselineskip}{13.0pt}
\setlength{\lineskip}{0.0pt}
\setlength{\lineskiplimit}{0.0pt}

%
%

\title{Shimura data and corners: topology
\forget{
\footnotemark
\footnotetext{To appear in ....}
}
}
\author{\footnotesize by\\ \\
\mbox{\hskip-2cm
\begin{minipage}{6cm} \begin{center} \begin{tabular}{c}
J\"org Wildeshaus \footnote{
Partially supported by the \emph{Agence Nationale de la Recherche}, project no.\ 
ANR–18–CE40–0017 ``P\'eroides en g\'eom\'etrie arithm\'etique et motivique''. }\\[0.2cm]
\footnotesize Universit\'e Sorbonne Paris Nord \\[-3pt]
\footnotesize LAGA, CNRS (UMR~7539)\\[-3pt]
\footnotesize F-93430 Villetaneuse\\[-3pt]
\footnotesize France\\
{\footnotesize \tt wildesh@math.univ-paris13.fr}
\end{tabular} \end{center} \end{minipage}
\hskip-2cm}
\\[2.5cm]
}
\maketitle
\begin{abstract}
\noindent
The purpose of this article is to give a new construction of the map relating the Borel--Serre
and the Baily--Borel compactifications of a Shimura variety, and to provide a close analysis of its
main properties. \\


\end{abstract}


\bigskip
\bigskip
\bigskip

\noindent {\footnotesize Math.\ Subj.\ Class.\ (2020) numbers: 
14G35 (11F23, 11F75, 14P05).
}

\eject

\tableofcontents

\bigskip


%
%

\setcounter{section}{-1}
\section{Introduction}
\label{Intro}



Let $M$ be a \emph{Shimura variety}, and $M^*$ its \emph{Baily--Borel compactification}
\cite{AMRT,BB}. It is widely accepted that interesting topological, geometric and arithmetic 
information can be obtained from the control of the \emph{cohomo\-logy} of $M$, of
$M^*$, of its \emph{boundary} $M^* - M$ and its \emph{strata}
(indexed by the conjugation classes of \emph{admissible parabolics}), 
and of the interplay, \emph{via} localization, of these data. 
A major difficulty for obtaining this control
comes from the fact that $M^* - M$ is in general highly singular. \\

In order to circumvent this difficulty, one is led to consider a proper map
with target $M^*$, satisfying three requirements: (A)~the pre-image of $M$ is
mapped isomorphically to $M$, (B)~the pre-image of the
boundary $M^* - M$ is ``less singular'' than $M^* - M$ itself, 
(C)~the higher direct images under the restriction of the proper map 
to the pre-image of each stratum of the boundary can be ``controlled''. \\

Given such a map, \emph{proper base change} relates cohomology of its source to 
cohomology of $M^*$; thanks to (B) and (C), 
it is reasonable to expect the control of the former to be easier to obtain. \\

The present article is concerned with the map $p$ from the \emph{Borel--Serre
compactification} $M^{BS}$ \cite{BS} to (the space of complex points of)
$M^*$, whose existence was established in \cite{Z}.
Our aim is to give a close analysis of $p$, aiming in particular at the aspects related to
points~(B) and (C) above. \\

More precisely, we first
(1)~make explicit the link between \emph{Shimura data} \cite{P1}
and \emph{spaces of type $S-\BQ$} (the latter being  the notion underlying all the constructions from \cite{BS}).
Part~(1) is indeed necessary: as we shall explain below, 
the space $\FX$ underlying Shimura data $(P,\FX)$ is \emph{not} in general
of type $S-\BQ$ under the group $P(\BR)$ of real points of $P$.
The ``link'' hinted at is therefore somewhat subtle to establish.
Once this foundational issue is clarified, we  
(2)~identify the \emph{geodesic action} (the main ingredient 
for the construction of $M^{BS}$ \cite{BS}), 
in terms of the Shimura data $(P,\FX)$,
(3)~give an adelic definition of $M^{BS}$,
(4)~provide a new, explicit, and direct construction of the map $p$.
Using this construction, we 
(5)~make explicit the interplay of the closures of the strata of $M^{BS}$
(indexed by the conjugation classes of \emph{all} parabolics),
and the pre-images under $p$ of the strata of $M^*$,
(6)~describe the fibres of $p$, together with their induced stratifications
(using again the language of spaces of type $S-\BQ$),
(7)~prove that if the ``level'' of $M$ 
is \emph{neat}, then the restriction of $p$ to the pre-image of each
stratum of $M^*$ is a locally trivial stratified fibration. 
Variants of (6) and (7) exist for the map $\FX^{BS} \to \FX^*$  
between the spaces ``covering'' $M^{BS}$ and $M^*$, and for
the quotient map $p^r: M^{rBS} \to M^*$ from the \emph{reductive
Borel--Serre compactification} \cite{Z1} to $M^*$; 
the non-stratified versions of these variants
are essentially known, at least when $\FX$ is connected,
thanks to \cite{Z} (concerning $\FX^{BS} \to \FX^*$) and \cite{G} (concerning $p^r$).\\ 

The numbering of this program corresponds to the sections of this paper.
Section~\ref{4a} relates to point~(B) above, and Sections~\ref{5} and \ref{6}, 
to point~(C). The main applications to the cohomology of $M^*$ will be spelled out 
in a separate article \cite{W1}. \\

As is apparent from the above discussion, there is a certain overlap with \cite{Z} and \cite{G}. 
Here is what we see as our main new contributions:
as far as the results are concerned, the direct construction of $p$ and the
detailed study of the interplay of the stratifications with its fibres and
local splittings; concerning the methods, the systematic use of the language of 
Shimura data (a corollary of which is the total absence of Langlands decompositions). 
Of course, the method is not without influence on the results:
we do think that their precise versions actually necessitate
their formulation in the language of Shimura data (this is particularly true for our description of the geodesic action).  \\

However, in order to have the constructions from \cite{BS} at one's disposal, 
one \emph{a priori} needs 
homogeneous spaces, that are connected
\cite[Rem.~2.4~(1)]{BS}. This is not necessarily the case for the space 
$\FX$ underlying $(P , \FX)$ --- not
even in those cases where $\FX$ is a conjugation class of morphisms from the
Deligne torus to $P_\BC$ (rather than just a finite covering of such a class). 
Note that this phenomenon already occurs for $P = \GL_2$ and $\FX$ equal
to the complement of the real line in $\BC \,$, with the standard action of 
$\GL_2(\BR)$.
The solution consists in considering the (finitely many) connected components $\FX^0$ of $\FX$, 
and proving that for each of them,
the action of the stabilizer $\Stab_{P(\BR)}(\FX^0)$ on $\FX^0$ extends cano\-ni\-cally
to an action of the whole of $P(\BR)$ on $\FX^0$ (\emph{sic!}); this extension coincides with the given 
action on $\FX$ if and only if
$\FX$ is connected. The difference of the two actions manifests itself 
in the stabilizers of points: those of the 
extended action always contain a maximal compact subgroup of $P(\BR)$. 
The corresponding statement 
for the stabilizers of the
action underlying the data $(P,\FX)$ is not necessarily true: in general,
maxi\-mal compact subgroups of the stabilizer of $x \in \FX^0$ 
are maximal compact only in $\Stab_{P(\BR)}(\FX^0)$, but not necessarily in $P(\BR)$. 
This phenomenon remains in general visible at the level of the derived group
$P^{der}$.
The difference of the two actions also concerns analyticity: for example, 
on the upper half plane, by the above, the action of the subgroup $\GL_2(\BR)^+$ of $\GL_2(\BR)$ 
of matrices of positive determinant can be canonically extended to an action of the whole of $\GL_2(\BR)$.
But only $\GL_2(\BR)^+$ acts by
$\BC$-analytic automorphisms (see our Example~\ref{1H}). \\ 

Let us now give a description of the individual sections. \\

Section~\ref{1} first gives a formal framework for extensions of the action from
a subgroup of finite index of a real Lie group $H$ to the whole of $H$.
We modify the main technical notion from \cite{BS}, whence our
definition of \emph{homogeneous space of type $S'$}
(Definition~\ref{1B}). Using enough information from \cite{P1}, we prove 
(Theorem~\ref{1F}) that every connected component $\FX^0$ of $\FX$ is
homogeneous of type $S'$, provided that the Shimura data $(P,\FX)$ sa\-tisfy two technical
conditions (which for the purpose of this introduction will be ignored). 
Using the extension principle for spaces of type $S'$ (Proposition~\ref{1C}),
we prove, as mentioned above, that the action of the stabilizer $\Stab_{P(\BR)}(\FX^0)$
on $\FX^0$ can be canonically extended to an action of $P(\BR)$. 
We then prove that with this extended action, the space $\FX^0$ is 
of type $S-\BQ$ in the sense of \cite{BS} (Theorem~\ref{1M}). \\

We are thus in a position to apply the program from \cite{BS} to $\FX$
(by applying it to every connected component), the essential
ingredient of this program being the \emph{geodesic action}.
The purpose of Section~\ref{2} is to give a description of this action, that
makes explicit use of the information coded within the Shimura data $(P,\FX)$, and
more precisely, of the \emph{weight co-characters} associated to the points of $\FX$.
According to \cite{BS}, the geodesic action is a collection of actions on $\FX$,
one for each parabolic subgroup $Q$ of $P$. 
If the intersection of two parabolics $Q$ and $Q'$ remains parabolic,
then the action indexed by $Q \cap Q'$ is generated by those indexed by $Q$ and $Q'$,
meaning that the geodesic action is known once the actions indexed by maximal proper
parabolics are known. This observation allows to reduce our explicit description
of the geodesic action
to the case of maximal proper parabolics of $P$. 
The action associated to a maximal proper parabolic is an action of the group $\BR^*_+$ of strictly
positive real numbers. In order to describe its effect, note that ma\-xi\-mal proper parabolics of $P$ 
are admissible. But for admissible parabolics $Q$, the full force of the theory of 
\emph{boundary components} of Shimura data \cite{P1} can be employed: 
the Shimura data $(P,\FX)$ contain in particular a rule assigning
to any point $x \in \FX$ a weight co-character of $P_\BC$, 
which for the purpose of this introduction will be denoted by $w_x \,$.
The same is true for any boundary component associated to $Q$; here
the assigned weight co-character will be denoted by $w_{x,Q} \,$. 
While the
technical conditions from Section~\ref{1} gua\-rantee that $w_x$ is defined over $\BR$,
this is not the case for $w_{x,Q}$ (unless $Q = P$). 
We therefore define the \emph{real part of $x$ with respect to $Q$} (Definition~\ref{2D}). 
This real part
$\re_Q(x)$ is a point of one of the boundary components; its weight-co-character
$w_{\re_Q(x)}$ is conjugate to $w_{x,Q}$, and defined over $\BR$.
The main insight of Section~\ref{2} is that if $Q$ is not only admissible,
but maximal proper, then the geodesic action of $r \in \BR^*_+$ maps 
$x$ to its image under $w_{\re_Q(x)}(r) \in P(\BR)$ (Theorem~\ref{2H}). 
This result will be at the heart of our construction, 
and proof of continuity of the map $p$ (Section~\ref{4}). 
The reader familiar with \cite{BS} will note that our approach to the parameterization
of the geodesic action is dual to the one of \loccit : we use co-characters (since
they can be ``read off'' the Shimura data), while \loccit \ uses characters. 
Theorem~\ref{2H} also contains the dictionary allowing to pass from one point of
view to the other. \\

In Section~\ref{3}, we follow \cite{BS}, to construct, starting from 
the geodesic action on $\FX$, the
\emph{manifold with corners} $\FX^{BS}$
(denoted $\bar{\FX}$ in \loccit), which contains $\FX$ as an open subset.
In order to obtain from $\FX^{BS}$ the Borel--Serre compactification
\[
M^K (P,\FX) (\BC)^{BS} := P (\BQ) \backslash \bigl( \FX^{BS} \times P (\BA_f) / K \bigr) 
\]
of level $K$ (Definition~\ref{3A}), care needs to be applied
as the quotient is formed with respect to the action of $P$ underlying
the Shimura data, while $\FX^{BS}$ is the result of a limit process involving
the geodesic action. One therefore needs to know that the two actions commute
(Proposition~\ref{2L}; \cite[Prop.~3.4]{BS} if $\FX$ is connected).  
The main results from \cite{BS} then carry over; in particular (Theorem~\ref{3B}),
the Borel--Serre compactification is indeed compact. \\

We are then in a position to construct the map $p^K$ relating the Borel--Serre compactification and
\[
M^K (P,\FX)^* (\BC) = P (\BQ) \backslash \bigl( \FX^* \times P (\BA_f) / K \bigr) \; ,
\] 
the space of complex points of the Baily--Borel compactification $M^K (P,\FX)^*$.
Here, the space $\FX^*$ is obtained as a partial compactification of $\FX$, by adding
one stratum for each boundary component of $(P,\FX)$. The map $p^K$ will be the map 
induced by $p \times \id_{P (\BA_f) / K} \,$, where $p$ is a
$P(\BQ)$-equivariant map from $\FX^{BS}$ to $\FX^*$.
We aim for a map respecting the stratifications of both $\FX^{BS}$ and $\FX^*$.
The space $\FX^{BS}$ is stratified by its \emph{faces} $e(Q)$, one for
each parabolic subgroup $Q$ of $P$, while the stratification of $\FX^*$ is indexed
(only) by the admissible parabolics of $P$. As a preliminary step
to the construction of $p$, one therefore needs
a rule, to be denoted $\adm$, associating to a parabolic of $P$ an admisssible one.
This rule is certainly known already (``$Q$ and $\adm(Q)$ have the same Hermitian factor'';
see e.g.\ \cite[Sect.~5.7, Figure~7]{G}), but its description using the language
of Shimura data is very useful. The result is Theorem~\ref{4D}, which 
for any parabolic $Q$ establishes existence and unicity of an admissible parabolic
subgroup $\adm(Q)$ containg $Q$, and such that the minimal normal subgroup $\adm_{Sh}(Q)$
of $\adm(Q)$ underlying the boundary component associated to $\adm(Q)$, is contained in $Q$:
\[
\adm_{Sh}(Q) \subset Q \subset \adm(Q) \; .
\]
Theorem~\ref{4D} also gives an important characterization of $\adm(Q)$ by a mini\-mality pro\-perty:
if $Q'$ is an admissible parabolic subgroup containing $Q$, then the closure
of the stratum of $\FX^*$ associated to $Q'$ contains the stratum associated to $\adm(Q)$. 
As a first application of minimality, we prove (Corollary~\ref{4I})
that the geodesic action indexed by $Q$ leaves invariant the real part 
$\re_{\adm(Q)}$ with respect to $\adm(Q)$.
In other words, the action only ma\-ni\-fests itself \emph{via} the \emph{imaginary
part with respect to $\adm(Q)$},
as defined in \cite[Sect.~4.14]{P1}. A second application of
minimality then gives an estimation on the variation of the imaginary part
under the geodesic action indexed by $Q$ (Corollary~\ref{4J}).
Corollary~\ref{4I} then allows a rather straightforward
definition of $p$ (Construction~\ref{4K}): each face $e(Q)$ is canonically identified with
the quotient of $\FX$ by the geodesic action; since points of $\FX$, that are conjugate under this
action, only differ by their imaginary parts, they are conjugate under the unipotent radical
of $\adm(Q)$. Therefore, their classes in $\FX^*$ are the same. 
Theorem~\ref{4L} then states that $p$ is continuous; the vital ingredient of the proof
is the estimation from Corollary~\ref{4J}. \\ 

Section~\ref{4a} gives an analysis of the canonical stratifications of the Borel--Serre, and the
Baily--Borel compactification, respectively. Its main result (Theorem~\ref{4aG}) decribes
the interplay of the closures of the canonical strata of the Borel--Serre compactification
$M^K (P,\FX) (\BC)^{BS}$,
and the pre-images under $p^K$ of the canonical strata of the Baily-Borel compactification
$M^K (P,\FX)^*$. 
One particular feature of
Theorem~\ref{4aG} is that these latter pre-images have much better separation
properties than the strata of $M^K (P,\FX)^*$ themselves. \\

Section~\ref{5} is then concerned with the description of the fibres of the maps 
$p: \FX^{BS} \to \FX^*$ 
(Theorem~\ref{5S}) and $p^K$, together with their stratifications (Theorem~\ref{5Z}).
Fix a point $z_0$ of $\FX^*$, and assume that it lies in a boundary component associated
to an admissible parabolic $Q_1$. Given Theorem~\ref{4D}, the parabolics $Q$ in the
pre-image of $Q_1$ under $\adm$ are precisely those satisfying
\[
P_1 \subset Q \subset Q_1 
\]
($P_1:=$ the normal subgroup of $Q_1$ underlying the boundary component containing $z_0$);
in particular, among the intersections of faces $e(Q)$ with the fibre $p^{-1}(z_0)$, the
intersection $p^{-1}(z_0) \cap e(Q_1)$ is of maximal dimension. We prove (Proposition~\ref{5P})
that $p^{-1}(z_0) \cap e(Q_1)$ is actually a space of type $S - \BQ$, 
in the sense of \cite{BS}, associated 
to an algebraic subgroup $C_1$ of $Q_1$, which is expli\-citly defined
in terms of our Shimura data. Theorem~\ref{5S} then states that the inclusion of
$p^{-1}(z_0) \cap e(Q_1)$ into $\FX^{BS}$ extends uniquely to a continuous map
from the manifold with corners $(p^{-1}(z_0) \cap e(Q_1))^{BS}$, identifying the
latter with the full fibre $p^{-1}(z_0)$. The identification is $C_1$-equivariant;
this allows to control the fibres of $p^K$ (Corollary~\ref{5Y}). 
Using the analysis from the previous section,
we then get the stratified versions of our results. \\  

In the final Section~\ref{6}, it is proved (Theorem~\ref{6F})
that over each of the canonical strata of
the Baily--Borel compactification $M^K (P,\FX)^*$, the map $p^K$ is a locally trivial fibration,
provided the level $K$ is neat; actually, the result still holds after restriction
to any canonical stratum of
the Borel--Serre compactification $M^K (P,\FX) (\BC)^{BS}$.
The strategy of proof follows the pattern from Section~\ref{5}: establish the analoguous
result for the map $p: \FX^{BS} \to \FX^*$, taking care to trace
the equivariant behaviour of our splittings, and of their interaction with
the stratifications (Theorem~\ref{6A} and its corollaries). 
Let again $Q_1$ be an admissible parabolic subgroup of $P$. 
Both the stratum $e(Q_1)$ of $\FX^{BS}$
and the stratum of $\FX^*$ associated to $Q_1$ are disjoint unions, indexed by
the boundary components $(P_1,\FX_1)$  
associated to $Q_1$. Fix one such. Its contribution to
$\FX^*$ equals the space $\FX_1/W_1$ in the pure Shimura data
$(P_1,\FX_1)/W_1$ underlying $(P_1,\FX_1)$ ($W_1$ := the unipotent radical
of $P_1$), while the contribution to $e(Q_1) \subset \FX^{BS}$ is the quotient by the geodesic
action on $\FX \cap \FX_1$. The map $p$ maps the class of $x \in \FX \cap \FX_1$ 
modulo the geodesic action to the image of $x$ in $\FX_1/W_1$ under the quotient
map from $\FX_1$ to $\FX_1/W_1$
(see the above discussion of the content of Section~\ref{4}). In order to get
the desired splittings, the language of mixed Shimura data turns out, yet again,
to provide the appropriate setting: in fact, we prove (Proposition~\ref{5B})
that the projection of any mixed Shimura data to the underlying pure Shimura data is
$\BR$-analytically split (this is then applied to the boundary component
$(P_1,\FX_1)$). The splitting is canonical (hence its equivariant
behaviour can be controlled) as soon as a fibre is fixed. A more detailed
analysis of the splitting in the context of boundary components then shows 
(Corollary~\ref{5G}) that it respects the subset $\FX \cap \FX_1$ of $\FX_1$;
therefore it can be composed with the quotient map to $e(Q_1)$,
thus providing the sought-for splitting of the restriction of $p$ to $e(Q_1)$. 
But by Construction~\ref{4K} (see Theorem~\ref{4D}),
the latter is dense in the pre-image by $p$ of the disjoint union $\coprod \FX_1/W_1$
(\emph{i.e.}, of the full stratum of $\FX^*$ associated to $Q_1$),
just as the intersection $p^{-1}(z_0) \cap e(Q_1)$ is dense in the fibre $p^{-1}(z_0)$.
It then turns out that our splitting of $p_{\tei e(Q_1)}$ is compatible
with the identification of $(p^{-1}(z_0) \cap e(Q_1))^{BS}$ with the fibre $p^{-1}(z_0)$
from Theorem~\ref{5S}, allowing for a relative version of this result. 
This yields the splitting of the restriction of $p$ to
the full pre-image $p^{-1}(\coprod \FX_1/W_1)$ (Theorem~\ref{6A}). \\

Even though our main interest lies in the situation where the unipotent radical $W$ 
of the group $P$ is trivial,
part of our results are valid in a more general context. 
This is true in particular for the material contained in Sections~\ref{1}--\ref{3},
the part of Section~\ref{4} concerning the map $\adm$, and the 
analysis of the canonical stratification
of the Borel--Serre compactification (the first half of Section~\ref{4a}). 
In particular, the reader might find it interesting that the Borel--Serre compactification  
of the universal Abelian scheme over Shimura varieties of Abelian type exists,
and can be studied. \\ 

I wish to thank G.~Ancona, M.~Cavicchi, M.~Goresky, S.~Morra 
and J.~Tilouine for useful comments and discussions. 
Two invitations to the \emph{Institut de recherche math\'ematique avanc\'ee} at Strasbourg
helped to speed up the process of writing up; these invitations, along with the very
friendly and stimulating atmosphere were greatly appreciated. Special thanks go to
D.~Blotti\`ere and S.~Morel for the organization of the \emph{Groupe de travail ``vari\'et\'es
de Shimura d'apr\`es R.~Pink''} (Chevaleret, spring term 2004-05), in the course of which 
I had the occasion to present Theorem~\ref{4D} (the definition of the map $\adm$), 
as well as an erroneous version of Theorem~\ref{2H} (the explicit description of
the geodesic action). 
I am indebted to the referee for her or his very meticulous reading of a first
version of the present text.
Last, but not least, it is difficult to appreciate, but impossible
not to acknowledge the role of Covid19 in the making of this article. \\

{\bf Conventions}: For an affine algebraic or real Lie group $H$, let us denote by $H^0$ 
(following \cite{BS}) 
its neutral connected component (in the Zariski topo\-logy if $H$ is algebraic; in the ordinary 
topology if $H$ is a real Lie group). Note that if $H$ is affine algebraic over $\BR$, then
$H^0(\BR)$ might be strictly larger than $H(\BR)^0$. In contrast to \cite{BS}, 
our homogeneous spaces will carry left, not right actions.
Even when the group is $P$, its elements will be denoted by $g$, in order to avoid  
confusion with the map $p: \FX^{BS} \to \FX^*$. Finally (and still following \cite{BS}), 
the reader needs to be prepared to find 
overlined symbols $\overline{\ast}$, whose meaning differs according to the nature
of the symbol $\ast$: when $Q$ is an affine algebraic group, then
$\bar{Q}$ denotes its maximal reductive quotient; when $A$ equals the Lie group
$( \BR_+^* )^r$, for some natural number $r$, then $\bar{A}$ denotes the partial
compactification $(0,+\infty]^r$ obtained by adding the point $+\infty$ in each coordinate;
when $e$ is a canonical stratum of $\FX^{BS}$ or of one of the three (Borel--Serre,
reductive Borel--Serre, 
or Baily--Borel) compactifications, then $\overline{e}$ denotes its closure in the ambient
space. \\


\bigskip
%
%

\section{Homogeneous spaces of type $S'$}
\label{1}



The aim of this section is to clarify the relation between mixed Shimura data \cite{P1} 
and homogeneous spaces of type $S - \BQ$ \cite{BS}. In order to do so, it turns out to 
be useful to first introduce and study a preliminary notion, that of \emph{homogeneous space
of type $S'$}. This study culminates in the first main result of this section (Theorem~\ref{1F}):
provided the Shimura data $(P, \FX)$
satisfy two hypotheses $(+)$ and $(U=0)$ (which will be specified),
any connected component $\FX^0$ of $\FX$ underlies in a canonical way 
a space of type $S'$.
This result implies in particular (Corollary~\ref{1G}) that the action of the stabilizer
$\Stab_{P(\BR)}(\FX^0)$ extends cano\-ni\-cally
to an action 
of the whole of $P(\BR)$ on $\FX^0$. The second main result (Theorem~\ref{1M}) then states
that for this extended action, the homogeneous space $\FX^0$ underlies a canonical structure
of space of type $S - \BQ$. \\

Let us prepare our modification of the main notion of \cite{BS}. We start by recalling a result
due to Mostow.

\begin{Thm}[{Mostow}] \label{Mos}
Let $N$ be a real Lie group with only finitely many connected components, and 
$K$ a ma\-xi\-mal compact subgroup of $N$. Then the canonical morphism 
(of finite groups)
\[
K/K^0 \longto N/N^0
\]
is an isomorphism. 
\end{Thm}

\begin{Proof}
It results from \cite[Thm.~3.2~(1)]{M} that the inclusion
of $K$ into $N$ induces indeed a bijection on the level of connected components.
\end{Proof}

\begin{Lem} \label{1A}
Let $N$ and $Z$ be real Lie groups, equipped with an action of
$Z$ on $N$ by automorphisms (of Lie groups). 
Denote by $Z \ltimes N$ the semi-direct product formed with respect to this action.
Assume $N$ to have only finitely many connected components.
Let $N'$ be a subgroup of $N$ containing $N^0$, and $K$ a ma\-xi\-mal compact subgroup of $N$.
Assume that $N'$ and $K$ are stable under the action of $Z$. 
Then the intersection $N' \cap K$ is maximal compact in $N'$,
and the canonical application of quotients
\[
Z \ltimes N' / Z \ltimes (N' \cap K) \longto Z \ltimes N / Z \ltimes K
\]
is a bijection.
\end{Lem}

\begin{Proof}
We may assume that $Z$ is trivial. 

All maximal compact subgroups of $N$ are conjugate to each other
under $N^0$ \cite[Thm.~3.1~(2)]{M}, hence under $N'$.
Therefore (cmp.~\cite[proof of Thm.~3.2]{M}),
the intersection $N' \cap K$ is maximal compact in $N'$. 

Applying Theorem~\ref{Mos} first to $N$, and then to $N'$, 
we have $K^0 = N^0 \cap K$ and 
\[
N' = N^0 \cdot (N' \cap K) \; .
\]
It follows that the canonical application of quotients
\[
N^0 / K^0 \longto N' / (N' \cap K)
\]
is bijective. The analogous statement being true also for $N$ instead of $N'$, our claim
follows.
\end{Proof}

For a real Lie group $H$ and a closed normal subgroup $N$ of $H$, 
a \emph{connected complement of $N$ in $H$} is a closed connected
subgroup $Z$ of $H$ such that 
\[
H = Z \ltimes N \; .
\]
If such a connected complement of $N$ in $H$ exists, then
we say that \emph{$N$ admits connected complements}.

\begin{Def} \label{1B}
Let $H$ be a real Lie group, and $N$ a closed normal subgroup of $H$. We assume 
$N$ to have only finitely many connected components, and to admit connected complements. 
A \emph{homogeneous space of type $S'$ under $(H,N)$} is a pair 
\[
\bigl( X^0,(K_x)_{x \in X^0} \bigr)
\]
consisting of a connected left 
homogeneous space $X^0$ under some closed subgroup $H'$ of $H$ containing $H^0$,
and a family
$(K_x)_{x \in X^0}$ of maximal compact subgroups $K_x$ of $N$ indexed by the points
of $X^0$, such that
\begin{enumerate}
\item[(i)] for all $x \in X^0$, we have
\[
\Stab_{H'}(x) = Z_x \ltimes (H' \cap K_x) \subset Z_x \ltimes N = H \; ,
\]
for some connected complement $Z_x$ of $N$ in $H$, which leaves
$K_x$ stable under conjugation,
\item[(ii)] for all $h \in H'$ and all $x \in X^0$, we have $K_{hx} = \Int(h) K_x$.
\end{enumerate}
\end{Def}

\begin{Rem} \label{1Ba}
Three observations concerning Definition~\ref{1B} are in order. \\[0.1cm]
(a)~Since $H'$ contains $H^0$, it also contains all connected complements of $N$ in $H$.
In particular, the intersection $H' \cap K_x$ \emph{is} stable under conjugation
by $Z_x$. \\[0.1cm]
(b)~If $\bigl( X^0,(K_x)_{x \in X^0} \bigr)$ is homogeneous of type $S'$ under $(H,N)$, 
then for all $x \in X^0$, we have $\Stab_N(x) = H' \cap K_x$. 
Thus, $\Stab_N(x)$ lies between $K_x^0$ and $K_x$. Thanks to 
Lemma~\ref{1A}, we can be more precise:
the intersection $H' \cap K_x = N' \cap K_x$
is maximal compact in $N'$. 
Hence so is $\Stab_N(x)$. \\[0.1cm]
(c)~Once a maximal compact subgroup $K_x$ of $N$ satisfies condition~\ref{1B}~(i),
for a point $x$ of a connected left 
homogeneous space $X^0$ under $H'$, 
the family $(K_x)_{x \in X^0}$ can be recovered, by \emph{defining} 
\[
K_{hx} := \Int(h) K_x
\]
for $h \in H'$. (We leave it to the reader to verify that this 
\emph{is} well-defined.) The resulting pair 
$\bigl( X^0,(K_x)_{x \in X^0} \bigr)$ is then a homogeneous space of type $S'$ under $(H,N)$.
\end{Rem}

\begin{Prop} \label{1C}
Let $H$ be a real Lie group, and $N$ a closed normal subgroup of $H$. We assume 
$N$ to have only finitely many connected components, and to admit connected complements.
Let $\bigl( X^0,(K_x)_{x \in X^0} \bigr)$ be a homogeneous space of type $S'$ under $(H,N)$. 
Then there is a unique structure of homogeneous space under $H$ on $X^0$,
extending the given action of $H' \subset H$ on $X^0$, and such that
for the induced action of $N$, we have
\[
\Stab_N (x) = K_x \, , \, \forall \, x \in X^0 \; .
\]
In particular, for this extended action,
the subgroups $\Stab_N (x)$ are maximal compact (in $N$). 
\end{Prop}

\begin{Rem} \label{1D}
Proposition~\ref{1C} states that a
homoge\-neous space of type $S'$ under $(H,N)$ is the same as a
homogeneous space under $H$, for which the stabilizers of the induced action of $N$
are maximal compact. In particular, the subgroup $H'$ in Definition~\ref{1B}
need not be specified (as long as its action on $X^0$ is transitive). 
\end{Rem}

\begin{Proofof}{Proposition~\ref{1C}}
Fix a point $x \in X^0$. Put $K := K_x$, $Z := Z_x$, and write
\[
\Stab_{H'}(x) = Z \ltimes (H' \cap K) 
\]
as in Definition~\ref{1B}. The map
\[
H' / \Stab_{H'}(x) \longto X^0 \; , \; [h] \longmapsto hx
\]
is bijective. The group $H'$ contains $Z$ (Remark~\ref{1Ba}~(a)); therefore,
\[
H' = Z \ltimes (H' \cap N) \subset Z \ltimes N = H \; .
\]
Thus, 
\[
H' / \Stab_{H'}(x) = Z \ltimes (H' \cap N) / Z \ltimes (H' \cap K) \; .
\]
By Lemma~\ref{1A} (with $N' := H' \cap N$), the canonical application
\[
H' / \Stab_{H'}(x) \longto
H / Z \ltimes K
\]
is a bijection. 

Therefore, the action of $H'$ on $X^0$ can be extended to $H$ in such a way that 
\[
\Stab_{N}(x) = \Stab_{H}(x) \cap N = K \; .
\]
As for uniqueness of the extension, it suffices to note that
\[
H = Z \ltimes N = Z \ltimes \bigl( (H' \cap N)K \bigr) = H'K \; ,
\]
the second equation being valid as $H' \cap N \supset N^0$, and $N^0K = N$ (Theorem~\ref{Mos}).
\end{Proofof}

\begin{Cor} \label{1Ca}
In the situation of Proposition~\ref{1C}, let $C \subset H$
a closed subgroup with only finitely many connected components, and $x \in X^0$. 
Assume that 
$C \cap K_x$ is maximal compact in $C$.
Then the $C^0$-orbit $C^0x$ of $x$ under the neutral connected component of $C$
is stable under the action of the restriction to $C$ of the
extended action of $H$ on $X^0$ from Proposition~\ref{1C}. 
\end{Cor}

\begin{Proof}
As $C^0$ is contained in  
$\Stab_C (C^0x)$, 
the latter is equal to the disjoint union of certain connected components
of $C$. But $\Stab_C (C^0x)$ 
contains $\Stab_C(x)$ (since $C^0$ is normal in $C$), 
hence $\Stab_{C\cap N}(x) = C \cap K_x$, 
which by our assumption
is maximal compact in $C$. 
Therefore (Theorem~\ref{Mos}) $\Stab_C (C^0x)$ meets
all connected components of $C$. 
\end{Proof}

Let us set up our main example of homoge\-neous spaces of type $S'$.
Let $(P, \FX)$ be \emph{mixed Shimura data} \cite[Def.~2.1]{P1}. 
In particular, $P$ is a connected algebraic linear group over $\BQ$, and
$P(\BR)$ acts on the complex ma\-ni\-fold $\FX$ by complex analytic automorphisms.
Denote by $W$ the unipotent radical of $P$. If $P$ is reductive, \emph{i.e.}, 
if $W=0$, then $(P, \FX)$ is called \emph{pure}. \\

We shall systematically denote by $G$ the maximal reductive quotient $P/W$ of $P$. 
Two additional hypotheses $(+)$ and $(U=0)$ will frequently be imposed on $(P, \FX)$: 
\begin{enumerate}
\item [$(+)$] The neutral connected component $Z (G)^0$ of the center $Z (G)$ of 
$G$ is, up to isogeny, a direct product of a $\BQ$-split torus with a torus 
$T$ of compact type (\emph{i.e.}, $T(\BR)$ is compact) defined over $\BQ$.
\end{enumerate}
\begin{enumerate}
\item [$(U=0)$] The weight $(-2)$-part of $P$ is trivial, in other words \cite[Def.~2.1~(v)]{P1},
the normal subgroup $U$ of $W$ (which is part of the data $(P, \FX)$) is trivial. 
\end{enumerate}
Note that condition $(U=0)$ implies that $\FX$ is actually left homogeneous under $P(\BR)$. \\

For an affine algebraic connected group $L$ defined over $\BQ$, 
denote by ${}^0 \! L$ the normal subgroup
\[
{}^0 \! L := \bigcap_{\chi} \ker \chi^2
\]
of $L$ \cite[Sect.~1.1]{BS}, where $\chi$ runs over all characters $L \to \GQm$.
We then have the following result.

\begin{Prop}[{\cite[Prop.~1.2]{BS}}] \label{1E}
Let $L$ be an affine algebraic connected group defined over $\BQ$.
Then the closed normal subgroup ${}^0 \! L(\BR)$ of $L(\BR)$
contains every compact subgroup of 
$L(\BR)$. It admits connected complements in $L(\BR)$. The set of
such connected complements equals the set of subgroups of the form $T(\BR)^0$, 
for the $\BR$-conjugates $T$ of maximal $\BQ$-split subtori
of the radical $RL$ of $L$. 
\end{Prop}

The case $L=P$ will be of interest for us.
Note that $RP$ is an extension of $RG$ by $W$. 
Under the projection from $P$ to $G$, any maximal $\BQ$-split 
subtorus of $RP$ therefore maps isomorphically to \emph{the} maximal $\BQ$-split subtorus
of $Z(G)$. The analogous statement holds for maximal $\BR$-split subtori.

\begin{Cor} \label{1Ea}
Assume hypothesis $(+)$. 
Let $T$ be a maximal $\BR$-split subtorus of $RP_\BR$. \\[0.1cm]
(a)~The subgroup $T(\BR)^0$ of $P(\BR)$ is a connected complement 
of ${}^0 \! P(\BR)$. \\[0.1cm]
(b)~Let $F$ be a subgroup of $P(\BR)$ containing $T(\BR)^0$. Then 
\[
F = T(\BR)^0 \ltimes \bigl(F \cap {}^0 \! P(\BR) \bigr) 
                        \subset T(\BR)^0 \ltimes {}^0 \! P(\BR) = P(\BR) \; .
\]
If in addition $T(\BR)^0$ lies in the center of $F$, then the semi-direct product is direct:
\[
F = T(\BR)^0 \times \bigl( F \cap {}^0 \! P(\BR) \bigr) \; . 
\]
\end{Cor}

The space $\FX$ has finitely many connected components. 
The stabilizer $\Stab_{P(\BR)}(\FX^0)$ of any such component $\FX^0$
contains $P(\BR)^0$. 
As part of our Shimura data $(P, \FX)$, there is a map
\[
h : \FX \longto \Hom ( \BS_\BC, P_\BC ) \; , \; x \longmapsto h_x \; , \;
\]
where $\BS$ denotes the Deligne torus. It lands in $\Hom ( \BS, P_\BR )$ if $(U=0)$
\cite[Def.~2.1~(ii)]{P1}. \\

Recall that $w: \GRm \to \BS$ denotes the weight co-character of $\BS$.
If the morphism $h_x \circ w$ is defined over $\BR$, then
according to \cite[Def.~2.1~(v)]{P1}, the subgroup $\Cent_{P_\BR}(h_x \circ w)$ of $P_\BR$
surjects onto $G_\BR$ under the canonical projection from $P_\BR$ to $G_\BR$,
with a finite kernel. 
But as the centralizer of a torus, it is connected. Therefore, $\Cent_{P_\BR}(h_x \circ w)$
is a Levi subgroup of $P_\BR$. By definition, it contains the image of $h_x$. 
Denote by 
\[
{}^0 \! \Cent_{P_\BR}(h_x \circ w) \subset \Cent_{P_\BR}(h_x \circ w)
\]
the pre-image of $({}^0 G)_\BR \subset G_\BR$ under the isomorphism
$\Cent_{P_\BR}(h_x \circ w) \isoto G_\BR$. Thus,
\[
{}^0 \! \Cent_{P_\BR}(h_x \circ w) = \Cent_{P_\BR}(h_x \circ w) \cap {}^0 \! P_\BR \; .
\]

\begin{Prop} \label{1Ipre}
Assume hypotheses $(+)$ and $(U=0)$. For all $x \in \FX$, the automorphism
$\Int(h_x(i))$ induces a Cartan involution on ${}^0 \! \Cent_{P_\BR}(h_x \circ w)$.
\end{Prop}

Here, we use the generalization of the notion of Cartan involution to reductive groups
over $\BR$ from \cite[Def.~1.7]{BS}. \\

For the definition of the quotient Shimura data $(P, \FX)/N$ by a normal subgroup $N$
of $P$, we refer to \cite[Prop.~2.9]{P1}. Denote by $P^{der}$ the derived group of $P$.

\medskip

\begin{Proofof}{Proposition~\ref{1Ipre}}
The automorphism $\Int(h_x(i))$ being algebraic, all we need to show is that 
the intersection of its fixed points with ${}^0 \! \Cent_{P(\BR)}(h_x \circ w)$
is maximal compact (in ${}^0 \! \Cent_{P(\BR)}(h_x \circ w)$).

Since $\Cent_{P_\BR}(h_x \circ w)$ is a Levi subgroup, we may prove the statement
on the level of the quotient $(P, \FX)/W$, \emph{i.e.}, we may assume that $(P, \FX) = (G, \FX)$
is pure, in which case $\Cent_{G_\BR}(h_x \circ w) = G_\BR$. 

According to \cite[Def.~2.1~(vi)]{P1}, the involution $\Int(h_x(i))$ induces a Cartan
involution on $G^{der}_\BR$ (which is contained in $({}^0 G)_\BR$). 
We may thus divide out $G^{der}$, \emph{i.e.}, we may
assume that $G$ is a torus. 

On the one hand, the involution $\Int(h_x(i))$ is then trivial. 
On the other hand, 
according to $(+)$, the neutral connected component ${}^0\! G^0$ of 
${}^0\! G$ is a torus of compact type. In other words, the whole of ${}^0\! G(\BR)$
is compact. 
\end{Proofof}

\begin{Def} \label{1Not}
Assume hypotheses $(+)$ and $(U=0)$. For $x \in \FX$, define
\[
L_x := \Cent_{P_\BR}(h_x \circ w) \; ,
\]
\[
{}^0 \! L_x := {}^0 \! \Cent_{P_\BR}(h_x \circ w) \; ,
\]
\[
Z_x := Z_d(L_x) \; , \;
\]
the maximal $\BR$-split subtorus of the center $Z(L_x)$ of $L_x$, and
\[
K_x := \Cent_{L_x(\BR)}(h_x(i)) \cap {}^0 \! L_x(\BR) \; .
\]
\end{Def}

Thus, $L_x$ is a Levi subgroup of $P_\BR$, and $K_x$ is maximal compact in 
${}^0 \! L_x(\BR)$, hence in ${}^0 \! P(\BR)$, hence in $P(\BR)$. 
Note that $Z_x(\BR)$ and $K_x$ centralize
each other; in particular, $K_x$ is stable under conjugation by $Z_x(\BR)^0$.
Our notation is compatible with the
homogeneous structure of $\FX$: for $x \in \FX$ and $g \in P(\BR)$, we have
\[
L_{gx} = \Int(g) L_x \quad \text{and} \quad K_{gx} = \Int(g) K_x \; .
\]
These choices being fixed, the first main result of this section reads as follows.

\begin{Thm}  \label{1F}
Assume hypotheses $(+)$ and $(U=0)$. Let $\FX^0$ be a connected component 
of $\FX$. Then $(\FX^0,(K_x)_{x \in \FX^0})$
is homogeneous of type $S'$ under $(P(\BR),{}^0 \! P(\BR))$.
\end{Thm}

The proof of Theorem~\ref{1F} will be given after Proposition~\ref{1K}.

\begin{Cor} \label{1G}
Under the hypotheses of Theorem~\ref{1F}, 
there is a uni\-que structure of homogeneous space under $P(\BR)$ on $\FX^0$,
extending the action of $\Stab_{P(\BR)}(\FX^0)$,
and such that for the induced action of ${}^0 \! P(\BR)$, we have
\[
\Stab_{{}^0 \! P(\BR)}^{\ext}(x) = K_x
\]
for all $x \in \FX^0$. In particular, for this extended action,
$\Stab_{{}^0 \! P(\BR)}^{\ext}(x)$ is maximal compact
in ${}^0 \! P(\BR)$ (hence, in $P(\BR)$), for all $x \in \FX^0$.
\end{Cor}

\begin{Rem} \label{1rem}
Care needs to be taken when considering actions of $P(\BR)$. There are indeed two of them:
(1)~the action on $\FX$ underlying the Shimura data $(P, \FX)$, which is one by 
complex analytic automorphisms, (2)~the action on $\FX^0$
resulting from Corollary~\ref{1G}, which is one by homeomorphisms. 
They induce the same action by $\Stab_{P(\BR)}(\FX^0)$
(hence by $P(\BR)^0$) on $\FX^0$. But they differ unless $\FX = \FX^0$ is connected. 
Confusion may arise when considering stabilizers of points. For this reason, let us
denote, as in the statement of Corollary~\ref{1G}, by $\Stab_{\argdot}^{\ext}(x)$
the stabilizer of the extended action~(2), and by $\Stab_{\argdot}(x)$ the stabilizer
of the original action~(1). We thus have 
\[
\Stab_{\argdot}(x) = \Stab_{\argdot}^{\ext}(x) \cap \Stab_{\argdot}(\FX^0) 
\subset \Stab_{\argdot}^{\ext}(x)
\]
if $\FX^0$ is the connected component of $\FX$ containing $x$;
in particular, the inclusion $\Stab_{\argdot}(x) 
\subset \Stab_{\argdot}^{\ext}(x)$ might be proper.
\end{Rem}

\begin{Proofof}{Corollary~\ref{1G}}
This results from Theorem~\ref{1F} and Proposition~\ref{1C}.
\end{Proofof}

\begin{Ex} \label{1H}
Consider the Shimura data $(\GL_{2,\BQ},\FH_2)$ \cite[Ex.~2.25, case $g=1$]{P1}, where
$\FH_2 := \BC - \BR$ is the union of the upper and the lower half plane in $\BC$,
equipped with the usual transitive action  $(A,\tau) \mapsto A(\tau)$
of $\GL_2(\BR)$, where for
\[
A = \begin{pmatrix} a & b\\ 
                       c & d
\end{pmatrix} \; ,
\]
we have 
\[
A(\tau) = \frac{a\tau+b}{c\tau+d} \; .
\]
The stabilizer of each of the two connected components is 
\[
\GL_2(\BR)^+ := \{ A \in \GL_2(\BR) \tei \det(A) > 0 \} \; .
\] 
Furthermore,
\[
{}^0 \! \GL_{2,\BQ} = \{ A \in \GL_{2,\BQ} \tei \det(A) = \pm 1 \}
\]
by definition (note that the group of characters of $\GL_{2,\BQ}$ is generated by $\det$).
Denote by $c$ the morphism from $\GL_2(\BR)$ to $\Aut(\FH_2)$ mapping a matrix $A$ to
the identity if $\det(A) > 0$, and to complex conjugation if $\det(A) < 0$. 
We leave it to the reader to show that the action of $\GL_2(\BR)$
from Corollary~\ref{1G} on each of the connected components of $\FH_2$ is given by
\[
(A,\tau) \longmapsto c(A) \bigl( A(\tau) \bigr) \; .
\] 
Note that the matrix
\[
\begin{pmatrix} -1 & 0\\ 
                 0 & 1
\end{pmatrix} 
\]
lies in $\Stab_{\GL_2(\BR)}^{\ext}(\tau)$, for
any purely imaginary element $\tau$ of $\FH_2$, 
while it does not belong to any $\Stab_{\GL_2(\BR)}(\tau)$.
Note also that while the action underlying the Shimura data is one by complex analytic automorphisms,
this is no longer true for the extended action of $\GL_2(\BR)$  
from Corollary~\ref{1G} on each of the connected components of $\FH_2$.
\end{Ex}

Let us prepare the proof of Theorem~\ref{1F}.

\begin{Prop} \label{1L}
Assume $(+)$ and $(U=0)$. Then for all $x \in \FX$,
the subgroup $Z_x(\BR)^0$ of $P(\BR)$ is a connected complement 
of ${}^0 \! P(\BR)$, and
\[
\Stab_{P(\BR)}(x) = Z_x(\BR)^0 \times \Stab_{\, {}^0\! P(\BR)}(x) \subset P(\BR) \; .
\]
\end{Prop}

\begin{Proof}
The Levi subgroup 
$L_x$ contains the image of $h_x$, hence its center $Z(L_x)$
centralizes $h_x$. Any (topologically) connected subgroup of $Z(L_x)(\BR)$
therefore stabilizes $x$ \cite[Cor.~2.12]{P1}; this is the case in particular
for $Z_x(\BR)^0$. 

As the maximal $\BR$-split subtorus of the center of a Levi subgroup,
$Z_x$ is maximal $\BR$-split in the radical $RP_\BR$. From $(+)$ and Corollary~\ref{1Ea}~(a),
we deduce that
\[
P(\BR) = Z_x(\BR)^0 \ltimes {}^0\! P(\BR) \; .
\]
It remains to apply Corollary~\ref{1Ea}~(b) to $F = \Stab_{P(\BR)}(x)$ (which is contained in
$\Cent_{P(\BR)}(h_x)$, hence in $L_x$).
\end{Proof} 

\begin{Prop} \label{1K}
Assume $(+)$ and $(U=0)$. Then for all $x \in \FX$, we have
\[
\Stab_{\, {}^0\! P(\BR)}(x) = \Stab_{\, {}^0\! P(\BR)}(\FX^0) \cap K_x 
= \Stab_{P(\BR)}(\FX^0) \cap K_x \; ,
\]
where $\FX^0$ is the connected component of $\FX$ containing $x$. 
\end{Prop}

\begin{Proof}
Since $K_x$ is contained in ${}^0\! P(\BR)$, the second of the two equations is obvious.
 
As for the first, the inclusion 
``$\Stab_{\, {}^0\! P(\BR)}(x) \subset \Stab_{\, {}^0\! P(\BR)}(\FX^0) \cap K_x$'' 
results from the definition of $K_x \ $: indeed,
an element $g$ of
${}^0\! P(\BR)$ stabilizing $x$ centralizes $h_x$, hence both 
$h_x \circ w$ (meaning that $g \in L_x$, hence $g \in {}^0 \! L_x$)
and $h_x(i)$. 

In order to get the equality we look for, note first that its right hand side
\[
\tilde{K_x} := \Stab_{\, {}^0\! P(\BR)}(\FX^0) \cap K_x
\]
is maximal compact in $\Stab_{\, {}^0\! P(\BR)}(\FX^0)$ (cmp.\ Remark~\ref{1Ba}~(b)).
According to Theorem~\ref{Mos}, the inclusion of $\tilde{K_x}$ into
$\Stab_{\, {}^0\! P(\BR)}(\FX^0)$ therefore induces a bijection on the level of
connected components.

Next, we claim that the neutral connected component $\tilde{K_x^0}$
is contained in the left hand side $\Stab_{\, {}^0\! P(\BR)}(x)$.
Since $\tilde{K_x^0}$ is connected, this claim is equivalent, thanks to \cite[Cor.~2.12]{P1}, to
\[
\tilde{K_x^0} \subset \Cent_{P(\BR)}(h_x) \; ,
\]
and therefore results from the inclusion
\[
\bigl( \Cent_{L_x(\BR)}(h_x(i)) \bigr)^0 \subset \Cent_{P(\BR)}(h_x) \; .
\]
In order to show this latter inclusion, consider the Lie algebra
\[
(\Lie L_x)^{\Ad(h_x(i))} 
\]
of $\bigl( \Cent_{L_x(\BR)}(h_x(i)) \bigr)^0$. The 
subspace of $\Lie L_x$ of invariants 
under $\Ad(h_x(i))$
equals the direct sum of Hodge types $(p,q)$, for $p-q$ divisible by $4$. 
But according to \cite[Def.~2.1~(iv)]{P1}, the only Hodge types possibly occurring
in $\Lie L_x$ are $(-1,1)$, $(0,0)$, and $(1,-1)$
(recall that the group $L_x$ is a Levi subgroup
of $P_\BR$). 
Therefore, $(\Lie L_x)^{\Ad(h_x(i))}$ equals
the Hodge type $(0,0)$, \emph{i.e.}, it equals $(\Lie L_x)^{\Ad(h_x(\BS))}$,
which is contained in $(\Lie P_\BR)^{\Ad(h_x(\BS))}$.

We have proved that
\[ 
\tilde{K_x^0} \subset \Stab_{\, {}^0\! P(\BR)}(x) \subset \tilde{K_x} \; ,
\]
meaning that $\Stab_{\, {}^0\! P(\BR)}(x)$ equals the disjoint union of certain
of the connected components of $\tilde{K_x}$. The desired equation
will therefore be proved once we show that $\Stab_{\, {}^0\! P(\BR)}(x)$ has at least as many connected components
as $\tilde{K_x}$, \emph{i.e.}, as $\Stab_{\, {}^0\! P(\BR)}(\FX^0)$. But the quotient 
\[
\Stab_{\, {}^0\! P(\BR)}(\FX^0) / \Stab_{\, {}^0\! P(\BR)}(x)
\]
is homeomorphic to $\FX^0$, which is connected. Therefore, $\Stab_{\, {}^0\! P(\BR)}(x)$
must meet every connected component of $\Stab_{\, {}^0\! P(\BR)}(\FX^0)$.
\end{Proof}

\forget{

\begin{Prop} \label{1I}
(a)~(\cite[Prop.~2.11]{P1}). The canonical morphism of Shimura data
\[
(P, \FX) \longto (P, \FX)/P^{der} \times (P, h(\FX))
\]
is an embedding (in the sense of \cite[Def.~2.3]{P1}). \\[0.1cm]
(b)~Assume $(U=0)$. Then for all $x \in \FX$, the stabilizer
\[
\Stab_{P^{der}(\BR)}(x) \subset P^{der}(\BR)
\]
of the induced action of $P^{der}(\BR)$ equals the connected component of 
a maximal compact subgroup of $P^{der}(\BR)$. 
\end{Prop}

\begin{Proof} (Cmp.\ \cite[proof of Prop.~2.11]{P1}.)
Put $K' := \Stab_{P^{der}(\BR)}(x)$. According to (a), we have $K' = \Cent_{P^{der}(\BR)}(h_x)$.
Note that $P^{der}$ contains $W$ since the latter is of weight $(-1)$, and the quotient
$P/P^{der}$ is (a torus, hence) reductive; therefore, 
the Shimura data $(P, \FX)/P^{der}$ are pure.
By \cite[Lemma~1.17~(a)]{P1} and hypothesis $(U=0)$, we may assume that $W=0$, \emph{i.e.},
that $(P, \FX) = (G, \FX)$ is pure.

Being the centralizer of a torus, the algebraic subgroup $\Cent_{G^{der}_\BR}(h_x)$
of $G^{der}_\BR$ is connected. Its Lie algebra equals

We thus have
\[
\Lie K' = (\Lie G^{der}_\BR)^{\Ad(h_x)(i)} \; .
\]
But  Thus, 
$K'$ is contained in a maximal compact subgroup $K$ of $G^{der}(\BR)$,
and contains $K^0$. Since every
algebraically connected compact linear group is topologocally connected, $K'$ is connected;
it therefore equals $K^0$. 
\end{Proof}
}
\medskip

\begin{Proofof}{Theorem~\ref{1F}}
In the notation of Definition~\ref{1B}, we have $H = P(\BR)$, $N = {}^0\! P(\BR)$,
$X^0 = \FX^0$, and $H' = \Stab_{P(\BR)}(\FX^0)$. Let $x \in \FX^0$. According to 
Proposition~\ref{1L},
\[
\Stab_{H'}(x) = \Stab_{P(\BR)}(x) = Z_x(\BR)^0 \times \Stab_{\, {}^0\! P(\BR)}(x) \; ,
\]
and $Z_x(\BR)^0$ is a connected complement of ${}^0\! P(\BR)$ in $P(\BR)$. 
We already know that $K_x$ is stable under conjugation by $Z_x(\BR)^0$. 
The validity of the criterion from Definition~\ref{1B}
thus follows from the equation
\[
\Stab_{\, {}^0\! P(\BR)}(x) = \Stab_{P(\BR)}(\FX^0) \cap K_x 
\]
(Proposition~\ref{1K}).
\end{Proofof}

\forget{
For the remark in parentheses, we refer to Proposition~\ref{1E}.

According to Proposition~\ref{1I}~(b), the intersection
\[
\Stab_{\, {}^0\! P(\BR)}(x) \cap P^{der}(\BR)
\]
is the connected component of a maximal compact subgroup of $P^{der}(\BR)$. 
Given the formation of the quotient $(P, \FX)/P^{der}$ \cite[proof of Prop.~2.9]{P1},
which induces surjections on stabilizers,
we may thus assume that $P^{der} = 0$, \emph{i.e.}, that $P$ is a torus.
But then, according to $(+)$, the group ${}^0\! P$ is a torus of compact type.
We conclude by noting that as $\FX$ is necessarily a finite set \cite[Ex.~2.6]{P1},
the subgroup ${}^0\! P(\BR)^0$ of ${}^0\! P(\BR)$ acts trivially on the whole of $\FX$.
}

\begin{Cor} \label{1J}
Assume hypotheses $(+)$ and $(U=0)$. Let $x \in \FX$. Then
\[
\Stab_{P(\BR)}^{\ext}(x) = Z_x(\BR)^0 \times K_x \; .
\]
The subgroup $\Stab_{P(\BR)}^{\ext}(x)$ contains the whole of $Z(L_x)(\BR)$, and
\[
\Stab_{P(\BR)}^{\ext}(x) = Z_x(\BR) \cdot K_x = Z(L_x)(\BR) \cdot K_x \; .
\]
\end{Cor}

\begin{Proof}
The subgroup $\Stab_{P(\BR)}^{\ext}(x)$ contains $\Stab_{P(\BR)}(x)$, hence $Z_x(\BR)^0$
(Proposition~\ref{1L}). According to Corollary~\ref{1Ea}~(b),
\[
\Stab_{P(\BR)}^{\ext}(x) = Z_x(\BR)^0 \ltimes \Stab_{\, {}^0\! P(\BR)}^{\ext}(x) \; ,
\]
which according to Corollary~\ref{1G} equals $Z_x(\BR)^0 \ltimes K_x$.
But $Z_x(\BR)^0$ and $K_x$ centralize each other. This proves the formula
\[
\Stab_{P(\BR)}^{\ext}(x) = Z_x(\BR)^0 \times K_x \; .
\]
Now $K_x$ is maximal compact in $L_x(\BR)$. 
Its intersection with the normal subgroup
$Z(L_x)(\BR)$ is therefore maximal compact, \emph{i.e.}, it equals \emph{the} 
maximal compact subgroup of $Z(L_x)(\BR)$. But this latter group
is a complement of $Z_x(\BR)^0$ in $Z(L_x)(\BR)$. 
Therefore, $Z(L_x)(\BR)$ is indeed contained in $Z_x(\BR)^0 \times K_x \, $, and hence, the
inclusion
\[
Z_x(\BR)^0 \times K_x \subset Z(L_x)(\BR) \cdot K_x 
\]
is an equality. Therefore, so are the inclusions 
\[
Z_x(\BR)^0 \times K_x \subset Z_x(\BR) \cdot K_x \subset Z(L_x)(\BR) \cdot K_x \; .
\]
\end{Proof}

For future reference, let us formulate Propositions~\ref{1L} and \ref{1K}, 
and Corollary~\ref{1J} in the pure case
$(P,\FX)=(G,\FX)$.
Denote by $Z_d:= Z_d(G_\BR)$ the maximal $\BR$-split subtorus of $Z(G_\BR)$. Under hypothesis
$(+)$, it coincides with the base change to $\BR$ of the maximal $\BQ$-split subtorus 
$Z_{d,\BQ}(G)$ of $Z(G)$. 

\begin{Cor} \label{1Lcor}
Assume that $(P,\FX)=(G,\FX)$ is pure. Also
assume $(+)$. \\[0.1cm]
(a)~The subgroup $Z_d(\BR)^0$ is a connected complement 
of ${}^0 G(\BR)$ in $G(\BR)$. \\[0.1cm]
(b)~For all $x \in \FX$,
\[
\Stab_{G(\BR)}(x) = Z_d(\BR)^0 \times \Stab_{\, {}^0 G(\BR)}(x) 
= Z_{d,\BQ}(G)(\BR)^0 \times \Stab_{\, {}^0 G(\BR)}(x) \; ,
\]
and 
\[
\Stab_{\, {}^0 G(\BR)}(x) = \Stab_{G(\BR)}(\FX^0) \cap K_x \; ,
\]
where $\FX^0$ is the connected component of $\FX$ containing $x$. \\[0.1cm]
(c)~For all $x \in \FX$,
\[
\Stab_{G(\BR)}^{\ext}(x) = Z_d(\BR)^0 \times K_x  = Z(G)(\BR) \cdot K_x \; .
\]
\end{Cor}

\begin{Rem}
The following will not be needed in the sequel. If $(P,\FX)=(G,\FX)$, 
and if $(+)$ is valid, then the Cartan
involution $\theta_{K_x}$ associated to $K_x$ by \cite[Def.~1.7]{BS} can be made explicit:
the group $G_\BR$ is an almost direct product of $Z_d$ and ${}^0 G_\BR$.
The restriction of $\theta_{K_x}$ to ${}^0 G(\BR)$ equals $\Int(h_x(i))$, while
the restriction of $\theta_{K_x}$ to $Z_d(\BR)$ equals inversion $x \mapsto x^{-1}$.
Indeed, according to \cite[Def.~1.7]{BS}, the involution $\theta_{K_x}$ \emph{exists},
and stabilizes each normal subgroup of $G_\BR$; this is therefore true in particular
for $Z_d$ and ${}^0 G_\BR$. On the two of them, the restrictions are again Cartan
involutions (as $K_x$ intersects maximally with each normal subgroup of $G(\BR)$).
Now use Proposition~\ref{1Ipre} and the definition of $K_x$ for ${}^0 G_\BR$, and  
the unique Cartan involution for the $\BR$-split torus $Z_d$.
\end{Rem}

Here is the second main result of this section. 

\begin{Thm} \label{1M}
Assume hypotheses $(+)$ and $(U=0)$. Let $\FX^0$ be a connected component 
of $\FX$, equipped with the extended action of $P(\BR)$ from Corollary~\ref{1G}. 
Then the pair $(\FX^0,(L_x)_{x \in \FX^0})$
is a \emph{space of type $S-\BQ$ under $P$} in the sense of \cite[Def.~2.3]{BS}.
\end{Thm}

Note that contrary to \loccit, we use group actions from the left,
not from the right. 
\medskip

\begin{Proofof}{Theorem~\ref{1M}}
The connected normal solvable subgroup $R_{\FX^0}$ of $P_\BR$ from \cite[Def.~2.3, SI]{BS}
is chosen to be equal to $R_d P$, the base change to $\BR$ of the $\BQ$-split radical of $P$.
Thus, $R_{\FX^0}$ is an extension of $Z_d$ by $W_\BR$, and all $Z_x$ are maximal tori
of $R_{\FX^0}$. 

According to Corollary~\ref{1J}, 
\[
\Stab_{P(\BR)}^{\ext}(x) = Z_x(\BR) \cdot K_x 
\]
for all $x \in \FX^0$, where $K_x$ is maximal compact, and normalizes $Z_x(\BR)$ (as it
actually centralizes $Z_x(\BR)$). This proves the validity of axiom \cite[Def.~2.3, SI]{BS}.

We also have $\Stab_{P(\BR)}^{\ext}(x) \subset L_x$ and $L_{gx} = \Int(g) L_x$ by definition,
for all $g \in \Stab_{P(\BR)}(\FX^0)$ and $x \in \FX^0$. In order to 
establish the validity of axiom 
\cite[Def.~2.3, SII]{BS}, it remains to show the formula
\[
L_{gx} = \Int(g) L_x
\]
for all $g \in P(\BR)$ and $x \in \FX^0$, with respect to the extended action of $P(\BR)$. 
But for all $x \in \FX^0$, 
\[
P(\BR) = \Stab_{P(\BR)}(\FX^0) K_x
\]
(Theorem~\ref{Mos}). For $g \in K_x$, we have $gx = x$ (Corollary~\ref{1G}) and
$\Int(g) L_x = L_x$ (as $K_x \subset L_x$).
\end{Proofof}

\begin{Cor} \label{1Mcor}
Assume that $(P,\FX)=(G,\FX)$ is pure. Also
assume $(+)$. Let $\FX^0$ be a connected component of $\FX$, equipped with the
extended action of $G(\BR)$. Then the map $x \mapsto G^{der}(\BR) \cap K_x$ 
identifies $\FX^0$ in a $G(\BR)$-equivariant manner 
with the space of maximal compact subgroups of $G^{der}(\BR)$.
\end{Cor}

\begin{Proof}
This is the content of \cite[Ex.~2.5~(2)]{BS} (observe that the group denoted $A$
in \loccit \ is trivial as $Z(G)(\BR)$ is contained in all stabilizers 
(Corollary~\ref{1Lcor}~(c))).
\end{Proof}

It turns out that the second component of the data $(\FX^0,(L_x)_{x \in \FX^0})$
is determined by the first. More generally, we have the following result.

\begin{Prop} \label{1Mgen}
Let $H$ be an affine algebraic group defined over $\BQ \,$, and
$(X,(L_x)_{x \in X})$ a space of type $S$ under $H$ in the sense of \cite[Def.~2.3]{BS}. Denote by
$R_u H$ the unipotent radical of $H$. 
Assume that for some ($\Longleftrightarrow$ any) $x \in X$, the only element of $R_u H(\BR)$
centralizing $\Stab_{H(\BR)}(x)$ is the neutral element.
Then for all points 
$x$ of $X$, the subgroup $L_x$ of $H_\BR$ is the unique Levi subgroup  
containing $\Stab_{H(\BR)}(x)$.
\end{Prop}

\begin{Proof}
Indeed, any Levi subgroup containing
$\Stab_{H(\BR)}(x)$ is conjugate to $L_x$ by an element $u$ of $R_u H(\BR)$
centralizing $\Stab_{H(\BR)}(x)$ \cite[Lemma~1.3]{BS}.
\end{Proof}

\begin{Rem} \label{1Mgenr}
In Proposition~\ref{1Mgen}, 
assume that $(X,(L_x)_{x \in X})$ is of type $S-\BQ$ under $H$. \\[0.1cm]
(a)~Denote by $R_d H$ the base change to $\BR$ of the $\BQ$-split radical of $H$.
Given the form the stabilizers of points can take \cite[Def.~2.3, SI]{BS}, the condition from Proposition~\ref{1Mgen}
is met under the following hypothesis: 
\begin{enumerate}
\item [$(uL)$] For some ($\Longleftrightarrow$ any) maximal torus $S$ of $R_d H$, the only element of $R_u H(\BR)$
centralizing $S$ is the neutral element.
\end{enumerate}
Note that among such $S$, there are the base changes to $\BR$ of maximal tori  
of the $\BQ$-split radical of $H$. \\[0.1cm]
(b)~The theory of \emph{standard parabolic subgroups} shows that hypothesis~$(uL)$ is met when $H$ is a parabolic
subgroup of a connected reductive group \cite[Thm.~4.15~b)]{BT}. \\[0.1cm]
(c)~Hypothesis~$(uL)$ is met for $H$ equal to $P$, the algebraic group underlying our Shimura data $(P, \FX)$
\cite[Def.~2.1~(v)]{P1}. According to (a), Proposition~\ref{1Mgen} therefore applies to the pair
$(\FX^0,(L_x)_{x \in \FX^0})$ from Theorem~\ref{1M}.
\end{Rem}

\forget{
Proposition~\ref{1Mgen} can be generalized to 
\emph{spaces of type $S$} under a connected affine algebraic group $H$ defined over
$\BQ$ \cite[Def.~2.3]{BS}, provided the connected normal solvable subgroup $R_X$
of $H_{\BR}$ from \cite[Def.~2.3, SI]{BS} can be chosen to contain $R_d H$. One then has
\[
L_x = \Cent_{H_{\BR}} \bigl( R_X \cap \overline{ \Stab_{H(\BR)}(x)} \bigr) \; ,
\]
for all points $x$ of $X$.

torus $S_x$ is 
conjugate to a maximal $\BQ$-split torus. Therefore, its centralizer
$\Cent_{H_{\BR}} (S_x)$ is a Levi subgroup of $H_{\BR}$. Since
it contains $L_x \,$, we necessarily have   
\[
L_x = \Cent_{H_{\BR}} ( S_x) \; .
\]

According to \cite[Rem.~2.2]{BS}, thanks to the hypothesis on connectivity,
the subgroup $L_x$ centralizes $R_d H \cap L_x \,$. We have 
\[
R_d H \cap L_x = S_x
\]
\cite[Lemma~2.1~(iv), (iii)]{BS}. Thus, $L_x$ is a subgroup of
\[
\Cent_{H_{\BR}} ( S_x ) 
\]
containing $\Stab_{H(\BR)}(x)$. Since $L_x$ is a Levi subgroup of $H_\BR$, it is a Levi subgroup
of $\Cent_{H_{\BR}} ( S_x )$. 
, and by 

For $x \in X$,
denote by $\overline{ \Stab_{H(\BR)}(x)}$ the Zariski closure of
$\Stab_{H(\BR)}(x)$.
Put
\[
S_x:= R_d H \cap \overline{ \Stab_{H(\BR)}(x)} \; .
\]
}
\begin{Prop} \label{1Mgenprop}
Let $H$ be an affine algebraic group defined over $\BQ \,$,
and $(X,(L_x)_{x \in X})$ and $(Y,(L_y)_{y \in Y})$ two spaces 
of type $S-\BQ$ under $H$. \\[0.1cm]
(a)~There is an isomorphism
\[
\bigl( X,(L_x)_{x \in X} \bigr) \isoto \bigl( Y,(L_y)_{y \in Y} \bigr)
\]
of spaces of type $S-\BQ$ under $H$, meaning that there is an $H(\BR)$-equivariant
real analytic
isomorphism $\psi: X \isoto Y$, such that $L_{\psi(x)} = L_x$ for all $x \in X$. \\[0.1cm]
(b)~If $H$ is connected, then any two isomorphisms $\psi$ and $\psi'$ as in (a) can be interpolated
by isomorphisms of spaces of type $S-\BQ$ under $H$, meaning that there
is a continuous map
\[
\Psi: X \times [0,1] \longto Y
\]
such that for any $t \in [0,1]$, the map
\[
\psi_t : X \longto Y \; , \; x \longmapsto \Psi(x,t)
\]
is an isomorphism of spaces of type $S-\BQ$, and $\psi_0 = \psi$ and $\psi_1 = \psi'$. \\[0.1cm]
(c)~If $H$ is connected, then the isomorphism in (a) is unique under the following condition:
for some (equivalently, any) point $x$ of $X$, the group of real points
of the center $Z(L_x)$ of $L_x$ is contained in the stabilizer $\Stab_{H(\BR)}(x)$,
and the only element of $R_u H(\BR)$
centralizing $L_x$ is the neutral element.
\end{Prop}

\begin{Proof}
(a): given the form the stabilizers of points can take \cite[Def.~2.3, SI]{BS},
and the fact that all subgroups of this form are conjugate 
to each other under $H(\BR)$ \cite[Lemma~2.1~(iii)]{BS}, the homogeneous
spaces $X$ and $Y$ are isomorphic. We may therefore assume that 
\[
\bigl( Y,(L_y)_{y \in Y} \bigr) = \bigl( X,(L_x')_{x \in X} \bigr) \; ,
\]
for choices of Levi subgroups $L_x'$ possibly different from $L_x \,$, $x \in X$. 
Choose and fix a point $x_0$ of $X$. The Levi subgroups $L_{x_0}$ and $L_{x_0}'$ both contain
$\Stab_{H(\BR)}(x_0)$; therefore \cite[Lemma~1.3]{BS} 
\[
L_{x_0} = u L_{x_0}' u^{-1} \; ,
\]
for some element $u$ of $R_u H(\BR)$ centralizing $\Stab_{H(\BR)}(x_0)$. But  
\[
u L_{x_0}' u^{-1} = L_{u x_0}' \; .
\]
Therefore, the unique $H(\BR)$-equivariant automorphism of $X$ mapping $x_0$ to $u x_0$
is an isomorphism of spaces of type $S-\BQ$.

\noindent (b), (c): given~(a), we may assume that 
\[
\bigl( Y,(L_y)_{y \in Y} \bigr) = \bigl( X,(L_x)_{x \in X} \bigr) \; .
\]
Choose and fix $x_0 \in X$.
An automorphism of $(X,(L_x)_{x \in X})$ maps $x_0$ to a point $x_1$ with the same stabilizer,
and such that $L_{x_1} = L_{x_0} \,$. Any such point $x_1$ determines a unique automorphism.
We therefore need to consider the intersection $I_{x_0}$ in $H(\BR)$ of
the normalizers of $\Stab_{H(\BR)}(x_0)$ and of $L_{x_0} \,$. 

The group $H(\BR)$ is the semi-direct product of $R_u H(\BR)$ and of $L_{x_0}(\BR)$.
Let $u \in R_u H(\BR)$ and $h \in L_{x_0}(\BR)$, and assume that their product
$hu$ normalizes a certain subgroup $D$ of $L_{x_0}(\BR)$. For any element $d$ of $D$,
the commutator
\[   
udu^{-1}d^{-1} = h^{-1}(hu)d(hu)^{-1}hd^{-1}
\]
lies in the intersection of $R_u H(\BR)$ (as $R_u H$ is normal in $H$)
and $L_{x_0}(\BR)$ (as $D \subset L_{x_0}(\BR)$), and is therefore trivial.
Thus, the element $u$ centralizes $D$, and
the element $h$ is contained in the normalizer of $D$ in $L_{x_0}(\BR)$.

Applying this to both $D = L_{x_0}(\BR)$ and $D = \Stab_{H(\BR)}(x_0)$, we see that $I_{x_0}$
is the (semi-)direct product of $\Cent_{R_u H} (L_{x_0})(\BR)$ 
and of the normalizer $N_{x_0}$ of $\Stab_{H(\BR)}(x_0)$ in $L_{x_0}(\BR)$.
Now
\[
\Stab_{H(\BR)}(x_0) = S_{x_0}(\BR)  \cdot K_{x_0} \; ,
\]
for some maximal compact subgroup $K_{x_0}$ of $L_{x_0}(\BR)$,
and some subgroup $S_{x_0}$ of the center $Z(L_{x_0})$ of $L_{x_0}$
\cite[Def.~2.3, Rem.~2.2, Lemma~2.1~(iv)]{BS}.
Given that in connected semi-simple groups, maximal
compact subgroups are equal to their normalizers, we see that
\[
N_{x_0} = Z(L_{x_0})(\BR)  \cdot K_{x_0}
\]
(see \cite[Ex.~2.5]{BS}).

Under the hypthesis made in (c), the group $I_{x_0}$
thus equals $\Stab_{H(\BR)}(x)$ itself, \emph{i.e.}, $x_1 = x_0$. This means
that there are no non-trivial automorphisms of $(X,(L_x)_{x \in X})$.

The more general statement from (b) holds since the quotient
\[
\left( \Cent_{R_u H} (L_{x_0})(\BR) \cdot Z(L_{x_0})(\BR)  \cdot K_{x_0} \right) / \left(S_{x_0}(\BR)  \cdot K_{x_0} \right)
\]
is connected (Theorem~\ref{Mos}).
\end{Proof}

\begin{Rem}
(a)~As the proof shows,
Proposition~\ref{1Mgenprop} holds more ge\-ne\-rally for spaces $(X,(L_x)_{x \in X})$ and $(Y,(L_y)_{y \in Y})$ 
of type $S$ under $H$, provided that they are of type $S$ with respect to the \emph{same} choice of connected normal 
solvable subgroups $R_X = R_Y$ of $H_\BR$ (see \cite[Def.~2.3, SI]{BS}). \\[0.1cm]
(b)~If $(X,(L_x)_{x \in X})$ is of type $S-\BQ$ under $H$, then
the condition ``$Z(L_x) \subset \Stab_{H(\BR)}(x)$''
from Proposition~\ref{1Mgenprop}~(c) can be rephrased in group theoretic terms,
given \cite[Def.~2.3]{BS}:
the neutral connected component of the center of
the maximal reductive quotient of $H$ is, up to isogeny,
a direct product of a $\BQ$-split torus with a torus of compact type
(cmp.\ hypothesis $(+)$ on Shimura data).
\end{Rem}

Let us return to our Shimura data $(P, \FX)$.

\begin{Rem} \label{1N}
It is important to distinguish the two actions of $P(\BR)$ (see Remark~\ref{1rem}): 
(1)~the action on $\FX$ underlying the Shimura data $(P, \FX)$, which is one by 
complex analytic automorphisms, 
(2)~the action on 
each connected component $\FX^0$ of $\FX$, extended from the action
of $\Stab_{P(\BR)}(\FX^0)$ by using Corollary~\ref{1G},
and which is one by real analytic automorphisms. \\[0.1cm]
(a)~The \emph{mixed Shimura varieties} associated to $(P,\FX)$ are
indexed by the
open compact subgroups of $P (\BA_f)$. If $K$ is one such, then the
analytic
space of $\BC$-valued points of the corresponding variety
$M^K := M^K (P,\FX)$ is given as
\[
M^K (\BC) := P (\BQ) \backslash ( \FX \times P (\BA_f) / K ) 
\]
\cite[Def.~3.1]{P1}. As is explained in \cite[Sect.~3.2]{P1}, the analytic space
$M^K (\BC)$ is a finite disjoint union of quotients of the form
\[
\Gamma(g)_+ \backslash \FX^0 \; ,
\]
for connected components $\FX^0$ of $\FX$ and elements $g \in P (\BA_f)$, where
\[
\Gamma(g) := P(\BQ) \cap g K g^{-1} 
\]
is a congruence subgroup of $P(\BQ)$, and 
\[
\Gamma(g)_+ := \Stab_{P(\BQ)}(\FX^0) \cap \Gamma(g)
\]
is of finite index in $\Gamma(g)$.
(The identification of $\Gamma(g)_+ \backslash \FX^0$
with an analytic subset of $M^K (\BC) = P (\BQ) \backslash ( \FX \times P (\BA_f) / K )$ 
maps the class of $x$ to the class of $(x,g)$.)

Since $\Gamma(g)_+$ is contained in $\Stab_{P(\BR)}(\FX^0)$, actions~(1) and (2)
coincide on $\Gamma(g)_+ \,$. Therefore, when forming the 
connected components of the Shimura varieties
associated to $(P, \FX)$, we need not worry about the distinction between (1) and (2). \\[0.1cm]
(b)~Each connected component $\FX^0$ of $\FX$
being a space of type $S-\BQ$ under $P$ according to Theorem~\ref{1M}, it admits
\emph{geodesic actions} \cite[Sect.~3]{BS}. These are indexed by parabolic subgroups $Q$
of $P$. For any such, the geodesic action is an action of 
the group of real points of the center $Z(\bar{Q})$  
of the maximal reductive quotient $\bar{Q}$ of $Q$. 
Its definition involves lifts (indexed by $x \in \FX^0$) of $Z(\bar{Q})(\BR)$ to subgroups
of $P(\BR)$, and the action by these subgroups induced from
the action~(2) giving rise to the $S-\BQ$-structure \cite[Sect.~3.2]{BS}. 

The construction of the partial compactification $\bar{\FX}^0$ of $\FX^0$
from \cite[Sect.~4, 5, 7]{BS} only uses the geodesic
actions of the neutral connected
components $Z(\bar{Q})(\BR)^0$ of the $Z(\bar{Q})(\BR)$. All these lift
to subgroups of $P(\BR)^0$, 
hence of $\Stab_{P(\BR)}(\FX^0)$; therefore, as far as the formation of 
$\bar{\FX}^0$ is concerned, there is no difference between actions~(1) and (2).
(However,
care needs to and will be applied when considering the equivariant behaviour of the result.) 
\end{Rem}

\begin{Cor} \label{1P}
Assume hypotheses $(+)$ and $(U=0)$. Let $\FX^{0,i}$ and $\FX^{0,k}$ 
be two connected components of $\FX$.
Then as spaces of type $S-\BQ\,$, the pairs 
$(\FX^{0,i},(L_x)_{x \in \FX^{0,i}})$ and $(\FX^{0,k},(L_x)_{x \in \FX^{0,k}})$ 
are canonically isomorphic. More precisely, 
there is a unique real analytic isomorphism
\[
\psi: \FX^{0,i} \isoto \FX^{0,k} \; ,
\]
that is equivariant for the extended action of $P(\BR)$ (action~(2)
in the notation used in Remark~\ref{1N}), and such that 
\[
L_{\psi(x)} = L_x \; , \; \forall \; x \in \FX^{0,i} \; .
\]
\end{Cor}

\begin{Proof}
Apply Proposition~\ref{1Mgenprop}~(c), Corollary~\ref{1J} and \cite[Def.~2.1~(v)]{P1}.
\end{Proof}

\begin{Rem}
The isomorphism $\psi$ from Corollary~\ref{1P} can be described explicitely.
For $g \in P(\BR)$ and $x \in \FX$, let us denote by $gx$ the image of $x$
under the action underlying the Shimura data (action~(1)
in the notation used in Remark~\ref{1N}), and by $g^{\ext}x$ the image of $x$
under the extended action on the connected component containing $x$
(action~(2)
in the notation used in Remark~\ref{1N}).
For $g \in P(\BR)$, define a real analytic
automorphism $\psi_g$ of $\FX$ by the formula $x \mapsto g (g^{\ext})^{-1} x$.  
We leave it to the reader to show: (i)~the map $\psi_g$ only depends on
the class of $g$ in $P(\BR)/P(\BR)^0$; more precisely,
its restriction to any connected component $\FX^0$ of $\FX$ only depends on 
the class of $g$ in $P(\BR)/\Stab_{P(\BR)}(\FX^0)$,
(ii)~for all $x \in \FX$, we have $K_{\psi_g(x)} = K_x$,
(iii)~the map $\psi_g$ is equivariant for the extended action of $P(\BR)$ 
(to prove~(iii), use
the formula $P(\BR) = P(\BR)^0 K_x$ (Theorem~\ref{Mos}), (i), and (ii)).

The isomorphism $\psi$ is then defined as the restriction of $\psi_g$
to $\FX^{0,i}$, for any $g \in P(\BR)$ satisfying $\FX^{0,k} = g \FX^{0,i}$.
\end{Rem}

\begin{Rem} \label{1O}
There are variants of all results of this section when hypothesis $(U=0)$
is not satisfied. They concern the subspace
\[ 
\re (\FX) := \{ x \in \FX \; , \; h_x : \BS_\BC \longto P_\BC 
\quad \text{is defined over} \; \BR \} 
\]
of $\FX$
(which is the best possible replacement as $\FX$ is not homogeneous under $P(\BR)$, 
but only
under $P(\BR)U(\BC)$).

By contrast, hypothesis $(+)$ cannot be avoided. 
\end{Rem}


\bigskip
%
%

\section{Geodesic action}
\label{2}



Throughout the section, we fix mixed Shimura data $(P,\FX)$. 
The aim of the present section is to give an
explicit description, in terms of the data $(P,\FX)$, of the geodesic action
(Theorem~\ref{2H}). \\

In order to do so, we first need to study the passage from $P$ to parabolic subgroups $Q$,
and the rule \cite[Sect.~2.7~(3)]{BS} of fixing Levi subgroups of the latter. 
This study will culminate in Theorem~\ref{2F}, and its Corollary~\ref{2G}.

\begin{Prop} [{\cite[Prop.~1.8]{BS}}] \label{2A}
Let $K$ a maximal compact subgroup of $P(\BR)$,
and $L$ a Levi subgroup of $P_\BR$, such that $K \subset L(\BR)$.
Then there is a canonical map $Q \mapsto L_{K,Q}$ 
on the set of those closed subgroups $Q$ of $P_\BR$,
for which $Q(\BR)$ acts transitively on the space of maximal compact subgroups of $P(\BR)$. 
The map associates
to any such $Q$ a Levi subgroup $L_{K,Q}$ of $Q$. The latter is characterized by the following requirements:
\begin{enumerate}
\item[(i)] $L_{K,Q}$ is contained in $L$,
\item[(ii)] $L_{K,Q}$ is stable under the Cartan involution 
$\theta_K$ of $L$ with respect to $K$.
\end{enumerate}
Explicitly, we have
\[
L_{K,Q} = (Q \cap L) \cap \theta_K(Q \cap L) \; .
\]
\end{Prop}

\begin{Cor} \label{2B}
Let $K$ a maximal compact subgroup of $P(\BR)$,
and $L$ a Levi subgroup of $P_\BR$, such that $K \subset L(\BR)$. \\[0.1cm]
(a)~Let $Q$ be a closed subgroup of $P_\BR$, for which $Q(\BR)$ acts transitively 
on the space of 
maximal compact subgroups of $P(\BR)$,
and $g \in P(\BR)$. Then writing $L' := \Int(g)(L)$, we have
\[
\Int(g)(L_{K,Q}) = L'_{\Int(g)(K),\Int(g)(Q)} \; .
\]
(b)~Let $Q_1$ and $Q_2$ be two closed subgroups of $P_\BR$, and assume that 
$Q_1 (\BR) \cap Q_2 (\BR)$
acts transitively on the space of maximal compact subgroups of $P(\BR)$.
Then
\[
L_{K,Q_1 \cap Q_2} = L_{K,Q_1} \cap L_{K,Q_2} \; .
\]
\end{Cor}

In the context of the geodesic action, Proposition~\ref{2A} 
and Corollary~\ref{2B} will only be applied
to the base change to $\BR$ of parabolic subgroups of $P$ (defined over $\BQ$). \\

If $(P,\FX)$ satisfy hypotheses
$(+)$ and $(U=0)$, then
for any point $x$ of $\FX$, we have the Levi subgroup 
\[
L_x = \Cent_{P_\BR}(h_x \circ w) 
\]
of $P_\BR$, and the maximal compact subgroup 
\[
K_x = \Cent_{L_x(\BR)}(h_x(i)) \cap {}^0 \! L_x(\BR) \; .
\]
of $P(\BR)$ (Definition~\ref{1Not}). 

\begin{Def} \label{2C}
Assume hypotheses $(+)$ and $(U=0)$, and
let $x \in \FX$. For any parabolic subgroup $Q$ of $P$, denote by $L_{x,Q}$ the
Levi subgroup $L_{K_x,Q_\BR}$ associated to the Levi subgroup $L = L_x$, 
the maximal
compact subgroup $K_x$ and the parabolic subgroup $Q_\BR$ by the map from
Proposition~\ref{2A}.
\end{Def}

We need to give an explicit formula for $L_{x,Q}$. Since any parabolic subgroup of $P$
is a finite intersection of maximal proper parabolic subgroups, the general formula
will result from Corollary~\ref{2B}~(b), once we know it for maxi\-mal proper 
$Q$. \\

We shall consider a slightly more general situation: let $Q_j$ be an \emph{admissible}
parabolic subgroup of $P$, meaning \cite[Def.~4.5]{P1} that under the projection, first
from $P$ to $G$, and then from $G$ to its adjoint group
$G^{ad}$, the group $Q_j$ is the pre-image
of a product of parabolic subgroups of 
the simple constituents of $G^{ad}$, each either maximal proper or equal to the 
respective simple
constituent. 
To $Q_j$ is associated
a canonical normal subgroup $P_j$ of $Q_j$ \cite[4.7]{P1}. (We added the subscript $j$ 
in order to have the same index for the parabolic subgroup and for its canonical
normal subgroup.) There is a finite collection
of \emph{rational boundary components} $(P_j , \FX_j)$ associated to
$P_j$, and indexed by the
$P_j (\BR)$-orbits in $\pi_0 (\FX)$ \cite[4.11]{P1}.
The $(P_j , \FX_j)$ are themselves mixed Shimura data. 
In particular, there is a distinguished closed subgroup $U_j$ of the unipotent
radical $W_j$ of $P_j$ (the weight $(-2)$-part of $P_j$),
which depends only on $Q_j$, but not on $\FX_j$ \cite[Lemma~4.8]{P1}. 
The proof of \cite[Cor.~4.10]{P1} shows that the
$(P_j, \FX_j)$ satisfy hypothesis $(+)$ since $(P, \FX)$ does. 
(By contrast, unless $Q_j=P$, they \emph{do not} satisfy hypothesis $(U=0) \,$!) \\

For a fixed $P_j (\BR)$-orbit $\CO$ in $\pi_0 (\FX)$, write $\FX^+$ for the disjoint union
of the connected components of $\FX$ in $\CO$, and $(P_j , \FX_j)$ for the Shimura data
associated to $\CO$; according to \cite[Sect.~4.11]{P1}, there is then a $P_j(\BR)$-equivariant map
\[
\FX^+ \longto \FX_j \; .
\]
Taking the disjoint union over all orbits in $\pi_0 (\FX)$, we get a map
\[
\FX \longto \coprod \FX_j \; ;
\]
this map is $Q_j(\BR)$-equivariant \cite[Sect.~4.11]{P1}, injective and complex analytic
\cite[Cor.~4.13]{P1}.

\begin{Def} \label{2D}
Let $x \in \FX$, and $Q_j$ an admissible parabolic subgroup of $P$. 
Denote by $x_j \in \FX_j$ the image of $x$ under the map $\FX^+ \to \FX_j$ 
associated to the $P_j (\BR)$-orbit in $\pi_0 (\FX)$ containing $\pi_0(x)$.
We denote by $\im_j(x) \in U_j(\BC)$ 
the \emph{imaginary part of $x_j$} \cite[Sect.~4.14]{P1}, and
define the \emph{real part of $x$ with respect to $Q_j$} as 
\[
\re_j(x) := \im_j(x)^{-1} x_j \in \FX_j \; .
\]
\end{Def}

Thus, the morphism 
\[
h_{\re_j(x)} : \BS_\BC \longto P_{j,\BC} \subset Q_{j,\BC} \subset P_\BC
\]
is defined over $\BR$. 

\begin{Prop} \label{2E}
Let $x \in \FX$, and $Q_j$ an admissible parabolic subgroup of $P$. 
Assume that the morphism 
\[
h_x : \BS_\BC \longto P_\BC 
\]
is defined over $\BR$ (this is automatically the case if $(U=0)$
\cite[Def.~2.1~(ii)]{P1}).\\[0.1cm]
(a)~Elements of $h_x \circ w(\GRm)$ and of $h_{\re_j(x)}(\BS)$ commute with each other.
In particular,
the co-characters $h_x \circ w$ and 
$h_{\re_j(x)} \circ w$ of $Q_{j,\BR}$ commute with each other. \\[0.1cm]
(b)~We have the relation
\[
\Int(h_x(i)) \bigl( h_{\re_j(x)} \circ w \bigr)
= \bigl( h_{\re_j(x)} \circ w \bigr)^{-1} \cdot (h_x \circ w)^2 
\]
between co-characters of $Q_{j,\BR}$. \\[0.1cm]
(c)~Modulo the unipotent radical of $Q_{j,\BR}$, the product
\[
\bigl( h_{\re_j(x)} \circ w \bigr) \cdot (h_x \circ w)^{-1} 
\]
equals the \emph{co-character canonically associated to the parabolic subgroup $Q_{j,\BR}$ of 
$P_\BR$} \cite[Sect.~4.1]{P1}. In particular, modulo the unipotent radical of $Q_j$,
this product is defined over $\BQ$. \\[0.1cm]
(d)~The element $\im_j(x)$ belongs to $\Cent_{P(\BC)}(h_x \circ w)$.
\end{Prop}

The reader may be used to parameterize parabolic subgroups using \emph{cha\-racters}
rather than co-characters, following \cite[Sect.~8]{BT}. The dictionary is the following:
denoting by $\lambda$ the product from Proposition~\ref{2E}~(c), 
and by $\lambda'$ its composition with the epimorphism from $Q_j$ to $Q_j' := Q_j/W \subset G$,
a simple character $\alpha$
in the root system associated to a fixed minimal parabolic of $G$ contained in $Q_j'$ belongs
to the parameters of $Q_j'$ if and only if $\alpha \circ \lambda' = 0$. Reciprocally,
if $Q_j'$ is parameterized by simple $\alpha_i$, for $i \in I$, then 
the co-character $\lambda'$ canonically associated to $Q_j'$
equals the sum of all fundamental weights associated to the simple roots
$\alpha_i$, for $i \not \in I$.  

\medskip

\begin{Proofof}{Proposition~\ref{2E}}
According to \cite[Prop.~4.6~(b)]{P1}, there is a morphism
\[
\omega = \omega_x: H_0 \longto P_\BR
\]
associated to $x$ ($H_0=$ the reductive subgroup of $\BS \times \GL_{2,\BR}$ 
defined in \cite[Sect.~4.3]{P1}), such that 
\[
h_x = \omega \circ h_0 \; ,
\]
where $h_0$ is the morphisms $\BS_\BC \to H_{0,\BC}$ defined by
\[
h_0 (z_1,z_2) := \Bigl( (z_1,z_2) , 
      \begin{pmatrix} \frac{1}{2}(z_1+z_2)  & -\frac{i}{2} (z_1-z_2)\\ 
                      \frac{i}{2} (z_1-z_2) & \frac{1}{2}(z_1+z_2)
      \end{pmatrix}   \Bigr) \; .
\]
Indeed, the unicity statement in \cite[Prop.~4.6]{P1}
and our hypothesis on $h_x$ being real, imply that $\omega$ \emph{is} defined over $\BR$.
Defining $h_\infty: \BS_\BC \to H_{0,\BC}$ by
\[
h_\infty (z_1,z_2) := \Bigl( (z_1,z_2) , 
      \begin{pmatrix} z_1z_2  & i (1-z_1z_2)\\ 
                      0 & 1
      \end{pmatrix}   \Bigr) \; ,
\]
one then has
\[
h_{x_j} = \omega \circ h_\infty 
\]
by definition of the point $x_j$ \cite[Sect.~4.11]{P1}.

Note that setting 
\[
u_\infty :=  \Bigl( (1,1), \begin{pmatrix} 1  & i \\ 
                                      0 & 1
                      \end{pmatrix}  \Bigr) \; ,
\]
we have 
\[
\Int(u_\infty)^{-1}(h_\infty) : (z_1,z_2) \longmapsto \Bigl( (z_1,z_2) , 
      \begin{pmatrix} z_1z_2  & 0\\ 
                      0 & 1
      \end{pmatrix}   \Bigr) \; ;
\]
in particular, $h_\infty' := \Int(u_\infty)^{-1}(h_\infty)$ is defined over $\BR$.
Hence so is
\[
\omega \circ h_\infty' =
\Int(\omega(u_\infty))^{-1}(\omega \circ h_\infty) = \Int(\omega(u_\infty))^{-1}(h_{x_j}) \; .
\]
We thus have $\im_j(x) = \omega(u_\infty)$, and 
\[
h_{\re_j(x)} = \omega \circ h_\infty' \; .
\]
Statements~(a), (b) and (d) then follow from the corresponding statements on the level 
of $H_0$: indeed,
\[
h_0 \circ w : z \longmapsto \Bigl( (z,z) , 
      \begin{pmatrix} z & 0\\ 
                      0 & z
      \end{pmatrix}   \Bigr) 
\]
is central, whence (a) and (d).
Furthermore,
\[ 
h_\infty' \circ w : z \longmapsto \Bigl( (z,z) , 
      \begin{pmatrix} z^2 & 0\\ 
                      0 & 1
      \end{pmatrix}   \Bigr) \; ,
\]
and
\[
\Int(h_0(i)) \bigl( h_\infty' \circ w \bigr)
= \bigl( h_\infty' \circ w \bigr)^{-1} \cdot (h_0 \circ w)^2 \; ,
\]
as both sides map $z$ to
\[
\Bigl( (z,z) , 
      \begin{pmatrix} 1 & 0\\ 
                      0 & z^2
      \end{pmatrix}   \Bigr) \; .
\]
In order to prove part~(c), note first that according to \cite[Prop.~4.6~(b)~(iii)]{P1},
the co-character canonically associated to $Q_j$ is congruent, modulo $W$, to
the product
\[
\bigl( h_{x_j} \circ w \bigr) \cdot (h_x \circ w)^{-1} \; .
\]
But mo\-du\-lo $U_j$, the latter equals 
$( h_{\re_j(x)} \circ w ) \cdot (h_x \circ w)^{-1}$.
\end{Proofof}

From now on, the Shimura data $(P,\FX)$ are assumed to 
satisfy hypotheses
$(+)$ and $(U=0)$. The results from Section~\ref{1} are then at our disposal; in particular
(Theorem~\ref{1M}), each connected component $\FX^0$ 
of $\FX$ is canonically equipped with a structure
of space of type $S-\BQ$ under $P$. 

\begin{Thm} \label{2F}
Let $x \in \FX$, and $Q_j$ an admissible parabolic subgroup of $P$. Then
\[
L_{x,Q_j} = L_x \cap \Cent_{P_\BR} \bigl( h_{\re_j(x)} \circ w \bigr) 
= \Cent_{P_\BR}(h_x \circ w) \cap \Cent_{P_\BR} \bigl( h_{\re_j(x)} \circ w \bigr) \; .
\]
\end{Thm}

\begin{Proof}
First, let us show that the group
\[
L_x^{Q_j} := L_x \cap \Cent_{P_\BR} \bigl( h_{\re_j(x)} \circ w \bigr)
\]
is indeed a Levi subgroup of $Q_{j,\BR}$. As $L_x$ is a Levi subgroup of $P_\BR$,
we may suppose that $(P,\FX) = (G,\FX)$ is pure, in which case $L_x = G_\BR$, and 
\[
L_x^{Q_j} = \Cent_{G_\BR} \bigl( h_{\re_j(x)} \circ w \bigr) \; .
\]
According to \cite[Prop.~4.6~(b)~(iii)]{P1}, $L_x^{Q_j}$ is contained in $Q_{j,\BR}$.
According to \cite[Def.~2.1~(v)]{P1} (applied to the Shimura data $(P_j , \FX_j)$), 
\[
P_{j,\BR} \cap L_x^{Q_j} = \Cent_{P_{j,\BR}} \bigl( h_{\re_j(x)} \circ w \bigr) 
\]
is a Levi subgroup of $P_{j,\BR}$. Furthermore, as $h_{\re_j(x)} \circ w$ factors over
$P_{j,\BR}$, the quotient $Q_{j,\BR} / P_{j,\BR}$ is pure of weight zero under
$h_{\re_j(x)} \circ w$. Thus, the subgroup $L_x^{Q_j}$ of $Q_{j,\BR}$
splits the projection from $Q_{j,\BR}$ to $Q_{j,\BR} / W_{j,\BR}$. But according
to \cite[proof of Lemma~4.8]{P1}, the unipotent radical of $Q_j$ equals $W_j$.

Next, in order to show that $L_x^{Q_j} = L_{x,Q_j}$, we need to establish that
$L_x^{Q_j}$ is stable under the Cartan involution $\theta_{K_x}$ of $L_x$ with respect to $K_x$.
In order to do so, we may again suppose that $(P,\FX) = (G,\FX)$ is pure. 
On the one hand, the
center of $G_\BR$ is obviously contained in $L_x^{Q_j}$. On the other hand, it is stable
under $\theta_{K_x}$ \cite[Prop.~1.6]{BS}. Therefore, we may even suppose that $G = G^{ad}$,
in which case $\theta_{K_x} = \Int(h_x(i))$.

Let then $q$ be an element of $L_x^{Q_j}$.
We have
\[
\Int \bigl( \theta_{K_x}(q) \bigr) \bigl( h_{\re_j(x)} \circ w \bigr)
= \Int \bigl( h_x(i) \bigr) \Int(q) \Int \bigl(h_x(i) \bigr) 
\bigl( h_{\re_j(x)} \circ w \bigr) \; ,
\]
which according to Proposition~\ref{2E}~(b) equals
\[
\Int \bigl( h_x(i) \bigr) \Int(q) 
\Bigl( \bigl( h_{\re_j(x)} \circ w \bigr)^{-1} \cdot (h_x \circ w)^2 \Bigr) 
= \Int \bigl( h_x(i) \bigr) \Int(q) 
\bigl( h_{\re_j(x)} \circ w \bigr)^{-1}  
\]
(note that $h_x \circ w$ is trivial as it lands in the center \cite[Def.~2.1~(iii)]{P1},
and $Z(G) = Z(G^{ad})$ is trivial). But $\Int(q)$ fixes 
$h_{\re_j(x)} \circ w$, hence its inverse, and
\[
\Int \bigl( \theta_{K_x}(q) \bigr) \bigl( h_{\re_j(x)} \circ w \bigr)
= \Int \bigl( h_x(i) \bigr) 
\bigl( h_{\re_j(x)} \circ w \bigr)^{-1} \; ,
\]
which by another application of Proposition~\ref{2E}~(b) equals $h_{\re_j(x)} \circ w$.
Altogether,
\[
\theta_{K_x}(q) \in \Cent_{G_\BR} \bigl( h_{\re_j(x)} \circ w \bigr) 
= L_x^{Q_j} \; .
\]
\end{Proof}

\begin{Cor} \label{2G}
Let $x \in \FX$, and $Q$ a parabolic subgroup of $P$. 
Write $Q = \bigcap_{j \in J} Q_j$, where the $Q_j$ are admissible parabolics. \\[0.1cm]
(a)~We have
\[
L_{x,Q} = L_x \cap \bigcap_{j \in J} 
           \Cent_{P_\BR} \bigl( h_{\re_j(x)} \circ w \bigr) 
= \Cent_{P_\BR}(h_x \circ w) \cap 
           \bigcap_{j \in J} \Cent_{P_\BR} \bigl( h_{\re_j(x)} \circ w \bigr) \; .
\]
(b)~The co-characters $h_x \circ w$ and $h_{\re_j(x)} \circ w$, $j \in J$, 
of $P_\BR$ all factor over the center $Z(L_{x,Q})$ of $L_{x,Q}$. In particular, they 
commute with each other.
\end{Cor}

\begin{Proof}
Part~(a) follows from Theorem~\ref{2F} and Corollary~\ref{2B}~(b).

As for part~(b), fix a maximal $\BR$-split subtorus $T$ of $P_\BR$ contained in $L_{x,Q}$.
Given the definition of $L_{x,Q}$, the torus $T$ centralizes 
$h_x \circ w$ and all $h_{\re_j(x)} \circ w$. Therefore, 
the torus $T$, together with the image of $h_x \circ w$, generates
an $\BR$-split subtorus of $P_\BR$, which by maximality must be
equal to $T$. In other words, the image of $h_x \circ w$ is contained in $T$.
One proceeds in the same way with each of the $h_{\re_j(x)} \circ w$, $j \in J$;
their images are thus all contained in $T$.  
In particular, the $h_x \circ w$ and $h_{\re_j(x)} \circ w$, $j \in J$, 
all factor over $L_{x,Q}$, and they commute with each other. 
Again by definition, elements of $L_{x,Q}$ commute with elements in the image
of each of the $h_x \circ w$ and $h_{\re_j(x)} \circ w$, $j \in J$; the latter
therefore belong to $Z(L_{x,Q})$. 
\end{Proof}

The formula in Corollary~\ref{2G} is valid in particular when $Q$ is admissible parabolic,
and the $Q_j$ are maximal proper. This gives a hint on ``how to count correctly'' the
number of co-characters which are necessary to determine the Levi subgroups, and more
precisely, to span the split part of their centers. \\

For a parabolic subgroup $Q$ of $P$, denote by $\bar{Q}$ its maximal reductive quotient,
\emph{i.e.}, the quotient by its unipotent radical
(for $Q=P$, we thus have $\bar{P} = G$), and by $\pi_Q$ the canonical
epimorphism from $Q$ to $\bar{Q}$. Denote by
$Z_{d,\BQ}(\bar{Q})$ the maximal $\BQ$-split subtorus of the center $Z(\bar{Q})$
of $\bar{Q}$. Note that an inclusion of parabolic subgroups $Q \subset Q'$ induces
a canonical monomorphism $Z(\bar{Q}') \into Z(\bar{Q})$.

\begin{Def} [{\cite[Sect.~4.2]{BS}}] \label{2Ga}
Let $Q$ be a parabolic subgroup of $P$. Define
\[
S_Q := Z_{d,\BQ}(\bar{Q}) / Z_{d,\BQ}(G) 
\]
and 
\[
A_Q := S_Q(\BR)^0 \; .
\]
\end{Def}

For the formation of the Borel--Serre compactification, we need to control the geode\-sic action
of the groups $A_Q$.
The main result of this section reads as follows.

\begin{Thm} \label{2H}
Let $Q$ be a parabolic subgroup of $P$. 
Write $Q = \bigcap_{j \in J} Q_j$, where the $Q_j$ are admissible parabolics. \\[0.1cm]
(a)~Let $j \in J$. The co-character 
\[ 
\pi_Q \circ h_{\re_j(x)} \circ w
\]
of $\bar{Q}$ does not depend on the choice of $x \in \FX$, is defined over 
$\BQ \,$, and factors over $Z_{d,\BQ}(\bar{Q})$. \\[0.1cm]
(b)~If the $Q_j$ are pairwise different maximal 
proper parabolics, then we have
\[
Z_{d,\BQ}(\bar{Q}) = Z_{d,\BQ}(G) \times 
          \prod_{j \in J} \bigl( \pi_Q \circ h_{\re_j(x)} \circ w (\GQm) \bigr) \; .
\]
(c)~Assume that the $Q_j$ are pairwise different maximal 
proper parabolics. 
Let $x \in \FX$. Then the subtorus of $Z(L_{x,Q})$ mapping
isomorphically to the base change to $\BR$ of $Z_{d,\BQ}(\bar{Q})$ under $\pi_Q$ equals
\[
Z_x \times 
          \prod_{j \in J} \bigl( h_{\re_j(x)} \circ w (\GRm) \bigr) \subset Z(L_{x,Q}) \; .
\]
\end{Thm} 

\begin{Proof}
The co-character
\[
\pi_P \circ h_x \circ w : \GRm \longto G_\BR
\]
is central \cite[Def.~2.1~(iii)]{P1}, hence by $(+)$ defined over $\BQ \,$.
Also, because of transitivity of the action of $P(\BR)$ on $\FX$, it does not depend
on the choice of $x$. \emph{A fortiori}, the same is true for 
$\pi_{Q_j} \circ h_x \circ w$, $j \in J$.

\noindent (a): for each $j$, the co-character 
$\pi_{Q_j} \circ h_{\re_j(x)} \circ w$ does not depend on the choice of 
$x \in \FX$, and is defined over $\BQ \,$, as follows from Proposition~\ref{2E}~(c),
and the above observation. \emph{A fortiori}, the same is true for 
$\pi_Q \circ h_{\re_j(x)} \circ w$. By Corollary~\ref{2G}~(b), the latter
lands in  $Z(\bar{Q})$, hence in $Z_{d,\BQ}(\bar{Q})$. 

\noindent (b): fix a minimal parabolic subgroup $B$ of $G$ contained in $Q' := Q/W$,
and a maximal $\BQ$-split torus $T$ of $B$. Let $\Delta$ be the basis of the root system
associated to the choices of $B$ and $T$. As the $Q_j$ are assumed maximal proper,
the co-character 
\[
\lambda_j' : \GQm \longto Q_j':= Q_j/W \longinto G
\]
canonically associated to $Q_j'$ is the fundamental weight associated to 
a unique element $\alpha_j$ of $\Delta$. In the normalization of \cite[Sect.~5.12]{BT},
the parabolic $Q_j'$ is thus parametrized by $\Delta - \{ \alpha_j \}$,
and $Q'$, by $\Delta - \{ \alpha_j \tei j \in J \}$. This means precisely that
$Z_{d,\BQ}(\bar{Q})$ equals the direct product of $Z_{d,\BQ}(G)$ and 
of the images of the co-characters induced by $\lambda_j'$, $j \in J$. Since by
Proposition~\ref{2E}~(c), the co-characters $\lambda_j'$ and $\pi_Q \circ h_{\re_j(x)} \circ w$
are congruent modulo $Z_{d,\BQ}(G)$, we have indeed
\[
Z_{d,\BQ}(\bar{Q}) = Z_{d,\BQ}(G) \times 
          \prod_{j \in J} \bigl( \pi_Q \circ h_{\re_j(x)} \circ w (\GQm) \bigr) \; .
\]
(c): use (b), and Corollary~\ref{2G}~(b).
\end{Proof}

According to Theorem~\ref{2H}~(b), the collection of co-characters 
$\pi_Q \circ h_{\re_j(x)} \circ w$, $j \in J$, provides a parameterization
\[
\prod_{j \in J} \bigl( \pi_Q \circ h_{\re_j(x)} \circ w \bigr) :
\bigl( \BR_+^* \bigr)^J \isoto A_Q
\]
provided that the $Q_j$, $j \in J$, are the pairwise different maximal proper parabolic
subgroups of $P$ containing $Q$. Actually, thanks to Theorem~\ref{2H}~(a),
this parameterization does not depend on $x$.

\begin{Def} \label{Ppar}
Let $Q$ be a parabolic subgroup of $P$. Let $Q_j$, $j \in J$ be the pairwise different maximal proper parabolic subgroups of $P$ containing $Q$.
Define
\[
par_Q : \bigl( \BR_+^* \bigr)^J \isoto A_Q
\]
as the product
\[
\prod_{j \in J} \bigl( \pi_Q \circ h_{\re_j(\argdot)} \circ w \bigr) :=
\prod_{j \in J} \bigl( \pi_Q \circ h_{\re_j(x)} \circ w \bigr) \; ,
\]
for any $x \in \FX$.
\end{Def}

\begin{Rem} \label{par}
The parameterization $par_Q$ from Definition~\ref{Ppar}
coincides with that used in \cite[Sect.~4.2]{BS} up to certain
positive powers $n_j \,$, $j \in J$.

More precisely, for $x \in \FX$ and
$q \in \BR^*_+ \,$,
consider the element
\[
\underline{q}_j := (1,1,\ldots,1,q,1,\ldots,1) \quad\quad \text{ ($q$ in position $j$)} \; ,
\]
of $(\BR^*_+)^J$, viewed as an element of $A_Q$ under the
parameterization of \cite[Sect.~4.2~(2)]{BS},
and the result $\underline{q}_j \cdot x$ of the geodesic action of $\underline{q}_j$ on $x$
\cite[Sect.~3.2]{BS}. Then
\[
\underline{q}_j \cdot x = \bigr( h_{\re_j(x)} \circ w(q^{1/n_j}) \bigr)x
\]
(for the action of $h_{\re_j(x)} \circ w(q^{1/n_j}) \in P(\BR)$ on $x$ given by the
Shimura data). In other words, we have the formula
\[
\underline{q}_j \cdot x =
\left( par_Q (1,1,\ldots,1,q^{1/n_j},1,\ldots,1) \right) \cdot x \; ,
\]
which gives the comparison of the effect of the parametrization
of \cite[Sect.~4.2~(2)]{BS} and the parameterization $par_Q$ from Definition~\ref{Ppar}.

Indeed, according to Theorem~\ref{2H}~(b),
the quotient $S_Q= Z_{d,\BQ}(\bar{Q}) / Z_{d,\BQ}(G)$ admits two morphisms: the isogeny
\[
\prod_{j \in J} \bigl( \pi_Q \circ h_{\re_j(x)} \circ w \bigr) :
\BG^J_{m,\BQ} \longto S_Q \; ,
\]
or equivalently (by Proposition~\ref{2E}~(c)), $\prod_{j \in J} \lambda_j'$,
with target $S_Q$, and the isogeny
\[
\prod_{j \in J} \alpha_j : S_Q \longto \BG^J_{m,\BQ}
\]
with source $S_Q$. Both induce isomorphisms between $A_Q = S_Q(\BR)^0$ and $(\BR^*_+)^J$.
In \cite[Sect.~4.2~(2)]{BS}, the inverse of $\prod_{j \in J} \alpha_j$
is used in order to parameterize the geodesic action. Thus, putting
$n_j:= \alpha_j \circ \lambda_j'$, we have
\[
\bigl( \prod_{j \in J} \alpha_j \bigr)^{-1} (q_j)_j
= \prod_{j \in J} \bigl( \pi_Q \circ h_{\re_j(x)} \circ w (q_j^{1/n_j})_j \bigr) \in A_Q
\]
for all $(q_j)^j \in (\BR^*_+)^J$.
In particular,
\[
\underline{q}_j = \pi_Q \circ h_{\re_j(x)} \circ w (q^{1/n_j}) \; .
\]
According to the definition of the geodesic action \cite[Sect.~3.2]{BS},
\[
\underline{q}_j \cdot x = a_x x \; ,
\]
where $a_x \in Z(L_{x,Q})(\BR) \subset P(\BR)$ is the unique lift of
$\underline{q}_j$ under $\pi_Q$. Our claim thus follows from Theorem~\ref{2H}~(c)
(and from Remark~\ref{1N}~(b)).

Note that $n_j$ depends only
on the $P(\BQ)$-conjugation class of $Q_j$ (not on $Q$).
\end{Rem}

Readers not willing to identify the set $J$ of indices
for the maximal proper parabolic subgroups $Q_j$ of $P$ containing $Q$ with a subset of the
natural numbers, by \emph{choosing} an order on the indices, may find comfort in
letting $J$ be \emph{equal} to the set of co-characters $\pi_{Q_j} \circ h_{\re_j(x)} \circ w$
(which determines the $Q_j$, thanks to \cite[Prop.~4.6~(b)~(iii)]{P1}).

\begin{Cor} \label{2Ha}
Let $Q$ be an admissible parabolic subgroup of $P$.
Write $Q = \bigcap_{j \in J} Q_j$, where the $Q_j$, $j \in J$, are pairwise different maximal 
proper parabolics. Let $x \in \FX$, and write $\re_0(x)$ for the real part of $x$
with respect to $Q$. \\[0.1cm]
(a)~Modulo $Z_x$, we have the relation
\[
h_{\re_0(x)} \circ w = \prod_{j \in J} \bigl( h_{\re_j(x)} \circ w \bigr) \circ \Delta_\BR
\]
of co-characters of $Z(L_{x,Q})$, 
where $\Delta_\BR$ denotes the diagonal inclusion of $\GRm$ into 
$\BG^J_{m,\BR}$. \\[0.1cm]
(b)~Modulo $Z_{d,\BQ}(G)$, we have the relation
\[
\pi_Q \circ h_{\re_0(x)} \circ w = 
\prod_{j \in J} \bigl( \pi_Q \circ h_{\re_j(x)} \circ w \bigr) \circ \Delta
\]
of co-characters of $Z_{d,\BQ}(\bar{Q})$, 
where $\Delta$ denotes the diagonal inclusion of $\GQm$ into 
$\BG^J_{m,\BQ}$.
\end{Cor}

\begin{Proof}
Apply Corollary~\ref{2G}~(b) twice, according to the representations $Q = Q$ and 
$Q = \bigcap_{j \in J} Q_j$ of $Q$ as intersection of admissible parabolics. 
We obtain that the images of all co-characters $h_{\re_0(x)} \circ w$ and 
$h_{\re_j(x)} \circ w$, $j \in J$, are indeed contained in $Z(L_{x,Q})$.

Given Theorem~\ref{2H}~(a), parts~(a) and (b) of the statement are equivalent.
As for (b), use Proposition~\ref{2E}~(c), the fact that $\pi_P \circ h_x \circ w$
is central \cite[Def.~2.1~(iii)]{P1}, and the relation 
$\lambda_0 = \prod_{j \in J} \lambda_j \circ \Delta$ between co-characters
canonically associated to $Q$, and to the $Q_j$, respectively.
\end{Proof}

\begin{Rem}
For the improper parabolic $Q=P$, we have $J = \emptyset$, $\re_0(x) = x$, 
and Corollary~\ref{2Ha} reduces to (a)~$h_x \circ w$ is central in $L_x$, 
(b)~$\pi_P \circ h_x \circ w$ is central in $G$ \cite[Def.~2.1~(iii)]{P1}.
\end{Rem}

Concerning the parameterization
\[
par_Q = \prod_{j \in J} \bigl( \pi_Q \circ h_{\re_j(\argdot)} \circ w \bigr) :
\bigl( \BR_+^* \bigr)^J \isoto A_Q \; ,
\]
the following will be of use in Section~\ref{3} (cmp.~\cite[proof of Lemma~3.2~(ii)]{BS}).

\begin{Lem} \label{2Ipre}
Let $Q$ be a parabolic subgroup of $P$.
Write $Q = \bigcap_{j \in J} Q_j$, where the $Q_j$, $j \in J$, are pairwise different maximal 
proper parabolics. Let $\gamma \in P(\BQ)$. 
Note that $\Int(\gamma) (Q) = \bigcap_{j \in J} \Int(\gamma) (Q_j)$, and
consider the para\-meterizations of $A_Q$ and of $A_{\Int(\gamma)(Q)}$,
\[
par_Q =  \prod_{j \in J} \bigl( \pi_Q \circ h_{\re_j(\argdot)} \circ w \bigr) :
\bigl( \BR_+^* \bigr)^J \isoto A_Q
\]
and 
\[
par_{\Int(\gamma) (Q)} = \prod_{j \in J} \bigl( \pi_{\Int(\gamma)(Q)} \circ h_{\re_j(\argdot)} 
                                                                        \circ w \bigr) :
\bigl( \BR_+^* \bigr)^J \isoto A_{\Int(\gamma)(Q)} \; ,
\]
respectively. Then $par_Q$ and $par_{\Int(\gamma) (Q)}$ satisfy the relation
\[
par_{\Int(\gamma) (Q)} = \Int(\gamma)(par_Q) \; .
\]
\end{Lem} 

\begin{Proof}
This results from the formula $h_{\re_j(\gamma x)} = \Int(\gamma) \circ h_{\re_j(x)}$, $j \in J$,
$x \in \FX$, whose proof we leave to the reader.
\end{Proof}

\begin{Rem} 
When hypothesis $(U=0)$ is not satisfied, then the results of this section
are valid only for $x \in \re (\FX)$ (cmp.~Remark~\ref{1O}); this principle was made explicit
in Proposition~\ref{2E}.
\end{Rem}


\bigskip
%
%

\section{The manifold with corners $\FX^{BS}$}
\label{3}



We continue to consider mixed Shimura data $(P,\FX)$ satisfying hypotheses
$(+)$ and $(U=0)$. In this section, we first carry out the program from 
\cite{BS}, making explicit the necessary modifications, in order to construct 
a partial com\-pactification $\FX^{BS}$ of $\FX$, using the geodesic action from
Section~\ref{2}. As before, particular care needs to be employed when distinguishing
the natural action of $P(\BR)$ on $\FX$ underlying our Shimura data,
and the action of $P(\BR)$ on every connected component of $\FX$ obtained
by extending the natural action of the stabilizer (Section~\ref{1}). Since the former turns out
to commute with the geodesic action (which is obtained from the latter)
(Proposition~\ref{2L}), 
we may and do define the \emph{Borel--Serre compactification} of the Shimura varieties
associated to $(P,\FX)$ (Definition~\ref{3A}). One of the main results of \cite{BS} then implies
that the Borel--Serre compactification,
as its name suggests, is indeed compact (Theorem~\ref{3B}). 
Along the lines of the construction, we introduce what will turn out to be a key notion
for the applications to cohomology we have in mind: the notion of \emph{contractible map}
between topological spaces (Definition~\ref{2Ja}). The notion is local on the target;
in particular, it is invariant under passage to the base change by a covering.     
The open immersion of $\FX$ into $\FX^{BS}$ is the main example of a contractible map
we have in mind; our Proposition~\ref{2Jc} gives a stratified version of this observation.
The second part
of the section will not be used in the rest of this paper, and is concerned with 
functoriality of $\FX^{BS}$ in $(P,\FX)$. \\

Let $Q$ be a parabolic subgroup, and $Q_j$, $j \in J$, the pairwise different maximal 
proper parabolic subgroups of $P$ containing $Q$. 
By Definition~\ref{Ppar}, we have
\[
par_Q = \prod_{j \in J} \bigl( \pi_Q \circ h_{\re_j(\argdot)} \circ w \bigr) :
\bigl( \BR_+^* \bigr)^J \isoto A_Q \; .
\]
Consider the ``partial compactification''
\[
\bigl( \BR_+^* \bigr)^J \longinto (0,+\infty]^J
\]
obtained by adding the point $+\infty$ in each coordinate.
Define $\bar{A}_Q$ as the quotient of $A_Q \times (0,+\infty]^J$
by the diagonal action of $( \BR_+^* )^J \cong A_Q$. We get
\[
A_Q \longinto \bar{A}_Q
\]
in a canonical way. Note that $A_Q$ acts on $\bar{A}_Q$,
and that $par_Q$ extends to give a parameterization, denoted by the same symbol 
\[
par_Q : (0,+\infty]^J \isoto \bar{A}_Q \; .
\]
Denote by $\infty_Q \in \bar{A}_Q$ the image under $par_Q$ of the point
$(+\infty,\ldots,+\infty)$. The singleton $\{ \infty_Q \}$ 
is the unique closed stratum of the stratification of
$\bar{A}_Q$ by orbits under the action of $A_Q$.

\begin{Def} \label{2I}
Let $Q$ be a parabolic subgroup of $P$. \\[0.1cm]
(a)~The \emph{corner of $\FX$ associated to $Q$} is 
\[
\FX(Q) := \bar{A}_Q \, {}^{A_Q} \! \times \FX \; .
\]
Explicitly, $\FX(Q)$ is the space of equivalence classes of pairs $(a,x)$, for $a \in \bar{A}_Q$
and $x \in \FX$, with respect to the relation 
\[
(a,x) \sim (b,y) \Longleftrightarrow \exists \, \alpha \in A_Q \; , \; 
b = \alpha^{-1} a \;\; \text{and} \;\; y = \alpha \cdot x
\]
(where $\alpha \cdot x$ denotes as before the geodesic action of $\alpha$ on $x$). \\[0.1cm]
(b)~The \emph{face of $\FX$ associated to $Q$} is
\[
e(Q) := A_Q \backslash \FX \; .
\]    
\end{Def} 

\begin{Rem}
Part~(a) of Definition~\ref{2I} 
\emph{is} compatible with the definition from \cite[Sect.~5.1]{BS}: 
first, $\FX$ is the disjoint union of its connected components, each of which
can be canonically endowed with the structure of a space of type $S-\BQ$ under $P$
(Theorem~\ref{1M}),
and therefore carries the geodesic actions of the $A_Q \,$. 
Hence so does $\FX$.
Next, as was already
pointed out, \cite{BS} is based on 
group actions from the right, not from the left. The rule transforming
one type of action into the other being
$a^{-1} \cdot x \leftrightarrow x \cdot a$, we need to add $+\infty$ where \cite[Sect.~5.1]{BS} 
adds the point $0^+$
(cmp.\ \cite[footnote~(4) on p.~320]{Z}). 
\end{Rem}

If $\FX = \coprod \FX^0$ is the decomposition of $\FX$ into 
its finitely many connected components,
then likewise $\FX(Q) = \coprod \FX^0(Q)$ for every parabolic subgroup $Q$ of $P$,
where $\FX^0(Q)$ is the corner of the space $\FX^0$ of type $S-\BQ$ associated to $Q$
by \cite[Sect.~5.1]{BS}. \\

If $Q_1 \subset Q_2$ is an inclusion of parabolic subgroups of $P$, then
$A_{Q_2}$ is in a canonical way a subgroup of $A_{Q_1}$, and the inclusion 
$A_{Q_2} \into A_{Q_1}$ extends (uniquely) to a continuous map 
$\bar{A}_{Q_2} \to \bar{A}_{Q_1}$.
This extension is injective and $A_{Q_2}$-equivariant (cmp.\ \cite[Sect.~4.3 and 4.5]{BS}).  
It follows that the corner $\FX(Q_2)$ admits a canonical continuous map
to the corner $\FX(Q_1)$, which extends the identity on $\FX$. 
Elements $a,b \in \bar{A}_{Q_2}$, and $\alpha \in A_{Q_1}$ 
satisfy the relation $b = \alpha^{-1} a$ only if $\alpha \in A_{Q_2} \,$.
Therefore, the continuous map $\FX(Q_2) \to \FX(Q_1)$ is injective. 
Actually, it is an open immersion of manifolds with corners \cite[Prop.~5.3]{BS}.
Thus, we can imitate the construction from \cite[Sect.~7.1]{BS}.

\begin{Def} \label{2J}
The manifold with corners $\FX^{BS}$ is defined as the disjoint union 
of the $\bar{\FX}^0$, where $\FX^0$ runs through the connected components of $\FX$,
and $\bar{\FX}^0$ denotes the manifold with corners defined for the 
space $\FX^0$ of type $S-\BQ$ under $P$ in \cite[Sect.~7.1]{BS}.
\end{Def}

Thus, the $\bar{\FX}^0$ are the connected components of $\FX^{BS}$,
the $\FX(Q)$ are open submanifolds with corners of $\FX^{BS}$, providing a covering
of $\FX^{BS}$, and there is a stratification
\[
\FX^{BS} = \coprod_R e(R) \; ,
\]
where $R$ runs over all parabolic subgroups of $P$. \\

Except for connectivity of the spaces, 
the results from
\cite[Sect.~5 and 7]{BS} 
essentially carry over \emph{verbatim} to $\FX$.
First example: the stratification of the corner $\FX(Q) = \bar{A}_Q \, {}^{A_Q} \! \times \FX$ 
induced by $\FX^{BS} = \coprod_R e(R)$ equals the stratification induced 
by the orbits under the action of $A_Q$ on $\bar{A}_Q$. 
We have 
\[
\FX(Q) = \coprod_{R \supset Q} e(R) \; ,
\]
where $R$ runs over all parabolics containing $Q$
\cite[5.1~(6)]{BS}. The unique open stratum is $e(P) = \FX$.
The unique closed stratum is $e(Q)$, which is identified with
the subspace of $\FX(Q) = \bar{A}_Q \, {}^{A_Q} \! \times \FX$ of equivalence classes 
$[(\infty_Q,x)]$ of
points of the form $(\infty_Q,x)$, \emph{via} the bijection
\[
A_Q x \longleftrightarrow [(\infty_Q,x)] \; , \; x \in \FX \; . 
\]
Second example: the closure $\overline{e(R)}$ of a face $e(R)$ in $\FX^{BS}$ equals
\[
\overline{e(R)} = \coprod_{R' \subset R} e(R') \; ,
\]
where $R'$ runs over all parabolics contained in $R$ \cite[Prop.~7.3~(i)]{BS}.
Third example:
$\FX^{BS}$ is Hausdorff \cite[Thm.~7.8]{BS}. \\

The open immersion of $\FX$ into $\FX^{BS}$ is of a very specific local nature.
This motivates the following.

\begin{Def} \label{2Ja}
A continuous map $f:X \to Y$ of topological spaces
is called \emph{contractible} if the topology on $Y$
admits a basis $(V_i)_i$, for which the pre-images
$f^{-1}(V_i)$ are contractible, for all $i$.
\end{Def}

The property of a continuous map being contractible is local on its target.
Also, if $f$ and $g$ are contractible, then so is their direct product $f \times g$.
If $f$ is contractible, then the image of $f$ is dense in $Y$.

\begin{Exs} \label{2Jb}
Let $X$ be a topological space. \\[0.1cm] 
(a)~The map $X \onto \{\ast\}$ to the singleton is contractible if and only if
$X$ is contractible. \\[0.1cm]
(b)~The identity $\id_X$ is contractible if and only of $X$ is locally
contractible.  
\end{Exs}

\begin{Prop} \label{2Jc}
Let $R$ be a parabolic subgroup of $P$, and 
$\FX' \subset \FX^{BS}$ a locally closed union of strata
containing $e(R)$, and contained in the closure $\overline{e(R)}$
of $e(R)$ in $\FX^{BS}$.
Then the open immersion $\FX' \into \overline{e(R)}$ is contractible. 
\end{Prop}

\begin{Proof}
We have
$\overline{e(R)} = \coprod_{R' \subset R} e(R')$,
where $R'$ runs over all parabolics contained in $R$.
Our claim being local on $\FX^{BS}$, it suffices to show that for any 
parabolic subgroup $Q$ of $P$ contained in $R$, the immersion
$\FX' \cap \FX(Q) \into \overline{e(R)} \cap \FX(Q)$ is contractible.

Recall the parameterization
\[
par_Q : (0,+\infty]^J \isoto \bar{A}_Q \; .
\]
According to \cite[5.4~(5)]{BS}, there is a homeomorphism
\[
\bar{A}_Q \times e(Q) \isoto \FX(Q)
\]
respecting the stratifications. 
Furthermore \cite[5.4~(6)]{BS}, the pre-image of $\overline{e(R)} \cap \FX(Q)$
under this homeomorphism equals
\[
\bar{A}_{Q,R} \times e(Q) \; , 
\]
for a closed subset $\bar{A}_{Q,R}$ of $\bar{A}_Q$, which under $par_Q$ corresponds to
the subset defined by the conditions ``$x_j = + \infty$'', for all indices $j$ belonging to
a certain subset $J_R$ of $J$. Whence an induced parameterization
\[
par_{Q,R} : (0,+\infty]^{J-J_R} \isoto \bar{A}_{Q,R} \; .
\]
We thus get a commutative diagram 
\[  
\vcenter{\xymatrix@R-10pt{
\bar{A}_{Q,R}' \times e(Q) \ar@{->>}[r]^-\cong \ar@{^{ (}->}[d]_{k \times \id_{e(Q)}} & 
                                                \FX' \cap \FX(Q) \ar@{^{ (}->}[d] \\
\bar{A}_{Q,R} \times e(Q) \ar@{->>}[r]^-\cong & \overline{e(R)} \cap \FX(Q) 
\\}}
\]
for some open immersion $k: \bar{A}_{Q,R}' \into \bar{A}_{Q,R}$
of a union of $A_Q$-orbits in $\bar{A}_{Q,R}$. Under $par_{Q,R}$, this
subset $\bar{A}_{Q,R}'$ corresponds to a finite union of subsets
defined by conditions ``$x_j \ne + \infty$'', for certain indices $j \in J-J_R$.
This shows that $k$ is contractible.
Our claim follows: as $e(Q)$ is locally contractible 
--- in fact, it is diffeomorphic to a Euclidean space 
\cite[Sect.~3.9, Rem.~2.4~(1)]{BS}) --- the map $\id_{e(Q)}$ is contractible.
\end{Proof}

The following principle will be used repeatedly.

\begin{Prop} \label{Par}
In a connected algebraic group, two conjugate parabolic
subgroups, whose intersection remains parabolic, are identical.
\end{Prop}

\begin{Proof}
We may clearly suppose that the connected algebraic group in question is reductive.
Now apply \cite[Sect.~4.3]{BT}.
\end{Proof}

\begin{Prop} \label{2Jd}
Let $R$ be a parabolic subgroup of $P$, and 
$Q \subset R$ a maximal proper parabolic of $R$ (\emph{i.e.}, a
parabolic containing $Q$ and contained in $R$ is equal either to $Q$ or to $R$).
Then the open immersion 
\[
\bigl( \overline{e(R)} - e(R) \bigr) - \bigcup_{\gamma \in R(\BQ)} \overline{e(\gamma Q \gamma^{-1})}
\longinto \bigl( \overline{e(R)} - e(R) \bigr) - \bigcup_{\gamma \in R(\BQ)} e(\gamma Q \gamma^{-1})
\] 
is contractible.   
\end{Prop}

\begin{Proof}
The claim being local on $\FX^{BS}$, it suffices to show that for any 
parabolic subgroup $Q'$ of $P$ contained in $R$, the immersion $k$ of
\[
\bigl( (\overline{e(R)} \cap \FX(Q')) - e(R) \bigr) 
- \bigcup_{\gamma \in R(\BQ)} \overline{e(\gamma Q \gamma^{-1})}
\]
into 
\[
\bigl( (\overline{e(R)} \cap \FX(Q')) - e(R) \bigr) 
- \bigcup_{\gamma \in R(\BQ)} e(\gamma Q \gamma^{-1})
\]
is contractible.

Consider as in the previous proof the parameterization
\[
par_{Q'} : (0,+\infty]^J \isoto \bar{A}_{Q'} 
\]
and the homeomorphism
\[
\bar{A}_{Q'} \times e(Q') \isoto \FX(Q')
\]
respecting the stratifications \cite[5.4~(5)]{BS}, and identifying
\[
\bar{A}_{Q',R} \times e(Q')  
\]
with $\overline{e(R)} \cap \FX(Q')$ \cite[5.4~(6)]{BS}, for
$\bar{A}_{Q',R} \subset \bar{A}_{Q'}$ closed, 
corresponding to the conditions ``$x_j = + \infty$'', for all indices $j$ belonging to
a certain subset $J_R$ of $J$ under $par_{Q'}$. The induced parameterization
of $\bar{A}_{Q',R}$ is
\[
par_{Q',R} : (0,+\infty]^{J-J_R} \isoto \bar{A}_{Q',R} \; .
\]
We get a commutative diagram 
\[  
\vcenter{\xymatrix@R-10pt{
\bar{A}_{Q',R}^\infty \times e(Q') \ar@{->>}[r]^-\cong \ar@{^{ (}->}[d] & 
                                     (\overline{e(R)} \cap \FX(Q')) - e(R) \ar@{^{ (}->}[d] \\
\bar{A}_{Q',R} \times e(Q') \ar@{->>}[r]^-\cong & \overline{e(R)} \cap \FX(Q') 
\\}}
\]
where $\bar{A}_{Q',R}^\infty \subset \bar{A}_{Q',R}$ is the closed subset corresponding to the
condition ``$x_j = + \infty$ for at least one $j \in J-J_R$''.

If $Q'$ is not contained in any $R(\BQ)$-conjugate of $Q$, then 
both source and target of $k$ are equal to 
\[
(\overline{e(R)} \cap \FX(Q')) - e(R) \; ,
\]
which is homeomorphic to $\bar{A}_{Q',R}^\infty \times e(Q')$,
and hence locally contractible \cite[Sect.~3.9, Rem.~2.4~(1)]{BS}.

Let us assume that $Q'$ is contained in $\gamma Q \gamma^{-1}$, for some $\gamma \in R(\BQ)$. 
Proposition~\ref{Par} tells us that there is no other conjugate of $Q$ than $\gamma Q \gamma^{-1}$
containing $Q'$. Replacing $Q$ by $\gamma^{-1} Q \gamma$, we thus get
\[
\bigl( (\overline{e(R)} \cap \FX(Q')) - e(R) \bigr) 
- \bigcup_{\gamma \in R(\BQ)} \overline{e(\gamma Q \gamma^{-1})} = 
\bigl( (\overline{e(R)} \cap \FX(Q')) - e(R) \bigr) - \overline{e(Q)}
\]
and 
\[
\bigl( (\overline{e(R)} \cap \FX(Q')) - e(R) \bigr) 
- \bigcup_{\gamma \in R(\BQ)} e(\gamma Q \gamma^{-1}) = 
\bigl( (\overline{e(R)} \cap \FX(Q')) - e(R) \bigr) - e(Q) \; .
\]
Since $Q$ is maximal proper in $R$, there is an index 
$j_Q \in J-J_R$, such that the above commutative diagram can be completed to
\[  
\vcenter{\xymatrix@R-10pt{
\bar{A}_{Q',R}^{\infty,Q,0} \times e(Q') \ar@{->>}[r]^-\cong \ar@{^{ (}->}[d] & 
                   \bigl( (\overline{e(R)} \cap \FX(Q')) - e(R) \bigr) \cap e(Q) \ar@{^{ (}->}[d] \\
\bar{A}_{Q',R}^{\infty,Q} \times e(Q') \ar@{->>}[r]^-\cong \ar@{^{ (}->}[d] & 
        \bigl( (\overline{e(R)} \cap \FX(Q')) - e(R) \bigr) \cap \overline{e(Q)} \ar@{^{ (}->}[d] \\
\bar{A}_{Q',R}^\infty \times e(Q') \ar@{->>}[r]^-\cong \ar@{^{ (}->}[d] & 
                                     (\overline{e(R)} \cap \FX(Q')) - e(R) \ar@{^{ (}->}[d] \\
\bar{A}_{Q',R} \times e(Q') \ar@{->>}[r]^-\cong & \overline{e(R)} \cap \FX(Q') 
\\}}
\]
where $\bar{A}_{Q',R}^{\infty,Q}$ and $\bar{A}_{Q',R}^{\infty,Q,0}$ are the subsets of
$\bar{A}_{Q',R}$ defined by 
\[
\bar{A}_{Q',R}^{\infty,Q} := \{ a \in \bar{A}_{Q',R}^\infty \; , \; x_{j_Q}(a) = + \infty \} 
\]
and 
\[
\bar{A}_{Q',R}^{\infty,Q,0} := 
      \{ a \in \bar{A}_{Q',R}^\infty \; , \; x_{j_Q}(a) = + \infty \, , \, x_j \ne + \infty
      \; \text{for} \; j \ne j_Q \} \; . 
\]
We leave it to the reader to show that the immersion of the complement of
$\bar{A}_{Q',R}^{\infty,Q}$ into the complement of $\bar{A}_{Q',R}^{\infty,Q,0}$ 
(in $\bar{A}_{Q',R}^\infty$)
is contractible.
Our claim follows as $e(Q')$ is locally contractible  
\cite[Sect.~3.9, Rem.~2.4~(1)]{BS}.
\end{Proof}

\begin{Rem}
(a)~The hypothesis on $Q$ being maximal proper in $R$ is essential in Proposition~\ref{2Jd}. \\[0.1cm]
(b)~As the proofs shows, Propositions~\ref{2Jc} and \ref{2Jd} are valid in the general context of 
spaces of type $S-k$ \cite[Def.~2.3]{BS}.
\end{Rem}

\begin{Rem} \label{2K}
Let $Q$ be a parabolic subgroup of $P$.
In analogy to Remark~\ref{1N}, there are two different actions of $Q(\BR)$ on 
the corner $\FX(Q)$, both of which respect the stratification 
\[
\FX(Q) = \coprod_{Q \subset R \subset P} e(R) \; :
\]
(1)~the action induced by and extending the one on $\FX$ underlying the Shimura data $(P, \FX)$, 
(2)~the action induced by and extending the one on $\FX$ coming from extension 
of the action of $\Stab_{P(\BR)}(\FX^0)$ on each connected component $\FX^0$ of 
$\FX$ \emph{via} Corollary~\ref{1G}. On each connected component
$\FX^0(Q)$, they induce the same action of the subgroup $\Stab_{Q(\BR)}(\FX^0)$
of $Q(\BR)$.
In particular, actions (1) and (2) coincide on $P(\BR)^0 \cap Q(\BR)$. 

For both (1) and (2), recall that by definition,
\[
\FX(Q) = \bar{A}_Q \, {}^{A_Q} \! \times \FX \; ,
\] 
and use the restriction to $Q(\BR)$ of the respective actions 
(1) and (2) on $\FX$. It results that both actions
are by automorphisms of manifolds with corners.

Caution: for these constructions to work,
it is necessary for the geodesic action of $A_Q$ on $\FX$ to
commute with both actions of $Q(\BR)$. For action (2), this is 
the first statement of \cite[Prop.~3.4]{BS}. It thus remains to prove the following result.
\end{Rem}

\begin{Prop} \label{2L}
Let $Q$ be a parabolic subgroup of $P$. Then
the geo\-de\-sic action of $A_Q$
commutes with action (1) of $Q(\BR)$ on $\FX$.
\end{Prop}

\begin{Proof}
Let $a \in A_Q$, $x \in \FX$, and $q \in Q(\BR)$. We have 
\[
a \cdot x = a_x x \quad \text{and} \quad a \cdot (qx) = a_{qx} qx \; ,
\]
where $a_x$ and $a_{qx}$ are elements of $P(\BR)^0$, which by Corollary~\ref{2B}~(a)
satisfy the relation $a_{qx} = qa_xq^{-1}$.
It follows that 
\[
q (a \cdot x) = qa_x x = a_{qx} qx = a \cdot (qx) \; .
\]
\end{Proof}

\begin{Rem}  [cmp.\ {\cite[Prop.~7.6]{BS}}] \label{2M}
Using Lemma~\ref{2Ipre}, one shows that
the restrictions to $P(\BQ)W(\BR)$ of the actions (1) and (2) on $\FX$
extend to 
$\FX^{BS}$. Both preserve the structure of manifold with corners.
On each connected component
$\bar{\FX}^0 = (\FX^{BS})^0$, they induce the same action of 
$\Stab_{P(\BQ)W(\BR)}(\FX^0)$. In parti\-cu\-lar, they coincide
on $P(\BR)^0 \cap P(\BQ)W(\BR)$. The induced actions on the set of faces are the same,
namely $g: e(Q) \mapsto e(Int(g)(Q))$.
\end{Rem}

Note that the extensions of the actions from $\FX$
to the corners and to the whole of $\FX^{BS}$ (Remarks~\ref{2K} and \ref{2M})
are necessarily unique as $\FX$ is dense
in $\FX^{BS}$, and $\FX^{BS}$ is Hausdorff \cite[Thm.~7.8]{BS}. 

\begin{Def} \label{3A}
Let $K$ be an open compact subgroup of $P (\BA_f)$. Define the \emph{Borel--Serre
compactification of the Shimura variety $M^K (P,\FX)$} as the quotient space
\[
M^K (P,\FX) (\BC)^{BS} := P (\BQ) \backslash \bigl( \FX^{BS} \times P (\BA_f) / K \bigr) 
\]
formed with respect to the action~(1) of $P (\BQ)$ extending the one on $\FX$
underlying the Shimura data.
\end{Def}

Thus, the Borel--Serre compactification of $M^K (P,\FX)$ contains the space of
complex points 
\[
M^K  (P,\FX) (\BC) = P (\BQ) \backslash \bigl( \FX \times P (\BA_f) / K \bigr)
\]
of the Shimura variety $M^K  (P,\FX)$
as an open dense subset. If 
\[
M^K  (P,\FX) (\BC) = \coprod \Gamma(g)_+ \backslash \FX^0
\]
is a representation of $M^K  (P,\FX) (\BC)$ as finite disjoint union of its
connected components as in \cite[Sect.~3.2]{P1} (cmp.\ Remark~\ref{1N}~(a)), then
\[
M^K  (P,\FX) (\BC)^{BS} = \coprod \Gamma(g)_+ \backslash \bar{\FX}^0 \; .
\]
Recall that by definition, 
\[
\Gamma(g)_+ = \Stab_{P(\BQ)}(\FX^0) \cap \Gamma(g) \; ,
\]
with
\[
\Gamma(g) = P(\BQ) \cap g K g^{-1} \; .
\]
In particular, the group $\Gamma(g)_+$ is contained in $\Stab_{P(\BR)}(\FX^0)$,
meaning that as far as the formation of the connected components of 
the Borel--Serre compactification is concerned,
we need not worry about the difference of actions (1) and (2). 

\begin{Thm} [Borel--Serre] \label{3B}
Let $K$ be an open compact subgroup of $P (\BA_f)$. Then 
the Borel--Serre compactification of $M^K (P,\FX)$ is compact.
\end{Thm}

\begin{Proof}
We identified $M^K (P,\FX) (\BC)^{BS}$ with a finite disjoint union
\[
\coprod \Gamma(g)_+ \backslash \bar{\FX}^0 \; .
\]
According to \cite[Thm.~9.3]{BS}, each of the components 
$\Gamma(g)_+ \backslash \bar{\FX}^0$ is compact.
\end{Proof}

\forget{
\begin{Prop} \label{3Ba}
Let $K$ be an open compact subgroup of $P (\BA_f)$, 
$R$ a parabolic subgroup of $P$, and $p \in P (\BA_f)$. Let
$U \subset M^K (P,\FX) (\BC)^{BS}$ be a locally closed union of strata
containing $e^K(R,p)$, and contained in the closure $\overline{e^K(R,p)}$
of $e^K(R,p)$ in $M^K (P,\FX) (\BC)^{BS}$.
If $K$ is \emph{neat} in the sense of \cite[Sect.~0.6]{P1}, 
then the open immersion $U \into \overline{e^K(R,p)}$ is contractible.  
\end{Prop}

\begin{Proof}
First, the closure $\overline{e^K(R,p)}$ equals the image of 
\[
\overline{e(R)} \times \{ pK \} \subset \FX^{BS} \times P (\BA_f) / K
\]
under the quotient map 
\[
\FX^{BS} \times P (\BA_f) / K  \longonto M^K  (P,\FX) (\BC)^{BS}
\]
\cite[Prop.~9.4~(ii)]{BS}. Second, the 
restriction of this map
to $\overline{e(R)} \times \{ pK \}$ equals the quotient by the action
of the stabilizer (in $P(\BQ)$) of $e(R) \times \{ pK \}$ \cite[Prop.~9.4~(ii)]{BS}.
But this action is free on the whole of $\FX^{BS} \times P (\BA_f) / K$ 
as the group $K$ is supposed neat \cite[Sect.~9.5]{BS}.
Our claim therefore follows from Proposition~\ref{2Jc}.
\end{Proof}
}

The rest of the present section is concerned with functoriality,
and will not be used in the sequel.
In the general context of homogeneous spaces $(X,(L_x)_{x \in X})$ of type $S-k$,
little seems to have been said in the litera\-ture
about the functorial behaviour of the manifold with corners $\bar{X}$. 
Clearly, a morphism $(H,X,(L_x)_{x \in X}) \to (H',Y,(L'_y)_{y \in Y})$ 
should be a pair consisting
of a morphism of algebraic groups $\varphi: H \to H'$ and a continuous map $\psi: X \to Y$
satisfying (i)~$\psi(hx) = \varphi(h) \psi(x)$ for all $h \in H$ and $x \in X$,
(ii)~for all $x \in X$,
the image under $\varphi$ of the Levi subgroup $L_x$ of $H$ associated to $x$
is contained in the Levi subgroup $L'_{\psi(x)}$ of $H'$ associated to $\psi(x)$. 
But in order to obtain a map $\bar{X} \to \bar{Y}$, the construction from \cite{BS},
which involves in particular very careful choices of Levi subgroups of parabolics
(cmp.\ Proposition~\ref{2A}), suggests that at least another requirement needs
to be imposed,
namely \emph{compatibility of Cartan involutions}: denote by
$K_x$ and $K'_y$ the unique maximal compact subgroups of $\Stab_{H^0(\BR)}(x)$ and
$\Stab_{(H')^0(\BR)}(y)$, respectively \cite[Rem.~2.2]{BS}, then the associated 
Cartan involutions $\theta_{K_x}$ and $\theta_{K'_y}$ satisfy 
(iii)~$\varphi \circ \theta_{K_x} = \theta_{K'_{\psi(x)}} \circ \varphi : 
L^0_x \to (L')^0_{\psi(x)}$ for all $x \in X$. 
In order to relate the morphisms $L_x \to L_{\psi(x)}'$ induced by $\varphi$ 
to the maximal reductive quotients of $H$ and $H'$, respectively, it also appears
reasonable to demand that (iv)~the image under $\varphi$ of the unipotent radical
of $H$ is contained in the unipotent radical of $H'$. \\

We decided not to pursue these general considerations concerning functoriality. 
In the context of morphisms of Shimura data (where (i)--(iv) are automatically
satisfied), our result reads as follows.

\begin{Thm} \label{3C}
The manifold with corners $\FX^{BS}$ is functorial in $(P,\FX)$.
\end{Thm}

For the rest of the section, we fix a morphism $(\varphi,\psi): (P,\FX) \to (P',\FY)$
of Shimura data \cite[Def.~2.3]{P1}, both of which satisfy hypotheses $(+)$ and $(U=0)$. 
In particular, we have a $P(\BR)$-equivariant analytic map $\psi: \FX \to \FY$.
Given that $\FX$ is dense
in $\FX^{BS}$, and $\FY^{BS}$ is Hausdorff \cite[Thm.~7.8]{BS}, Theorem~\ref{3C}
is equivalent to the following: the composition
\[
\FX \stackrel{\psi}{\longto} \FY \longinto \FY^{BS}
\]
can be extended to a morphism of manifolds with corners $\psi^{BS}: \FX^{BS} \to \FY^{BS}$. \\

Our strategy has three steps: (A)~using $\varphi$,
define a map $Q \mapsto Q'$ associating to every parabolic
subgroup of $P$ a parabolic subgroup of $P'$, (B)~show that for every parabolic subgroup
$Q$ of $P$, the morphism $\varphi$ induces $\varphi_Q: A_Q \to A_{Q'} \,$, admitting a
(necessarily unique) continuous extension $\bar{A}_Q \to \bar{A}_{Q'} \,$, 
which in addition respects
the stratifications of source and target, 
(C)~show that for every $Q$, the analytic map $\psi$ is compatible
with the geodesic action of $A_Q$ on the source $\FX$, and of $A_{Q'}$ on the target 
$\FY$. \\

We denote by $W'$ the unipotent radical of $P'$. \\

As for step~(A), here is the most concise description of the map $Q \mapsto Q'$: choose a 
co-character $\lambda$ of $P$, such that $\Lie Q$ is the sum of $\Lie W$ and of all
non-negative weight spaces in $\Lie P$ under $\Ad \circ \lambda$. Then define 
$Q'$ as the unique connected subgroup of $P'$, whose Lie algebra is   
the sum of $\Lie W'$ and all
non-negative weight spaces in $\Lie P'$ under $\Ad \circ \varphi \circ \lambda$.
With an eye towards steps~(B) and (C), we need to be more explicit. \\

First, compatibly with what is suggested by our notation, the parabolic subgroup
of $P'$ associated by our map to $P$ is $P'$ itself. Next, assume that $Q_j$ is
admissible. The construction of $Q'_j$ is then spelled out in \cite[Sect.~4.16]{P1}:
choose any point $x \in \FX$, and 
let $Q'_{j,\BR}$ be the unique connected subgroup of $P'_\BR$, whose Lie algebra is   
the direct sum of all
non-negative weight spaces in $\Lie P'_\BR$ under 
$\Ad \circ \varphi \circ h_{\re_j(x)} \circ w$.
It is proved \cite[Sect.~4.16]{P1} that $Q'_{j,\BR}$ does not depend on the choice of $x$,
and is defined over $\BQ \,$, \emph{i.e.}, it is the base change to $\BR$ of a closed
connected subgroup $Q'_j$ of $P'$. According to \loccit, the subgroup $Q'_j$ is
admissible parabolic; actually, the morphism $(\varphi,\psi)$ induces morphisms
$(P_j,\FX_j) \to (P'_j,\FY_j)$ between the boundary components associated to 
the canonical normal subgroups $P_j$ of $Q_j$ and $P'_j$ of $Q'_j \,$, respectively. 
We therefore have $\varphi \circ h_{\re_j(x)} = h_{\re_j(\psi(x))}$
for all $x \in \FX$.
Finally, for an arbitrary parabolic subgroup $Q$ of $P$, write
$Q = \bigcap_{j \in J} Q_j$, where the $Q_j \,$, $j \in J$, are admissible, and define
$Q'$ as the intersection $\bigcap_{j \in J} Q'_j$.
With this definition, it is obvious that for two parabolics $Q$ and $R$ of $P$, whose
intersection is still parabolic, we have $(Q \cap R)' = Q' \cap R'$. In particular,
if $R \subset Q$, then $R' \subset Q'$.

\begin{Prop} \label{3D}
Let $Q$ be a parabolic subgroup of $P$. 
Write $Q = \bigcap_{j \in J} Q_j$, where the $Q_j$ are admissible parabolics. \\[0.1cm]
(a)~For all $x \in \FX$,
the Lie algebra $\Lie Q'_\BR$ is the direct sum of all 
weight spaces in $\Lie P'_\BR$, which are simultaneously non-negative under 
the commuting co-characters 
\[
\Ad \circ \varphi \circ h_{\re_j(x)} \circ w = \Ad \circ h_{\re_j(\psi(x))} \circ w \; ,
\]
$j \in J$. \\[0.1cm]
(b)~The subgroup $Q'$ of $P'$ is parabolic.  \\[0.1cm]
(c)~The pre-image of $Q'$ under $\varphi$ satisfies
\[
\varphi^{-1}(Q') = Q \ker(\varphi) \; .
\]
(d)~If $\varphi$ is an epimorphism, then $Q' = \varphi(Q)$. 
\end{Prop}

\begin{Proof}
The $h_{\re_j(x)} \circ w$, $j \in J$, 
commute with each other according to Corollary~\ref{2G}~(b).
Hence so do their compositions with $\varphi$.

The description of $\Lie Q'_\BR = \bigcap_{j \in J} \Lie Q'_{j,\BR}$ 
from part~(a) then follows from that of the individual $\Lie Q'_{j,\BR}$, $j \in J$, 
given before.

Part~(d) is implied by (a). 

As for part~(b), observe that the formation of $Q'$ is compatible with composition
of morphisms of Shimura data. By (d), statement (b) is true if $\varphi$ is an epimorphism.
In the general case,
note that each of the $Q'_j$ being parabolic, it contains the unipotent radical
$W'$ of $P'$. Hence so does $Q'$. Therefore, and by (d), 
statements (b) for $(\varphi,\psi)$ and for the composition of $(\varphi,\psi)$
with the canonical projection from $(P',\FY)$ to the quotient $(P',\FY)/W'$
are equivalent.
We may therefore assume that $(P',\FY)$ is pure,
\emph{i.e.}, that $P'$ is reductive. 
The morphism $(\varphi,\psi)$ factors over the quotient $(P,\FX)/ \ker \varphi$
\cite[Prop.~2.9]{P1}.
The Lie algebra of $W$ being equal to
the weight $(-1)$-part of $\Lie P$ \cite[Def.~2.1~(v)]{P1}, the kernel $\ker \varphi$
contains $W$. Thus, the quotient $(P,\FX)/ \ker \varphi$ is pure. 
Part~(b) being known for the epimorphism $(P,\FX) \longonto (P,\FX)/ \ker \varphi$,
we may therefore assume that $(P,\FX)$ is pure, too, and that $\varphi$ is a monomorphism.

It clearly suffices to show that $Q'_\BR$ is parabolic. 
Choose a maximal
$\BR$-split subtorus $T$ of $P'_\BR$ containing the images of the
$\varphi \circ h_{\re_j(x)} \circ w$, $j \in J$. On the $\BR$-vector space
$X_*(T) \otimes_\BZ \BR$, choose a notion of positivity, for which the
$\varphi \circ h_{\re_j(x)} \circ w$, $j \in J$, are all positive (this is possible
since the $h_{\re_j(x)} \circ w$, $j \in J$, are positive with respect to a minimal
parabolic subgroup of $P_\BR$ contained in $Q_\BR$). Then the direct sum of all weight
spaces of $\Lie P'_\BR$, which are simultaneously non-negative under all 
positive co-characters of $X_*(T)$, is the Lie algebra
of a minimal parabolic of $P'_\BR$ contained
in $Q'_\BR$ \cite[Sect.~5.12]{BT}.

Part~(a) implies that
\[
\bigl( \varphi^{-1}(Q') \bigr)^0 = Q \bigl( \ker(\varphi) \bigr)^0 \; :
\]
check on the level of Lie algebras, and use the fact that $Q$ is connected. 
Since $Q'$ is connected, we have
\[
\varphi^{-1}(Q') = \bigl( \varphi^{-1}(Q') \bigr)^0 \ker(\varphi) \; .
\]
This proves part~(c).
\end{Proof}

\begin{Cor} \label{3Da}
Let $Q$ be a parabolic subgroup of $P$, and $x \in \FX$. Then
$\varphi$ maps the Levi subgroup $L_{x,Q}$ of $Q$ to
the Levi subgroup $L_{\psi(x),Q'}$ of $Q'$.
\end{Cor}

\begin{Proof}
This follows from Proposition~\ref{3D}~(a) and Corollary~\ref{2G}.
\end{Proof}

Now for step~(B). Recall (Definition~\ref{2Ga}) that $A_Q$ is the 
neutral connected component of the group of $\BR$-points of the $\BQ$-split torus 
\[
S_Q = Z_{d,\BQ}(\bar{Q}) / Z_{d,\BQ}(G) \; . 
\]
The morphism $(\varphi,\psi): (P,\FX) \to (P',\FY)$ remains fixed.

\begin{Prop} \label{3E}
Let $Q$ be a parabolic subgroup of $P$. \\[0.1cm]
(a)~The image of the unipotent radical $\Rad^u(Q)$ of $Q$ under $\varphi$ is contained 
in the unipotent radical $\Rad^u(Q')$ of $Q'$. 
In particular, $\varphi$ induces a morphism, equally denoted
by $\varphi$, from $\bar{Q}$ to $\bar{Q}'$. \\[0.1cm]
(b)~The image of $Z_{d,\BQ}(\bar{Q})^0$ under $\varphi: \bar{Q} \to \bar{Q}'$
is contained in $Z_{d,\BQ}(\bar{Q}')^0$. In particular, $\varphi$ 
induces a morphism, still denoted by $\varphi$, from $S_Q^0$ to $S_{Q'}^0$.
\end{Prop}

\begin{Proof}
(a): according to \cite[Def.~2.1~(v)]{P1}, the Lie algebra of $W$ is equal to
the weight $(-1)$-part of $\Lie P$, and similarly for $\Lie W'$.
It follows that $\varphi$ maps $W$ to $W'$. Therefore, the statement is true for $Q = P$.

If $Q = Q_j$ is admissible, then $(\varphi,\psi)$ induces morphisms 
$(P_j,\FX_j) \to (P'_j,\FY_j)$ between the boundary components associated to 
the canonical normal subgroups $P_j$ of $Q_j$ and $P'_j$ of $Q'_j$, respectively
\cite[Sect.~4.16]{P1}. Again by \cite[Def.~2.1~(v)]{P1}, we have
$\varphi(\Rad^u(P_j)) \subset \Rad^u(P_j')$. But 
\[
\Rad^u(Q_j) = \Rad^u(P_j) W
\]
\cite[proof of Lemma~4.8]{P1}, and likewise for $\Rad^u(Q'_j)$. 
Therefore, the statement is true if $Q$ is admissible.

In the general case, write $Q = \bigcap_{j \in J} Q_j$, where the $Q_j$, $j \in J$, are admissible, hence $Q' = \bigcap_{j \in J} Q'_j$, use Proposition~\ref{3D}~(b)
and the fact that if an intersection of parabolic subgroups $R_j$ is parabolic, then 
the unipotent radical of that intersection is the product of the unipotent radicals of
the $R_j$. 

\noindent (b): we may assume that both $(P,\FX)=(G,\FX)$ and $(P',\FY)=(G',\FY)$ are pure. 
Choose a point $x \in \FX$.
The morphism $\varphi$ maps $\Stab_{G(\BR)}(x)$ to $\Stab_{G'(\BR)}(\psi(x))$. According
to Corollary~\ref{1Lcor}~(b), 
\[
Z_{d,\BQ}(G)(\BR)^0 \subset \Stab_{G(\BR)}(x) \; ,
\]
and $\Stab_{G'(\BR)}(\psi(x))$ is a direct product of $Z_{d,\BQ}(G')(\BR)^0$
and a compact group. 
Therefore, $\varphi$ maps $Z_{d,\BQ}(G)(\BR)^0$ to $Z_{d,\BQ}(G')(\BR)^0$, hence
$Z_{d,\BQ}(G)^0$ to $Z_{d,\BQ}(G')^0$. Thus, the statement is true for $Q = P$.

In the general case, write $Q = \bigcap_{j \in J} Q_j$, where the $Q_j$, $j \in J$, are 
pairwise different maximal proper parabolics. We have $Q' = \bigcap_{j \in J} Q'_j$.
According to Theorem~\ref{2H}~(b), 
\[
Z_{d,\BQ}(\bar{Q}) = Z_{d,\BQ}(G) \times 
          \prod_{j \in J} \bigl( \pi_Q \circ h_{\re_j(x)} \circ w (\GQm) \bigr) \; .
\]
But for all $j \in J$,
\[
\varphi \circ h_{\re_j(x)} \circ w = h_{\re_j(\psi(x))} \circ w \; .
\]
The parabolic subgroup $Q'_j$ of $P'$ is admissible \cite[Sect.~4.16]{P1};
according to Theorem~\ref{2H}~(a), the image under $\pi_{Q'}$ 
of $h_{\re_j(\psi(x))} \circ w$ is therefore contained
in $Z_{d,\BQ}(\bar{Q}')$. It remains to recall
that $\varphi(Z_{d,\BQ}(G)^0) \subset Z_{d,\BQ}(G')^0$
by the above.
\end{Proof}

Note that for $x \in \FX$, the diagram
\[
\vcenter{\xymatrix@R-10pt{
        L_{x,Q} \ar[r]^-{\pi_Q}_-\cong \ar[d]_\varphi & \bar{Q}_\BR \ar[d]^\varphi \\
        L_{\psi(x),Q'} \ar[r]^-{\pi_{Q'}}_-\cong & \bar{Q}'_\BR 
\\}}
\]
is commutative. 

\begin{Cor} \label{3F}
Let $Q$ be a parabolic subgroup of $P$. \\[0.1cm]
(a)~The morphism $\varphi$ induces a morphism
\[
\varphi: A_Q \longto A_{Q'} \; .
\]
(b)~Let $R$ be a parabolic subgroup of $P$ contained in $Q$. Then the diagram
\[
\vcenter{\xymatrix@R-10pt{
        A_Q \ar@{^{ (}->}[r] \ar[d]_\varphi & A_R \ar[d]^\varphi \\
        A_{Q'} \ar@{^{ (}->}[r] & A_{R'} 
\\}}
\]
is commutative.
\end{Cor}

In order to show that $\varphi: A_Q \to A_{Q'}$ extends to a map
$\bar{A}_Q \to \bar{A}_{Q'} \,$ respecting
the stratifications, we now give an explicit description in terms of the parameterizations
of $A_Q$ and $A_{Q'}$. Write $Q = \bigcap_{j \in J} Q_j$, where the $Q_j$, $j \in J$, 
are pairwise different maximal proper parabolics of $P$.  
Write $Q' = \bigcap_{i \in I} R_i$, where the $R_i$, $i \in I$, 
are pairwise different maximal proper parabolics of $P'$.
Thus, for every $i \in I$, there exists $j \in J$ such that
$Q'_j \subset R_i$. For $i \in I$ and $j \in J$, define $e_{i,j}$ to be equal to $1$
if $Q'_j \subset R_i$,
and to be equal to $0$ if $Q'_j \not \subset R_i$. Thus, for all $i \in I$, we have
\[
\sum_{j \in J} e_{i,j} \ge 1 \; .
\]
Consider the parameterizations 
\[
par_Q =  \prod_{j \in J} \bigl( \pi_Q \circ h_{\re_j(\argdot)} \circ w \bigr) :
\bigl( \BR_+^* \bigr)^J \isoto A_Q
\]
and
\[
par_{Q'} =  \prod_{i \in I} \bigl( \pi_{Q'} \circ h_{\re_i(\argdot)} \circ w \bigr) :
\bigl( \BR_+^* \bigr)^I \isoto A_{Q'} \; .
\] 

\begin{Prop} \label{3G}
Define the morphism
\[
inc_{Q,Q'} : \bigl( \BR_+^* \bigr)^J \longto \bigl( \BR_+^* \bigr)^I
\]
by sending 
\[
(q_j)_{j \in J} \quad \text {to} \quad 
             \bigl( \prod_{j \in J} q_j^{e_{i,j}} \bigr)_{i \in I} \; .
\]
Then the diagram
\[
\vcenter{\xymatrix@R-10pt{
        \bigl( \BR_+^* \bigr)^J \ar[d]_{inc_{Q,Q'}} \ar[r]^-{par_Q}_-\cong & 
                                                                   A_Q \ar[d]^\varphi \\
        \bigl( \BR_+^* \bigr)^I \ar[r]^-{par_{Q'}}_-\cong & A_{Q'} 
\\}}
\]
is commutative.
\end{Prop}

\begin{Proof}
Fix $x \in \FX$, and write $y := \psi(x)$. We have
\[ 
\varphi \circ par_Q = \prod_{j \in J} \bigl( \pi_{Q'} \circ h_{\re_j(y)} \circ w \bigr) :
\bigl( \BR_+^* \bigr)^J \longto A_{Q'} \; ,
\]
while 
\[
par_{Q'} =  \prod_{i \in I} \bigl( \pi_{Q'} \circ h_{\re_i(y)} \circ w \bigr) :
\bigl( \BR_+^* \bigr)^I \isoto A_{Q'} \; .
\]
But by Corollary~\ref{2Ha}~(b), for each $j \in J$, we have
\[
\pi_{Q'} \circ h_{\re_j(y)} \circ w = 
\prod_{\{i \in I, e_{i,j}=1\}} \bigl( \pi_{Q'} \circ h_{\re_i(y)} \circ w \bigr) 
                                                                        \circ \Delta \; ,
\]
where $\Delta$ denotes the diagonal inclusion of $\BR_+^*$ into 
$(\BR_+^*)^{\{i \in I, e_{i,j}=1\}}$.
\end{Proof}

\begin{Cor} \label{3H}
Let $Q$ be a parabolic subgroup of $P$. \\[0.1cm]
(a)~The morphism $\varphi: A_Q \to A_{Q'}$ extends (uniquely) to a continuous map
\[
\varphi: \bar{A}_Q \longto \bar{A}_{Q'} \; .
\]
This extension is $A_Q$-equivariant. \\[0.1cm]
(b)~The extension $\varphi: \bar{A}_Q \to \bar{A}_{Q'}$ 
respects the stratifications of $\bar{A}_Q$ and $\bar{A}_{Q'}$. \\[0.1cm]
(c)~We have $\varphi(\infty_Q) = \infty_{Q'}$.
\end{Cor}

\begin{Proof}
(a): for the first part of the statement,
it suffices to observe that the formula
\[
(q_j)_{j \in J} \longmapsto \bigl( \prod_{j \in J} q_j^{e_{i,j}} \bigr)_{i \in I}
\]
for $inc_{Q,Q'}$ from Proposition~\ref{3G} still makes sense if the coordinates
$q_j$ are allowed to take the value $+ \infty$ (with the usual conventions
for multiplication by $+ \infty$). The second part follows from unicity of
this extension (or from the explicit formula for $inc_{Q,Q'}$). 

\noindent (b): since the extension is $A_Q$-equivariant, it maps orbits under $A_Q$
to orbits under $A_{Q'}$. But these orbits constitute the 
respective stratifications.

\noindent (c): according to Proposition~\ref{3G}, 
\[
\varphi(\infty_Q) = par_{Q'} \bigl( \prod_{j \in J} (+ \infty)^{e_{i,j}} \bigr)_{i \in I}  = par_{Q'} \bigl( (+ \infty)^{\sum_{j \in J} e_{i,j}} \bigr)_{i \in I} \; .
\]
As we observed already, the natural number $\sum_{j \in J} e_{i,j}$
is stricly positive for every $i \in I$. 
\end{Proof}

Given the results already obtained, not much remains to be done for step~(C).
We keep the morphism $(\varphi,\psi): (P,\FX) \to (P',\FY)$. 

\begin{Prop} \label{3I}
The analytic map $\psi: \FX \to \FY$ is compatible with the geodesic actions. 
More precisely, let $Q$ be a parabolic subgroup of $P$, $x \in \FX$, $a \in A_Q$,
and $\varphi(a)$ the image of $a$ under $\varphi: A_Q \to A_{Q'}$.
Then
\[
\psi(a \cdot x) = \varphi(a) \cdot \psi(x) \; .
\] 
\end{Prop}

\begin{Proof}
The element $a$ is of the form $\pi_Q(a_x)$, for some $a_x \in Q(\BR)$ belonging
to the neutral connected component of the center of $L_{x,Q}$. According to
the definition of the geodesic action \cite[Sect.~3.2]{BS},
\[
a \cdot x = a_x x \; .
\]
By Corollary~\ref{3Da} and Proposition~\ref{3E}~(b), the element $\varphi(a_x) \in Q'(\BR)$
belongs
to the center of $L_{\psi(x),Q'}$ and lifts $\varphi(a)$. Therefore,
\[
\varphi(a) \cdot \psi(x) = \varphi(a_x) \psi(x) \; .
\]
The claim thus follows from the relation
\[
\psi(a_x x) = \varphi(a_x) \psi(x) \; ,
\]
which is satisfied as $(\varphi,\psi): (P,\FX) \to (P',\FY)$
is a morphism of Shimura data.
\end{Proof}

Theorem~\ref{3C} is implied by part~(a) of the following result.

\begin{Prop} \label{3K}
Let $(\varphi,\psi): (P,\FX) \to (P',\FY)$ be a morphism of Shimura data, 
both of which satisfy hypotheses $(+)$ and $(U=0)$. \\[0.1cm]
(a)~The analytic map $\psi: \FX \to \FY$ extends uniquely to a continuous map
$\psi^{BS}: \FX^{BS} \to \FY^{BS}$. This map is a morphism of manifolds 
with corners. \\[0.1cm] 
(b)~Let $Q$ be a parabolic subgroup of $P$. Then $\psi^{BS}$ maps the corner
$\FX(Q)$ to the corner $\FY(Q')$, and the face $e(Q)$ of $\FX^{BS}$  to the face $e(Q')$
of $\FY^{BS}$. \\[0.1cm]
(c)~Let $Q$ be a parabolic subgroup of $P$. Then the diagram
\[  
\vcenter{\xymatrix@R-10pt{
        \FX \ar@{->>}[r] \ar[d]_\psi & 
        A_Q \backslash \FX = e(Q) \ar@{^{ (}->}[r] \ar[d]^{\psi^{BS}_{\tei e(Q)}} &
        \FX^{BS} \ar[d]^{\psi^{BS}}                            \\
        \FY \ar@{->>}[r] & 
        A_{Q'} \backslash \FY = e(Q') \ar@{^{ (}->}[r] &
        \FY^{BS}   
\\}}
\]
is commutative.
\end{Prop} 

\begin{Proof}
(a): as mentioned already, unicity of $\psi^{BS}$ follows from the fact that $\FX$ is dense
in $\FX^{BS}$, and $\FY^{BS}$ is Hausdorff \cite[Thm.~7.8]{BS}. 
The same principle shows that continuous extensions of $\psi$ to open subsets of 
$\FX^{BS}$ containing $\FX$ necessarily coincide on the intersection, and hence glue
to give an extension of the union. It is therefore sufficient to show that 
for every parabolic $Q$ of $P$, the composition
\[
\FX \stackrel{\psi}{\longto} \FY \longinto \FY^{BS}
\]
can be continuously extended to the corner $\FX(Q)$, and that this 
extension is a morphism of manifolds with corners. Consider the continuous maps
\[
\varphi: \bar{A}_Q \longto \bar{A}_{Q'} \quad \text{and} \quad 
\psi: \FX \longto \FY \; ;
\]
according to Corollary~\ref{3H}~(a) and Proposition~\ref{3I}, both are $A_Q$-equivariant.
Therefore, they induce a continuous map
\[
\FX(Q) = \bar{A}_Q \, {}^{A_Q} \! \times \FX \longto 
\bar{A}_{Q'} \, {}^{A_{Q'}} \! \times \FY = \FY(Q') \; ,
\]
which is the continuous extension of $\psi$ we are looking for.
It is a morphism of manifolds with corners thanks to Corollary~\ref{3H}~(b).

\noindent (b): the inclusion $\psi^{BS} (\FX(Q)) \subset \FY(Q')$ follows from
our construction. As for $\psi^{BS} (e(Q)) \subset e(Q')$, use Corollary~\ref{3H}~(c),
and the identification of $e(Q)$ and $e(Q')$ with the subspaces of 
equivalence classes $[(\infty_Q,x)]$ and $[(\infty_{Q'},y)]$,
respectively.
 
\noindent (c): the identification 
\[
A_Q \backslash \FX = e(Q) \longleftrightarrow \{ [(\infty_Q,x)] \, , \, x \in \FX \}
\]
maps $A_Q x$ to $[(\infty_Q,x)]$, and likewise for $e(Q')$. But by definition
the extension $\psi^{BS}$ maps $[(\infty_Q,x)]$ to $[(\infty_{Q'},\psi(x))]$.
\end{Proof}

\begin{Rem}
Unicity of the extension $\psi^{BS}: \FX^{BS} \to \FY^{BS}$ implies that
it is equivariant for all continuous group actions, for which
$\psi: \FX \to \FY$ is equivariant. This is the case in particular for
the group $P(\BQ)W(\BR)$, and its action~(1) on $\FX$ and (\emph{via} 
$\varphi: P(\BQ)W(\BR) \to P'(\BQ)W'(\BR))$ on $\FY$.
\end{Rem} 

\begin{Cor} \label{3L}
Let $(\varphi,\psi): (P,\FX) \to (P',\FY)$ be a morphism of Shimura data, 
both of which satisfy hypotheses $(+)$ and $(U=0)$. 
Let $K \subset P(\BA_f)$ and $K' \subset P'(\BA_f)$ be open compact subgroups,
such that $\varphi(K) \subset K'$. 
Then the morphism
\[
[(\varphi,\psi)]_{K,K'} : M^K  (P,\FX) (\BC) \longto M^{K'}  (P',\FY) (\BC)
\]
of complex analytic spaces \cite[Sect.~3.4~(b)]{P1} extends uniquely to a 
continuous map
\[
M^K (P,\FX) (\BC)^{BS} \longto M^{K'} (P',\FY) (\BC)^{BS} 
\]
between the Borel--Serre compactifications.
It is a morphism of manifolds with corners if both $K$ and $K'$ are
neat.
\end{Cor}

\begin{Proof}
Uniqueness of the continuous extension results from Theorem~\ref{3B}.
It is given by the map induced by 
\[
\psi^{BS} \times \varphi (\BA_f) : 
\FX^{BS} \times P (\BA_f) \longto \FY^{BS} \times P' (\BA_f)
\]
on the quotients $M^K (P,\FX) (\BC)^{BS}$ and $M^{K'} (P',\FY) (\BC)^{BS}$, respectively.

If both $K$ and $K'$ are neat, then the quotient maps 
\[
\FX^{BS} \times P (\BA_f) / K  \longonto M^K  (P,\FX) (\BC)^{BS}
\]
and
\[
\FY^{BS} \times P' (\BA_f) / K' \longonto M^{K'} (P',\FY) (\BC)^{BS}
\]
are local homeomorphisms, and $M^K (P,\FX) (\BC)^{BS}$ and 
$M^{K'} (P',\FY) (\BC)^{BS}$ inherit the structure of manifold with corners
\cite[Sect.~9.5]{BS}.
The statement therefore follows from Proposition~\ref{3K}~(a).
\end{Proof}

\begin{Rem} 
Except for Definition~\ref{3A}, Theorem~\ref{3B} and Corollary~\ref{3L}, 
all constructions and results of this section admit variants for $\re (\FX)$
instead of $\FX$ when hypothesis $(U=0)$ is not satisfied. The manifold with corners
is $\re (\FX)^{BS}$, and for every open compact subgroup $K$ of $P(\BA_f)$, one
can form the quotient space
\[
M^K \bigl( P,\re(\FX) \bigr) (\BC)^{BS} := 
             P (\BQ) \backslash \bigl( \re (\FX)^{BS} \times P (\BA_f) / K \bigr) \; . 
\]
It is compact, and contains the ``real locus''
\[
M^K \bigl( P,\re(\FX) \bigr) (\BC) := 
             P (\BQ) \backslash \bigl( \re (\FX) \times P (\BA_f) / K \bigr) 
\]
of the space of complex points of the Shimura variety $M^K ( P,\FX)$
as an open dense subset. But it is not a compactification of $M^K ( P,\FX)(\BC)$. 
\end{Rem}


\bigskip
%
%

\section{The continuous map $p$ from the Borel--Serre to the Baily--Borel compactification}
\label{4}



The present section aims at an alternative proof of a result of Zucker's \cite{Z}:
the space $\FX^*$ domi\-nating the \emph{Baily--Borel compactification} of our Shimura varieties
is in a canonical way a quotient of $\FX^{BS}$. Actually,
\cite{Z} contains the analogous statement for all Satake compactifications
instead of only $\FX^*$. Our construction of the quotient map $p: \FX^{BS} \to \FX^*$
is explicit in terms of the Shimura data. 
Given the combinatorics of the strata of the
boundaries of $\FX^{BS}$ and of $\FX^*$, it
is necessary first to
define a rule $Q \mapsto \adm(Q)$ associating to each parabolic subgroup of $P$ 
an admissible one (Theorem~\ref{4D}).
Using this rule, we are then in a position to give the definition of $p$
(Construction~\ref{4K}). 
Our analysis of the geodesic action
with respect to the data involved in 
$Q \mapsto \adm(Q)$ then allows (Theorem~\ref{4L}) to prove that $p$ is
continuous and is, indeed, the unique continuous extension of $\id_\FX$. 
In particular, it coincides with the quotient map from \cite{Z}. \\

We fix mixed Shimura data $(P,\FX)$. Until Corollary~\ref{4I} , no further hypotheses  
on $(P,\FX)$ will be required.
Recall that to any admissible parabolic subgroup $Q_j$ of $P$ is associated
a canonical normal subgroup $P_j \subset Q_j$ \cite[4.7]{P1}, underlying 
finitely many boundary components of $(P,\FX)$.

\begin{Def} \label{4A}
Let $Q_1$ and $Q_2$ be two admissible parabolic subgroups
of $P$, with associated canonical normal subgroups $P_1 \subset Q_1$ and $P_2 \subset Q_2$.
We define the relation
\[
Q_2 \preceq Q_1
\]
to hold if boundary components $(P_1,\FX_1)$ and $(P_2,\FX_2)$ of $(P,\FX)$
can be chosen such that $(P_2,\FX_2)$ is a boundary component of $(P_1,\FX_1)$. 
\end{Def}

As for transitivity of the notion of boundary component, we refer to \cite[Sect.~4.17]{P1}.
It follows from \cite[Lemma~4.19~(b)]{P1} that $Q_2 \preceq Q_1$ only if $Q_1 \cap Q_2$
remains parabolic. Note also that if $Q_2 \preceq Q_1$, then \emph{any} boundary component
$(P_1,\FX_1)$, admits a boundary component of the form $(P_2,\FX_2)$, and
\emph{any} boundary component
$(P_2,\FX_2)$ is a boundary component of some $(P_1,\FX_1)$ (use the fact that both
\[ 
\FX \longto \coprod \FX_1 \quad \text{and} \quad \FX \longto \coprod \FX_2 
\]
are equivariant under $(Q_1 \cap Q_2)(\BR)$ \cite[4.11]{P1}).

\begin{Prop} \label{4B}
Let $Q_1$ and $Q_2$ be two admissible parabolic subgroups
of $P$, with associated canonical normal subgroups $P_1 \subset Q_1$ and $P_2 \subset Q_2$.
Assume that $Q_2 \preceq Q_1$.
Then the following are equivalent.
\begin{enumerate}
\item[(i)] $Q_1 = Q_2$,
\item[(ii)] $P_1 = P_2$,
\item[(iii)] the weight $(-2)$-part $U_2$ of $P_2$ is contained in the unipotent radical
$W_1$ of $P_1$.
\end{enumerate}
\end{Prop}

\begin{Proof}
Clearly (i) implies (ii), and (ii) implies (iii). 

$(P_2,\FX_2)$ being a boundary component of $(P_1,\FX_1)$, 
for a fixed choice of $x \in \FX_1$, there is 
a morphism
\[
\omega = \omega_{x}: H_{0,\BC} \longto P_{1,\BC}
\]
associated to $x$, such that 
\[
h_x = \omega \circ h_0 \quad \text{and} \quad h_{x_2} = \omega \circ h_\infty 
\]
\cite[Prop.~4.6~(b)]{P1}. Here, $x_2$ is the image of $x$
under the map
\[
\FX_1 \longto \coprod \FX_2 \; ,
\]
and $h_0$ and $h_\infty$ are the morphisms $\BS_\BC \to H_{0,\BC}$ 
from \cite[Sect.~4.3]{P1} (cmp.~the proof of Proposition~\ref{2E}); for our purposes,
it will be sufficient to recall that $h_0 \circ w$ 
and $h_\infty \circ w$ are of the shape
\[
h_0 \circ w : z \longmapsto \Bigl( (z,z) , 
      \begin{pmatrix} z & 0\\ 
                      0 & z
      \end{pmatrix}   \Bigr) 
\]
and
\[
h_\infty \circ w : z \longmapsto \Bigl( (z,z) , 
      \begin{pmatrix} z^2  & i (1-z^2)\\ 
                      0 & 1
      \end{pmatrix}   \Bigr) \; .
\]
Setting 
\[
u_\infty :=  \Bigl( (1,1), \begin{pmatrix} 1  & i \\ 
                                      0 & 1
                      \end{pmatrix}  \Bigr) \; ,
\]
we see that we have the relations
\[
\Int( h_0 \circ w (z))(u_\infty) = u_\infty \quad \text{and} \quad
\Int( h_\infty \circ w (z))(u_\infty) = u_\infty^{z^2} 
\]
for every $z \in \BR^*$. Thus, the element
\[
u := \omega(u_\infty)
\]
of $P_1(\BC)$
is of weight $0$ under $h_x \circ w$ and of weight $-2$ under $h_{x_2} \circ w$.
The latter property means that $u$ belongs to $U_2(\BC)$ \cite[Def.~2.1~(v)]{P1}.

Assuming hypothesis (iii), we have $u \in W_1(\BC)$. But then \cite[Def.~2.1~(v)]{P1},
it must be of weight $-1$ under $h_x \circ w$. Therefore, $u$ belongs to
the kernel of $\omega$. This kernel being normal in $H_{0,\BC} \,$, it contains
$H_{0,\BC}^{der} \,$. In other words, the morphism $\omega$ factors
over the quotient $\BS_\BC \times \BG_{m,\BC}$ of 
$H_{0,\BC} \subset \BS_\BC \times \GL_{2,\BC} \,$.

Therefore, the co-characters $h_x \circ w$ and $h_{x_2} \circ w$ are in fact equal.
By \cite[Prop.~4.6~(b)~(iii)]{P1}, the Lie algebra $\Lie Q_1$ is the sum of all 
non-negative weight spaces in $\Lie P$ under $\Ad \circ h_x \circ w$,
while $\Lie Q_2$ is the sum of all 
non-negative weight spaces in $\Lie P$ under $\Ad \circ h_{x_2} \circ w$.
\end{Proof}

The implication ``(ii) $\Rightarrow$ (i)'' from Proposition~\ref{4B}
shows in particular that the relations $Q_1 \preceq Q_2$ and $Q_2 \preceq Q_1$
imply equality $Q_1 = Q_2$.

\begin{Prop} \label{4C}
Let $Q_1$ and $Q_2$ be two admissible parabolic subgroups
of $P$, with associated canonical normal subgroups $P_1 \subset Q_1$ and $P_2 \subset Q_2$.
Assume that $Q_1 \cap Q_2$ remains parabolic, and that $P_1$ and $P_2$
normalize each other.
Then $Q_1 = Q_2$.
\end{Prop}

\begin{Proof}
Consider the quotient $(G^{ad},\FH^{ad}) := ((P,\FX)/W)/Z(G)$ of $(P,\FX)$.
Any parabolic subgroup $Q$ of $P$ is the pre-image of a unique parabolic
subgroup of $G^{ad}$, namely of $Q/WZ(G)$. The latter is admissible if and only if $Q$ is.
We may therefore show the claim on the level of $(G^{ad},\FH^{ad})$, \emph{i.e.}, we
may assume that $P$ is semi-simple, with trivial center.

Next, observe that the definition of the canonical normal subgroups $P_j$ 
only depends on the Shimura data $(P,h(\FX))$ \cite[Sect.4.7]{P1}. Thus, we may assume
that 
\[
h : \FX \longto \Hom ( \BS, P_\BR ) \; , \; x \longmapsto h_x 
\] 
is injective. As a consequence, under the decomposition $P = \prod_n P_n$ of $P$
into its simple factors, the space $\FX$ decomposes as well, in other words, we have
$(P,\FX) = \prod_n (P_n,\FX_n)$ in the sense of \cite[Def.~2.5]{P1}. 
The relation $\preceq$ being compatible with products, we may assume that $P$ is simple.

Fix a minimal parabolic subgroup $Q$ of $P$ contained in $Q_1 \cap Q_2$. 
Since $P$ is simple, the restriction of the relation $\preceq$ to the (finitely many)
admissible parabolic subgroups containing $Q$ is a total order
\cite[Rem.~(ii) on p.~91, after Cor.~4.20]{P1}. In particular,
we have $Q_1 \preceq Q_2$ or $Q_2 \preceq Q_1$. 

Assume that $Q_2 \preceq Q_1$, which implies $P_2 \subset P_1$. 
By assumption, $P_1$ normalizes $P_2$. Therefore, it
normalizes the weight $(-2)$-part $U_2$ of $P_2$. It being unipotent, this implies that
$U_2$ is contained in the unipotent radical $W_1$ of $P_1$. 
By Proposition~\ref{4B}, we have $Q_1 = Q_2$.

One argues symmetrically when $Q_1 \preceq Q_2$.   
\end{Proof}

\begin{Rem}
\emph{Via} \cite[Rem.~(ii) on p.~91, after Cor.~4.20]{P1}, 
the proof of Proposition~\ref{4C}
makes use of the explicit description 
of roots systems possibly occurring in the 
simple factors of $G^{ad}$. We tried in vain to find a proof not using this description. 
Observe that the same type of information is used in \cite{Z},
in order to establish the relation between certain Satake compactifications
and the Baily--Borel compactification \cite[Sect.~(3.11)]{Z}.
\end{Rem}

Here is the first main result of this section.

\begin{Thm} \label{4D}
Let $Q$ be a parabolic subgroup of $P$. \\[0.1cm]
(a)~There is a unique
admissible parabolic subgroup $Q_\infty$
of $P$, with associated canonical normal subgroup $P_\infty \subset Q_\infty$, such that
\[
P_\infty \subset Q \subset Q_\infty \; .
\]
(b)~The admissible parabolic $Q_\infty$ from (a) 
is the unique minimal element, with respect to the relation
$\preceq \,$, in the (finite) set of admissible parabolics containing $Q$.
\end{Thm}

\begin{Proof}
Let us first establish unicity in (a). If 
\[
P_j \subset Q \subset Q_j
\]
for $j = 1,2$, then $Q_1 \cap Q_2$ contains $Q$, and is therefore parabolic.
$P_1$ is contained in $Q_2$, and hence normalizes $P_2$. Similarly, $P_2$
normalizes $P_1$. According to Proposition~\ref{4C}, we have indeed $Q_1 = Q_2$.

Next, assume that existence has been proved 
for all proper boundary components of $(P,\FX)$. 
Let $Q_1$ be an admissible parabolic subgroup of $P$,
with associated canonical normal subgroup $P_1 \subset Q_1$,
and containing $Q$.
Then $Q \cap P_1$ is a parabolic subgroup of $P_1$
(as $P_1$ is normal in $Q_1$). If $Q_1 \ne P$, meaning that $P_1$
underlies proper boundary components of $(P,\FX)$, then by our induction
hypothesis, there exists an admissible parabolic subgroup $Q_{12}$
of $P_1$, with associated canonical normal subgroup $P_2 \subset Q_{12}$, such that
\[
P_2 \subset Q \cap P_1 \subset Q_{12} \; .
\]
According to \cite[Lemma~4.19]{P1}, there is a unique admissible parabolic subgroup $Q_2$
of $P$, whose associated canonical normal subgroup equals $P_2$, and such that
\[
Q_{12} = Q_2 \cap P_1 \; .
\]
Now $Q$, being contained in $Q_1$, normalizes $P_1$, hence $Q \cap P_1$,
hence (by unicity in (a)) $P_2$, hence (by unicity in \cite[Lemma~4.19]{P1}) $Q_2$.
The latter being parabolic, we have $Q \subset Q_2$. Thus,
\[
P_2 \subset Q \cap P_1 \subset Q \subset Q_2 \; .
\]
This construction establishes existence in (a) as long as $Q_1$ can be chosen to be proper. If
the only admissible parabolic is equal to $P$, then $Q = P$, and $\adm(P) = P$.
Thus, existence in (a) is proved in any case.

The above construction also shows (b), again as long as the admissible parabolic $P$
is excluded. But $P$ is the unique maximal element, with respect to $\preceq \,$,
in the set of all admissible parabolics.
\end{Proof}

\begin{Def} \label{4E}
Let $Q$ be a parabolic subgroup of $P$. Define $\adm(Q)$ as
the unique admissible parabolic from Theorem~\ref{4D}~(a), and $\adm_{Sh}(Q)$
as the normal subgroup canonically associated to $\adm(Q)$.
\end{Def}

\begin{Rem}
The map $Q \mapsto \adm(Q)$ can be described using 
Dynkin diagrams. We refer to \cite[Sect.~5.7, Figure~7]{G}
in the case when $G^{ad}$ is simple (noting that $P$ from \loccit \ is our $Q$,
and $Q$ in \loccit \ is our $\adm(Q)$). The principle ``$L_P$ and $L_Q$
have the same Hermitian factor'' from \loccit \ translates into the double inclusion 
\[
\adm_{Sh}(Q) \subset Q \subset \adm(Q) \; .
\]
\end{Rem}

\begin{Prop} \label{4F}
Let $Q_j$ be an admissible parabolic subgroup of $P$,
with associated canonical normal subgroup $P_j \subset Q_j$. 
Then the canonical epimorphism $pr_j: Q_j \onto Q_j/P_j$ induces an inclusion-preserving
bijection $R \mapsto pr_j^{-1}(R)$ between the set of parabolics of $Q_j/P_j$
and the set $\adm^{-1}(Q_j)$ of parabolics $Q$ of $P$ satisfying $\adm(Q) = Q_j \,$.
\end{Prop}

\begin{Proof}
This is obvious from the definition.
\end{Proof}

It will be important to have an alternative description of the fibres of $\adm$
at our disposal, using subgroups rather than quotients of the admissible parabolics.

\begin{Def} \label{4G}
Let $Q_j$ be an admissible parabolic subgroup of $P$,
with associated canonical normal subgroup $P_j \subset Q_j$. 
Define a closed connected normal subgroup $C_j$ of $Q_j$ as
\[ 
C_j := \{ q \in Q_j \; , \; \pi_{Q_j}(q) \in \Cent_{\bar{Q}_j}(\pi_{Q_j}(P_j)) \}^0 \; .
\]
\end{Def}

Here as before, the morphism $\pi_Q : Q \onto \bar{Q}$ is the canonical epimorphism
of a parabolic $Q$ to its maximal reductive quotient.

\begin{Cor} \label{4H}
Let $Q_j$ be an admissible parabolic subgroup of $P$,
with associated canonical normal subgroup $P_j \subset Q_j$. 
Then the monomorphism $C_j \into Q_j$ induces an inclusion-preserving
bijection $R \mapsto C_j \cap R$ between the set $\adm^{-1}(Q_j)$ and the set 
of parabolics of $C_j$. Its inverse is given by $S \mapsto SP_j$. 
\end{Cor}

\begin{Proof}
Write $\bar{P}_j := \pi_{Q_j}(P_j)$, and
choose a connected almost direct complement $H$ of $\bar{P}_j$ in $\bar{Q}_j$
(which is possible since both are reductive, and $\bar{P}_j$ is normal in $\bar{Q}_j$).
Then $\Cent_{\bar{Q}_j}(\bar{P}_j)^0 = Z( \bar{P}_j )^0 H$. The composition
\[
H \longinto \bar{Q}_j \longonto \bar{Q}_j/ \bar{P}_j
\]
is surjective, with finite kernel. Thus, the composition
\[
C_j \longinto Q_j \longonto Q_j/P_j \longonto \bar{Q}_j/ \bar{P}_j
\]
is surjective, with solvable kernel. It therefore induces a bijection of the sets
of parabolics. The same is true for the canonical epimorphism from $Q_j/P_j$
to $\bar{Q}_j/ \bar{P}_j$ (as its kernel is unipotent). Altogether, the composition
\[
C_j \longinto Q_j \longonto Q_j/P_j 
\]
induces a bijection of the sets
of parabolics of $Q_j/ P_j$ and of $C_j$. 
Now apply Proposition~\ref{4F}.
\end{Proof}

Theorem~\ref{4D} allows to analyze the geodesic action more closely.  

\begin{Cor} \label{4I}
Assume hypotheses $(+)$ and $(U=0)$.
Let $Q$ be a parabolic subgroup of $P$. 
For $x \in \FX$, denote by $\re_\infty(x)$ the real part of $x$ with respect to
$\adm(Q)$ (Definition~\ref{2D}).
Then 
\[
\re_\infty(a \cdot x) = \re_\infty(x) \; , \; \forall \; x \in \FX \; , \; a \in A_Q \; .
\]
\end{Cor}

\begin{Proof}
By definition, 
\[
a  \cdot x = a_x x \; ,
\]
for some element $a_x$ of $Q(\BR) \subset \adm(Q)(\BR)$. The group
$\adm(Q) =: Q_\infty$ normalizes the weight $(-2)$-part $U_\infty$ of 
$P_\infty := \adm_{Sh}(Q)$; therefore,
\[
\re_\infty(a  \cdot x) = \re_\infty(a_x x) = a_x \re_\infty(x) \; .
\]
Our claim is thus equivalent to $a_x$ belonging to the stabilizer of the point
$\re_\infty(x) \in \FX_\infty \,$.

Recall \cite[Sect.~3.2]{BS}
that $a_x$ belongs more precisely to
the neutral connected component of
\[
Z(L_{x,Q}) \subset L_{x,Q} \subset Q \; .
\]
Here, $L_{x,Q}$ is the Levi subgroup
\[
L_{x,Q} = \bigcap_j \Cent_{P_\BR} \bigl( h_{\re_j(x)} \circ w \bigr) \; ,
\]
where the intersection is over all admissible parabolics $Q_j$ containing $Q$, including $P$
(Corollary~\ref{2G}~(a)). According to Theorem~\ref{4D}~(b),
$Q_\infty \preceq Q_j$ for all $j$. By Proposition~\ref{2E}~(a) (applied to the admissible
parabolic $Q_\infty \cap P_j$ of $P_j$; cmp.\ \cite[Sect.~4.17]{P1}),
the image of $h_{\re_\infty(x)}$ centralizes all co-characters
$h_{\re_j(x)} \circ w$. It is thus contained
in $L_{x,Q}$. Thus, it commutes with any element in $Z(L_{x,Q})$.
In particular,
\[
h_{a_x \re_\infty(x)} = a_x h_{\re_\infty(x)} a_x^{-1} = h_{\re_\infty(x)} \; .
\]
But $a_x \in Q(\BR)^0 \subset Q_\infty(\BR)^0$, hence $a_x \re_\infty(x)$
belongs to the same connected component of $\FX_\infty$ as $\re_\infty(x)$.
\cite[Cor.~2.12]{P1}, applied to the Shimura data
$(P_\infty,\FX_\infty)$, thus allows to conclude that
\[
 a_x \re_\infty(x) = \re_\infty(x) \; ,
\]
as desired.
\end{Proof}

Corollary~\ref{4I} says that ``the geodesic action of $A_Q$ 
is visible only in the imaginary part, formed with respect to $\adm(Q)$''. 
It remains to obtain sufficient control on the effect 
of the action on the imaginary part. For an admissible 
parabolic subgroup $Q_j$ of $P$, consider the complex analytic map
\[
\FX \longto \coprod \FX_j
\]
of $\FX$ to the rational boundary components associated to $Q_j \,$. 
For any connected component $\FX^0$ of $\FX$,
denote by $C(\FX^0,P_j)$ the image of  the map $\im_j$ from Definition~\ref{2D}.
According to \cite[Prop.~4.15~(b)]{P1}, using additive
notation, $C(\FX^0,P_j)$ is an open convex cone 
in the real vector space $U_j(\BR)(-1)$ of purely imaginary elements of $U_j(\BC)$.
\emph{Via} the parameterization 
\[
par_Q =  \prod_{j \in J} \bigl( \pi_Q \circ h_{\re_j(\argdot)} \circ w \bigr) :
\bigl( \BR_+^* \bigr)^J \isoto A_Q
\]
from Definition~\ref{Ppar}, we identify $A_Q$ and $( \BR_+^* )^J$.

\begin{Cor} \label{4J}
Assume hypotheses $(+)$ and $(U=0)$.
Let $Q$ be a parabolic subgroup of $P$. 
For $x \in \FX$, denote by $\im_\infty(x)$ the imaginary part of $x$ with respect to
$\adm(Q)$. Then for all $a \in [1,+\infty)^J \subset A_Q$, 
there is an element $u(a)$ in the closure 
$\overline{C(\FX^0,\adm_{Sh}(Q))}$ of $C(\FX^0,\adm_{Sh}(Q))$, such that
\[
\im_\infty(a \cdot x) = \im_\infty(x) + u(a) \; , \; \forall \; x \in \FX^0 \; . 
\]
\end{Cor}

\begin{Proof}
Write $Q_\infty := \adm(Q)$, hence $P_\infty = \adm_{Sh}(Q)$.
The geodesic action respects $\FX^0$. Therefore, if the claim is true for two
elements $a$ and $b$ of the monoid $[1,+\infty)^J$, it is true for their product
(as $\overline{C(\FX^0,P_\infty)}$ is closed under sums).
We may therefore suppose that $a$ is of the form
\[
a = (1,1,\ldots,1,q,1,\ldots,1) \quad\quad \text{ ($q$ in position $j$)} \; , \;
\]
for some $j \in J$, and $q \ge 1$.
The index $j$ corresponds to one of the maximal proper parabolics $Q_j$ containing $Q$. Given that we used $par_Q$ to identify $A_Q$ and $( \BR_+^* )^J$,
we have
\[
a \cdot x = a_x x \; ,
\]
for $a_x := h_{\re_j(x)} \circ w(q)$ (Theorem~\ref{2H}~(c)).
The group $U_j$ equals the weight $(-2)$-part of $P_j$. Therefore, 
\[ 
a_x \im_j(x) = q^2 \im_j(x) \; .
\]
According to Theorem~\ref{4D}~(b), $Q_\infty \preceq Q_j$. This means that
the boundary component $(P_j,\FX_j)$ containing the image $x_j$ of $x$ under
$\FX \to \coprod \FX_j$ admits a boundary component of the shape $(P_\infty,\FX_\infty)$,
containing the image $x_\infty$ of $x_j$ under $\FX_j \to \coprod \FX_\infty$.
But $x_\infty$ is also the image of $x$ under $\FX \to \coprod \FX_\infty$,
given the transitivity of formation of boundary components \cite[Sect.~4.17]{P1}. 
We thus have
\[
\im_\infty(x) = \im_\infty(x_j) = \im_\infty \bigl( \im_j(x)\re_j(x) \bigr)
= \im_\infty \bigl( \re_j(x) \bigr) + \im_j(x)
\]
(note that as $(P_\infty,\FX_\infty)$ is a boundary component of $(P_j,\FX_j)$,
the group $U_j$ is contained in $U_\infty$ \cite[Lemma~4.4~(a)]{P1}).
According to Proposition~\ref{2E}~(d), applied to the admissible parabolic
$Q_\infty \cap P_j$ of $P_j$ and $\re_j(x) \in \FX_j$, 
the element $a_x = h_{\re_j(x)} \circ w(q)$
fixes $\im_\infty ( \re_j(x) )$. 

Altogether, 
\[
\im_\infty(a \cdot x) = \im_\infty(a_x x) = a_x \im_\infty(x)
= a_x \bigl( \im_\infty \bigl( \re_j(x) \bigr) + \im_j(x) \bigr)
\]
equals
\[
\im_\infty \bigl( \re_j(x) \bigr) + q^2 \im_j(x) =
\im_\infty(x) + (q^2-1) \im_j(x) \; ,
\] 
where $q^2-1 \ge 0$ by assumption. But $\im_j(x) \in C(\FX^0,P_j)$, which according 
to \cite[Prop.~4.21~(b)]{P1} is contained in $\overline{C(\FX^0,P_\infty)}$. 
\end{Proof}

For the rest of the section, the Shimura data $(P,\FX) = (G,\FX)$ are assumed to be pure
(hence $(U=0)$ is satisfied). We also assume hypothesis $(+)$. 
Let 
\[
\FX^* := \coprod \FX_j/W_j \; ,
\]
where the disjoint union is over all rational boundary components $(P_j,\FX_j)$
of $(G,\FX)$ ($W_j :=$ the unipotent radical of $P_j$). 
The set $\FX^*$ is equipped with the \emph{Satake topology}
\cite[Chap.~III, Sect.~6.1]{AMRT} (see \cite[Sect.~6.2]{P1} for an equivalent
description), with respect to which the inclusion $\FX \into \FX^*$ becomes
an open immersion with dense image. 
The space $\FX^*$ is Hausdorff \cite[Thm.~4.9~(iii), (iv)]{BB}. 
There is an action of $G(\BQ)$
by continuous automorphisms on $\FX^*$
\cite[Sect.~6.2, Sect.~4.16]{P1}, which is uniquely characterized by the requirement
of extending the action of $G(\BQ) \subset G(\BR)$ 
on $\FX$ underlying the Shimura data $(G,\FX)$. 

\begin{Cons} \label{4K}
We propose ourselves to extend the identity on $\FX$ to a map $p: \FX^{BS} \to \FX^*$. 
As $\FX^{BS} = \coprod_Q e(Q)$
($Q$ running over all parabolic subgroups of $G$), we need to define 
the restriction $p_Q$ of $p$ to the stratum $e(Q)$, for each  
parabolic subgroup $Q$ of $G$.
In order to do so, write $Q_\infty := \adm(Q)$, hence $P_\infty = \adm_{Sh}(Q)$. The 
unipotent radical of $P_\infty$ is denoted by $W_\infty \,$.
Consider the map 
\[
\FX \longto \coprod \FX_\infty 
\]
to the disjoint union of the finitely many spaces $\FX_\infty$ underlying
rational boundary components associated to $P_\infty$, and its composition
\[
\tilde{\pi}_\infty : \FX \longto \coprod \FX_\infty / W_\infty \subset \FX^*
\]
with the projections from $\FX_\infty$ to $\FX_\infty / W_\infty \,$.
Define 
\[
p_Q : e(Q) = A_Q \backslash \FX \longto \FX^* \; , \;
A_Q x \longmapsto \tilde{\pi}_\infty(x) \; .
\]
According to Corollary~\ref{4I}, $p_Q$ is well-defined.
It is continuous, as the Satake topology induces the quotient topology
on each $\FX_\infty / W_\infty$. By construction, it is the unique map
that makes the diagram
\[
\vcenter{\xymatrix@R-10pt{
        \FX \ar@{->>}[r] \ar@{=}[d] & e(Q) \ar[d]^{p_Q} \\
        \FX \ar[r]^-{\tilde{\pi}_\infty} & \coprod \FX_\infty / W_\infty 
\\}}
\]
commute. 

We claim that every element of $\coprod \FX_\infty / W_\infty$ 
is of the form
$\tilde{\pi}_\infty(x)$, for some $x \in \FX$.
Indeed, each $\FX_\infty$ is homogeneous under $P_\infty(\BR)U_\infty(\BC)$ 
as $(P_\infty,\FX_\infty)$ are mixed Shimura data. 
Therefore, the quotient $\FX_\infty / W_\infty$ is homogeneous under $P_\infty(\BR)$.
But by \cite[Sect.~4.11]{P1},
the map $\FX \to \coprod \FX_\infty$ is $P_\infty(\BR)$-equivariant,
and its image meets every connected component of $\coprod \FX_\infty \,$.

Thus,
the image of the map $p_Q$ equals the whole of $\coprod \FX_\infty / W_\infty$.
Note that for $Q = G$, the map $\tilde{\pi}_\infty$ is the identity \cite[Sect.~2.13]{P1}.

By definition, the map $p: \FX^{BS} \to \FX^*$ equals $p_Q$ on $e(Q)$, for
each parabolic subgroup $Q$ of $G$. Given the definition of the set
$\FX^*$, $p$ is surjective.
\end{Cons}

\forget{
\begin{Prop}
Let $Q$ be a parabolic subgroup of $G$. Consider the open submanifold
$\FX(Q)$ of $\FX^{BS}$.  Then
\[
p \bigl( \FX(Q) \bigr) = \coprod_{Q \subset Q_j} \FX_j/W_j \subset \FX^* \; ,
\]
where the union is over all rational boundary components 
$(P_j,\FX_j)$, for which $Q_j$ contains $Q$. 
\end{Prop}

\begin{Proof}
We have $\FX(Q) = \coprod_{Q \subset R} e(R)$, and for each $R$ occurring
in this union, the image $p(e(R))$ equals the union of those $\FX_j / W_j$, 
for which $Q_j = \adm(R)$.

On the one hand, since $Q \subset R$, any $Q_j$ containing $R$ also
contains $Q$; this is true in particular for $\adm(R)$. Therefore, $\adm(R)$ is one
of the admissible parabolics $Q_j$ satisfying $Q \subset Q_j$, and hence
\[
p \bigl( \FX(Q) \bigr) \subset \coprod_{Q \subset Q_j} \FX_j/W_j \; .
\]
On the other hand, for $R = Q_j$, we have $Q_j = \adm(R)$, hence
\[
p \bigl( \FX(Q) \bigr) \supset \coprod_{Q \subset Q_j} \FX_j/W_j \; .
\]
\end{Proof}
}
Here is the second main result of this section.

\begin{Thm} \label{4L} 
The map $p: \FX^{BS} \to \FX^*$ is the unique continuous extension of
the identity on $\FX$.
\end{Thm}

Note that as $\FX$ is dense in $\FX^{BS}$, and $\FX^*$ is Hausdorff \cite[Thm.~4.9~(iii), (iv)]{BB}, there is at most one continuous extension of $\id_\FX$. 
This establishes the uniqueness assertion. 
Before giving the proof of continuity of $p$, let us indicate how to
deduce existence of the continuous extension from  
the main results from \cite{Z} (which are formulated only 
for spaces of type $S - \BQ$ under semi-simple groups). 
It suffices to show that for every connected component $\FX^0$ of $\FX$,
the identity on $\FX^0$ extends to a continuous map $\bar{\FX}^0 \to \FX^*$.
According to Theorem~\ref{1M}, $\FX^0$ is a space of type $S - \BQ$ under $G$.
It follows \cite[Ex.~2.5~(2)]{BS} that a certain quotient $A \backslash \FX$ is of type
$S - \BQ$ under $G^{der}$. Here, $A$ is a complement of $\Stab_{Z(G)(\BR)}^{\ext}(x)$
in $Z(G)(\BR)$, for one, hence any point $x \in \FX^0$. But according to Corollary~\ref{1J},
we have 
\[
\Stab_{Z(G)(\BR)}^{\ext}(x) = Z(G)(\BR) \; ,
\]
hence $A$ is trivial, and $\FX^0$ itself is of type $S - \BQ$ under $G^{der}$.
Note that $\bar{\FX}^0$, defined with respect to the $S - \BQ$-structure under $G^{der}$,
coincides with $\bar{\FX}^0$, defined with respect to the $S - \BQ$-structure under $G$:
indeed, the map $Q \mapsto Q \cap G^{der}$ induces a bijection of sets
of parabolic subgroups, and the formation of the subgroups $S_Q$ of $G^{ad}$
\cite[Sect.~4.2]{BS} is invariant under this bijection. \\

By \cite[Sect.~(3.7)~(2)]{Z}, the identity $\id_\FX$ first extends to give a continuous map,
denoted $p^*$ in \loccit ,
from $\bar{\FX}^0$ to a space denoted ${}_\BQ \tilde{\FX}^* \,$; actually, the latter
is \emph{constructed} as a quotient of $\bar{\FX}^0$. Implicit in this construction
is the choice of an irreducible representation $\tau$ of $G^{der}$. According to
\cite[Theorem~(3.10)]{Z}, there is then a continuous bijection
from ${}_\BQ \tilde{\FX}^*$ to the Satake compactification,
denoted ${}_\BQ \FX^*$. 
Finally \cite[Sect.~(3.11)]{Z}, for certain choices of $\tau$, the spaces
${}_\BQ \FX^*$ and $(\FX^0)^* \subset \FX^*$ are homeomorphic.

\medskip

\begin{Proofof}{Theorem~\ref{4L}}
We need to show that for every point $\bar{x}$ of $\FX^{BS}$, the pre-image
under $p$ of every member of a fundamental system of neighbourhoods of $p(\bar{x}) \in \FX^*$
contains a neighbourhood of $\bar{x}$.

Let $Q$ be a parabolic subgroup of $G$, and 
\[
\bar{x} = A_Qx \in e(Q) \subset \FX^{BS} \; .
\]
Write $Q_\infty := \adm(Q)$, hence $P_\infty = \adm_{Sh}(Q)$.
The unipotent radical is denoted by $W_\infty \,$, and $U_\infty \,$
is its weight $(-2)$-part.
By definition, we have $p(\bar{x}) = \tilde{\pi}_\infty(x)$, where 
\[
\tilde{\pi}_\infty : \FX \longto \coprod \FX_\infty / W_\infty \subset \FX^* \; .
\]
Denote by $\FX^0$ the connected component of $\FX$ containing $x$,
and by $\FX_\infty$ the space for which $\tilde{\pi}_\infty(x) \in \FX_\infty / W_\infty$.
Let us describe a fundamental system of neighbourhoods of $\tilde{\pi}_\infty(x)$, following
\cite[Sect.~6.2]{P1}: fix a \emph{convex core} $D \subset C(\FX^0,P_\infty)$.
The precise definition of this notion will not be important; however, our choice of $D$
needs to be done in a way stable under addition by elements
of $\overline{C(\FX^0,P_\infty)}$:
\[
D + \overline{C(\FX^0,P_\infty)} \subset D
\]
\cite[Sect.~6.1]{P1}. The fundamental system has two parameters: first, the
neighbourhoods $\FY$ of $\tilde{\pi}_\infty(x)$ in $\FX_\infty / W_\infty \,$; second,
strictly positive real numbers $\mu$. To such a pair $(\FY,\mu)$, one
associates the open subset
\[
\FU_{(\FY,\mu)} :=
\coprod \tilde{\pi}_j \bigl( \tilde{\pi}_\infty^{-1}(\FY) \cap \im_\infty^{-1}(\mu D) \bigr)
\subset \FX^* \; , 
\]
where the disjoint union is over all rational boundary components $(P_j,\FX_j)$
between $(G,\FX)$ and $(P_\infty,\FX_\infty)$, and
\[
\tilde{\pi}_j : \FX \longto \coprod \FX_j / W_j \subset \FX^* 
\]
is defined in the same way as $\tilde{\pi}_\infty$. Set 
\[
\FV_{(\FY,\mu)} := \FU_{(\FY,\mu)} \cap \FX 
= \tilde{\pi}_\infty^{-1}(\FY) \cap \im_\infty^{-1}(\mu D) \; .
\]
The point $\bar{x} = A_Qx$ belongs to $e(Q)$, hence to the open
subset $\FX(Q) = \bar{A}_Q \, {}^{A_Q} \! \times \FX$ of $\FX^{BS}$.
We identify $A_Q$ and $( \BR_+^* )^J$, hence $\bar{A}_Q$ and $(0,+\infty]^J$
(Theorem~\ref{2H}~(b)).
Given the definition of the quotient topology, the image $\FW_{(\FY,\mu)}$
of the product $(1,+\infty]^J \times \FV_{(\FY,\mu)}$ under the projection
\[
\bar{A}_Q \times \FX \longonto \bar{A}_Q \, {}^{A_Q} \! \times \FX = \FX(Q)
\]
is open in $\FX(Q)$, hence in $\FX^{BS}$. It contains  
$[(\infty_Q,x)] = A_Q x$, \emph{i.e.}, it is a neighbourhood of $\bar{x}$.
Our proof will be complete once we have established that 
$p(\FW_{(\FY,\mu)}) \subset \FU_{(\FY,\mu)}$.

First, note that $\tilde{\pi}_\infty^{-1}(\FY)$ is stable under
the geodesic action of $A_Q$ (Corollary~\ref{4I}), and $\im_\infty^{-1}(\mu D)$
is stable under the geodesic action of $[1,+\infty)^J \subset A_Q$ (Corollary~\ref{4J}),
given our choice of $D$ (indeed, $\overline{C(\FX^0,P_\infty)}$ being a cone,
we have $\mu \overline{C(\FX^0,P_\infty)} = \overline{C(\FX^0,P_\infty)}$).
Therefore, $\FV_{(\FY,\mu)}$ is stable under the geodesic action of $[1,+\infty)^J$.
This means that for any parabolic $R$ of $G$ containing $Q$, the intersection  
$\FW_{(\FY,\mu)} \cap e(R)$ is contained in the image of $\FV_{(\FY,\mu)}$ 
under the projection $\FX \onto e(R)$.

Next, fix one of the rational boundary components $(P_j,\FX_j)$
between $(G,\FX)$ and $(P_\infty,\FX_\infty)$, and a parabolic $R$ of $G$ containing $Q$,
and such that $\adm(R) = Q_j$. The contribution of $e(R)$ to the intersection
\[
p(\FW_{(\FY,\mu)}) \cap \FX_j/W_j 
\]
equals $p_R (\FW_{(\FY,\mu)} \cap e(R))$, which by the above is contained in
the image of $\FV_{(\FY,\mu)}$ under the composition of $\FX \onto e(R)$ and $p_R$.

The commutativity of the diagram
\[
\vcenter{\xymatrix@R-10pt{
        \FX \ar@{->>}[r] \ar@{=}[d] & e(R) \ar[d]^{p_R} \\
        \FX \ar[r]^-{\tilde{\pi}_j} & \coprod \FX_j / W_j 
\\}}
\]
then establishes the desired inclusion 
\[
p_R (\FW_{(\FY,\mu)} \cap e(R))
\subset \tilde{\pi}_j \bigl( \FV_{(\FY,\mu)} \bigr) \; .
\]
\end{Proofof}

Let us note the following aspect of our construction.

\begin{Comp} \label{4Lb}
Fix an admissible parabolic subgroup $Q_j$ of $G$, and consider the disjoint union 
$\coprod \FX_j$ of the finitely many spaces $\FX_j$ underlying
rational boundary components associated to $Q_j \,$, as well as its quotient
$\coprod \FX_j / W_j$ by the action of the unipotent radical $W_j \,$, 
considered as a locally closed subset of $\FX^*$. Then
under the map $p: \FX^{BS} \to \FX^*$ from Theorem~\ref{4L},
\[
p^{-1} \bigl( \coprod \FX_j / W_j \bigr) = \coprod_{P_j \subset R \subset Q_j} e(R) \subset \FX^{BS}
\; ,
\]
where the disjoint union runs over the 
parabolics $R$ of $P$ contained in $Q_j$ and containing
$P_j \,$, the canonical normal subgroup of $Q_j \,$. 
\forget{\\[0.1cm]
(b)~Let $Q$ be a parabolic of $G$. If $P_j \not \subset Q$ or if $Q_j \cap Q$
is not parabolic, then the intersection of $p^{-1} ( \coprod \FX_j / W_j )$
and the closure $\overline{e(Q)}$ of $e(Q)$ in $\FX^{BS}$ satisfies
\[
p^{-1} \bigl( \coprod \FX_j / W_j \bigr) \cap \overline{e(Q)} = \emptyset \; .
\]
If $P_j \subset Q$ and $Q_j \cap Q$
is parabolic, then 
\[
p^{-1} ( \coprod \FX_j / W_j ) \cap \overline{e(Q)} 
= \coprod_{P_j \subset R \subset Q_j \cap Q} e(R)
\]
is open in $\overline{e(Q_j \cap Q)}$, and the immersion
\[
p^{-1} \bigl( \coprod \FX_j / W_j \bigr) \cap \overline{e(Q)} \longinto \overline{e(Q_j \cap Q)} 
\]
is contractible.
}  
\end{Comp}

\begin{Proof}
This follows directly from Construction~\ref{4K}
(as $P_j \subset R \subset Q_j$ if and only if $Q_j = \adm(R)$).
\forget{
As for (b), recall first that
\[
\overline{e(Q)} = \coprod_{R \subset Q} e(R) \; ,
\]
where $R$ runs over all parabolics contained in $Q$. 
Part~(a) then implies that
\[
p^{-1} ( \coprod \FX_j / W_j ) \cap \overline{e(Q)} 
= \coprod_{P_j \subset R \subset Q_j \cap Q} e(R) \; .
\]
Our claim 
thus follows from the definition of the topology on $\FX^{BS}$
(observe that $p^{-1} ( \coprod \FX_j / W_j ) \cap \overline{e(Q)}$ 
is the intersection of $\overline{e(Q_j \cap Q)}$
with the union of the corners $\FX(Q')$, for all $Q'$ containing $P_j$),
and from Proposition~\ref{2Jc}.
}
\end{Proof}

\begin{Rem}
Complement~\ref{4Lb} should be compared to the first statement of
\cite[Prop.~(3.8)~(ii)]{Z}.
\end{Rem}

As far as equivariance of $p$ is concerned, we have the following result.

\begin{Comp} \label{4La} 
The continuous map $p: \FX^{BS} \to \FX^*$ from Theorem~\ref{4L} is $G(\BQ)$-equivariant 
(with respect to the action~(1) on $\FX^{BS}$). 
\end{Comp}

\begin{Proof}
$\FX$ is dense in $\FX^{BS}$, and $\FX^*$ is Hausdorff \cite[Thm.~4.9~(iii), (iv)]{BB}.
\end{Proof}

\begin{Cor} \label{4M}
Let $K \subset G(\BA_f)$ be an open compact subgroup. Then the identity on
the Shimura variety $M^K (G,\FX)$ extends uniquely to a conti\-nuous map
\[
p^K: M^K  (G,\FX) (\BC)^{BS} \longto M^K (G,\FX)^* (\BC)   
\]
between the Borel--Serre compactification and the space of complex points
of the Baily--Borel compactification of $M^K (G,\FX)$.
\end{Cor}

\begin{Proof}
The space 
\[
M^K (G,\FX)^* (\BC) = G (\BQ) \backslash \bigl( \FX^* \times G (\BA_f) / K \bigr) 
\]
is compact \cite[Chap.~II, Thm.~2]{AMRT}, hence Hausdorff; therefore there is at
most one continuous extension of $\id_{M^K (G,\FX)}$. 
In order to define it, use the map
\[
p \times \id_{G(\BA_f)} : \FX^{BS} \times G (\BA_f) \longto \FX^* \times G (\BA_f)
\]
and Complement~\ref{4La}.
\end{Proof}


\bigskip
%
%

\section{The canonical stratifications}
\label{4a}



The aim of the present section is to define the canonical stratifications of both
the Borel-Serre (Definition~\ref{4aA}) and the Baily--Borel compactification
(Definition~\ref{4aF}). The main result (Theorem~\ref{4aG}) then gives
a description of the intersections of the closures of the
canonical strata of the Borel-Serre compactification
on the one hand, and the pre-images under the map $p^K$ of the canonical strata
of the Baily--Borel compactification on the other. An important feature of
these intersections is that the immersion into their closures is contractible
provided the ``level'' $K$ is neat (Theorem~\ref{4aG}~(c)). \\

We fix mixed Shimura data $(P,\FX)$ satisfying hypotheses $(+)$ and $(U=0)$,
and an open compact subgroup $K$ of $P(\BA_f)$. 
The stratification of $\FX^{BS}$ by the faces $e(R)$ induces a stratification of
the Borel--Serre compactification $M^K  (P,\FX) (\BC)^{BS}$
(cmp.\ \cite[Prop.~9.4~(i)]{BS}). Every stratum
is the image $e^K(R,g)$ of a copy 
\[
e(R) \times \{ gK \} \subset \FX^{BS} \times P (\BA_f) / K
\]
of a face $e(R)$ under the projection from $\FX^{BS} \times P (\BA_f) / K$ 
to 
\[
M^K (P,\FX) (\BC)^{BS} = P (\BQ) \backslash \bigl( \FX^{BS} \times P (\BA_f) / K \bigr) \; ,
\] 
for some $g \in P(\BA_f)$. 
Denoting by $\overline{e^K(R,g)}$ the closure of $e^K(R,g)$ in $M^K  (P,\FX) (\BC)^{BS}$,
one checks that $\overline{e^K(R,g)}$ equals the image of 
\[
\overline{e(R)} \times \{ gK \} \subset \FX^{BS} \times P (\BA_f) / K
\]
under the projection from $\FX^{BS} \times P (\BA_f) / K$ 
(cmp.\ \cite[Prop.~9.4~(ii)]{BS}), \emph{i.e.},
\[
\overline{e^K(R,g)} = \bigcup_{R' \subset R} e^K(R',g) \; ,
\]
where $R'$ runs over all parabolics contained in $R$.
The stratification of $M^K  (P,\FX) (\BC)^{BS}$ that will be of interest for us,
is coarser than the one by the $e^K(R,g)$.

\begin{Def} \label{4aA}
(a)~Let $R$ be a parabolic subgroup of $P$. Define 
\[
e^K \bigl( R,P(\BA_f) \bigr) \subset M^K  (P,\FX) (\BC)^{BS}
\]
to be the image of $e(R) \times P (\BA_f) / K \subset \FX^{BS} \times P (\BA_f) / K$
under the projection from $\FX^{BS} \times P (\BA_f) / K$. In other words,
\[
e^K \bigl( R,P(\BA_f) \bigr) = \bigcup_{g \in P (\BA_f)} e^K(R,g) \; .
\]
(b)~The \emph{canonical stratification of
the Borel--Serre compactification}
is the stratification by the $e^K ( R,P(\BA_f) )$, where $R$ runs through the
parabolic subgroups of $P$.
\end{Def}

We leave it to the reader to check that $e^K ( R,P(\BA_f) )$ and $e^K ( R',P(\BA_f) )$
have a non-empty intersection if and only if they are equal, which in turn is equivalent
to $R$ and $R'$ being $P(\BQ)$-conjugate to each other. In particular, the $e^K ( R,P(\BA_f) )$
\emph{do} form a stratification.
Also,
\[
e^K \bigl( R,P(\BA_f) \bigr) = \bigcup_{g \in P (\BA_f)} e^K(R,g)
= \coprod_{g \in I} e^K(R,g)
\]
if the index set $I$ on the right hand side is chosen to be a set of representatives
of the (finite) double quotient $R(\BQ) \backslash P(\BA_f) / K$. 
Denoting by $\overline{e^K ( R,P(\BA_f) )}$
the closure of $e^K(R,P(\BA_f))$, we thus have
\[
\overline{e^K \bigl( R,P(\BA_f) \bigr)}
= \bigcup_{g \in P (\BA_f)} \overline{e^K(R,g)}
= \bigcup_{g \in I} \overline{e^K(R,g)} \; ,
\]
and the latter union is disjoint:
\[
\overline{e^K \bigl( R,P(\BA_f) \bigr)}
= \coprod_{g \in I} \overline{e^K(R,g)} \; ,
\]
as follows from \cite[Cor.~7.7~(3)]{BS}. The closure
$\overline{e^K ( R,P(\BA_f) )}$
equals the image of $\overline{e(R)} \times P (\BA_f) / K$
under the projection from $\FX^{BS} \times P (\BA_f) / K$ to $M^K  (P,\FX) (\BC)^{BS}$. Thus,
\[
\overline{e^K \bigl( R,P(\BA_f) \bigr)} = \bigcup_{R' \subset R} e^K \bigl( R',P(\BA_f) \bigr) \; ,
\]
where $R'$ runs over all parabolics contained in $R$.

\begin{Prop} \label{4aB}
Let $Q$ and $R$ be parabolic subgroups of $P$. \\[0.1cm]
(a)~The stratum $e^K(R,P(\BA_f))$ is contained in $\overline{e^K ( Q,P(\BA_f) )}$
if and only if $R$ is contained in a $P(\BQ)$-conjugate of $Q$. \\[0.1cm]
(b)~We have
\[
\overline{e^K \bigl( R,P(\BA_f) \bigr)} \cap \overline{e^K \bigl( Q,P(\BA_f) \bigr)}
= \bigcup_\gamma \overline{e^K \bigl( R \cap \gamma Q \gamma^{-1},P(\BA_f) \bigr)} \; ,
\]
where the union runs over all $\gamma \in P(\BQ)$ such that $R \cap \gamma Q \gamma^{-1}$
is parabolic.
\end{Prop}

\begin{Proof}
(a): if $R \subset \gamma Q \gamma^{-1}$, then 
\[
e^K \bigl( R,P(\BA_f) \bigr) = e^K \bigl( \gamma^{-1} R \gamma,P(\BA_f) \bigr)
\subset \overline{e^K \bigl( Q,P(\BA_f) \bigr)} \; .
\]
Conversely, assume $e^K(R,P(\BA_f)) \subset \overline{e^K ( Q,P(\BA_f) )}$. Then
\[
e^K \bigl( R,P(\BA_f) \bigr) = e^K \bigl( R',P(\BA_f) \bigr) \; ,
\]
for some parabolic $R' \subset Q$. But then, $R$ and $R'$ are $P(\BQ)$-conjugate to each other.

\noindent (b): let $R'$ be a parabolic contained in $R$. Then according to (a),
\[
e^K \bigl( R' ,P(\BA_f) \bigr) \subset \overline{e^K \bigl( Q,P(\BA_f) \bigr)}
\]
if and only if $R' \subset \gamma Q \gamma^{-1}$, for some $\gamma \in P(\BQ)$.
But then the intersection $R \cap \gamma Q \gamma^{-1}$, containing a parabolic,
is itself parabolic.
\end{Proof}

Recall the map $\adm_{Sh}$ from Definition~\ref{4E} between the set of parabolics of $P$
and the set of subgroups of $P$ underlying boundary components of $(P,\FX)$.

\begin{Def} \label{4aC}
Let $Q$ be a parabolic subgroup of $P$. Define
\[
e(Q)' := \coprod_{\adm_{Sh}(Q) \subset R \subset Q} e(R) \subset \FX^{BS}
\]
and
\[
e^K \bigl( Q,P(\BA_f) \bigr)' := 
       \bigcup_{\adm_{Sh}(Q) \subset R \subset Q} e^K \bigl( R,P(\BA_f) \bigr) 
\subset M^K  (P,\FX) (\BC)^{BS} \; ,
\]
where the unions run over all parabolic subgroups $R$ contained in $Q \,$,
and containing $\adm_{Sh}(Q)$.
\end{Def}

We have the double inclusion
\[
e^K \bigl( Q,P(\BA_f) \bigr) \subset e^K \bigl( Q,P(\BA_f) \bigr)' 
\subset \overline{e^K \bigl( Q,P(\BA_f) \bigr)} \; .
\]
In particular, the closure of $e^K ( Q,P(\BA_f) )'$ equals 
$\overline{e^K ( Q,P(\BA_f) )}$. \\

\begin{Prop} \label{4aD}
Let $Q$ be a parabolic subgroup of $P$. \\[0.1cm]
(a)~The set $e^K ( Q,P(\BA_f) )'$ is the image of 
\[
e(Q)' \times P (\BA_f) / K \subset \FX^{BS} \times P (\BA_f) / K
\]
under the projection to $M^K  (P,\FX) (\BC)^{BS}$. \\[0.1cm]
(b)~The set $e^K ( Q,P(\BA_f) )'$ is open in $\overline{e^K ( Q,P(\BA_f) )}$. 
In particular, it is locally closed in $M^K  (P,\FX) (\BC)^{BS}$. \\[0.1cm]
(c)~If $K$ is neat, then the inclusion
\[
e^K \bigl( Q,P(\BA_f) \bigr)' \longinto \overline{e^K \bigl( Q,P(\BA_f) \bigr)}  
\]
is contractible.
\end{Prop}

\begin{Proof}
Let us first show the following claim $(*)$: if $R \subset Q$ is a parabolic such that 
\[
e^K \bigl( R,P(\BA_f) \bigr) \subset e^K \bigl( Q,P(\BA_f) \bigr)' \; ,
\]
then $\adm_{Sh}(Q) \subset R \subset Q$. 

The assumption on $R$ is equivalent to the existence of $\gamma \in P(\BQ)$ such that 
$\adm_{Sh}(Q) \subset \gamma R \gamma^{-1} \subset Q$. This implies that $R$ is contained
in both $Q$ and $\gamma^{-1} Q \gamma$. Therefore, the intersecion $Q \cap \gamma^{-1} Q \gamma$
remains parabolic. According to Proposition~\ref{Par}, this implies 
$Q = \gamma^{-1} Q \gamma$, \emph{i.e.}, 
that $\gamma \in Q$. But $\adm_{Sh}(Q)$ is normal in $\adm(Q)$, hence in $Q$.
Altogether, the inclusion $\adm_{Sh}(Q) \subset \gamma R \gamma^{-1} \subset Q$ implies
that $\adm_{Sh}(Q) \subset R \subset Q$.

Claim~(a) is a direct consequence of $(*)$.

As for (b), we have $e^K ( Q,P(\BA_f))' \subset \overline{e^K ( Q,P(\BA_f) )}$, and
$\overline{e^K ( Q,P(\BA_f) )}$ is the image of $\overline{e(Q)} \times P (\BA_f) / K$
under the projection from $\FX^{BS} \times P (\BA_f) / K$; actually, this projection
identifies $\overline{e^K ( Q,P(\BA_f) )}$ with a topological quotient of
$\overline{e(Q)} \times P (\BA_f) / K$. Now according to $(*)$, the pre-image of
$e^K ( Q,P(\BA_f))'$ under the projection equals
\[
e(Q)' \times P (\BA_f) / K \subset \FX^{BS} \times P (\BA_f) / K \; .
\]
But $e(Q)' \subset \overline{e(Q)}$ is open (it equals the intersection of 
$\overline{e(Q)}$ with the union of the corners $\FX(R)$, for all parabolics $R$
containing $\adm_{Sh}(Q)$). 

Part~(c) is a special case of the following result.
\end{Proof}

\begin{Prop} \label{4aDa}
Let $R$ be a parabolic subgroup of $P$. \\[0.1cm]
(a)~Let $e^K ( R,P(\BA_f) )^o \subset M^K  (P,\FX) (\BC)^{BS}$ be
a locally closed union of strata containing $e^K ( R,P(\BA_f) )$, 
and contained in the closure $\overline{e^K ( R,P(\BA_f) )}$.
If $K$ is neat, then the open immersion 
\[
e^K \bigl( R,P(\BA_f) \bigr)^o \longinto \overline{e^K \bigl( R,P(\BA_f) \bigr)}  
\]
is contractible. \\[0.1cm]
(b)~Define
\[
\partial \overline{e^K \bigl( R,P(\BA_f) \bigr)}
:= \overline{e^K \bigl( R,P(\BA_f) \bigr)} - e^K \bigl( R,P(\BA_f) \bigr) \; .
\]
Let $Q \subset R$ be a maximal proper parabolic of $R$ (\emph{i.e.}, a
parabolic containing $Q$ and contained in $R$ is equal either to $Q$ or to $R$).
If $K$ is neat, then the open immersion 
\[
\partial \overline{e^K \bigl( R,P(\BA_f) \bigr)} - \overline{e^K \bigl( Q,P(\BA_f) \bigr)}
\longinto \partial \overline{e^K \bigl( R,P(\BA_f) \bigr)} - e^K \bigl( Q,P(\BA_f) \bigr)
\] 
is contractible.   
\end{Prop}

\begin{Proof}
The restriction of the projection from $\FX^{BS} \times P (\BA_f) / K$
to each of the $\overline{e(R)} \times \{ gK \}$, $g \in P (\BA_f)$,
equals the quotient by the action
of the stabilizer in $P(\BQ)$ of $e(R) \times \{ gK \}$ 
\cite[Prop.~9.4~(ii)]{BS}.
But this action is free on the whole of $\FX^{BS} \times P (\BA_f) / K$ 
as the group $K$ is supposed neat \cite[Sect.~9.5]{BS}.
It follows that both (a) and (b) can be proved after pull-back to
$\overline{e(R)} \times \{ gK \}$, $g \in P (\BA_f)$.

\noindent (a): apply Proposition~\ref{2Jc}.

\noindent (b): note that for a parabolic subgroup $Q'$ of $R$, 
the set $e^K ( Q',g )$ is contained in $\overline{e^K ( Q,P(\BA_f) )}$
if and only if $Q' \subset \gamma Q \gamma^{-1}$, for some $\gamma \in P(\BQ)$.  
The intersection $R \cap \gamma R \gamma^{-1}$ then being parabolic
(as it contains $Q'$), we conclude
from Proposition~\ref{Par} that $\gamma \in R(\BQ)$. 
Now apply Proposition~\ref{2Jd}.
\end{Proof}

\forget{
\begin{Lem} \label{4aLem}
Let $Q$ be a parabolic subgroup of $P$, and $\gamma \in P(\BQ)$.
Then the following are equivalent.
\begin{enumerate}
\item[(i)] $Q = \gamma Q \gamma^{-1}$, \emph{i.e.}, $\gamma \in Q(\BQ)$,
\item[(ii)] $e(Q)' = e(\gamma Q \gamma^{-1})'$,
\item[(iii)] the intersection $e(Q)' \cap e(\gamma Q \gamma^{-1})'$
is not empty.
\end{enumerate}
\end{Lem}

\begin{Proof}
Clearly (i)$\Rightarrow$(ii)$\Rightarrow$(iii).

It remains to show that (iii) implies (i). Thus, suppose 
\[
e(Q)' \cap e(\gamma Q \gamma^{-1})' \ne \emptyset \; .
\]
This intersection then necessarily contains a whole stratum, which is of the form
$e(R)$, for some parabolic $R$ contained in both $Q$ and $\gamma Q \gamma^{-1}$. 
This implies in particular that $Q \cap \gamma Q \gamma^{-1}$ is parabolic. 
Now apply Proposition~\ref{Par}.
\end{Proof}
}
Here is the key property of the sets $e^K ( Q_j,P(\BA_f) )'$.

\begin{Prop} \label{4aE}
Let $Q_j$ be an admissible parabolic subgroup of $P$, 
with canonical normal subgroup $P_j \,$. 
Let $Q$ be a parabolic subgroup of $P$ containing $P_j \,$, 
and such that $Q_j \cap Q$ remains parabolic. Then
\[
e^K \bigl( Q_j,P(\BA_f) \bigr)' \cap \overline{e^K \bigl( Q,P(\BA_f) \bigr)}
= e^K \bigl( Q_j \cap Q ,P(\BA_f) \bigr)' \; .
\]
\end{Prop}

This statement should be compared to Proposition~\ref{4aB}~(b).

\medskip

\begin{Proofof}{Proposition~\ref{4aE}}
Write $Q_k := \adm(Q)$ and $P_k := \adm_{Sh}(Q)$.
Thus,
\[
P_k \subset Q \subset Q_k \; .
\]
We have $P_j \subset Q_j \cap Q \subset Q_j \,$, meaning that $\adm_{Sh}(Q_j \cap Q) = P_j$
and $\adm(Q_j \cap Q) = Q_j$.
This implies that
\[
e^K \bigl( Q_j,P(\BA_f) \bigr)' \cap \overline{e^K \bigl( Q,P(\BA_f) \bigr)}
\supset e^K \bigl( Q_j \cap Q ,P(\BA_f) \bigr)' \; .
\]
In order to show the reverse inclusion, 
note first that the parabolic $Q_j \cap Q$ is contained in $Q_k$. 
Therefore (Theorem~\ref{4D}~(b)), the relation $Q_j = \adm(Q_j \cap Q) \preceq Q_k$ is satisfied.

Next,
let $R \subset Q_j$ be a parabolic, and 
$\gamma \in P(\BQ)$ such that $P_j \subset R \subset Q_j \cap  \gamma Q \gamma^{-1}$.
In other words,
\[
e^K \bigl( R,P(\BA_f) \bigr) 
\subset e^K \bigl( Q_j,P(\BA_f) \bigr)' \cap \overline{e^K \bigl( Q,P(\BA_f) \bigr)} \; .
\]
We have 
$P_j \subset R \subset Q_j$, meaning that $\adm(R) = Q_j$.
But $R$ is also contained in the admissible parabolic $\gamma Q_k \gamma^{-1}$.
Therefore (Theorem~\ref{4D}~(b)), the relation $Q_j \preceq \gamma Q_k \gamma^{-1}$,
or equivalently,
$\gamma^{-1} Q_j \gamma \preceq Q_k$ holds.
Thus, the conjugate admissible parabolics $Q_j$ and $\gamma^{-1} Q_j \gamma$ both
give rise to Shimura data, which are boundary components of Shimura data associated to $Q_k$.
According to \cite[Rem.~(iii) on p.~91, after Cor.~4.20]{P1},  
the parabolics $Q_j$ and $\gamma^{-1} Q_j \gamma$
are conjugate under $P_k(\BQ)$, \emph{i.e.}, 
\[
\gamma^{-1} Q_j \gamma = \delta^{-1} Q_j \delta \; , 
\]
for some $\delta \in P_k$. 
The product $\delta \gamma^{-1}$ normalizes $Q_j$, and hence belongs to $Q_j$. The subgroup
$P_j$ of $Q_j$ being normal, we have $P_j = (\delta \gamma^{-1}) P_j (\delta \gamma^{-1})^{-1}$.
Hence the inclusion $P_j \subset R \subset Q_j \cap \gamma Q_k \gamma^{-1}$ implies
\[
P_j \subset (\delta \gamma^{-1}) R (\delta \gamma^{-1})^{-1}
\subset (\delta \gamma^{-1}) (Q_j \cap \gamma Q \gamma^{-1}) (\delta \gamma^{-1})^{-1}
= Q_j \cap \delta Q \delta^{-1} \; .
\] 
But $\delta \in P_k \subset Q$, therefore $\delta Q \delta^{-1} = Q$. Altogether,
\[
P_j \subset (\delta \gamma^{-1}) R (\delta \gamma^{-1})^{-1} \subset Q_j \cap Q \; ,
\]
hence
\[
e^K \bigl( R,P(\BA_f) \bigr) 
= e^K \bigl( (\delta \gamma^{-1}) R (\delta \gamma^{-1})^{-1},P(\BA_f) \bigr)
\subset e^K \bigl( Q_j \cap Q ,P(\BA_f) \bigr)' \; .
\]
\end{Proofof}

The hypothesis ``$P_j \subset Q$'' is essential in Proposition~\ref{4aE},
as is illustrated by the following result.

\begin{Prop} \label{4aEa}
Let $Q_j$ and $Q_k$ be admissible parabolic subgroups of $P$.
Assume that $Q_j \preceq Q_k$. Then 
\[
\overline{e^K \bigl( Q_j,P(\BA_f) \bigr)} \cap e^K \bigl( Q_k,P(\BA_f) \bigr)' 
\]
is non-empty if and only if $Q_j = Q_k$.
\end{Prop}

\begin{Proof}
Let $R$ be a parabolic of $P$, such that $\adm(R) = Q_k$, \emph{i.e.}, such that
\[
e^K \bigl( R,P(\BA_f) \bigr) \subset e^K \bigl( Q_k,P(\BA_f) \bigr)' \; . 
\]
If in addition we have
\[
e^K \bigl( R,P(\BA_f) \bigr) \subset \overline{e^K \bigl( Q_j,P(\BA_f) \bigr)} \; , 
\]
then according to Proposition~\ref{4aB}~(a), there exists $\gamma \in P(Q)$ such that
$R \subset \gamma Q_j \gamma^{-1}$. By Theorem~\ref{4D}~(b), we get 
$Q_k \preceq \gamma Q_j \gamma^{-1}$. Together with our assumption $Q_j \preceq Q_k$,
this implies that $Q_j \preceq \gamma Q_j \gamma^{-1}$. In particular, the intersection
$Q_j \cap \gamma Q_j \gamma^{-1}$ remains parabolic. From Proposition~\ref{Par},
we conclude that $\gamma Q_j \gamma^{-1} = Q_j$. Therefore, $Q_k \preceq Q_j$.
\end{Proof}

\begin{Def} \label{4aF}
Assume that the Shimura data $(P,\FX) = (G,\FX)$ are pure. \\[0.1cm]
(a)~Let $Q_j$ be an admissible parabolic subgroup of $G$. Consider the disjoint union 
$\coprod \FX_j$ of the finitely many spaces $\FX_j$ underlying
rational boundary components associated to $Q_j \,$, as well as its quotient
$\coprod \FX_j / W_j$ by the action of the unipotent radical $W_j \,$, 
considered as a locally closed subset of $\FX^*$. Define
\[
M^K (Q_j,\FX) \subset M^K (G,\FX)^*
\]
to be the image of $\coprod \FX_j / W_j \times G(\BA_f)/K$ under the projection
\[
\FX^* \times G (\BA_f) / K \longonto M^K (G,\FX)^* (\BC) \; . 
\]
(b)~The \emph{canonical stratification of the Baily--Borel compactification}
is the stratification by the $M^K (Q_j,\FX)$, where $Q_j$ runs through the
admissible para\-bolic subgroups of $G$.
\end{Def}

We leave it to the reader to check that $M^K (Q_j,\FX)$ and $M^K (Q_k,\FX)$
have a non-empty intersection if and only if they are equal, which in turn is equi\-valent
to $Q_j$ and $Q_k$ being $G(\BQ)$-conjugate to each other. 

\begin{Rem} \label{4aRem}
(a)~Let us connect the $M^K (Q_j,\FX)$ to the notation introduced in \cite[Sect.~6.3]{P1}:
in \loccit, quotients $\Delta_1 \backslash M^{\pi_1(K_f^1)} ((P_1,\FX_1)/W_1)$ 
of ``smaller'' Shimura varieties $M^{\pi_1(K_f^1)} ((P_1,\FX_1)/W_1)$ are considered,
indexed by rational boundary components $(P_1,\FX_1)$ and $g \in G(\BA_f)$.    
Then $M^K (Q_j,\FX)$ equals the union of the 
$\Delta_1 \backslash M^{\pi_1(K_f^1)} ((P_1,\FX_1)/W_1)$, for all $g \in G(\BA_f)$,
and all boundary components $(P_1,\FX_1)$ associated to $Q_j$ (\emph{i.e.},
satisfying $P_1 = P_j$, the canonical normal subgroup of $Q_j$). 
In particular, the canonical stratification of $M^K (G,\FX)$
is coarser than the stratification considered in \cite{P1}. \\[0.1cm]
(b)~\emph{A priori}, Definition~\ref{4aF} concerns a locally closed
subset $M^K (Q_j,\FX)(\BC)$ of the space of $\BC$-valued
points $M^K (G,\FX)^*(\BC)$ of $M^K (G,\FX)^*$.
It follows from \cite[Main Theorem~12.3 for the Baily--Borel compactification]{P1} that
this subset is indeed identified with the set of $\BC$-valued points of a locally
closed subscheme $M^K (Q_j,\FX)$ of $M^K (G,\FX)^*$.  
\end{Rem}

Putting everything together, we obtain the main result of this section.

\begin{Thm} \label{4aG}
Assume that $(P,\FX) = (G,\FX)$ are pure, and that they satisfy hypothesis $(+)$. 
Let $Q_j$ be an admissible parabolic subgroup of $G$. \\[0.1cm]
(a)~Under the map $p^K: M^K  (G,\FX) (\BC)^{BS} \to M^K (G,\FX)^* (\BC)$, we have
\[
(p^K)^{-1} \bigl( M^K (Q_j,\FX)(\BC) \bigr) = 
         e^K \bigl( Q_j, G(\BA_f) \bigr)' \subset M^K  (G,\FX) (\BC)^{BS} \; .
\]
(b)~Let $Q$ be a parabolic of $G$ containing $P_j \,$, and such that $Q_j \cap Q$ remains
parabolic. Then 
\[ 
(p^K)^{-1} \bigl( M^K (Q_j,\FX) (\BC) \bigr) \cap \overline{e^K \bigl( Q,G(\BA_f) \bigr)}
= e^K \bigl( Q_j \cap Q, G(\BA_f) \bigr)' \; .
\]
In particular,
the intersection $(p^K)^{-1} ( M^K (Q_j,\FX)(\BC) ) \cap \overline{e^K(Q,G(\BA_f))}$
contains $e^K(Q_j \cap Q,G(\BA_f))$, and is open in 
$\overline{e^K(Q_j \cap Q,G(\BA_f))}$. \\[0.1cm]
(c)~Assume $K$ to be neat.
Let $Q$ be a parabolic of $G$ containing $P_j \,$, and such that $Q_j \cap Q$ remains
parabolic. Then the immersion
\[
(p^K)^{-1} \bigl( M^K (Q_j,\FX)(\BC) \bigr) \cap \overline{e^K \bigl( Q,G(\BA_f)} 
\longinto \overline{e^K \bigl( Q_j \cap Q ,G(\BA_f) \bigr)}
\]
is contractible. In particular, the immersion 
\[
(p^K)^{-1} \bigl( M^K (Q_j,\FX)(\BC) \bigr) \longinto \overline{e^K \bigl( Q_j,G(\BA_f) \bigr)} 
\]
is contractible. 
\end{Thm}

\begin{Proof}
Part~(a) follows from Complement~\ref{4Lb}, part~(b) from (a),
Proposition~\ref{4aD}~(b) and Proposition~\ref{4aE},
and part~(c) from (b) and Proposition~\ref{4aD}~(c).
\end{Proof}

\begin{Rem} \label{4aH}
(a)~According to Theorem~\ref{4aG}~(a), the map $p^K$
is a \emph{morphism of stratifications}: indeed, the pre-image under $p^K$ of any stratum
of the canonical stratification of $M^K (G,\FX)^* (\BC)$ is a union of strata
of the canonical stratification of $M^K  (G,\FX) (\BC)^{BS}$. \\[0.1cm]
(b)~Theorem~\ref{4aG} states that the pre-images $e^K ( Q_j, G(\BA_f) )'$ 
of the canonical strata of the Baily--Borel compactification have much better separation
properties than the canonical strata $M^K (Q_j,\FX)$ themselves. Let us illustrate
what we mean: let $Q_j$ and $Q_k$ be admissible parabolics of $G$, and assume that
$Q_j \preceq Q_k$. Given the definition of the Satake topology, this implies that
$M^K (Q_j,\FX)$ is contained in the closure of $M^K (Q_k,\FX)$. But for the pre-images
$e^K ( Q_j, G(\BA_f) )'$ and $e^K ( Q_k, G(\BA_f) )'$, the situation is quite different:
the closure of $e^K ( Q_k, G(\BA_f) )'$ equals $\overline{e^K ( Q_k, G(\BA_f) )}$, and
\[
e^K \bigl( Q_j, G(\BA_f) \bigr)' \cap \overline{e^K \bigl( Q_k, G(\BA_f) \bigr)}
= e^K \bigl( Q_j \cap Q_k, G(\BA_f) \bigr)' \; ,
\]
which equals $e^K ( Q_j, G(\BA_f) )'$ if and only if $Q_j = Q_k$. Furthermore,
if $K$ is neat, then the inclusion of $e^K ( Q_j, G(\BA_f) )' \cap \overline{e^K ( Q_k, G(\BA_f))}$
into its closure is contractible --- a property which is (very) false in general
for the inclusion of $M^K (Q_j,\FX)(\BC)$ into its closure in $M^K (G,\FX)^*(\BC)$!
\end{Rem}

\forget{
The following principle will turn out to be important.

\begin{Lem} \label{4Ma}
Let $Q$ and $Q'$ be $G(\BQ)$-conjugate parabolic subgroups of $G$, and
$R$ a third parabolic. Assume that both $R \cap Q$ and $R \cap Q'$ remain parabolic.
Let $g \in G (\BA_f)$.
Then the following are equivalent.
\begin{enumerate}
\item[(i)] $e^K(R \cap Q,g) = e^K(R \cap Q',g)$,
\item[(ii)] $\overline{e^K(R \cap Q,g)} = \overline{e^K(R \cap Q',g)}$,
\item[(iii)] the intersection $\overline{e^K(R \cap Q,g)} \cap \overline{e^K(R \cap Q',g)}$
is not empty,
\item[(iv)] there exists $\gamma \in R(\BQ) \cap gKg^{-1}$, such that $Q' = \gamma Q\gamma^{-1}$.
\end{enumerate}
\end{Lem}

\begin{Proof}
Clearly (iv)$\Rightarrow$(i)$\Rightarrow$(ii)$\Rightarrow$(iii).

It remains to show that (iii) implies (iv). Thus, suppose 
\[
\overline{e^K(R \cap Q,g)} \cap \overline{e^K(R \cap Q',g)} \ne \emptyset \; .
\]
This intersection then necessarily contains a whole stratum, which is of the form
$e^K(R',g)$, for some parabolic $R'$ contained in $R \cap Q'$. As 
$e^K(R',g) \subset \overline{e^K(R \cap Q,g)}$,
the parabolic $R'$ is also contained in $\gamma (R \cap Q) \gamma^{-1}$, for some
$\gamma \in G(\BQ) \cap gKg^{-1}$. 
This implies in particular that both $R \cap \gamma R \gamma^{-1}$ and
$Q' \cap \gamma Q \gamma^{-1}$ are parabolic. 
But $\gamma Q \gamma^{-1}$ is conjugate to $Q'$ by assumption. Now
in a connected algebraic group, two conjugate parabolic
subgroups, whose intersection remains parabolic, are identical \?  
It follows that $Q' = \gamma Q \gamma^{-1}$, while $\gamma R \gamma^{-1} = R$, \emph{i.e.},
$\gamma \in R(\BQ)$.
\end{Proof}

In order to be able to formulate our next main result, we need to fix some notation.

\begin{Def} 
Let $Q_j$ be an admissible parabolic subgroup of $G$, 
with canonical normal subgroup $P_j \,$, $g \in G(\BA_f)$,
and $K$ an open compact subgroup of $G (\BA_f)$. 
Let $Q$ be a parabolic subgroup of $G$ containing $P_j \,$, and such that $Q_j \cap Q$ remains
parabolic. \\[0.1cm]
(a)~Define $( G(\BQ) \cap gKg^{-1} )_{Q_j,Q}'$ as the subset
of $G(\BQ) \cap gKg^{-1}$ consisting of 
those $\gamma$, for which $Q_j \cap \gamma Q \gamma^{-1}$ 
remains parabolic and contains $P_j \,$. (We thus have 
$Q_j \cap \gamma Q \gamma^{-1} \in \adm^{-1}(Q_j)$.) \\[0.1cm]
(b)~For $\gamma \in ( G(\BQ) \cap gKg^{-1} )_{Q_j,Q}'$, define
\[
e^K(Q_j \cap \gamma Q \gamma^{-1},g)' 
         := \bigcup_{P_j \subset R \subset Q_j \cap \gamma Q \gamma^{-1}} e^K(R,g) 
         \subset M^K  (G,\FX) (\BC)^{BS} \; .
\]
\end{Def}

Note that $( G(\BQ) \cap gKg^{-1} )_{Q_j,Q}'$ 
contains the neutral element of $G(\BQ)$ and is stable under both multiplication
by $Q(\BQ) \cap gKg^{-1}$ from the right, and multiplication by $Q_j(\BQ) \cap gKg^{-1}$ 
from the left. 

\begin{Ex} \label{4Maa}
If $Q \in  \adm^{-1}(Q_j)$, then
\[
( G(\BQ) \cap gKg^{-1} )_{Q_j,Q}' = Q_j(\BQ) \cap gKg^{-1} \; .
\]
Indeed, if $Q_j \cap \gamma Q \gamma^{-1}$ is parabolic, then so is $Q_j \cap \gamma Q_j \gamma^{-1}$
(as $Q$ is contained in $Q_j$ by assumption), and hence $\gamma Q_j \gamma^{-1} = Q_j$
\? .

In particular, we have
\[
( G(\BQ) \cap gKg^{-1} )_{Q_j,Q_j}' = Q_j(\BQ) \cap gKg^{-1} \; .
\]
\end{Ex}

\begin{Lem} \label{4Mb}
Let $Q_j$ be an admissible parabolic subgroup of $G$, $g \in G(\BA_f)$,
and $K$ an open compact subgroup of $G (\BA_f)$. 
Let $Q$ be a parabolic subgroup of $G$ containing $P_j \,$, and such that $Q_j \cap Q$ remains
parabolic. \\[0.1cm]
(a)~Let $\gamma_1, \gamma_2 \in (G(\BQ) \cap gKg^{-1})_{Q_j,Q}'$. The following are equivalent. 
\begin{enumerate}
\item[(i)] $e^K(Q_j \cap \gamma_1 Q \gamma_1^{-1},g) = e^K(Q_j \cap \gamma_2 Q \gamma_2^{-1},g)$,
\item[(ii)] the intersection 
$e^K(Q_j \cap \gamma_1 Q \gamma_1^{-1},g)' \cap e^K(Q_j \cap \gamma_2 Q \gamma_2^{-1},g)'$
is not empty, 
\item[(iii)] $\overline{e^K(Q_j \cap \gamma_1 Q \gamma_1^{-1},g)} 
= \overline{e^K(Q_j \cap \gamma_2 Q \gamma_2^{-1},g)}$,
\item[(iv)] the classes of $\gamma_1$ and $\gamma_2$ in the double quotient
\[
\bigl( Q_j(\BQ) \cap gKg^{-1} \bigr) \backslash
\bigl( G(\BQ) \cap gKg^{-1} \bigr)_{Q_j,Q}'\, / \bigl( Q(\BQ) \cap gKg^{-1} \bigr) 
\]
are the same. 
\end{enumerate} 
(b)~The double quotient
\[
\bigl( Q_j(\BQ) \cap gKg^{-1} \bigr) \backslash
\bigl( G(\BQ) \cap gKg^{-1} \bigr)_{Q_j,Q}' \, / \bigl( Q(\BQ) \cap gKg^{-1} \bigr) 
\]
is finite. \\[0.1cm]
(c)~Let $\gamma \in (G(\BQ) \cap gKg^{-1})_{Q_j,Q}'$. Then
$e^K(Q_j \cap \gamma Q \gamma^{-1},g)'$ equals 
the image of 
\[
\bigl( p^{-1} ( \coprod \FX_j / W_j ) \cap \overline{e(\gamma Q \gamma^{-1})}  \bigr)  \times \{gK\}
\]
under the projection 
\[
\FX^* \times G (\BA_f) / K \longonto M^K (G,\FX)^* (\BC) \; . 
\]
\end{Lem}

\begin{Proof}
(a): clearly (iv)$\Rightarrow$(i)$\Rightarrow$(ii).

The intersection 
$e^K(Q_j \cap \gamma_1 Q \gamma_1^{-1},g)' \cap e^K(Q_j \cap \gamma_2 Q \gamma_2^{-1},g)'$
is contained in 
$\overline{e^K(Q_j \cap \gamma_1 Q \gamma_1^{-1},g)} 
\cap \overline{e^K(Q_j \cap \gamma_2 Q \gamma_2^{-1},g)}$. The implications
(ii)$\Rightarrow$(iii)$\Rightarrow$(iv) therefore follow from Lemma~\ref{4Ma}. 

(b): let $\Omega$ be a set of representatives of 
\[
\bigl( Q_j(\BQ) \cap gKg^{-1} \bigr) \backslash
\bigl( G(\BQ) \cap gKg^{-1} \bigr)_{Q_j,Q}' \, / \bigl( Q(\BQ) \cap gKg^{-1} \bigr) \; .
\] 
According to (a), the union 
$\cup_{\gamma \in \Omega} e^K(Q_j \cap \gamma Q \gamma^{-1},g)$ of strata
is disjoint in $M^K  (G,\FX) (\BC)^{BS}$. But according to \cite[Prop.~9.4~(i)]{BS},
the stratification of $M^K  (G,\FX) (\BC)^{BS}$ is finite.

(c): this follows from Complement~\ref{4Lb}.
\end{Proof}

Putting everything together, we obtain the third main result of this section.

\begin{Thm} \label{4N}
Fix an admissible parabolic subgroup $Q_j$ of $G$, 
with canonical normal subgroup $P_j \,$, and consider the disjoint union 
$\coprod \FX_j$ of the finitely many spaces $\FX_j$ underlying
rational boundary components associated to $Q_j \,$, as well as its quotient
$\coprod \FX_j / W_j$ by the action of the unipotent radical $W_j \,$.
Let $K$ be an open compact subgroup of $G (\BA_f)$, 
$g \in G(\BA_f)$, and define $pr_{gK} (\coprod \FX_j/W_j) \subset M^K (G,\FX)^* (\BC)$ as 
the image of $\coprod \FX_j/W_j \times \{gK\}$ under the projection 
\[
\FX^* \times G (\BA_f) / K \longonto M^K (G,\FX)^* (\BC) \; . 
\]
(a)~Under the map $p^K: M^K  (G,\FX) (\BC)^{BS} \longto M^K (G,\FX)^* (\BC)$, we have
\[
(p^K)^{-1} \bigl( pr_{gK} (\coprod \FX_j/W_j) \bigr) = 
         \bigcup_{P_j \subset R \subset Q_j} e^K(R,g) \subset M^K  (G,\FX) (\BC)^{BS}
\; ,
\]
where the union runs over the 
parabolics $R$ of $G$ contained in $Q_j$ and containing $P_j \,$. 
In particular, $(p^K)^{-1} ( pr_{gK} (\coprod \FX_j/W_j) )$ 
contains $e^K(Q_j,g)$, and is open in its closure $\overline{e^K(Q_j,g)}$. \\[0.1cm]
(b)~Let $Q$ be a parabolic of $G$ containing $P_j \,$, and such that $Q_j \cap Q$ remains
parabolic. Choose a set $\Omega$ of representatives of 
\[
\bigl( Q_j(\BQ) \cap gKg^{-1} \bigr) \backslash
\bigl( G(\BQ) \cap gKg^{-1} \bigr)_{Q_j,Q}' \, / \bigl( Q(\BQ) \cap gKg^{-1} \bigr) \; .
\] 
Then 
the intersection $(p^K)^{-1} ( pr_{gK} (\coprod \FX_j/W_j) ) \cap \overline{e^K(Q,g)}$
contains the disjoint union $\coprod_{\gamma \in \Omega} e^K(Q_j \cap \gamma Q \gamma^{-1},g)$,  
and is open in the disjoint union of closures
$\coprod_{\gamma \in \Omega} \overline{e^K(Q_j \cap \gamma Q \gamma^{-1},g)}$. 
In particular, the closure of 
the intersection $(p^K)^{-1} ( pr_{gK} (\coprod \FX_j/W_j) ) \cap \overline{e^K(Q,g)}$
equals
$\coprod_{\gamma \in \Omega} \overline{e^K(Q_j \cap \gamma Q \gamma^{-1},g)}$. \\[0.1cm]
(c)~Assume $K$ to be neat.
Let $Q$ be a parabolic of $G$ containing $P_j \,$, and such that $Q_j \cap Q$ remains
parabolic. Choose a set $\Omega$ of representatives of 
\[
\bigl( Q_j(\BQ) \cap gKg^{-1} \bigr) \backslash
\bigl( G(\BQ) \cap gKg^{-1} \bigr)_{Q_j,Q}' \, / \bigl( Q(\BQ) \cap gKg^{-1} \bigr) \; .
\] 
Then the immersion
\[
(p^K)^{-1} \bigl( pr_{gK} (\coprod \FX_j/W_j) \bigr) \cap \overline{e^K(Q,g)} 
\longinto \coprod_{\gamma \in \Omega} \overline{e^K(Q_j \cap \gamma Q \gamma^{-1},g)}
\]
is contractible. In particular, the immersion 
\[
(p^K)^{-1} \bigl( pr_{gK} (\coprod \FX_j/W_j) \bigr) \longinto \overline{e^K(Q_j,g)} 
\]
is contractible. 
\end{Thm}

\begin{Proof}
According to Complement~\ref{4Lb}, 
\[
(p^K)^{-1} \bigl( pr_{gK} (\coprod \FX_j/W_j) \bigr) \subset M^K (G,\FX) (\BC)^{BS}
\]
equals the image of $\coprod_{P_j \subset R \subset Q_j} e(R) \times \{ gK \}$ under
the projection from $\FX^{BS} \times G (\BA_f) / K$.
This implies (a).

As for (b), note first that as the pre-image of a locally closed subset of 
$M^K (G,\FX)^* (\BC)$, 
\[
(p^K)^{-1} \bigl( pr_{gK} (\coprod \FX_j/W_j) \bigr) \subset M^K (G,\FX) (\BC)^{BS}
\]
is locally closed. Therefore, the intersection
$(p^K)^{-1} ( pr_{gK} (\coprod \FX_j/W_j) ) \cap \overline{e^K(Q,g)}$
is open in its closure. Next, we have
\[
\overline{e^K(Q,g)} = \bigcup_{R \subset Q} e^K(R,g) \; ,
\]
hence according to (a),
\[
(p^K)^{-1} \bigl( pr_{gK} (\coprod \FX_j/W_j) \bigr) \cap \overline{e^K(Q,g)} 
= \bigcup_{\gamma \in (G(\BQ) \cap gKg^{-1})} 
                       \bigcup_{P_j \subset R \subset Q_j \cap \gamma Q \gamma^{-1}} e^K(R,g) \; , 
\]
which equals 
\[
\bigcup_{\gamma \in ( G(\BQ) \cap gKg^{-1} )_{Q_j,Q}'}e^K(Q_j \cap \gamma Q \gamma^{-1},g)' \; . 
\]
By Lemma~\ref{4Mb}~(a),
\[
(p^K)^{-1} \bigl( pr_{gK} (\coprod \FX_j/W_j) \bigr) \cap \overline{e^K(Q,g)} 
= \coprod_{\gamma \in \Omega} e^K(Q_j \cap \gamma Q \gamma^{-1},g)' \; ,
\]
and the union
\[
\bigcup_{\gamma \in \Omega} \overline{e^K(Q_j \cap \gamma Q \gamma^{-1},g)} 
\]
remains disjoint.

To show (c), note that the restriction of the projection from $\FX^{BS} \times P (\BA_f) / K$
to each of the $\overline{e(Q_j \cap \gamma Q \gamma^{-1})} \times \{ gK \}$ 
equals the quotient by the action
of the stabilizer in $G(\BQ)$ of $e(Q_j \cap \gamma Q \gamma^{-1}) \times \{ gK \}$ 
\cite[Prop.~9.4~(ii)]{BS}.
But this action is free on the whole of $\FX^{BS} \times G (\BA_f) / K$ 
as the group $K$ is supposed neat \cite[Sect.~9.5]{BS}.
Our claim therefore follows from Lemma~\ref{4Mb}~(c) and Complement~\ref{4Lb}.
\end{Proof}
}


\bigskip
%
%

\section{The fibres of the map $p$}
\label{5}



The first main result of this section (Theorem~\ref{5S}) gives a description of the
fibres of the map $p: \FX^{BS} \to \FX^*$ constructed in Section~\ref{4},
by identifying them with manifolds with corners $\FZ^{BS}$ associated to 
spaces $\FZ$ of type $S - \BQ$ associated to certain algebraic subgroups of the group
$P$ underlying our Shimura data. Theorem~\ref{5S} will then be used to describe
the fibres of the map $p^K$ between the Borel--Serre and the Baily--Borel 
compactification of the Shimura variety (Corollary~\ref{5Y}~(b)).
Combining this with the results from the previous section, 
we get the second main result (Theorem~\ref{5Z}),
which yields a precise analysis of the intersection of the fibres of $p^K$ with the 
(closures of the) canonical strata of the Borel--Serre compactification. 
For the sake of completeness, we also give the corresponding results concerning the
\emph{reductive Borel--Serre compactification}; here, we recover in particular 
the part of \cite[Thm.~5.8]{G} concerning the fibres of the map from
that compactification to the Baily--Borel compactification (Corollary~\ref{5Y}~(c)). \\

We fix mixed Shimura data $(P,\FX)$, and an admissible 
parabolic subgroup $Q_j$ of $P$,
with canonical normal subgroup $P_j \,$. Consider the complex analytic map
\[
\iota_j: \FX \longto \coprod \FX_j
\]
of $\FX$ to the disjoint union of the finitely many spaces $\FX_j$
underlying rational boundary components associated to $Q_j$,
and its composition
\[
\tilde{\pi}_j : \FX \longto \coprod \FX_j / W_j \; .
\]
with the canonical epimorphisms $\pi_j: \FX_j \to \FX_j/W_j$.
The map $\iota_j$ is equivariant with respect to $Q_j(\BR)U(\BC)$
\cite[4.11]{P1}, hence so is
$\tilde{\pi}_j$. 
Recall (Definition~\ref{4G}) the closed
connected normal subgroup
\[ 
C_j = \{ q \in Q_j \; , \; \pi_{Q_j}(q) \in \Cent_{\bar{Q}_j}(\pi_{Q_j}(P_j)) \}^0 
\]
of $Q_j$. By its very definition, the group $C_j(\BR)U(\BC)$ acts trivially
on all $h(\FX_j/W_j)$. Therefore, according to
\cite[Cor.~2.12]{P1}, the neutral connected component $C_j(\BR)^0U(\BC)$
acts trivially on the target of $\tilde{\pi}_j$. In other words, it respects
all fibres of $\tilde{\pi}_j$.

\begin{Prop} \label{5K}
The induced action of $C_j(\BR)^0U(\BC)$ on each of the fibres of $\tilde{\pi}_j$
is transitive. In particular, the fibres of $\tilde{\pi}_j$ are connected.
\end{Prop}

\begin{Proof}
Let $(P_j,\FX_j)$ be one of the rational boundary components
associated to $Q_j$, and $z_0 \in \FX_j / W_j$. For $x$, $y \in \tilde{\pi}_j^{-1}(z_0)$,
write $x_j := \iota_j(x)$ and $y_j := \iota_j(y)$. 
These two elements of $\FX_j$
being mapped to $z_0$ under $\pi_j$, there are $w \in W_j(\BR)$
and $u \in U_j(\BC)$ such that $y_j = uw x_j$.
We may assume that $u$ is purely imaginary. 
In additive notation, the element $u$ thus equals $\im(y_j) - \im(x_j)$.

The relation $y_j = uw x_j$ shows that $x_j$ and $y_j$ are in the same connected
component of $\FX_j \,$. 
According to \cite[Sect.~4.11]{P1}, the points $x$ and $y$ are
in the same connected component, say
$\FX^0$, of $\FX$. 
The action of $Q_j(\BR)^0U(\BC)$ on the cone
$C(\FX^0,P_j)$ is transitive \cite[Prop.~4.15]{P1}. 
We have $Q_j = C_j P_j^{der}$, and $P_j^{der}(\BR)$
acts trivially on $U_j$ \cite[Prop.~2.14~(a)]{P1}
(note that this result can be applied since by \cite[Rem.~(ii) on p.~82,
after Sect.~4.11]{P1}, the Shimura data $(P_j,\FX_j)$ are irreducible). 
Therefore, the induced action of $C_j(\BR)^0U(\BC)$ on $C(\FX^0,P_j)$ is transitive.
We may thus suppose that $x_j$ and $y_j$ have the same imaginary part.

It follows that $u$ is trivial, \emph{i.e.}, that $y_j = w x_j$.
But the map $\iota_j$ is
$Q_j(\BR)$-equivariant and injective, 
and the element $w$ belongs to $W_j(\BR) \subset C_j(\BR)^0$.
\end{Proof}

\begin{Prop} \label{5L}
Assume $(+)$ and $(U=0)$. 
Then each of the fibres of $\tilde{\pi}_j$ is stable under
the restriction to $C_j(\BR)$ of the
extended action of $P(\BR)$ on $\FX$ from Corollary~\ref{1G}. 
\end{Prop}

\begin{Proof}
If $K$ is maximal compact in $P(\BR)$, then
$C_j(\BR) \cap K$ is maximal compact in $C_j(\BR)$ 
as $C_j$ is normal in $Q_j$, and $Q_j$ is parabolic (cmp.~\cite[Sect.~1.4]{BS}).
Now apply Corollary~\ref{1Ca}.
\end{Proof}

\begin{Cor} \label{5M}
Assume $(+)$ and $(U=0)$. 
Then each of the fibres of $\tilde{\pi}_j$
is a homogeneous space under the restriction to $C_j(\BR)$ of the
extended action of $P(\BR)$ on $\FX$ from Corollary~\ref{1G}.
\end{Cor}

\begin{Proof}
Apply Propositions~\ref{5K} and \ref{5L}, noting that 
the restrictions to $C_j(\BR)^0 \subset P(\BR)^0$ of the extended action~(2)
and the action~(1) underlying the Shimura data coincide.
\end{Proof}

\begin{Rem} \label{5Ma}
(a)~An alternative proof of Corollary~\ref{5M} could be given by first verifying that
$(\tilde{\pi}_j^{-1}(z_0),(C_j(\BR) \cap K_x)_{x \in \tilde{\pi}_j^{-1}(z_0)})$
is homogeneous of type $S'$ under $(C_j(\BR),C_j(\BR) \cap {}^0 \! P(\BR))$,
and then using Proposition~\ref{1C} in order to get an extended action 
of $C_j(\BR)$. Unicity
guarantees that the result is indeed equal to the
restriction of the extended action of $P(\BR)$. \\[0.1cm]
(b)~If $z_0 \in \FX_j/W_j \,$, then
the restrictions of~(1) the action of $P(\BR)$ underlying the Shimura data and
(2)~the extended action of $P(\BR)$, coincide on the subgroup
\[
\Cent_{C_j(\BR)} \bigl( (\FX_j/W_j)^0 \bigr) 
\]
of finite index of $C_j(\BR)$, consisting of the elements acting trivially on
(equi\-valently \cite[Cor.~2.12]{P1}, stabilizing)
the connected component $(\FX_j/W_j)^0$ of $\FX_j/W_j$ containing $z_0$. \\[0.1cm] 
(c)~It is straightforward to show that the pair 
\[
\bigl( \tilde{\pi}_j^{-1}(z_0) , 
             \bigl( C_{j,\BR} \cap L_{x,Q_j} \bigr)_{x \in \tilde{\pi}_j^{-1}(z_0)} \bigr)
\]
is a space of type $S$ under $C_j$ in the sense of \cite[Def.~2.3]{BS}.
Note: the
connected normal solvable subgroup of $C_{j,\BR}$ 
from \cite[Def.~2.3, SI]{BS} needs to be defined as 
the base change to $\BR$ of the subgroup of
the $\BQ$-split radical of $C_j$ 
generated by the unipotent radical of $W_j$ (hence of $C_j)$
and the $\BQ$-split radical of $P$ (cmp.~\cite[Lemma~2.6]{BS}). \\[0.1cm]
(d)~Unless $Q_j$ equals $P$, 
\[
\bigl( \tilde{\pi}_j^{-1}(z_0) , 
             \bigl( C_{j,\BR} \cap L_{x,Q_j} \bigr)_{x \in \tilde{\pi}_j^{-1}(z_0)} \bigr)
\]
is \emph{not} of type $S-\BQ$ under $C_j \,$, the reason being that the stabilizers
\[
\Stab_{C_j(\BR)}^{\ext}(x) = Z_x(\BR)^0 \times \bigl( C_j(\BR) \cap K_x \bigr)
= Z_x(\BR) \cdot (C_j(\BR) \cap K_x) 
\]
are too small. Indeed, the image of $\Stab_{C_j(\BR)}^{\ext}(x)^0$
under $\pi_{Q_j}$
has trivial intersection with 
$\pi_{Q_j} \circ h_{\re_j(x)} \circ w (\GRm) \subset Z_{d,\BQ}(\bar{Q}_j)_\BR$.
\end{Rem}

Recall (Corollary~\ref{4I}) that the geodesic action 
of $A_{Q_j}$ respects each of the fibres of 
$\tilde{\pi}_j$. Furthermore \cite[Prop.~3.4]{BS}, it commutes with the action of 
$C_j(\BR)$.

\begin{Prop} \label{5P}
Assume $(+)$ and $(U=0)$. Let $z_0 \in \coprod \FX_j / W_j \,$. 
Let $\tilde{\pi}_j^{-1}(z_0)$ be equipped with the action of $C_j(\BR)$ from
Proposition~\ref{5L}. 
Then the pair
\[
\bigl( A_{Q_j} \! \backslash \tilde{\pi}_j^{-1}(z_0) , 
     \bigl( C_{j,\BR} \cap L_{x,Q_j} 
                   \bigr)_{x \in A_{Q_j} \! \backslash \tilde{\pi}_j^{-1}(z_0)} \bigr)
\]
is a space of type $S-\BQ$ under $C_j \,$.
\end{Prop}

\begin{Proof}
Denote by $\FX^0$ the connected component of $\FX$ containing $\tilde{\pi}_j^{-1}(z_0)$
(the latter \emph{is} connected thanks to Proposition~\ref{5K}). 
According to \cite[Sect.~3.9]{BS}, 
\[
\bigl( A_{Q_j} \! \backslash \FX^0 , 
            \bigl( L_{x,Q_j} \bigr)_{x \in A_{Q_j} \! \backslash \FX^0} \bigr)
\]
is a space of type $S-\BQ$ under $Q_j \,$.
Given the definition of $C_j \,$, the neutral connected component $Z(L_{x,Q_j})^0$ 
of the center $Z(L_{x,Q_j})$ is contained in $C_{j,\BR}$, for every $x \in \FX$,
and hence, for every $x \in \tilde{\pi}_j^{-1}(z_0)$.

Our claim then follows from the following abstract principle concerning
spaces of type $S-\BQ$: 
let $(X,(L_x)_{x \in X})$ be of type $S-\BQ$ under $G$, 
and $N \subset G$ a normal subgroup. Both $G$ and $N$ are assumed to be connected. 
Let $Y \subset X$ be an orbit under $N(\BR)$,
such that for every $x \in Y$, we have $Z(L_x)^0 \subset N_\BR$. Then
$(Y,(N_\BR \cap L_x)_{x \in Y})$ is of type $S-\BQ$ under $N$. 
We leave the proof of this principle to the reader (but see \cite[Sect.~2.7~(1)]{BS}
for the case $Y = X$).
\end{Proof}

\begin{Prop} \label{5Pa}
Assume $(+)$ and $(U=0)$. Let $z_0 \in \coprod \FX_j / W_j$ be a point,
and $(\FX_j / W_j)^0 \subset \coprod \FX_j / W_j$ a connected component of
$\coprod \FX_j / W_j \,$. Then there is a canonical choice of
$z_1 \in (\FX_j / W_j)^0$, such that the pairs 
\[
\bigl( A_{Q_j} \! \backslash \tilde{\pi}_j^{-1}(z_0) , 
     \bigl( C_{j,\BR} \cap L_{x,Q_j} 
                   \bigr)_{x \in A_{Q_j} \! \backslash \tilde{\pi}_j^{-1}(z_0)} \bigr)
\]
and
\[
\bigl( A_{Q_j} \! \backslash \tilde{\pi}_j^{-1}(z_1) , 
     \bigl( C_{j,\BR} \cap L_{x,Q_j} 
                   \bigr)_{x \in A_{Q_j} \! \backslash \tilde{\pi}_j^{-1}(z_1)} \bigr)
\]
are canonically isomorphic as spaces of type $S-\BQ\,$: 
there is a canonical real analytic isomorphism
\[
\psi: A_{Q_j} \! \backslash \tilde{\pi}_j^{-1}(z_0) 
\isoto \; A_{Q_j} \! \backslash \tilde{\pi}_j^{-1}(z_1) \; ,
\]
that is equivariant for the action of $C_j(\BR)$ from
Proposition~\ref{5L}, and such that 
\[
C_{j,\BR} \cap L_{\psi(x),Q_j} = C_{j,\BR} \cap L_{x,Q_j} \; , \; 
\forall \; x \in A_{Q_j} \! \backslash \tilde{\pi}_j^{-1}(z_0) \; .
\]
\end{Prop}

Note that the existence of \emph{some} isomorphism of spaces of type $S-\BQ$ is guaranteed
by Proposition~\ref{1Mgenprop}~(a).

\medskip

\begin{Proofof}{Proposition~\ref{5Pa}}
The continuous map $\tilde{\pi}_j$ is surjective
with connected fibres (Proposition~\ref{5K}); 
denote by $\FX^{0,i}$ the connected component of $\FX$
contai\-ning $\tilde{\pi}_j^{-1} (z_0)$, 
and by $\FX^{0,k}$ the connected component of $\FX$
containing $\tilde{\pi}_j^{-1} ((\FX_j / W_j)^0)$. 
According to Corollary~\ref{1P},
there is a canonical isomorphism 
\[
\psi: \bigl( \FX^{0,i},(L_x)_{x \in \FX^{0,i}} \bigr) 
\isoto \bigl( \FX^{0,k},(L_x)_{x \in \FX^{0,k}} \bigr) 
\]
of spaces of type $S-\BQ$ under $P$. It induces an isomorphism
\[
\bigl( A_{Q_j} \! \backslash \FX^{0,i} , 
            \bigl( L_{x,Q_j} \bigr)_{x \in A_{Q_j} \! \backslash \FX^{0,i}} \bigr)
\isoto \bigl( A_{Q_j} \! \backslash \FX^{0,k} , 
            \bigl( L_{x,Q_j} \bigr)_{x \in A_{Q_j} \! \backslash \FX^{0,k}} \bigr)            
\]
of spaces of type $S-\BQ$ under $Q_j$ \cite[Sect.~3.9]{BS}.

It remains to define $z_1$ as the image of $\psi(x_0)$ under $\tilde{\pi}_j$, 
for any $x_0 \in \tilde{\pi}_j^{-1} (z_0)$.
\end{Proofof}

Given Proposition~\ref{5P}, the manifolds with corners 
$(A_{Q_j} \! \backslash \tilde{\pi}_j^{-1}(z_0))^{BS}$, 
for\-med with respect to the action of $C_j(\BR)$ from
Proposition~\ref{5L}, are defined \cite[Sect.~7.1]{BS}
(in the notation of \loccit, they would be denoted 
$\overline{A_{Q_j} \! \backslash \tilde{\pi}_j^{-1}(z_0)}$), for every 
$z_0 \in \coprod \FX_j / W_j \,$. The space 
$(A_{Q_j} \! \backslash \tilde{\pi}_j^{-1}(z_0))^{BS}$
carries an action of $\Stab_{Q_j(\BQ)}(z_0)$, induced by the $\Stab_{Q_j(\BQ)}(z_0)$-action
on $A_{Q_j} \! \backslash \tilde{\pi}_j^{-1}(z_0) \subset A_{Q_j} \! \backslash \FX^0$.
Its restriction to $\Stab_{C_j(\BQ)}(z_0) \subset \Stab_{Q_j(\BQ)}(z_0)$ coincides
with the restriction to $\Stab_{C_j(\BQ)}(z_0) \subset C_j(\BQ)$ of the
$C_j(\BQ)$-action on $(A_{Q_j} \! \backslash \tilde{\pi}_j^{-1}(z_0))^{BS}$ from
\cite[Prop.~7.6]{BS}. \\

Recall from Corollary~\ref{4H} that the map $Q \mapsto C_j \cap Q$ is
a bijection between $\adm^{-1}(Q_j)$ and the set of parabolics of $C_j$.
In order to relate the manifolds with corners
$(A_{Q_j} \! \backslash \tilde{\pi}_j^{-1}(z_0))^{BS}$, $z_0 \in \coprod \FX_j / W_j \,$, 
to $\FX^{BS}$, we need to compare geodesic actions,
formed with respect to parabolics of $P$, and with respect to parabolics of $C_j \,$.

\begin{Prop} \label{5Q}
Let $Q$ be a parabolic contained in $Q_j$ and containing
$P_j$ (in other words, we have $Q \in \adm^{-1}(Q_j)$). 
Write $\overline{C_j \cap Q}$ for the maximal reductive quotient of $C_j \cap Q$.
Then the inclusion $C_j \cap Q \into Q$ 
induces an equality of unipotent radicals
\[
\Rad^u(C_j \cap Q) = \Rad^u(Q)
\]
and an isomorphism of centers
\[
Z \bigl( \overline{C_j \cap Q} \bigr)^0 \isoto Z \bigl( \bar{Q} \bigr)^0 \; .
\]
\end{Prop}

\begin{Proof}
By definition,
\[
C_j \cap Q = \{ q \in Q \; , \; \pi_{Q_j}(q) \in \Cent_{\bar{Q}_j}(\pi_{Q_j}(P_j)) \}^0 \; ,
\]
which equals 
\[
\{ q \in Q \; , \; qgq^{-1}g^{-1} \in W_j \, , \, \forall \, g \in P_j \}^0 
\]
(since the unipotent radical of $Q_j$ equals $W_j$
\cite[proof of Lemma~4.8]{P1}).

The parabolic $Q$ being contained in $Q_j$, its unipotent radical
$\Rad^u(Q)$ contains $W_j$. Let $q \in Q$ such that 
\[
qgq^{-1}g^{-1} \in \Rad^u(Q) \, , \, \forall \, g \in P_j \; . 
\]
The subgroup $P_j$ being normal in $Q$, the commutators $qgq^{-1}g^{-1}$ all lie in 
$P_j \,$. 
Hence, if in addition they belong to $\Rad^u(Q)$, 
then they are in $\Rad^u(Q) \cap P_j = W_j \,$. Therefore,
\[
C_j \cap Q = \{ q \in Q \; , \; \pi_Q(q) \in \Cent_{\bar{Q}}(\pi_Q(P_j)) \}^0 \; . 
\]
In particular, this equation shows that $\Rad^u(Q)$ is contained in $C_j \cap Q$.
It also shows that on the one hand,
$Z(\bar{Q})^0$ is contained in $\overline{C_j \cap Q}$, hence
in $Z(\overline{C_j \cap Q})$. On the other hand, we have 
$\bar{Q} = (\overline{C_j \cap Q}) \pi_Q(P_j)$ (Corollary~\ref{4H}), and hence 
$Z(\overline{C_j \cap Q}) \subset Z(\bar{Q})$.
\end{Proof}

It follows that the inclusion 
``$Z(L_{x,Q})^0 \subset C_{j,\BR}$'' used in the proof of Proposition~\ref{5P}
is valid for more general choices of $Q$.

\begin{Cor} \label{5R}
Assume $(+)$ and $(U=0)$. 
Let $Q \in \adm^{-1}(Q_j)$. \\[0.1cm]
(a)~The inclusion $C_j \cap Q \into Q$ 
induces an isomorphism
\[
A_{C_j \cap Q} \isoto A_{Q_j} \! \backslash A_Q \; .
\]
(b)~For all $x \in \FX$, we have
\[
Z \bigl( C_{j,\BR} \cap L_{x,Q} \bigr)^0 = Z \bigl( L_{x,Q} \bigr)^0 \; . 
\]
In particular, $Z(L_{x,Q})^0 \subset C_{j,\BR}$.
\end{Cor}

\begin{Proof}
Claim (a) follows from Proposition~\ref{5Q}, Definition~\ref{2Ga},
and the fact that the neutral connected component of $Z(G)$ is contained in 
$\overline{C_j \cap Q}$, hence in $Z(\overline{C_j \cap Q})$. 

As for claim~(b), the intersection $C_{j,\BR} \cap L_{x,Q}$ is a Levi subgroup
of $(C_j \cap Q)_\BR$ as $C_j \cap Q$ is a normal subgroup of $Q$, 
that according to the first statement
of Proposition~\ref{5Q} contains $\Rad^u(Q)$.
Now use the second statement of Proposition~\ref{5Q}.
\end{Proof}

Therefore, ``the geodesic action lies in $C_j$'' as long as it is formed with
respect to parabolics $Q$ in $\adm^{-1}(Q_j)$. 
In the sequel, we shall identify $A_{C_j \cap Q}$ and $A_{Q_j} \! \backslash A_Q$, for every
such $Q$. \\

Consider the inclusion $\tilde{\pi}_j^{-1}(z_0) \longinto \FX$
of one of the fibres of $\tilde{\pi}_j$. It induces an inclusion
\[
A_{Q_j} \! \backslash \tilde{\pi}_j^{-1}(z_0) \longinto A_{Q_j} \! \backslash \FX = e(Q_j)
\]
(which was exploited in the proof of Proposition~\ref{5P}).
The space $e(Q_j)$ is one of the faces of $\FX^{BS}$.
Here is the first main result
of this section.

\begin{Thm} \label{5S}
Assume hypotheses $(+)$ and $(U=0)$. Let $z_0 \in \FX_j / W_j \,$. \\[0.1cm]
(a)~The inclusion 
\[
A_{Q_j} \! \backslash \tilde{\pi}_j^{-1}(z_0) \longinto e(Q_j) \longinto \FX^{BS}
\]
extends uniquely to a continuous map
\[
\kappa_{z_0,j}: 
\bigl( A_{Q_j} \! \backslash \tilde{\pi}_j^{-1}(z_0) \bigr)^{BS} \longto \FX^{BS} \; .
\]
The map $\kappa_{z_0,j}$ is a morphism of manifolds with corners. \\[0.1cm]
(b)~The morphism $\kappa_{z_0,j}$ is injective. 
It is a morphism of stratifications; more precisely,
for any parabolic $Q$ of $P$, we have 
\[
\kappa_{z_0,j}^{-1} (e(Q)) = \emptyset 
\]
if $Q \not \in \adm^{-1}(Q_j)$, and
\[
\kappa_{z_0,j}^{-1} (e(Q)) = e(C_j \cap Q) 
\subset \bigl( A_{Q_j} \! \backslash \tilde{\pi}_j^{-1}(z_0) \bigr)^{BS}
\]
if $Q \in \adm^{-1}(Q_j)$. In particular,
\[
\kappa_{z_0,j}^{-1} (e(Q_j)) = e(C_j) = A_{Q_j} \! \backslash \tilde{\pi}_j^{-1}(z_0) \; .
\]
(c)~The morphism $\kappa_{z_0,j}$ is $C_j(\BQ)W(\BR)$-equivariant 
with respect to the action from Proposition~\ref{5L}.
It is therefore
$\Cent_{C_j(\BQ)W(\BR)} ((\FX_j/W_j)^0)$-equivariant 
with respect to the action underlying the Shimura data, $(\FX_j/W_j)^0$
denoting the connected component of $\FX_j/W_j$ containing $z_0$. \\[0.1cm]
(d)~If the Shimura data $(P,\FX) = (G,\FX)$ are pure, then 
$\kappa_{z_0,j}$ yields an $\BR$-analytic
identification of $( A_{Q_j} \! \backslash \tilde{\pi}_j^{-1}(z_0))^{BS}$
with the fibre of $p : \FX^{BS} \to \FX^*$ over 
$z_0 \in \FX_j / W_j \subset \FX^*$. The map
$\kappa_{z_0,j}$ is the unique continuous extension
of the identification of $A_{Q_j} \! \backslash \tilde{\pi}_j^{-1}(z_0)$
with the subspace $p^{-1}(z_0) \cap e(Q_j)$ of $\FX^{BS}$.
\end{Thm}

\begin{Proof}
(a), (b): as for the unicity statement in (a), we argue as usual 
($A_{Q_j} \! \backslash \tilde{\pi}_j^{-1}(z_0)$ is dense in 
$( A_{Q_j} \! \backslash \tilde{\pi}_j^{-1}(z_0))^{BS}$, and $\FX^{BS}$
is Hausdorff \cite[Thm.~7.8]{BS}).
By Corollary~\ref{4H}, the map $Q \mapsto C_j \cap Q$ is
a bijection between $\adm^{-1}(Q_j)$ and the set of parabolics of $C_j$.
By Corollary~\ref{5R}, this bijection is compatible with the geodesic action.
The remaining claims therefore follow from the definition of the spaces 
$( A_{Q_j} \! \backslash \tilde{\pi}_j^{-1}(z_0))^{BS}$ and
$\FX^{BS}$ \cite[Sect.~7.1]{BS}.

\noindent (c): the inclusion
\[
A_{Q_j} \! \backslash \tilde{\pi}_j^{-1}(z_0) \longinto e(Q_j) \longinto \FX^{BS}
\]
is $C_j(\BQ)W(\BR)$-equivariant. Therefore, the first part of the claim follows from unicity of
$\kappa_{z_0,j}$ (part~(a)). As for the second, apply Remark~\ref{5Ma}~(b).

\noindent (d): recall (Complement~\ref{4Lb}) that the pre-image $p^{-1}(\coprod \FX_j / W_j)$
equals the disjoint union of those faces $e(Q)$ for which $\adm(Q) = Q_j \,$.
By Construction~\ref{4K}, the restriction $p_Q$ of $p$ to each such face $e(Q)$ fits into a commutative
diagram
\[
\vcenter{\xymatrix@R-10pt{
        \FX \ar@{->>}[r] \ar@{=}[d] & e(Q) = A_Q \backslash \FX \ar[d]^{p_Q} \\
        \FX \ar[r]^-{\tilde{\pi}_j} & \coprod \FX_j / W_j 
\\}}
\]
But the base change to $z_0 \in \FX_j / W_j$ of this diagram is
\[
\vcenter{\xymatrix@R-10pt{
        \tilde{\pi}_j^{-1}(z_0) \ar@{->>}[r] \ar@{=}[d] & 
        A_Q \backslash \tilde{\pi}_j^{-1}(z_0) \ar[d]^{p_Q} \\
        \tilde{\pi}_j^{-1}(z_0) \ar[r]^-{\tilde{\pi}_j} & \{ z_0 \}
\\}}
\]
\end{Proof}

\begin{Cor} \label{5Sa}
Assume that the Shimura data $(P,\FX) = (G,\FX)$ are pure, and that they satisfy 
hypothesis $(+)$. 
Let $z_0 \in \coprod \FX_j / W_j$ be a point,
and $(\FX_j / W_j)^0 \subset \coprod \FX_j / W_j$ a connected component of
$\coprod \FX_j / W_j \,$. Then there is a canonical choice of
$z_1 \in (\FX_j / W_j)^0$, such that 
(with respect to the action from Proposition~\ref{5L},)
the fibres 
$p^{-1}(z_0)$ and $p^{-1}(z_1)$ are $C_j(\BQ)$-equivariantly canoni\-cally
isomorphic as $\BR$-analytic stratified spaces : there is a cano\-ni\-cal
$C_j(\BQ)$-equivariant isomorphism $p^{-1}(z_0) \isoto p^{-1}(z_1)$, 
that restricts to give 
$p^{-1}(z_0) \cap e(Q) \isoto p^{-1}(z_1) \cap e(Q)$, for every 
parabolic $Q \in \adm^{-1}(Q_j)$.
\end{Cor}

\begin{Proof}
Apply Proposition~\ref{5Pa} and Theorem~\ref{5S}.
\end{Proof}

\begin{Rem}
The isomorphism
\[
\kappa_{z_0,j}: \bigl( A_{Q_j} \! \backslash \tilde{\pi}_j^{-1}(z_0) \bigr)^{BS}
\isoto p^{-1}(z_0) 
\]
from Theorem~\ref{5S}~(d) is \emph{a priori} compatible with \cite[Cor.~(3.8)]{Z},
but we were unable to deduce either result from the other.
\end{Rem}

In the situation of Theorem~\ref{5S}, the quotient
\[
A_Q \backslash \tilde{\pi}_j^{-1}(z_0) = 
A_Q \backslash \bigl( A_{Q_j} \! \backslash \tilde{\pi}_j^{-1}(z_0) \bigr)
\]
is of type $S-\BQ$ under $C_j \cap Q$,
for every $Q \in \adm^{-1}(Q_j)$ (Proposition~\ref{5P}, \cite[Sect.~3.9]{BS}).
In particular, the manifold with corners $(A_Q \! \backslash \tilde{\pi}_j^{-1}(z_0))^{BS}$
is defined.

\begin{Cor} \label{5T}
Assume that the Shimura data $(P,\FX) = (G,\FX)$ are pure, and that they satisfy 
hypothesis $(+)$. Let $z_0 \in (\FX_j / W_j)^0 \subset \FX^*$.
Let $Q \in \adm^{-1}(Q_j)$. \\[0.1cm]
(a)~The morphism $\kappa_{z_0,j}$ yields a 
$\Cent_{(C_j \cap Q)(\BR)} ((\FX_j/W_j)^0)$-equivariant $\BR$-analy\-tic isomorphism
\[
A_Q \! \backslash \tilde{\pi}_j^{-1}(z_0) \isoto p^{-1}(z_0) \cap e(Q) \; ,
\]
and a $\Cent_{(C_j \cap Q)(\BQ)} ((\FX_j/W_j)^0)$-equivariant isomorphism of manifolds with corners 
\[
\bigl( A_Q \! \backslash \tilde{\pi}_j^{-1}(z_0) \bigr)^{BS}
\isoto p^{-1}(z_0) \cap \overline{e(Q)} \; .
\]
(b)~Both $p^{-1}(z_0) \cap e(Q)$ and 
$p^{-1}(z_0) \cap \overline{e(Q)}$ are contractible. 
In particular, they are connected. \\[0.1cm]
(c)~Let $\Gamma_C$ be an arithmetic subgroup of $(C_j \cap Q)(\BQ)$. Then the action
of $\Cent_{\Gamma_C} ((\FX_j/W_j)^0)$ on $p^{-1}(z_0) \cap \overline{e(Q)}$
is properly discontinuous. It is free if $\Gamma_C$ is 
\emph{neat} in the sense of \cite[Sect.~17.1]{B}.
\end{Cor}

\begin{Proof}
Except for the statement on $\Cent_{(C_j \cap Q)(\BR)} ((\FX_j/W_j)^0)$-equivarian\-ce
(rather than just $\Cent_{(C_j \cap Q)(\BQ)} ((\FX_j/W_j)^0)$-equivariance),
part~(a) follows from \cite[Prop.~7.3~(i)]{BS} and Theorem~\ref{5S}.
Use the diagram
\[
\vcenter{\xymatrix@R-10pt{
        \tilde{\pi}_j^{-1}(z_0) \ar@{->>}[r] \ar@{=}[d] & 
        A_Q \backslash \tilde{\pi}_j^{-1}(z_0) \subset e(Q) \ar[d]^{p_Q} \\
        \tilde{\pi}_j^{-1}(z_0) \ar[r]^-{\tilde{\pi}_j} & \{ z_0 \}
\\}}
\]
for the action of $\Cent_{(C_j \cap Q)(\BR)} ((\FX_j/W_j)^0)$,
recalling \cite[Prop.~3.4]{BS} that the geodesic action of 
$A_Q$ commutes with the action of $Q(\BR)$.

In order to prove (b), apply (a) and \cite[Lemma~8.3.1 and its proof]{BS}. 

As for (c), use (a) and \cite[Thm.~9.3, Sect.~9.5]{BS}.
\end{Proof}

Recall (\cite[Def.~5.2]{G}, cmp.~\cite[(4.1)]{Z1}) 
the definition of the space $\FX^{rBS}$ 
(denoted $\bar{D}^{RBS}$ in \cite{G})
as a quotient of $\FX^{BS}$: two elements $x$ and $y$ of $\FX^{BS}$
map to the same point of $\FX^{rBS}$ only
if they belong to the same face. If $x$, $y \in e(Q)$, 
then $x$ and $y$ map to the same point of $\FX^{rBS}$ if and only if there is
an element $w$ belonging to the real points of the unipotent radical $\Rad^u(Q)$ of $Q$, 
such that $y = wx$. Set theoretically, we thus have
\[
\FX^{rBS} = \coprod_Q e^r(Q) \; ,
\]
where $e^r(Q) := \Rad^u(Q)(\BR) \backslash e(Q)$. By construction, 
the epimorphism from $\FX^{BS}$ to $\FX^{rBS}$ respects the stratifications.
It induces on $\FX^{rBS}$ a continuous action of $G(\BQ)$, with respect to which it is
equivariant.  

\begin{Cor} \label{5U}
Assume that the Shimura data $(P,\FX) = (G,\FX)$ are pure, and that they satisfy 
hypothesis $(+)$. \\[0.1cm]
(a)~The map $p: \FX^{BS} \to \FX^*$ factors over a continuous, $G(\BQ)$-equivariant
map $p^r: \FX^{rBS} \to \FX^*$. \\[0.1cm]
(b)~Let $z_0 \in (\FX_j / W_j)^0 \subset \FX^*$. Then 
$\kappa_{z_0,j}$ induces a $\Cent_{\bar{C}_j(\BQ)} ((\FX_j/W_j)^0)$-equivariant
identification of 
\[
\bigl( (A_{Q_j} W_j(\BR)) \backslash \tilde{\pi}_j^{-1}(z_0) \bigr)^{rBS}
\]
with the fibre of $p^r : \FX^{rBS} \to \FX^*$ over 
$z_0$. In particular, the fibre $(p^r)^{-1}(z_0)$ is connected.
\end{Cor}

\begin{Proof}
Let $Q$ be a parabolic of $G$.
Write $Q_\infty := \adm(Q)$, hence $P_\infty = \adm_{Sh}(Q)$. 
According to Corollary~\ref{4H}, 
\[
Q = (C_\infty \cap Q) P_\infty \; .
\]
The group $C_\infty$ contains $W_\infty$; therefore, 
\[
\Rad^u(Q) = \Rad^u(C_\infty \cap Q) \subset C_\infty \; ,
\]
and $\Rad^u(Q)(\BR) \subset C_\infty(\BR)^0$.

But $C_\infty(\BR)^0$ acts trivially on $\FX_\infty / W_\infty$. Together
with Complement~\ref{4La}, this shows part~(a).

If $Q = Q_j$ (hence $Q_\infty = Q_j$), 
then by the above, $\Rad^u(C_j) = \Rad^u(Q_j)$, which equals $W_j$
\cite[proof of Lemma~4.8]{P1}. Thus, the quotient $\bar{C}_j$
is contained in $\bar{Q}_j$, and
$(A_{Q_j} W_j(\BR)) \backslash \tilde{\pi}_j^{-1}(z_0)$ is identified
with $(p^r)^{-1}(z_0) \cap e^r(Q_j)$. 
We leave it to the reader to deduce part~(b)
from the construction of $\FX^{rBS}$, and from Theorem~\ref{5S}~(d).
\end{Proof}

\begin{Def} \label{5V}
Assume that the Shimura data $(P,\FX) = (G,\FX)$ are pure, and that they satisfy 
hypothesis $(+)$. 
Let $K$ be an open compact subgroup of $G (\BA_f)$. Define the \emph{reductive Borel--Serre
compactification of the Shimura variety $M^K (G,\FX)$} as the quotient space
\[
M^K (G,\FX) (\BC)^{rBS} 
                := G (\BQ) \backslash \bigl( \FX^{rBS} \times G (\BA_f) / K \bigr) \; .
\]
\end{Def}

\begin{Cor} \label{5W}
Assume that the Shimura data $(P,\FX) = (G,\FX)$ are pure, and that they satisfy 
hypothesis $(+)$. 
Let $K \subset G(\BA_f)$ be an open compact subgroup. Then 
\[
p^K: M^K  (G,\FX) (\BC)^{BS} \longto M^K (G,\FX)^* (\BC)   
\]
factors over a continuous map
\[
p^{r,K}: M^K  (G,\FX) (\BC)^{rBS} \longto M^K (G,\FX)^* (\BC)   
\]
between the reductive Borel--Serre compactification and the space of complex points
of the Baily--Borel compactification of $M^K (G,\FX)$.
\end{Cor}

\begin{Proof}
By Corollary~\ref{5U}~(a), the map $p^r \times \id_{G(\BA_f)}$ 
is $G(\BQ)$-equivariant.
\end{Proof}

\begin{Rem}
With the exception of Proposition~\ref{4aDa} and part~(c) of Theorem~\ref{4aG},
all constructions and results of Section~\ref{4a} admit analogues for  
the reductive Borel--Serre compactification, and for the maps
\[
p^{r,K}: M^K  (G,\FX) (\BC)^{rBS} \longto M^K (G,\FX)^* (\BC)   
\]
from Corollary~\ref{5Y}. We leave the details to the reader.
\end{Rem}

In order to identify the fibres of $p^K$ and $p^{r,K}$, we need to analyze the nature
of the action of arithmetic subgroups of $Q_j$ on $\FX^*$. Recall that 
\[
\tilde{\pi}_j : \FX \longto \coprod \FX_j / W_j 
\]
is equivariant with respect to $Q_j(\BR)U(\BC)$.

\begin{Prop} \label{5X}
Let $z_0 \in \FX_j / W_j$, and 
denote by $(\FX_j / W_j)^0$ the connected component of $\FX_j / W_j$
containing $z_0 \,$. Let $\Gamma \subset Q_j(\BQ)$ an arithmetic subgroup. 
Define $\Gamma_C:= C_j(\BQ) \cap \Gamma$. \\[0.1cm]
(a)~The inclusion 
$\Cent_{\Gamma_C}((\FX_j / W_j)^0) \subset \Stab_{\Gamma_C}(z_0)$ is an equality. 
In parti\-cular, the subgroup $\Stab_{\Gamma_C}(z_0)$ of $\Gamma_C$ does not
change when replacing $z_0$ by another point in $(\FX_j / W_j)^0$. \\[0.1cm]
(b)~The group $\Cent_{\Gamma_C}((\FX_j / W_j)^0) = \Stab_{\Gamma_C}(z_0)$ 
is an arithmetic subgroup of $C_j(\BQ)$. \\[0.1cm]
(c)~If the Shimura data $(P,\FX)$ satisfy hypothe\-sis $(+)$, then the index
\[
\bigl[ \Stab_\Gamma(z_0) : \Stab_{\Gamma_C}(z_0) \bigr]
\]   
is finite. \\[0.1cm]
(d)~Assume that the Shimura data $(P,\FX)$ satisfy hypothe\-sis $(+)$. 
If $\Gamma$ is neat,
then (so is $\Gamma_C$, and) both inclusions
\[
\Stab_{\Gamma_C}(z_0) \subset \Stab_\Gamma(z_0) \quad \text{and} \quad
\Stab_{\Gamma_C}(z_0) \subset \Gamma_C
\]
are equalities.
Therefore, we have $\Cent_{\Gamma_C}((\FX_j / W_j)^0) = \Stab_\Gamma(z_0) = \Gamma_C$. 
In parti\-cular, the subgroup $\Stab_{\Gamma}(z_0)$ of $\Gamma$ does not
change when replacing $z_0$ by another point in $\FX_j / W_j \,$,
or $\FX_j / W_j$ by another rational boundary component associated to 
the same admissible parabolic $Q_j$. \\[0.1cm]
(e)~Assume that the Shimura data $(P,\FX)$ satisfy hypothe\-sis $(+)$. 
If $\Gamma$ is neat,
then its action on $\coprod \FX_j / W_j$ factors through an action of $\Gamma / \Gamma_C$.
This latter action is free.
\end{Prop}

\begin{Proof}
First, let us show that 
\[
C_j(\BR)^0 U(\BC) \subset 
           \Cent_{C_j(\BR)U(\BC)} \bigl( (\FX_j / W_j)^0 \bigr) = \Stab_{C_j(\BR)U(\BC)}(z_0) \; .
\quad\quad\quad (\ast)
\]
For this, recall that 
by definition, the group $C_j(\BR)U(\BC)$ acts trivially
on all $h(\FX_j/W_j)$. Therefore
\cite[Cor.~2.12]{P1}, an element of $C_j(\BR)U(\BC)$ stabilizes
a connected component of $\FX_j/W_j$ 
if and only if it centralizes the connected component in question.
This principle proves $(\ast)$, but also the following: let $c \in C_j(\BR) U(\BC)$. Then
\[
c \; \text{acts trivially on} \; \pi_0(\FX_j/W_j) \Longleftrightarrow
c \; \text{acts trivially on} \; \FX_j/W_j \; .
\quad\quad\quad (\ast \ast)
\]
Observation~$(\ast)$ clearly implies claims~(a) and (b).

Next, write $\bar{C}_j := \pi_j(C_j) \subset \bar{Q}_j$. 
We have $\bar{Q}_j = \bar{C}_j (P_j/W_j)$, hence the quotient $Q_j/C_j$ is canonically
isomorphic to $(P_j/W_j) / Z$, where $Z:= \bar{C}_j \cap (P_j/W_j)$ is contained
in the center of $P_j/W_j$. The Shimura data $(P_j,\FX_j)$ satisfy $(+)$
\cite[proof of Cor.~4.10]{P1}, hence so does $((P_j,\FX_j)/W_j)/Z$.

The image $\Gamma/\Gamma_C$
of $\Gamma$ under the epimorphism $Q_j \onto Q_j/C_j \cong (P_j/W_j) / Z$ is arithmetic.
It follows from $(+)$, and from \cite[Lemma~1.3]{BW} that 
the action of $\Gamma/\Gamma_C$ on $((P_j,\FX_j)/W_j)/Z$ is properly discontinuous.
In particular, its stabilizers are finite.
If $\Gamma$ is neat, then they are trivial. 

This shows claim~(c), and the equality
\[
\Stab_{\Gamma_C}(z_0) = \Stab_\Gamma(z_0) 
\]
from (d) (under the neatness assumption).
Consider the algebraic representation of $C_j$ on the set of connected components
of $(\FX_j / W_j)$. Its restriction to any neat subgroup of $C_j(\BQ)$ is trivial. 
The equality
\[
\Stab_{\Gamma_C}(z_0) = \Gamma_C 
\]
therefore follows from observation $(\ast \ast)$.

Up to freeness of the action of $\Gamma / \Gamma_C$,
claim~(e) follows from (d). Claim~(d) also implies that the action of $\Gamma / \Gamma_C$
is without fixed points. But according to \cite[Sect.~6.3]{P1}, the action of 
$\Gamma$ on $\coprod \FX_j/W_j$ 
is properly discontinuous in the sense of \cite[Sect.~0.4]{P1}. Therefore, 
the action of $\Gamma / \Gamma_C$ is properly discontinuous in the usual sense. 
\end{Proof}

\begin{Cor} \label{5Y}
Assume that the Shimura data $(P,\FX) = (G,\FX)$ are pure, and that they satisfy 
hypothesis $(+)$. 
Let $K$ be an open compact subgroup of $G (\BA_f)$, $z_0 \in \coprod \FX_j / W_j \,$, and
$g \in G(\BA_f)$. 
Consider the point 
\[
[(z_0,gK)] \in M^K (Q_j,\FX) (\BC) \subset M^K (G,\FX)^* (\BC) \; ,
\]
(\emph{i.e.}, the $G(\BQ)$-class of $(z_0,gK) \in \FX^* \times G(\BA_f)/K)$,)
and define
\[
H:= H(z_0,gK) := \Stab_{Q_j(\BQ)}(z_0) \cap gKg^{-1} \subset Q_j(\BQ)
\]
and
\[
H_C:= H_C(z_0,gK) := C_j(\BQ) \cap H = \Stab_{C_j(\BQ)}(z_0) \cap gKg^{-1} 
\subset C_j(\BQ) \; .
\] 
(a)~The group $H_C$ is an arithmetic subgroup of $C_j(\BQ)$, which is of finite index in
$H$. If $K$ is neat, then so is $H_C$, and both inclusions
\[
H_C \subset H \quad \text{and} \quad H_C \subset C_j(\BQ) \cap gKg^{-1}
\]
are equalities; in particular, if $K$ is neat, then so is $H = C_j(\BQ) \cap gKg^{-1}$, 
and $H$ does not depend on the choice of
$z_0$ in $\coprod \FX_j / W_j \,$. \\[0.1cm]
(b)~The isomorphism
\[
\kappa_{z_0,j} : \bigl( A_{Q_j} \! \backslash \tilde{\pi}_j^{-1}(z_0) \bigr)^{BS}
\isoto p^{-1}(z_0)
\]
from Theorem~\ref{5S}~(d) induces an isomorphism
\[
H \backslash \bigl( A_{Q_j} \! \backslash \tilde{\pi}_j^{-1}(z_0) \bigr)^{BS}
\isoto H \backslash p^{-1}(z_0) = (p^K)^{-1} \bigl( [(z_0,gK)] \bigr) \; .
\]
In particular, the fibres of $p^K$ are connected.
If $K$ is neat, then the above isomorphism gives
\[
\bigl( C_j(\BQ) \cap gKg^{-1} \bigr) \backslash 
                  \bigl( A_{Q_j} \! \backslash \tilde{\pi}_j^{-1}(z_0) \bigr)^{BS}
\isoto (p^K)^{-1} \bigl( [(z_0,gK)] \bigr) \: .
\]
It identifies the fibre $(p^K)^{-1} ( [(z_0,gK)] )$  
with the quotient by the free action of $C_j(\BQ) \cap gKg^{-1}$ 
on the contractible space 
$( A_{Q_j} \! \backslash \tilde{\pi}_j^{-1}(z_0) )^{BS}$. \\[0.1cm]
(c)~The isomorphism $\kappa_{z_0,j}$ from Theorem~\ref{5S}~(d) induces an isomorphism
\[
\pi_{Q_j}(H) \backslash \bigl( (A_{Q_j} W_j(\BR)) \backslash \tilde{\pi}_j^{-1}(z_0) \bigr)^{rBS}
\isoto (p^{r,K})^{-1} \bigl( [(z_0,gK)] \bigr) \; .
\]
In particular, the fibre $(p^{r,K})^{-1} ([(z_0,gK)])$ is connected.
\end{Cor}

\begin{Proof}
(a): put $\Gamma := Q_j(\BQ) \cap gKg^{-1}$, and apply Proposition~\ref{5X}~(b)--(d).

\noindent (b): we leave it to the reader to prove first that under 
\[
p^K : G (\BQ) \backslash \bigl( \FX^{BS} \times G (\BA_f) / K \bigr)
\longto G (\BQ) \backslash \bigl( \FX^* \times G (\BA_f) / K \bigr) \; ,  
\]
the pre-image of the point $[(z_0,gK)]$ is indeed identified, \emph{via} 
$[(x,gK)] \mapsto Hx$, to the quotient by
the action of $H$ on $p^{-1}(z_0)$. Then use Theorem~\ref{5S}~(d) and part~(a).
The claim concerning freeness of the action of $H_C$ (provided the latter is neat)
is Corollary~\ref{5T}~(c).

\noindent (c): the strategy is formally identical to the one employed for (b)
(replace Theorem~\ref{5S}~(d) by
Corollary~\ref{5U}~(b)).
\end{Proof}

\begin{Rem} \label{5Gor}
(a)~According to \cite[Thm.~9.3]{BS}, the action of $H_C$ on 
$( A_{Q_j} \! \backslash \tilde{\pi}_j^{-1}(z_0) )^{BS}$ is properly discontinuous.
Therefore (Corollary~\ref{5Y}~(a)), the same is true for the action of $H$.
It follows (Corollary~\ref{5Y}~(b)) that the fibre $(p^K)^{-1} ( [(z_0,gK)] )$ 
is locally isomorphic to a quotient of 
$( A_{Q_j} \! \backslash \tilde{\pi}_j^{-1}(z_0) )^{BS}$ by a finite group,
which is trivial if $K$ is neat. \\[0.1cm]
(b)~The isomorphism
\[
\pi_{Q_j}(H_C) \backslash \bigl( (A_{Q_j} W_j(\BR)) \backslash \tilde{\pi}_j^{-1}(z_0) \bigr)^{rBS}
\isoto (p^{r,K})^{-1} \bigl( [(z_0,gK)] \bigr)
\]
from Corollary~\ref{5Y}~(c) is already known, at least if $\FX$ is connected: 
see \cite[Thm.~5.8]{G} (which is
stated without the assumption on neatness). It can be safely supposed that
Corollary~\ref{5Y}~(b) is known as well. Its not being stated in \loccit
\ seems to be solely due to the wish to avoid the use of
homogeneous spaces under non-reductive algebraic groups. \\[0.1cm]
(c)~Let $\FZ$ be a subspace of $M^K (G,\FX)^* (\BC)$.
The statement on connectedness of the fibres from Corollary~\ref{5Y}~(b) has 
the following implication: if $q: \FZ' \onto \FZ$ is a covering of $\FZ$,
then it is possible to factor the base change of $p^K$ to $\FZ$ through a
surjection with target $\FZ'$ only if $q$ is a homeomorphism
(consider the base change of $q$ to the universal cover of $\FZ$...). \\[0.1cm]
(d)~The preceding observation~(c) applies in particular to the quotients
$\FZ = \Delta_1 \backslash M^{\pi_1(K_f^1)} ((P_1,\FX_1)/W_1)$ occurring in the
stratification considered in
\cite[Sect.~6.3]{P1} (cmp.\ Remark~\ref{4aRem}~(a)): 
if $K$ is neat, then the ``smaller'' Shimura variety
$M^{\pi_1(K_f^1)} ((P_1,\FX_1)/W_1)$ is a (finite) covering of $\FZ$
\cite[Prop.~1.1~(d)]{BW}. The latter can thus be dominated by $p^K$ only if
the action of $\Delta_1$ is trivial, \emph{i.e.}, the covering
$M^{\pi_1(K_f^1)} ((P_1,\FX_1)/W_1) \onto \FZ$ is a homeomorphism. 
\end{Rem}

Let us prepare the stratified version of Corollary~\ref{5Y} (Theorem~\ref{5Z}).

\begin{Def} \label{5Ya}
Assume that the Shimura data $(P,\FX) = (G,\FX)$ are pure, and that they satisfy 
hypothesis $(+)$. 
Let $K$ be an open compact subgroup of $G (\BA_f)$, $z_0 \in \coprod \FX_j / W_j \,$,
$g \in G(\BA_f)$, and $Q \in \adm^{-1}(Q_j)$. \\[0.1cm]
(a)~Define
\[
e^K(Q,g)_{z_0} \subset (p^K)^{-1} \bigl( [(z_0,gK)] \bigr)
\]
as the image of $p^{-1}(z_0) \cap e(Q)$ under the projection 
\[
p^{-1}(z_0) \onto (p^K)^{-1} \bigl( [(z_0,gK)] \bigr) \; , \;
x \longmapsto [(x,gK)] \; .
\]
(b)~Define
\[
\overline{e^K(Q,g)_{z_0}} \subset (p^K)^{-1} \bigl( [(z_0,gK)] \bigr)
\]
as the closure of $e^K(Q,g)_{z_0}$ (in $(p^K)^{-1} ( [(z_0,gK)] )$
or in $M^K (G,\FX) (\BC)^{BS}$).
\end{Def}

Recall that by definition,
$e^K(Q,g)$ equals the image of $e(Q)$ under $x \mapsto [(x,gK)]$. We thus have
\[
e^K(Q,g)_{z_0} \subset (p^K)^{-1} \bigl( [(z_0,gK)] \bigr) \cap e^K(Q,g) 
\]
and
\[
\overline{e^K(Q,g)_{z_0}} \subset (p^K)^{-1} \bigl( [(z_0,gK)] \bigr) \cap \overline{e^K(Q,g)} \; .
\]
Recall also that $p^{-1}(z_0) \cap e(Q_j)$ is dense in $p^{-1}(z_0)$
(Theorem~\ref{5S}~(d)). 
It follows that
\[
\overline{e^K(Q_j,g)_{z_0}} = (p^K)^{-1} \bigl( [(z_0,gK)] \bigr) \; .
\]

\begin{Prop} \label{5Yb}
We keep the setting of Definition~\ref{5Ya}. Define
\[
H(z_0,gK) := \Stab_{Q_j(\BQ)}(z_0) \cap gKg^{-1}
\]
as in Corollary~\ref{5Y}. \\[0.1cm]
(a)~$\overline{e^K(Q,g)_{z_0}}$ is the
image of $p^{-1}(z_0) \cap \overline{e(Q)}$ under $x \mapsto [(x,gK)]$.
Both $e^K(Q,g)_{z_0}$ and
$\overline{e^K(Q,g)_{z_0}}$ are connected. \\[0.1cm]
(b)~We have 
\[
(p^K)^{-1} \bigl( [(z_0,gK)] \bigr) \cap \overline{e^K \bigl( Q,G(\BA_f) \bigr)}
= (p^K)^{-1} \bigl( [(z_0,gK)] \bigr) \cap e^K \bigl( Q,G(\BA_f) \bigr)' \; ,
\]
where $e^K(Q,G(\BA_f))'$ is the subset of $\overline{e^K(Q,G(\BA_f))}$
from Definition~\ref{4aC}. \\[0.1cm]
(c)~We have
\[
(p^K)^{-1} \bigl( [(z_0,gK)] \bigr) \cap e^K \bigl( Q,G(\BA_f) \bigr)
= \bigcup_{\gamma \in Q_j(\BQ)} e^K \bigl( \gamma Q \gamma^{-1},g \bigr)_{z_0}
\]
and
\[
(p^K)^{-1} \bigl( [(z_0,gK)] \bigr) \cap \overline{e^K \bigl( Q,G(\BA_f) \bigr)}
= \bigcup_{\gamma \in Q_j(\BQ)} \overline{e^K \bigl( \gamma Q \gamma^{-1},g \bigr)_{z_0}} \; .
\]
(d)~Let $\gamma_1, \gamma_2 \in Q_j(\BQ)$. The following are equivalent. 
\begin{enumerate}
\item[(i)] $e^K(\gamma_1 Q \gamma_1^{-1},g)_{z_0} 
= e^K(\gamma_2 Q \gamma_2^{-1},g)_{z_0}$,
\item[(ii)] $\overline{e^K(\gamma_1 Q \gamma_1^{-1},g)_{z_0}} = 
\overline{e^K(\gamma_2 Q \gamma_2^{-1},g)_{z_0}}$,
\item[(iii)] the intersection 
$\overline{e^K(\gamma_1 Q \gamma_1^{-1},g)_{z_0}} \cap 
\overline{e^K(\gamma_2 Q \gamma_2^{-1},g)_{z_0}}$
is not emp\-ty, 
\item[(iv)] the classes of $\gamma_1$ and $\gamma_2$ in the double quotient
\[
H(z_0,gK) \backslash Q_j(\BQ) / Q(\BQ)  
\]
are the same. 
\end{enumerate} 
(e)~We have
\[
(p^K)^{-1} \bigl( [(z_0,gK)] \bigr) \cap e^K \bigl( Q,G(\BA_f) \bigr)
= \coprod_{\gamma \in H(z_0,gK) \backslash Q_j(\BQ) / Q(\BQ)} e^K \bigl( \gamma Q \gamma^{-1},g \bigr)_{z_0} 
\]
and
\[
(p^K)^{-1} \bigl( [(z_0,gK)] \bigr) \cap \overline{e^K \bigl( Q,G(\BA_f) \bigr)}
= \coprod_{\gamma \in H(z_0,gK) \backslash Q_j(\BQ) / Q(\BQ)} 
\overline{e^K \bigl( \gamma Q \gamma^{-1},g \bigr)_{z_0}} \; .
\]
(f)~The map $\gamma \mapsto \overline{e^K \bigl( \gamma Q \gamma^{-1},g \bigr)_{z_0}}$
induces a canonical bijection between the
double quotient $H(z_0,gK) \backslash Q_j(\BQ) / Q(\BQ)$ and 
the set of connected components of 
\[
(p^K)^{-1} \bigl( [(z_0,gK)] \bigr) \cap \overline{e^K \bigl( Q,G(\BA_f) \bigr)} \; .
\]
In particular, the double quotient $H(z_0,gK) \backslash Q_j(\BQ) / Q(\BQ)$ is finite. 
\end{Prop}

\begin{Proof}
(a): the image of the closed subset 
$p^{-1}(z_0) \cap \overline{e(Q)}$ under the quotient map $x \mapsto [(x,gK)]$
is closed.
Now use Corollary~\ref{5T}~(a) and (b).

\noindent (b): apply Theorem~\ref{4aG}~(b).

\noindent (c): observe that for any $\gamma \in Q_j(\BQ)$, 
\[
\overline{e^K \bigl( \gamma Q \gamma^{-1},g \bigr)_{z_0}}
\subset \overline{e^K \bigl( \gamma Q \gamma^{-1}, G(\BA_f) \bigr)}
= \overline{e^K \bigl( Q, G(\BA_f) \bigr)} \; .
\]
Therefore,
\[
(p^K)^{-1} \bigl( [(z_0,gK)] \bigr) \cap \overline{e^K \bigl( Q,G(\BA_f) \bigr)}
\supset \bigcup_{\gamma \in Q_j(\BQ)} \overline{e^K \bigl( \gamma Q \gamma^{-1},g \bigr)_{z_0}} \; .
\]
In order to show the reverse inclusion, let $x \in p^{-1}(z_0) \cap e(R)$, for
some parabolic $R$ of $G$, and assume that the point $[(x,gK)]$ belongs to
\[
(p^K)^{-1} \bigl( [(z_0,gK)] \bigr) \cap \overline{e^K \bigl( Q,G(\BA_f) \bigr)}
\stackrel{(b)}{=} (p^K)^{-1} \bigl( [(z_0,gK)] \bigr) \cap e^K \bigl( Q,G(\BA_f) \bigr)' \; .
\]
Since $[(x,gK)]$ obviously belongs to $e^K(R,G(\BA_f))$, we have $R = \gamma R' \gamma^{-1}$,
for some parabolic $R'$ containing $P_j$ and contained in $Q$,
and some $\gamma \in G(\BQ)$. In particular,
$\adm(R') = Q_j$. But since $p^{-1}(z_0) \cap e(R)$ is not empty, $\adm(R) = Q_j$ 
as well (Complement~\ref{4Lb}). The map $\adm$ being $G(\BQ)$-equivariant, we conclude
that $\gamma$ belongs in fact to $Q_j(\BQ)$. Therefore, 
$x \in \overline{e(\gamma Q \gamma^{-1})}$, and
\[
[(x,gK)] \in \overline{e^K \bigl( \gamma Q \gamma^{-1},g \bigr)_{z_0}} \; .
\]
The proof of the equality
\[
(p^K)^{-1} \bigl( [(z_0,gK)] \bigr) \cap e^K \bigl( Q,G(\BA_f) \bigr)
= \bigcup_{\gamma \in Q_j(\BQ)} e^K \bigl( \gamma Q \gamma^{-1},g \bigr)_{z_0}
\]
is left to the reader.

\noindent (d): clearly (iv)$\Rightarrow$(i)$\Rightarrow$(ii)$\Rightarrow$(iii).
It remains to show that (iii) implies (iv). Thus, assume that
the intersection 
\[
\overline{e^K \bigl( \gamma_1 Q \gamma_1^{-1},g \bigr)_{z_0}} \cap 
\overline{e^K \bigl( \gamma_2 Q \gamma_2^{-1},g \bigr)_{z_0}}
\]
is not empty. According to (a), there exist parabolics $R_1 \subset \gamma_1 Q \gamma_1^{-1}$
and $R_2 \subset \gamma_2 Q \gamma_2^{-1}$, points
\[
x_i \in p^{-1}(z_0) \cap e(R_i) \; , \; i = 1, 2 \; ,
\]   
and $h \in G(\BQ) \cap gKg^{-1}$, 
such that $x_2 = h x_1$. We have $\adm(R_i) = Q_j$, $i = 1,2$  (Complement~\ref{4Lb}). Since 
\[ 
R_2 = h R_1 h^{-1} \; , 
\]
we find as in (c) that $h \in Q_j(\BQ)$. 
But $h$ stabilizes $z_0$; therefore, $h \in H(z_0,gK)$. 
The intersection $\gamma_2 Q \gamma_2^{-1} \cap (h\gamma_1) Q (h\gamma_1)^{-1}$
contains $R_2$ and is therefore parabolic. According to Propostion~\ref{Par},
\[
\gamma_2 Q \gamma_2^{-1} = (h\gamma_1) Q (h\gamma_1)^{-1} \; .
\]
In other words, we have
$\gamma_2 = (h\gamma_1) q$, for some element $q$ of $Q$.

\noindent (e): this follows from (c) and (d).

\noindent (f): each $\overline{e^K ( \gamma Q \gamma^{-1},g )_{z_0}}$ is closed. By (a), it is
connected. 
Using (e), we see that it is the complement of a union of 
$\overline{e^K ( \gamma' Q {\gamma'}^{-1},g )_{z_0}}$, the union in question being 
indexed by the complement of the class of $\gamma$
in $H(z_0,gK) \backslash Q_j(\BQ) / Q(\BQ)$. Therefore, if this index set is finite,
then $\overline{e^K ( \gamma Q \gamma^{-1},g )_{z_0}}$ is open, and our
claim is proved. 

As for finiteness of $H(z_0,gK) \backslash Q_j(\BQ) / Q(\BQ)$,
observe first that by Corollary~\ref{5Y}~(a), we may replace $H(z_0,gK)$ by
an arithmetic subgroup $H'$ of $C_j(\BQ)$. Next, write $S := C_j \cap Q$.
This subgroup of $C_j$ is parabolic, and $Q = S P_j$ (Corollary~\ref{4H}).  
The inclusion of $C_j(\BQ) P_j(\BQ)$ into $Q_j(\BQ)$ being of finite index 
(as $Q_j = C_j P_j$), our claim follows from the finiteness of the double quotient
\[
H' \backslash C_j(\BQ) / S(\BQ)
\]
\cite[Cor.~15.7]{B}
\end{Proof}

\begin{Rem} \label{5Yc}
(a)~\emph{Via} $\gamma \mapsto \gamma Q \gamma^{-1}$, the double quotient 
\[
H(z_0,gK) \backslash Q_j(\BQ) / Q(\BQ)
\] 
is in bijection with the set of orbits
under the action of $H(z_0,gK)$ on the parabolic subgroups of $Q_j$ that are
conjugate to $Q$. 

It will be important to reformulate this observation, using parabolic subgroups
of $C_j$. The bijective correspondence $Q' \mapsto C_j \cap Q'$, $S \mapsto S P_j$
from Corollary~\ref{4H} between parabolics of $Q_j$ containing $P_j$ on the one hand,
and parabolics of $C_j$ on the other, 
respects $Q_j$-conjugation (as $C_j$ and $P_j$ are normal in $Q_j$).
However, a $Q_j(\BQ)$-orbit may consist of several $C_j(\BQ)$-orbits (the number of
such $C_j(\BQ)$-orbits equals the index of $C_j(\BQ) Q(\BQ)$ in $Q_j(\BQ)$). 

Here is the precise statement: \emph{via} $\gamma \mapsto C_j \cap \gamma Q \gamma^{-1}$, 
the double quotient 
\[
H(z_0,gK) \backslash Q_j(\BQ) / Q(\BQ)
\] 
is in bijection with the set of orbits
under the action of $H(z_0,gK)$ on the parabolic subgroups of $C_j$ that are
$Q_j(\BQ)$-conjugate to $C_j \cap Q$. \\[0.1cm]
(b)~If $K$ is neat, then according to Corollary~\ref{5Y}~(a), 
the group
\[
H(z_0,gK) = C_j(\BQ) \cap gKg^{-1}
\]
does not depend on $z_0 \in \coprod \FX_j / W_j \,$.
Hence, neither does the double quotient
\[
H(z_0,gK) \backslash Q_j(\BQ) / Q(\BQ) \; .
\]
\end{Rem}

Putting everything together, we obtain the second main result of this section. 

\begin{Thm} \label{5Z}
Assume that the Shimura data $(P,\FX) = (G,\FX)$ are pure, and that they satisfy 
hypothesis $(+)$. 
Let $K$ be an open compact subgroup of $G (\BA_f)$, $z_0 \in \coprod \FX_j / W_j \,$,
$g \in G(\BA_f)$, and $Q \in \adm^{-1}(Q_j)$. 
Consider the point $[(z_0,gK)] \in M^K (G,\FX)^* (\BC)$,
the canonical stratum $e^K(Q,G(\BA_f))$ of $M^K (G,\FX) (\BC)^{BS}$,
and its closure $\overline{e^K(Q,G(\BA_f))}$. 
Define
\[
H(z_0,gK) := \Stab_{Q_j(\BQ)}(z_0) \cap gKg^{-1}
\]
as in Corollary~\ref{5Y}. \\[0.1cm]
(a)~The projection 
\[
p^{-1}(z_0) \onto (p^K)^{-1} \bigl( [(z_0,gK)] \bigr) \; , \;
x \longmapsto [(x,gK)]
\]
induces isomorphisms between
\[
H(z_0,gK) \backslash 
    \Bigr( \coprod_{\gamma \in Q_j(\BQ) / Q(\BQ)} 
                 \bigl( p^{-1}(z_0) \cap e ( \gamma Q \gamma^{-1}) \bigr) \Bigr) 
\subset H(z_0,gK) \backslash p^{-1}(z_0)
\]
and
\[
(p^K)^{-1} \bigl( [(z_0,gK)] \bigr) \cap e^K \bigl( Q,G(\BA_f) \bigr) \; ,
\]
as well as between
\[
H(z_0,gK) \backslash 
    \Bigr( \coprod_{\gamma \in Q_j(\BQ) / Q(\BQ)} 
                 \bigl( p^{-1}(z_0) \cap \overline{e ( \gamma Q \gamma^{-1})} \bigr) \Bigr) 
\subset H(z_0,gK) \backslash p^{-1}(z_0)
\]
and
\[
(p^K)^{-1} \bigl( [(z_0,gK)] \bigr) \cap \overline{e^K \bigl( Q,G(\BA_f) \bigr)} 
= (p^K)^{-1} \bigl( [(z_0,gK)] \bigr) \cap e^K \bigl( Q,G(\BA_f) \bigr)' \; .
\]
In particular, the projection $x \mapsto [(x,gK)]$ induces an isomorphism
\[
H(z_0,gK) \backslash p^{-1}(z_0) \isoto (p^K)^{-1} \bigl( [(z_0,gK)] \bigr)  \; .
\]
(b)~The choice of a (necessarily finite) set $\Omega \subset Q_j(\BQ)$ of representatives of 
$H(z_0,gK) \backslash Q_j(\BQ) / Q(\BQ)$ induces identifications of
\[
\coprod_{\gamma \in \Omega} \Bigl( \bigl(H(z_0,gK) \cap \gamma Q(\BQ) \gamma^{-1} \bigr) 
         \backslash \bigl( p^{-1}(z_0) \cap e ( \gamma Q \gamma^{-1}) \bigr) \Bigr) 
\]
and 
\[
H(z_0,gK) \backslash 
    \Bigr( \coprod_{\gamma \in Q_j(\BQ) / Q(\BQ)} 
                 \bigl( p^{-1}(z_0) \cap e ( \gamma Q \gamma^{-1}) \bigr) \Bigr) \; ,
\]
as well as of 
\[
\coprod_{\gamma \in \Omega} \Bigl( \bigl(H(z_0,gK) \cap \gamma Q(\BQ) \gamma^{-1} \bigr) 
         \backslash \bigl( p^{-1}(z_0) \cap \overline{e ( \gamma Q \gamma^{-1})} \bigr) \Bigr) 
\]
and 
\[
H(z_0,gK) \backslash 
    \Bigr( \coprod_{\gamma \in Q_j(\BQ) / Q(\BQ)} 
            \bigl( p^{-1}(z_0) \cap \overline{e ( \gamma Q \gamma^{-1})} \bigr) \Bigr) \; .
\]
(c)~If $K$ is neat, then for every $\gamma \in Q_j(\BQ)$, we have
\[
H(z_0,gK) \cap \gamma Q(\BQ) \gamma^{-1} 
= (C_j(\BQ) \cap gKg^{-1}) \cap \gamma Q(\BQ)  \gamma^{-1} \; .
\]
The component 
\[
\Bigl( \bigl(H(z_0,gK) \cap \gamma Q(\BQ) \gamma^{-1} \bigr) 
         \backslash \bigl( p^{-1}(z_0) \cap \overline{e ( \gamma Q \gamma^{-1})} \bigr) \Bigr)
\]
is equal to the quotient by the free action of 
the group $(C_j(\BQ) \cap gKg^{-1}) \cap \gamma Q(\BQ)  \gamma^{-1}$ 
on the contractible space $p^{-1}(z_0) \cap \overline{e(\gamma Q \gamma^{-1})}$. \\[0.1cm]
(d)~If $K$ is neat, and $(\FX_j / W_j)^0 \subset \coprod \FX_j / W_j$ 
a connected component of
$\coprod \FX_j / W_j \,$, then there is a canonical choice of
$z_1 \in (\FX_j / W_j)^0$, and a canonical isomorphism
\[
(p^K)^{-1} \bigl( [(z_0,gK)] \bigr) \isoto
(p^K)^{-1} \bigl( [(z_1,gK)] \bigr) 
\]
of $\BR$-analytic stratified spaces.
In particular, the latter induces a canonical isomorphism between
\[
(p^K)^{-1} \bigl( [(z_0,gK)] \bigr) \cap e^K \bigl( Q,G(\BA_f) \bigr)' 
\]
and
\[
(p^K)^{-1} \bigl( [(z_1,gK)] \bigr) \cap e^K \bigl( Q,G(\BA_f) \bigr)' \; .
\]
\end{Thm}

\begin{Proof}
We leave it to the reader to show, using Proposition~\ref{Par}, that the union
\[
\coprod_{\gamma \in Q_j(\BQ) / Q(\BQ)} 
                 \bigl( p^{-1}(z_0) \cap \overline{e ( \gamma Q \gamma^{-1})} \bigr) 
\subset p^{-1}(z_0)
\]
is indeed disjoint. Part~(b) follows as $\gamma Q(\BQ) \gamma^{-1}$ 
equals the stabilizer in $G(\BQ)$ of 
$e(\gamma Q \gamma^{-1})$, and also of $\overline{e(\gamma Q \gamma^{-1})}$,
for any $\gamma \in Q_j(\BQ)$.

According to Proposition~\ref{5Yb}~(e), we have
\[
(p^K)^{-1} \bigl( [(z_0,gK)] \bigr) \cap e^K \bigl( Q,G(\BA_f) \bigr)
= \coprod_{\gamma \in \Omega} e^K \bigl( \gamma Q \gamma^{-1},g \bigr)_{z_0} 
\]
and
\[
(p^K)^{-1} \bigl( [(z_0,gK)] \bigr) \cap \overline{e^K \bigl( Q,G(\BA_f) \bigr)}
= \coprod_{\gamma \in \Omega} 
\overline{e^K \bigl( \gamma Q \gamma^{-1},g \bigr)_{z_0}} \; ,
\]
while 
\[
(p^K)^{-1} \bigl( [(z_0,gK)] \bigr) \cap \overline{e^K \bigl( Q,G(\BA_f) \bigr)}
= (p^K)^{-1} \bigl( [(z_0,gK)] \bigr) \cap e^K \bigl( Q,G(\BA_f) \bigr)'
\]
(Proposition~\ref{5Yb}~(b)). We thus may treat every component 
$e^K ( \gamma Q \gamma^{-1},g )_{z_0}$, resp.\
$\overline{e^K ( \gamma Q \gamma^{-1},g )_{z_0}}$,
separately. By Definition~\ref{5Ya}~(a), the former is the image of
$p^{-1}(z_0) \cap e(\gamma Q \gamma^{-1})$ under the projection 
\[
p^{-1}(z_0) \onto (p^K)^{-1} \bigl( [(z_0,gK)] \bigr) \; , \;
x \longmapsto [(x,gK)] \; .
\]
According to Proposition~\ref{5Yb}~(a), the latter is the image of
$p^{-1}(z_0) \cap \overline{e(\gamma Q \gamma^{-1})}$ under the same projection. 

Then, part~(a) follows from (b), and 
from Corollary~\ref{5Y}~(b).

The part of claim~(c) concerning the group by which we divide, follows from 
Corollary~\ref{5Y}~(a).
As for freeness of the action, and contractibility of 
$p^{-1}(z_0) \cap \overline{e(\gamma Q \gamma^{-1})}$,
use Corollary~\ref{5Y}~(b) and Corollary~\ref{5T}~(b).

Part~(d) follows from (a), (c), and from Corollary~\ref{5Sa}.
\end{Proof}

\forget{
\begin{Cor} \label{5Zb}
Assume that the Shimura data $(P,\FX) = (G,\FX)$ are pure, and that they satisfy 
hypothesis $(+)$. 
Let $K$ be a neat open compact subgroup of $G (\BA_f)$, 
$g \in G(\BA_f)$, and $Q \in \adm^{-1}(Q_j)$. 
Fix a (necessarily finite) set $\Omega \subset Q_j(\BQ)$ of representatives of 
$\bigl( C_j(\BQ) \cap gKg^{-1} \bigr) \backslash Q_j(\BQ) / Q(\BQ)$.
Then the map 
\[
\FX^{BS} \longto M^K  (G,\FX) (\BC)^{BS} \; , \;
x \longmapsto [(x,gK)]
\]
induces isomorphisms 
\[
\coprod_{\gamma \in \Omega} 
\Bigl( \bigl((C_j \cap \gamma Q \gamma^{-1})(\BQ) \cap gKg^{-1} \bigr) 
                                                        \backslash e ( \gamma Q \gamma^{-1}) \Bigr) 
\isoto e^K \bigl( Q,Q_j(\BQ)gK \bigr)   
\]
and 
\[
\coprod_{\gamma \in \Omega} 
\Bigl( \bigl((C_j \cap \gamma Q \gamma^{-1})(\BQ) \cap gKg^{-1} \bigr) 
                                                        \backslash e ( \gamma Q \gamma^{-1})' \Bigr) 
\isoto e^K \bigl( Q,Q_j(\BQ)gK \bigr)' \; .
\]
{\bf CAUTION: this statement is wrong. Le right hand sides need to be replaced by
their base changes to $H_C \backslash \FX^{BS}$!}

Here, the $e ( \gamma Q \gamma^{-1})'$ are the open subsets of 
$\overline{e ( \gamma Q \gamma^{-1})}$ from Definition~\ref{4aC}, while
\[
e^K \bigl( Q,Q_j(\BQ)gK \bigr) := \bigcup_{h \in Q_j(\BQ)gK} e^K ( Q,h) 
\subset e^K \bigl( Q,G(\BA_f) \bigr)
\]
and 
\[
e^K \bigl( Q,Q_j(\BQ)gK \bigr)' := 
\bigcup_{h \in Q_j(\BQ)gK} \bigcup_{P_j \subset R \subset Q} e^K ( R,h) 
\subset e^K \bigl( Q,G(\BA_f) \bigr)' \; .
\]
(b)~For every $\gamma \in Q_j(\BQ)$,
the action of $(C_j \cap \gamma Q \gamma^{-1})(\BQ) \cap gKg^{-1}$ 
on $e(\gamma Q \gamma^{-1})'$ is free. 
\end{Cor} 

\begin{Proof}
Given that the fibres of $p$ over $\coprod \FX_j / W_j$
are stable under $C_j(\BQ) \cap gKg^{-1}$ (Corollary~\ref{5Y}~(a)), both (a) and (b) follow from
Theorem~\ref{5Z}: indeed, Complement~\ref{4Lb} tells us that
\[
p^{-1} \bigl( \coprod \FX_j / W_j \bigr) \cap e ( \gamma Q \gamma^{-1})
= e ( \gamma Q \gamma^{-1})
\]
and
\[
p^{-1} \bigl( \coprod \FX_j / W_j \bigr) \cap \overline{e ( \gamma Q \gamma^{-1})} 
= e ( \gamma Q \gamma^{-1})'
\]
for all $\gamma \in Q_j(\BQ)$, which takes care of the left hand sides of our two isomorphisms.
As for the right hand sides, let $x \in \FX^{BS}$. 
Then $[(x,gK)]$ lies in $e^K ( R,G(\BA_f) )$,
for some parabolic $R$ of $G$, if and only if there is $\tau \in G(\BQ)$ such that
$x \in e(\tau^{-1} R \tau)$. 
If we apply this principle to $x$ belonging to $p^{-1} ( \coprod \FX_j / W_j ) \subset \FX^{BS}$,
\emph{i.e.} (Complement~\ref{4Lb}), to $e(Q_j)'$, and to $R \subset Q_j$, then we see
that
\[
x \in e^K \bigl( R,G(\BA_f) \bigr) \Longleftrightarrow
\exists \, \tau \in Q_j(\BQ) \, , \, x \in e(\tau^{-1} R \tau) 
\]
($R$ is contained in $Q_j \cap \tau Q_j \tau^{-1}$; apply Proposition~\ref{Par}),
and this condition implies that $P_j \subset R$ (as $P_j$ is normal in $Q_j$).  
We conclude that
\[
\bigcup_{z_0 \in \coprod \FX_j / W_j} (p^K)^{-1} \bigl( [(z_0,gK)] \bigr) \cap 
                                                                                                        e^K \bigl( Q,G(\BA_f) \bigr) 
=  e^K \bigl( Q,Q_j(\BQ)gK \bigr)  \; ,
\]
(use the above for $R = Q$, noting that $e^K ( \tau^{-1} R \tau,g) = e^K ( R,\tau g)$), while
\[
\bigcup_{z_0 \in \coprod \FX_j / W_j} (p^K)^{-1} \bigl( [(z_0,gK)] \bigr) \cap 
                                                                                  \overline{e^K \bigl( Q,G(\BA_f) \bigr)}
=  e^K \bigl( Q,Q_j(\BQ)gK \bigr)'  \; .
\]
(use the above for $R \subset Q$).
\end{Proof}
}
We leave it to the reader to formulate and prove the variants of Proposition~\ref{5Yb} and Theorem~\ref{5Z}
valid for the fibres of 
\[
p^{r,K}: M^K  (G,\FX) (\BC)^{rBS} \longto M^K (G,\FX)^* (\BC) \; .
\]
\forget{
It will turn out to be important to control the dependence of the group
\[
(C_j \cap Q)(\BQ) \cap gKg^{-1}
\]
occurring in Theorem~\ref{5Z}~(b) on the class $gK$. 

\begin{Prop} \label{5Za}
Assume that the Shimura data $(P,\FX) = (G,\FX)$ are pure.
Denote as usual by $P_j$ the canonical normal subgroup of $Q_j \,$, and by
$W_j$ the unipotent radical.
Let $K$ be an open compact subgroup of $G (\BA_f)$, $g \in G(\BA_f)$, and $Q \in \adm^{-1}(Q_j)$. 
Define 
\[
\Gamma(Q,gK) := (C_j \cap Q)(\BQ) \cap gKg^{-1}
\]
and
\[
\Lambda(Q,gK) := (C_j \cap Q)(\BQ) \cap W_j(\BQ) \cdot gKg^{-1} \; .
\] 
(a)~We have
\[
\Lambda(Q,gK) = W_j(\BQ) \cdot \Gamma(Q,gK) \; .
\] 
(b)~For any element $p_j$ of $P_j(\BA_f)$, we have
\[
\Lambda(Q,gK) = \Lambda(Q,p_jgK) \; .
\] 
\end{Prop}

\begin{Proof}
(a): the group $(C_j \cap Q)(\BQ)$ contains $W_j(\BQ)$. 

\noindent (b): by symmetry, and (a), it suffices to show the inclusion
\[
\Gamma(Q,p_jgK) \subset \Lambda(Q,gK) \; .
\] 
Thus, let $c$ be an element of $\Gamma(Q,p_jgK) \subset (C_j \cap Q)(\BQ)$
By definition of $C_j$, the elements $c$ and $p_j$ commute modulo $W_j$, \emph{i.e.,}
\[
c= w p_j^{-1}cp_j
\]
for some element $w$ of $W_j(\BA_f)$. As $\Gamma(Q,p_jgK) \subset p_jgKg^{-1}p_j^{-1}$,
we see that $c = w p_j^{-1}cp_j \in W_j(\BA_f) \cdot gKg^{-1} = W_j(\BQ) \cdot gKg^{-1}$
(by strong approximation for unipotent groups).
\end{Proof}
}


\bigskip

%
%

\section{Triviality of $p$ over each stratum}
\label{6}



In this last section, we prove (Theorems~\ref{6A} and \ref{6F}~(a)), that 
over any stratum of their target, both maps $p$ and $p^K$ 
are locally trivial fibrations (whose fibres where described in Section~\ref{5});
as far as $p^K$ is concerned, this result necessitates the level $K$ to be neat.
As for $p$, our result is known, since implied by \cite[Prop.~(3.8)~(ii)]{Z}. 
The analogue of the result concerning $p^K$ for the projection from the reductive Borel--Serre
to the Baily--Borel compactification holds, too (Theorems~\ref{6F}~(b)); we thus
recover \cite[Thm.~5.8]{G}.
For the cohomological applications we have in mind, it will be important
to sharpen the local triviality result on the map $p$: 
according to Theorem~\ref{6Ca}, it is actually (globally) trivial over any stratum
(the point being that the stratum may have several connected components). \\

We fix mixed Shimura data $(P,\FX)$. Until Complement~\ref{5J}, no further hypotheses  
on $(P,\FX)$ will be required. We write $(G,\FX/W)$ for the quotient $(P,\FX)/W$,
and $\pi$ for the canonical epimorphism from $(P,\FX)$ to $(G,\FX/W)$. We need to
study sections of the map $\FX \onto \FX/W$ induced by $\pi$; the latter map
will be denoted by the same symbol. 

\begin{Def} \label{5A}
Let $z_0 \in \FX/W$. Define a map
\[
\Phi_{z_0} : \pi^{-1}(z_0) \times \FX/W \longto \FX \; , \; 
(x,z) \longmapsto \ell x \; ,
\]
where $\ell \in \Cent_{P(\BR)U(\BC)}(h_x \circ w)$ is such that 
$\pi(\ell) z_0 = z$.
\end{Def}

\begin{Prop} \label{5B}
Let $z_0 \in \FX/W$. \\[0.1cm]
(a)~The map $\Phi_{z_0}$ is well defined. \\[0.1cm]
(b)~The map $\Phi_{z_0}$ is an $\BR$-analytic isomorphism. \\[0.1cm]
(c)~There is a commutative diagram
\[
\vcenter{\xymatrix@R-10pt{
\pi^{-1}(z_0) \times \FX/W \ar[rr]^-{\Phi_{z_0}} \ar@{->>}[dr]_-{pr_2} && 
\FX \ar@{->>}[dl]^-{\pi} \\
& \FX/W  &
\\}}
\]
where $pr_2$ denotes the projection onto the second factor. \\[0.1cm]
(d)~The map $\Phi_{z_0}$ is $P(\BR)U(\BC)$-equivariant
in the following sense: 
\[
P(\BR)U(\BC) = W(\BR)U(\BC) \rtimes \Cent_{P(\BR)U(\BC)}(h_x \circ w) 
\]
for all $x \in \FX$; for
$w \in W(\BR)U(\BC)$, we have
\[
\Phi_{z_0} (wx,z) = w \Phi_{z_0} (x,z) \; , \; 
\forall \, (x,z) \in \pi^{-1}(z_0) \times \FX/W \; ,
\] 
and for $x \in \pi^{-1}(z_0)$ and 
$\ell \in \Cent_{P(\BR)U(\BC)}(h_x \circ w)$, we have
\[
\Phi_{z_0} (x,\pi(\ell)z) = \ell \Phi_{z_0} (x,z) \; , \; 
\forall \, z \in \FX/W \; .
\]
(e)~Let $z_1$ be another point of $\FX/W$. Then
for all $(y,z) \in \pi^{-1}(z_1) \times \FX/W$, we have
\[
\Phi_{z_1} (y,z) = \Phi_{z_0} \bigl( \Phi_{z_1}(y,z_0),z \bigr) \; .
\]
\end{Prop}

\begin{Proof}
We first repeat the remark preceding Proposition~\ref{1Ipre},
showing that for any $x \in \FX$, the subgroup $\Cent_{P_\BC}(h_x \circ w)$
is a Levi subgroup of $P_\BC \,$. 
According to \cite[Def.~2.1~(ii)]{P1}, its image in $(P/U)_\BC$
is defined over $\BR$, whence the first equality in (d),
\[
P(\BR)U(\BC) = W(\BR)U(\BC) \rtimes \Cent_{P(\BR)U(\BC)}(h_x \circ w) \; .
\]
Fix $x_0 \in \pi^{-1}(z_0)$.
We identify $\FX$ in a $P(\BR)U(\BC)$-equivariant way with
\[
W(\BR)U(\BC) \times 
\bigl( \Cent_{P(\BR)U(\BC)}(h_{x_0} \circ w) / \Stab_{P(\BR)U(\BC)}(x_0) \bigr) \; ,
\]
by sending $gx_0$ to the class of $g$ modulo $\Stab_{P(\BR)U(\BC)}(x_0)$, for
every $g$ in $P(\BR)U(\BC) = W(\BR)U(\BC) \rtimes \Cent_{P(\BR)U(\BC)}(h_{x_0} \circ w)$;
likewise, we identify $\FX/W$ with
\[
G(\BR) / \Stab_{G(\BR)}(z_0) \; .
\]
These identifications are $\BR$-analytic, and $\pi: \FX \onto \FX/W$ corresponds to
the projection onto the second factor
$\Cent_{P(\BR)U(\BC)}(h_{x_0} \circ w) / \Stab_{P(\BR)U(\BC)}(x_0)$, 
followed by the map induced by
$\Cent_{P(\BR)U(\BC)}(h_{x_0} \circ w) \isoto G(\BR)$.

\noindent (a): the isomorphism $\Cent_{P(\BR)U(\BC)}(h_x \circ w) \isoto G(\BR)$
induces an isomorphism $\Stab_{P(\BR)U(\BC)}(x) \isoto \Stab_{G(\BR)}(z_0)$
\cite[Lemma~1.17~(a), Cor.~2.12]{P1}.

\noindent (b)--(d): note first that for $g \in P(\BR)U(\BC)$ and $x \in \FX$, we have
\[
\Cent_{P(\BR)U(\BC)}(h_{gx} \circ w) = 
\Int(g) \bigl( \Cent_{P(\BR)U(\BC)}(h_x \circ w) \bigr)
\]
and 
\[
\Stab_{P(\BR)U(\BC)}(g x) =
\Int(g) \bigl( \Stab_{P(\BR)U(\BC)}(x) \bigr) \; .
\]
It follows that
in our $\BR$-analytic identifications
\[
\FX =  W(\BR)U(\BC) \times 
\bigl( \Cent_{P(\BR)U(\BC)}(h_{x_0} \circ w) / \Stab_{P(\BR)U(\BC)}(x_0) \bigr) 
\]
and
\[
\FX/W = G(\BR) / \Stab_{G(\BR)}(z_0) \; ,
\]
the map $\Phi_{z_0}$ associates to
\[
(u,[g]) \in W(\BR)U(\BC) \times G(\BR) / \Stab_{G(\BR)}(z_0) 
\]
the pair
\[
(u,[\ell]) \in 
W(\BR)U(\BC) \times 
\bigl( \Cent_{P(\BR)U(\BC)}(h_{x_0} \circ w) / \Stab_{P(\BR)U(\BC)}(x_0) \bigr) \; , 
\]
where $\ell \in \Cent_{P(\BR)U(\BC)}(h_{x_0} \circ w)$ maps to $g$ under $\pi$.

\noindent (e): let $\ell' \in \Cent_{P(\BR)U(\BC)}(h_y \circ w)$ such that
$\pi(\ell') z_1 = z_0$. Then by definition,
\[
\Phi_{z_1}(y,z_0) = \ell' y \; ,
\]
and
\[
\Phi_{z_0} \bigl( \Phi_{z_1}(y,z_0),z \bigr) = \ell_0 \ell' y \; ,
\]
with $\ell_0 \in \Cent_{P(\BR)U(\BC)}(h_{\ell' y} \circ w)
= \Int(\ell') ( \Cent_{P(\BR)U(\BC)}(h_y \circ w) )$ such that 
$\pi(\ell_0) z_0 = z$. Thus, the element $\ell_0$ is of the form $\Int(\ell') (\ell_1)$,
with $\ell_1 \in \Cent_{P(\BR)U(\BC)}(h_y \circ w)$. But then, the product
$\ell_0 \ell' = \ell' \ell_1$ also belongs to $\Cent_{P(\BR)U(\BC)}(h_y \circ w)$, and
$\pi(\ell_0 \ell')(z_1) = \pi(\ell_0) z_0 = z$. 
\end{Proof}

\begin{Rem}
The following will not be used in the sequel. 
If hypotheses $(+)$ and $(U=0)$ are satisfied, then we have the geodesic actions
on $\FX$ and on $\FX/W$. The formula
\[
\Phi_{z_0} (x,\pi(\ell)z) = \ell \Phi_{z_0} (x,z) \; , \; 
\forall \, z \in \FX/W \; ,
\]
for $\ell \in \Cent_{P(\BR)}(h_x \circ w) = L_x(\BR)$
(Proposition~\ref{5B}~(d)), and Proposition~\ref{3I} together imply that
\[
\Phi_{z_0} (x,a \cdot z) = a \cdot \Phi_{z_0} (x,z) \; , \; 
\forall \, z \in \FX/W \; ,
\]
for any parabolic $Q$ of $P$ and $a \in A_Q$ (which we identify with $A_{\pi(Q)}$).
It follows that $\Phi_{z_0}$ extends to an isomorphism of manifolds with corners
\[
\Phi_{z_0}^{BS} : \pi^{-1}(z_0) \times (\FX/W)^{BS} \isoto \FX^{BS} \; .
\]
There is a commutative diagram
\[
\vcenter{\xymatrix@R-10pt{
\pi^{-1}(z_0) \times (\FX/W)^{BS} \ar[rr]^-{\Phi_{z_0}^{BS}} \ar@{->>}[dr]_-{pr_2} && 
\FX^{BS} \ar@{->>}[dl]^-{\pi^{BS}} \\
& (\FX/W)^{BS}  &
\\}}
\]
where $pr_2$ denotes the projection onto the second factor, and $\pi^{BS}$
is the map from Proposition~\ref{3K}~(a).
In particular, the latter is a trivial fibration.
\end{Rem}

\begin{Def} \label{5C}
Let $x \in \FX$. Define a map
\[
s_x : \FX/W \longto \FX \; , \; 
z \longmapsto \Phi_{\pi(x)} (x,z) \; .
\]
\end{Def}

\begin{Prop} \label{5D}
Let $x \in \FX$. \\[0.1cm]
(a)~The map $s_x$ is an $\BR$-analytic section of $\pi$. \\[0.1cm]
(b)~The map $s_x$ is $G(\BR)$-equivariant
in the following sense: 
for
$g \in G(\BR)$, we have
\[
s_x(gz) = \ell s_x (z) \; , \; 
\forall \, z \in \FX/W \; , 
\]
where $\ell \in \Cent_{P(\BR)U(\BC)}(h_x \circ w)$ is such that
$\pi(\ell) = g$. 
In particular, the image of $s_x$ equals the orbit of $x$ under 
$\Cent_{P(\BR)U(\BC)}(h_x \circ w)$. \\[0.1cm]
(c)~Let $y$ be another point of the image of $s_x$. 
Then $s_y = s_x$. \\[0.1cm]
(d)~Let $z \in \FX/W$. 
Then the imaginary part \cite[Sect.~4.14]{P1}
of $s_x(z)$ equals that of $x$:
\[
\im \bigl( s_x(z) \bigr) = \im ( x ) \in U(\BC) \; .
\]
In particular, the composition $\im \circ s_x$ is constant.
\end{Prop}

The proof of Proposition~\ref{5D} makes use of the following result.

\begin{Lem} \label{5Da}
Let $x \in \FX$, and $\ell \in \Cent_{P(\BR)U(\BC)}(h_x \circ w)$. Then
\[
\im \bigl( \ell x \bigr) = \im ( x ) \in U(\BC) \; .
\]
\end{Lem}

\begin{Proof}
Put $u := \im(x)$. 
Thus, the point $x_0 := u^{-1}x$ belongs to $\re (\FX)$, meaning that 
the morphism $h_{x_0} : \BS_\BC \to P_\BC$ is defined over $\BR$.
Hence so
is the Levi subgroup 
\[
\Cent_{P_\BC}(h_{x_0} \circ w)
\] 
(it actually coincides
with the base change from $\BR$ to $\BC$ of the group $L_{x_0}$ from
Definition~\ref{1Not}).
As $\Cent_{P(\BR)U(\BC)}(h_{x_0} \circ w)$ maps isomorphically to $G(\BR)$
under $\pi$, it follows that
\[
\Cent_{P(\BR)U(\BC)}(h_{x_0} \circ w) = \Cent_{P(\BR)}(h_{x_0} \circ w) \subset P(\BR) \; .
\]
The product $u^{-1} \ell u$ is an element of 
\[
\Cent_{P(\BR)U(\BC)}(h_{x_0} \circ w) = \Cent_{P(\BR)}(h_{x_0} \circ w) 
\]
and is therefore defined over $\BR$.
Therefore, the point $u^{-1} (\ell x) = (u^{-1} \ell u) x_0$
belongs to $\re (\FX)$. By definition, $u = \im(\ell x)$.
\end{Proof}

\medskip

\begin{Proofof}{Proposition~\ref{5D}}
Use the definition, and 
Proposition~\ref{5B}~(b)--(e) and Lemma~\ref{5Da}.
\end{Proofof}

The following is immediate from the definition.

\begin{Prop} \label{5E}
Let $z_0 \in \FX/W$. Then 
\[
\Phi_{z_0}(x,z) = s_x(z) \; , \; 
\forall \, (x,z) \in \pi^{-1}(z_0) \times \FX/W \; .
\]
\end{Prop}

As in \cite[Chap.~4]{P1}, let us denote by $U(\BR)(-1)$ 
the real vector space  of purely imaginary elements of $U(\BC)$.

\begin{Cor} \label{5F}
Let $\FH' \subset \FX/W$ and $C \subset U(\BR)(-1)$ be open subsets. Consider the open subset
\[
\FV := \pi^{-1}(\FH') \cap \im^{-1}(C)
\]
of $\FX$. Then the induced complex analytic map $\pi_{\tei \FV} :\FV \to \FH'$ is 
$\BR$-analytically split.
More precisely, 
for every $z_0 \in \FH'$, the map $\Phi_{z_0}$ induces an $\BR$-analytic
isomorphism
\[
\pi_{\tei \FV}^{-1}(z_0) \times \FH' \isoto \FV \; .
\]
\end{Cor}

\begin{Proof}
Apply Propositions~\ref{5D}~(d) and \ref{5E}.
\end{Proof}

For an admissible 
parabolic subgroup $Q_j$ of $P$, consider as in Section~\ref{5}
the complex analytic, $Q_j(\BR)U(\BC)$-equivariant map
\[
\tilde{\pi}_j : \FX \longto \coprod \FX_j / W_j \; .
\]

\begin{Cor} \label{5G}
Let $Q_j$ be an admissible parabolic subgroup of $P$. Then
the complex analytic map $\tilde{\pi}_j : \FX \to \coprod \FX_j / W_j$ is 
$\BR$-analytically split.
More precisely, let $(P_j,\FX_j)$ be a rational boundary component
associated to $Q_j \,$, and $z_0 \in \FX_j / W_j$.
Denote by $(\FX_j / W_j)^0$ the connected component of $\FX_j / W_j$
containing $z_0 \,$, and define 
$\FX^0 := \tilde{\pi}_j^{-1}((\FX_j / W_j)^0)$. \\[0.1cm]
(a)~The map 
\[
\Phi_{z_0} : \pi_j^{-1}(z_0) \times \FX_j / W_j \isoto \FX_j 
\]
associated to the mixed Shimura data $(P_j,\FX_j)$
induces an $\BR$-analytic isomorphism
\[
\tilde{\Phi}_{z_0,j} : \tilde{\pi}_j^{-1}(z_0) \times (\FX_j / W_j)^0 \isoto \FX^0 \; .
\] 
(b)~There is a commutative diagram
\[
\vcenter{\xymatrix@R-10pt{
\tilde{\pi}_j^{-1}(z_0) \times (\FX_j / W_j)^0 
                             \ar[rr]^-{\tilde{\Phi}_{z_0,j}} \ar@{->>}[dr]_-{pr_2} && 
\FX^0 \ar@{->>}[dl]^-{\tilde{\pi}_{j \tei \FX^0}} \\
& (\FX_j / W_j)^0  &
\\}}
\]
where $pr_2$ denotes the projection onto the second factor. \\[0.1cm]
(c)~For all $x \in \tilde{\pi}_j^{-1}(z_0)$, we have
\[
\tilde{\Phi}_{z_0,j}(x,z_0) = x \; .
\]
\end{Cor}

\begin{Proof}
Put $\FH := (\FX_j / W_j)^0$.
According to \cite[Sect.~4.11]{P1}, its pre-image $\FX^0 = \tilde{\pi}_j^{-1}(\FH)$
is a connected component of $\FX$. Define
$C$ as the image $C(\FX^0,P_j)$ of the composition of 
\[
\iota_j: \FX \longto \coprod \FX_j
\]
and
the imaginary part $\im_j$ associated to $(P_j,\FX_j)$.
Then \cite[Prop.~4.15~(a), (b)]{P1}, the subset $C$ of $U_j(\BR)(-1)$ is open,
and the map $\iota_j$ identifies $\FX^0$ and
\[
\FV := \pi_j^{-1}(\FH) \cap \im_j^{-1}(C) \; .
\]
Now apply Corollary~\ref{5F} and Proposition~\ref{5B}~(c)
to $(P_j,\FX_j)$, $\FH$ and $C$.
\end{Proof} 

\begin{Rem} \label{5H}
In the situation of Corollary~\ref{5G}, the map $\Phi_{z_0}$ actually
induces an isomorphism of the larger open subsets
\[
\tilde{\pi}_j^{-1}(z_0) \times \FX_j / W_j \isoto \im_j^{-1}(C) \; .
\]
But unless $\FX$ is connected,
the right hand side is not in general contained in the subset $\FX$ of $\FX_j$.
\end{Rem}

\begin{Ex} \label{5I}
Consider the Shimura data $(\GL_{2,\BQ},\FH_2)$ of Example~\ref{1H}, 
and the admissible parabolic subgroup $Q_j$ of $\GL_{2,\BQ}$ of upper
triangular matrices. According to \cite[Sect.~4.25]{P1},
there is a unique rational boundary component $(P_j,\FX_j)$
associated to $Q_j$. It is isomorphic to the mixed Shimura data
denoted $(P_0,\FX_0)$ in \cite[Ex.~2.24]{P1}. We have  
\[
\FH_2 = \FH^+ \coprod \FH^- \; ,
\]
where $\FH^+$ and $\FH^-$ denote the upper and the lower half plane in $\BC$,
respectively, 
\[
\FX_j = \BC \coprod \BC 
\]
(two copies of $\BC$), and  
\[
\FX_j/W_j = \{+\} \coprod \{-\} \; .
\]
The map $\iota_j : \FH_2 \to \FX_j$ equals the inclusion of $\FH^+$ 
into the first component $\BC$ and the inclusion of $\FH^-$ into the second
\cite[Sect.~4.26]{P1}, and 
the canonical epimorphism $\FX_j \onto \FX_j/W_j$ is the projection to the respective
base points, $+$ for the first, and $-$ for the second component $\BC \,$.

This illustrates Remark~\ref{5H}, even in a situation where Proposition~\ref{5B}
is trivial: put $z_0 := +$. The $\BR$-analytic isomorphism 
\[
\Phi_+ : \BC \times \FX_j/W_j \isoto \FX_j
\]
($\BC$ being the first component of $\FX_j$) maps $(x,+)$ to the point
$x$ in the first component $\BC$ of $\FX_j$, and $(x,-)$ to the point
$x$ in the second component $\BC \,$. Therefore,
\[
\Phi_+ \bigl( \FH^+ \times \FX_j/W_j \bigr) = \FH^+ \coprod \FH^+ \not \subset \FH_2 \; .
\]
Indeed, the splitting of $\tilde{\pi}_j$ provided by Corollary~\ref{5G}
equals $\Phi_+$ (only) on $\BC \times \{+\} \subset \BC \times \FX_j/W_j$,
and $\Phi_-$ on $\BC \times \{-\}$, for the second component $\BC$ of $\FX_j$.
\end{Ex}

Recall (Definition~\ref{4G}) the closed
connected normal subgroup
\[ 
C_j = \{ q \in Q_j \; , \; \pi_{Q_j}(q) \in \Cent_{\bar{Q}_j}(\pi_{Q_j}(P_j)) \}^0 
\]
of $Q_j$. 

\begin{Comp} \label{5J}
Let $(P_j,\FX_j)$ be a rational boundary component
associated to $Q_j \,$, and $z_0 \in \FX_j / W_j$.
Denote by $(\FX_j / W_j)^0$ the connected component of $\FX_j / W_j$
containing $z_0 \,$, and define 
$\FX^0 := \tilde{\pi}_j^{-1}((\FX_j / W_j)^0)$. Then the map
\[
\tilde{\Phi}_{z_0,j} : \tilde{\pi}_j^{-1}(z_0) \times (\FX_j / W_j)^0 \isoto \FX^0 
\]
is $\Cent_{C_j(\BR)U(\BC)} ( (\FX_j/W_j)^0 )$-equivariant
in the following sense: 
for all $q \in \Cent_{C_j(\BR)U(\BC)} ( (\FX_j/W_j)^0 )$,
we have
\[
\tilde{\Phi}_{z_0,j}(qx,z) = q \tilde{\Phi}_{z_0,j}(x,z) \; , \;
\forall \, (x,z) \in \tilde{\pi}_j^{-1}(z_0) \times (\FX_j / W_j)^0 \; .
\]
\end{Comp}

\begin{Proof}
The desired equality can be checked in $\FX_j$ as $\iota_j$ is injective. 
Write $x_j$ for the image of $x$ under $\iota_j$.
By Definition~\ref{5A},
\[
\tilde{\Phi}_{z_0,j}(x,z) = \ell x_j \; ,
\]
where $\ell \in \Cent_{P_j(\BR)U_j(\BC)}(h_{x_j} \circ w)$ is such that 
$\pi_j(\ell) z_0 = z$. But then, the element $q\ell q^{-1}$ belongs to
$\Cent_{P_j(\BR)U_j(\BC)}(h_{qx_j} \circ w)$.
Since $q \in C_j$, it satisfies 
$\pi_j(q \ell q^{-1}) = \pi_{Q_j}(q) \pi_j(\ell) \pi_{Q_j}(q)^{-1} = \pi_j(\ell)$. 
Thus, again by Definition~\ref{5A},
\[
\tilde{\Phi}_{z_0,j}(qx,z) = \bigl( q \ell q^{-1} \bigr) qx_j \; ,
\]
which equals $q \tilde{\Phi}_{z_0,j}(x,z)$.
\end{Proof}

We are ready to analyze the compatibility between the splittings $\tilde{\Phi}_{z_0,j}$
and the geodesic action. Recall the Levi subgroups $L_{x,Q}$ 
(Definition~\ref{2C}, Corollary~\ref{2G}).

\begin{Lem} \label{5N}
Assume hypotheses $(+)$ and $(U=0)$. 
Let $Q_j$ be an admissible parabolic subgroup of $P$,
and $Q$ a parabolic contained in $Q_j$ and containing
$P_j$ (in other words, we have $Q \in \adm^{-1}(Q_j)$). 
Let $x \in \FX$. The inclusion 
\[
\Cent_{{}^0 \! P_j^0(\BR)}(h_{\iota_j(x)} \circ w) 
\subset \Cent_{({}^0 \! P_j^0(\BR))U_j(\BC)}(h_{\iota_j(x)} \circ w) 
\]
is an equality, and 
\[
\Cent_{{}^0 \! P_{j,\BR}^0}(h_{\iota_j(x)} \circ w) \subset L_{x,Q} \; .
\]
\end{Lem}

\begin{Proof}
Let $q \in {}^0 \! P_j^0(\BR)$ and $u \in U_j(\BR)(-1)$. 
For any $y \in \FX_j$, we then have
\[
\im ((qu)y) = \Int(q) (u + \im(y)) 
\]
\cite[Sect.~4.14]{P1}. The Shimura data $(P_j,\FX_j)$, are irreducible
\cite[Rem.~(ii) on p.~82, after Sect.~4.11]{P1}.
Therefore \cite[Prop.~2.14~(a)]{P1}, the group $P_j$ acts on $\Lie U_j$ through
a character. But all characters of $P_j$ are trivial on ${}^0 \! P_j^0$
\cite[Sect.~1.1]{BS}. It follows that
\[
\im ((qu)y) = u + \im(y) \; .   
\]
This shows that if $(qu)y$ and $y$ have the same imaginary part, then $u$ is trivial.
Using Lemma~\ref{5Da}, we conclude that
\[
\Cent_{({}^0 \! P_j^0(\BR))U_j(\BC)} \bigl(h_{\iota_j(x)} \circ w \bigr) 
\subset {}^0 \! P_j^0(\BR) \; ,
\]
which establishes the first of the two claims. Next,
\[
\Cent_{{}^0 \! P_{j,\BR}^0} \bigl(h_{\iota_j(x)} \circ w \bigr) 
= \Cent_{{}^0 \! P_{j,\BR}^0} \bigl(h_{\re_j(x)} \circ w \bigr) \; ,
\]
again since ${}^0 \! P_{j,\BR}^0$ centralizes $U_j$.

In order to establish the second claim, it therefore suffices to show that
\[
\Cent_{P_{j,\BR}} \bigl(h_{\re_j(x)} \circ w) \subset L_{x,Q} \; .
\]
Recall (Corollary~\ref{2G}~(a)) that 
\[
L_{x,Q} \subset \Cent_{P_\BR} \bigl( h_{\re_j(x)} \circ w \bigr) \; .
\]
It follows that there is an inclusion
\[
P_{j,\BR} \cap L_{x,Q} \subset \Cent_{P_{j,\BR}} \bigl(h_{\re_j(x)} \circ w)  \; .
\]
But both sides are Levi subgroups of $P_{j,\BR}$ (recall that $P_j$ is normal in $Q$). 
Therefore, the inclusion is an equality.
\end{Proof}

\begin{Cor} \label{5Na}
Assume hypotheses $(+)$ and $(U=0)$. 
Let $Q_j$ be an admissible parabolic subgroup of $P$.
Let $(P_j,\FX_j)$ be a rational boundary component
associated to $Q_j \,$, and $z_0$ and $z$ be points in the same connected component
of $\FX_j / W_j$. Then for all $x \in \tilde{\pi}_j^{-1}(z_0)$, there exists
\[
\ell \in  \Cent_{P_j(\BR)U_j(\BC)}(h_{\iota_j(x)} \circ w) 
\cap \bigcap_{Q \in \adm^{-1}(Q_j)} L_{x,Q}(\BR)
\]
such that $\pi_j(\ell) z_0 = z$.
\end{Cor}

\begin{Proof}
Write $G_j := P_j/W_j$.
By Corollary~\ref{1Lcor}~(a), (b), the stabilizer of $z_0$
in $G_j(\BR)$ contains a connected complement of ${}^0 G_j(\BR)$ in $G_j(\BR)$.
Therefore, the subgroup ${}^0 G_j(\BR)$ in $G_j(\BR)$ acts transitively on $\FX_j / W_j$.
Choose $g \in {}^0 G_j(\BR)^0$ such that $gz_0 = z$, and 
$\ell \in \Cent_{P_j(\BR)U_j(\BC)}(h_{\iota_j(x)} \circ w)$
such that $\pi_j(\ell) = g$.
Given the choice of $g$, the element $\ell$ actually belongs to
$\Cent_{({}^0 \! P_j(\BR)^0)U_j(\BC)}(h_{\iota_j(x)} \circ w)$.
Now apply Lemma~\ref{5N}.
\end{Proof}

Corollary~\ref{5Na} and the relation 
``$\Int(\ell) (L_{x,Q}) = L_{\ell x, \Int(\ell)(Q)}$'' together imply the following
result.  

\begin{Cor} \label{5Nb}
Assume hypotheses $(+)$ and $(U=0)$. 
Let $Q_j$ be an admissible parabolic subgroup of $P$.
Let $(P_j,\FX_j)$ be a rational boundary component
associated to $Q_j \,$, and $z_0 \in \FX_j / W_j$. 
Denote by $(\FX_j / W_j)^0$ the connected component of $\FX_j / W_j$
containing $z_0 \,$, and define 
$\FX^0 := \tilde{\pi}_j^{-1}((\FX_j / W_j)^0)$.
Then the map
\[
\tilde{\Phi}_{z_0,j} : \tilde{\pi}_j^{-1}(z_0) \times (\FX_j / W_j)^0 \isoto \FX^0 
\]
respects the Levi subgroups $L_{\argdot,Q} \,$,
for $Q \in \adm^{-1}(Q_j)$, in the following sense: for all
points $x$ in $\tilde{\pi}_j^{-1}(z_0)$, and $z \in (\FX_j / W_j)^0$, we have
\[
L_{\tilde{\Phi}_{z_0,j}(x,z),Q} = L_{x,Q} 
\]
for all $Q \in \adm^{-1}(Q_j)$.
\end{Cor}

\begin{Cor} \label{5O}
Assume hypotheses $(+)$ and $(U=0)$. 
Let $Q_j$ be an admissible parabolic subgroup of $P$,
$(P_j,\FX_j)$ a rational boundary component
associated to $Q_j \,$, and $z_0 \in \FX_j / W_j$.
Denote by $(\FX_j / W_j)^0$ the connected component of $\FX_j / W_j$
containing $z_0 \,$, and define 
$\FX^0 := \tilde{\pi}_j^{-1}((\FX_j / W_j)^0)$. Then the $\BR$-analytic isomorphism
\[
\tilde{\Phi}_{z_0,j} : \tilde{\pi}_j^{-1}(z_0) \times (\FX_j / W_j)^0 \isoto \FX^0 
\]
is equivariant for the geodesic action: for all $Q \in \adm^{-1}(Q_j)$, and
for all $a \in A_Q$, we have
\[
\tilde{\Phi}_{z_0,j}(a \cdot x,z) = a \cdot \tilde{\Phi}_{z_0,j}(x,z) \; , \;
\forall \, (x,z) \in \tilde{\pi}_j^{-1}(z_0) \times (\FX_j / W_j)^0 \; .
\]
\end{Cor}

\begin{Proof}
Recall that
\[
a \cdot x = a_x x \; ,
\]
where $a_x \in Z(L_{x,Q})(\BR)$ maps to $a$ under $\pi_Q \,$. 
By Corollary~\ref{5R}~(b) and Complement~\ref{5J},
\[
\tilde{\Phi}_{z_0,j}(a \cdot x,z) = \tilde{\Phi}_{z_0,j}(a_x x,z)
= a_x \tilde{\Phi}_{z_0,j}(x,z) \; .
\]
Our statement is therefore proved once we establish the equality
\[
a \cdot \tilde{\Phi}_{z_0,j}(x,z) = a_x \tilde{\Phi}_{z_0,j}(x,z) \; .
\] 
By Definition~\ref{5A},
\[
\tilde{\Phi}_{z_0,j}(x,z) = \ell x \; ,
\]
where $\ell \in \Cent_{P_j(\BR)U_j(\BC)}(h_{\iota_j(x)} \circ w)$ is such that 
$\pi_j(\ell) z_0 = z$. 
According to Corollary~\ref{5Na}, $\ell$ can be chosen in $L_{x,Q}(\BR)$.
But then, $\ell$ and $a_x$ commute with each other.
Therefore,
\[
a_x = \ell a_x \ell^{-1} \in Z(L_{\ell x,Q})(\BR) \; ,
\]
and 
\[
a \cdot \tilde{\Phi}_{z_0,j}(x,z) = 
a \cdot (\ell x) = a_x \ell x = a_x \tilde{\Phi}_{z_0,j}(x,z) \; .
\]
\end{Proof}

For the rest of the section, 
assume that the Shimura data $(P,\FX) = (G,\FX)$ are pure, and that they satisfy 
hypothesis $(+)$. Also, fix an admissible parabolic subgroup $Q_j$ of $G$.
Consider as before the composition $\tilde{\pi}_j: \FX \to \coprod \FX_j / W_j$ of the inclusion
of $\FX$ into the disjoint union of the spaces $\FX_j$
underlying rational boundary components associated to $Q_j$,
and the canonical epimorphisms $\pi_j: \FX_j \to \FX_j/W_j$. \\

Recall the map $p: \FX^{BS} \to \FX^*$. By Construction~\ref{4K}, its restriction
$p_{Q_j}$ to the face $e(Q_j)$ is defined by the commutativity of the diagram
\[
\vcenter{\xymatrix@R-10pt{
        \FX \ar@{->>}[r] \ar@{=}[d] & e(Q_j) = A_{Q_j} \! \backslash \FX \ar[d]^{p_{Q_j}} \\
        \FX \ar[r]^-{\tilde{\pi}_j} & \coprod \FX_j / W_j 
\\}}
\]
Corollary~\ref{5O} implies in particular that $\tilde{\Phi}_{z_0,j}$ induces
an isomorphism
\[
\bigl( p^{-1}(z_0) \cap e(Q_j) \bigr) \times (\FX_j / W_j)^0 
= A_{Q_j} \backslash \tilde{\pi}_j^{-1}(z_0) \times (\FX_j / W_j)^0 
                                                          \isoto A_{Q_j} \backslash \FX^0 \; , 
\]
which we shall denote by $A_{Q_j} \backslash \tilde{\Phi}_{z_0,j}$. 
Its composition with the inclusions
\[
p^{-1}(z_0) \cap e(Q_j) 
\longinto \bigl( p^{-1}(z_0) \cap e(Q_j) \bigr) \times (\FX_j / W_j)^0 \; , \;
x \longmapsto (x,z_0)
\]
and 
\[
A_{Q_j} \backslash \FX^0 \longinto A_{Q_j} \backslash \FX = e(Q_j)
\]
equals the inclusion of $p^{-1}(z_0) \cap e(Q_j)$
into $e(Q_j)$, as follows from Corollary~\ref{5G}~(c). \\

The following is the first main result of this section.  

\begin{Thm} \label{6A}
Let $(P_j,\FX_j)$ be a rational boundary component
associated to $Q_j \,$, and $z_0 \in \FX_j / W_j$.
Denote by $(\FX_j / W_j)^0$ the connected component of $\FX_j / W_j$
containing $z_0 \,$, and define 
$\FX^0 := \tilde{\pi}_j^{-1}((\FX_j / W_j)^0)$. \\[0.1cm]
(a)~The composition of the isomorphism
\[
A_{Q_j} \backslash \tilde{\Phi}_{z_0,j}: 
\bigl( p^{-1}(z_0) \cap e(Q_j) \bigr) \times (\FX_j / W_j)^0 
\isoto A_{Q_j} \backslash \FX^0 
\]
and the inclusion
\[
A_{Q_j} \backslash \FX^0 
\longinto e(Q_j) \longinto \FX^{BS} 
\]
extends uniquely to a continuous map
\[
k_{z_0,j}: p^{-1}(z_0) \times (\FX_j / W_j)^0 \longto \FX^{BS} \; .
\]
The map $k_{z_0,j}$ is a morphism of manifolds with corners. \\[0.1cm]
(b)~We have 
\[
k_{z_0,j}(x,z_0) = x \; , \; \forall \, x \in p^{-1}(z_0) \; .
\] 
(c)~The morphism $k_{z_0,j}$ is injective. 
It respects the stratifications; more precisely,
for any parabolic $Q$ of $G$, we have 
\[
k_{z_0,j}^{-1} (e(Q)) = \bigl( p^{-1}(z_0) \cap e(Q) \bigr) \times (\FX_j / W_j)^0 
\]
if $Q \in \adm^{-1}(Q_j)$, and
\[
k_{z_0,j}^{-1} (e(Q)) = \emptyset
\]
if $Q \not \in \adm^{-1}(Q_j)$. \\[0.1cm]
(d)~The morphism $k_{z_0,j}$ is 
$\Cent_{C_j(\BQ)} ((\FX_j/W_j)^0)$-equivariant:
we have
\[
k_{z_0,j}(qx,z) = q k_{z_0,j}(x,z) \; , \;
\forall \, (x,z) \in p^{-1}(z_0) \times (\FX_j / W_j)^0 \; ,
\] 
for all $q \in \Cent_{C_j(\BQ)} ( (\FX_j/W_j)^0 )$. \\[0.1cm]
(e)~The diagram
\[
\vcenter{\xymatrix@R-10pt{
p^{-1}(z_0) \times (\FX_j / W_j)^0 \ar@{^{ (}->}[r]^-{k_{z_0,j}} \ar@{->>}[d]_-{pr_2} & 
                                   \FX^{BS} \ar@{->>}[d]^-p \\
(\FX_j / W_j)^0 \ar@{^{ (}->}[r] & \FX^*  &
\\}}
\]
is commutative 
($pr_2:=$ the projection onto the second factor). \\[0.1cm]
(f)~The diagram from (e) is Cartesian. In other words, the morphism 
$k_{z_0,j}$ yields an
identification of $p^{-1}(z_0) \times (\FX_j / W_j)^0$
with the pre-image under $p: \FX^{BS} \to \FX^*$ of 
$(\FX_j / W_j)^0 \subset \FX^*$. 
\end{Thm}

\begin{Proof}
Using the isomorphism
\[
\kappa_{z_0,j}: 
\bigl( A_{Q_j} \! \backslash \tilde{\pi}_j^{-1}(z_0) \bigr)^{BS} \isoto p^{-1}(z_0) 
\]
from Theorem~\ref{5S}, the required map $k_{z_0,j}$ is identified with
\[
\bigl( A_{Q_j} \! \backslash \tilde{\pi}_j^{-1}(z_0) \bigr)^{BS} \times (\FX_j / W_j)^0
\longto \FX^{BS} \; ,
\]
extending the inclusion
\[
A_{Q_j} \backslash \tilde{\pi}_j^{-1}(z_0) \times (\FX_j / W_j)^0 
\stackrel{A_{Q_j} \backslash \tilde{\Phi}_{z_0,j}}{\longto}
A_{Q_j} \backslash \FX^0 \longinto e(Q_j) \longinto \FX^{BS} \; .
\]
(a), (c), (d): imitate the proof of Theorem~\ref{5S}~(a), (b), (c),
using Corollary~\ref{5O} and Complement~\ref{5J}.
\forget{

as for the unicity statement in (a), we argue as usual 
($A_{Q_j} \! \backslash \tilde{\pi}_j^{-1}(z_0) \times (\FX_j / W_j)^0$ is dense in 
$( A_{Q_j} \! \backslash \tilde{\pi}_j^{-1}(z_0))^{BS} \times (\FX_j / W_j)^0$, 
and $\FX^{BS}$ is Hausdorff \cite[Thm.~7.8]{BS}).
Recall from Corollary~\ref{4H} that the map $Q \mapsto C_j \cap Q$ is
a bijection between $\adm^{-1}(Q_j)$ and the set of parabolics of $C_j$.
By Corollary~\ref{5R}, this bijection is compatible with the geodesic action.
The remaining claims therefore follow from the definition of the spaces 
$( A_{Q_j} \! \backslash \tilde{\pi}_j^{-1}(z_0))^{BS}$ and
$\FX^{BS}$ \cite[Sect.~7.1]{BS}.
}

\noindent (b): as $p^{-1}(z_0) \cap e(Q_j)$ is dense in $p^{-1}(z_0)$, and $\FX^{BS}$
is Hausdorff \cite[Thm.~7.8]{BS}, the inclusion of $p^{-1}(z_0)$ into $\FX^{BS}$ is
the only continuous extension to $p^{-1}(z_0)$ of the inclusion of $p^{-1}(z_0) \cap e(Q_j)$.

\forget{
(d): the inclusion
\[
A_{Q_j} \! \backslash \tilde{\pi}_j^{-1}(z_0) \times (\FX_j / W_j)^0 \longinto \FX^{BS}
\]
is $\Cent_{C_j(\BQ)W(\BR)} ((\FX_j/W_j)^0)$-equivariant (see Remark~\ref{5Ma}~(b)).
The claim thus follows from unicity of $k_{z_0,j}$ (part~(a)). 
}

\noindent (e): as $A_{Q_j} \backslash \tilde{\pi}_j^{-1}(z_0) \times (\FX_j / W_j)^0$ is dense in 
$( A_{Q_j} \! \backslash \tilde{\pi}_j^{-1}(z_0) )^{BS} \times (\FX_j / W_j)^0$, 
and $\FX^*$ is Hausdorff \cite[Thm.~4.9~(iii), (iv)]{BB}, it suffices to show that
\[
\vcenter{\xymatrix@R-10pt{
A_{Q_j} \backslash \tilde{\pi}_j^{-1}(z_0) \times (\FX_j / W_j)^0 
                             \ar@{^{ (}->}[r]^-{k_{z_0,j}} \ar@{->>}[d]_-{pr_2} & 
\FX^{BS} \ar@{->>}[d]^-p \\
(\FX_j / W_j)^0 \ar@{^{ (}->}[r] & \FX^*  &
\\}}
\]
commutes. But 
the restriction of $k_{z_0,j}$ to 
$A_{Q_j} \backslash \tilde{\pi}_j^{-1}(z_0) \times (\FX_j / W_j)^0$
equals 
\[
A_{Q_j} \backslash \tilde{\pi}_j^{-1}(z_0) \times (\FX_j / W_j)^0 
\stackrel{A_{Q_j} \backslash \tilde{\Phi}_{z_0,j}}{\longto}
A_{Q_j} \backslash \FX^0 \longinto e(Q_j) \longinto \FX^{BS} \; .
\]
Our claim therefore follows from Construction~\ref{4K} and Corollary~\ref{5G}~(b).

\noindent (f): by Construction~\ref{4K}, the pre-image $p^{-1}(\coprod \FX_j / W_j)$
equals the disjoint union of those faces $e(Q)$ for which $\adm(Q) = Q_j$.
Therefore, 
\[
p^{-1} \bigl( (\FX_j / W_j)^0 \bigr) = \coprod_{\adm(Q) = Q_j} A_Q \backslash \FX^0 \; .
\]
It thus suffices to observe that the right hand side of this equation equals the image of
$k_{z_0,j}$.
\end{Proof}

\begin{Rem}
Parts~(e) and (f) of Theorem~\ref{6A} imply in particular that
over any connected component $(\FX_j / W_j)^0$ of a stratum
of $\FX^*$,
the map $p: \FX^{BS} \to \FX^*$ is a trivial fibration with 
contractible fibres. This result is already known: \cite[Prop.~(3.8)~(ii)]{Z}.
\end{Rem}

\begin{Cor} \label{6C}
Let $(P_j,\FX_j)$ be a rational boundary component
associated to $Q_j \,$, and $z_0 \in \FX_j / W_j$.
Denote by $(\FX_j / W_j)^0$ the connected component of $\FX_j / W_j$
containing $z_0 \,$. Let $Q \in \adm^{-1}(Q_j)$. 
Then 
\[
k_{z_0,j}: p^{-1}(z_0) \times (\FX_j / W_j)^0 \isoto 
                                      p^{-1} \bigl( (\FX_j/W_j)^0 \bigr) 
\]
restricts to a $\Cent_{(C_j \cap Q)(\BR)} ((\FX_j/W_j)^0)$-equivariant $\BR$-analy\-tic isomorphism
\[
\bigl( p^{-1}(z_0) \cap e(Q) \bigr) \times (\FX_j / W_j)^0
                    \isoto p^{-1} \bigl( (\FX_j/W_j)^0 \bigr) \cap e(Q) \; ,
\]
and a $\Cent_{(C_j \cap Q)(\BQ)} ((\FX_j/W_j)^0)$-equivariant isomorphism
of manifolds with corners
\[
\bigl( p^{-1}(z_0) \cap \overline{ e(Q) } \bigr) \times (\FX_j / W_j)^0
                    \isoto p^{-1} \bigl( (\FX_j/W_j)^0 \bigr) \cap \overline{e(Q)} \; .
\]
\end{Cor}

\begin{Proof}
For all claims except $\Cent_{(C_j \cap Q)(\BR)} ((\FX_j/W_j)^0)$-equivariance of
\[
\bigl( p^{-1}(z_0) \cap e(Q) \bigr) \times (\FX_j / W_j)^0
                    \isoto p^{-1} \bigl( (\FX_j/W_j)^0 \bigr) \cap e(Q) \; ,
\]
apply Theorem~\ref{6A}~(c) and (d) (recall that by Theorem~\ref{6A}~(c), the only 
faces $e(R)$ having non-empty intersection with $p^{-1} ( (\FX_j/W_j)^0 )$,
are those for which $R \in \adm^{-1}(Q_j)$).
Consider the diagram
\[
\vcenter{\xymatrix@R-10pt{
        \FX \ar@{->>}[r] \ar@{=}[d] & A_Q \backslash \FX = e(Q) \ar[d]^{p_Q} \\
        \FX \ar[r]^-{\tilde{\pi}_j} & \coprod \FX_j / W_j 
\\}}
\]
as well as its base changes
\[
\vcenter{\xymatrix@R-10pt{
        \FX^0 \ar@{->>}[r] \ar@{=}[d] & 
        A_Q \backslash \FX^0 = p^{-1} \bigl( (\FX_j/W_j)^0 \bigr) \cap e(Q) \ar[d]^{p_Q}\\
        \FX^0 \ar[r]^-{\tilde{\pi}_j} & (\FX_j / W_j)^0 
\\}}
\]
to $(\FX_j / W_j)^0$ (where $\FX^0 := \tilde{\pi}_j^{-1}((\FX_j / W_j)^0)$), and
\[
\vcenter{\xymatrix@R-10pt{
        \tilde{\pi}_j^{-1}(z_0) \ar@{->>}[r] \ar@{=}[d] & 
        A_Q \backslash \tilde{\pi}_j^{-1}(z_0) = p^{-1}(z_0) \cap e(Q) \ar[d]^{p_Q} \\
        \tilde{\pi}_j^{-1}(z_0) \ar[r]^-{\tilde{\pi}_j} & \{ z_0 \}
\\}}
\]
to $\{ z_0 \}$,
recalling \cite[Prop.~3.4]{BS} that the geodesic action of 
$A_Q$ commutes with the action of $Q(\BR)$. The diagrams are therefore
$\Cent_{(C_j \cap Q)(\BR)} ((\FX_j/W_j)^0)$-equivariant. The restriction 
\[
\bigl( p^{-1}(z_0) \cap e(Q) \bigr) \times (\FX_j / W_j)^0
                    \longto p^{-1} \bigl( (\FX_j/W_j)^0 \bigr) \cap e(Q) \; ,
\]
of $k_{z_0,j}$ to $( p^{-1}(z_0) \cap e(Q) ) \times (\FX_j / W_j)^0$ 
is thus identified with the map
\[
A_Q \backslash \tilde{\pi}_j^{-1}(z_0) \times (\FX_j / W_j)^0
                    \longto A_Q \backslash \FX^0
\]
induced by $\tilde{\Phi}_{z_0,j}$, which 
is indeed $\Cent_{(C_j \cap Q)(\BR)} ((\FX_j/W_j)^0)$-equivariant (Com\-ple\-ment~\ref{5J}).
\end{Proof}

Theorem~\ref{6A} 
and Corollary~\ref{6C} can be sharpened, to show that the restriction of $p$ to the whole of 
$p^{-1} ( \coprod \FX_j/W_j )$ (not just its connected components) is (globally) trivial.

\begin{Thm} \label{6Ca}
Let $z_0 \in \coprod \FX_j / W_j \,$. 
Then the isomorphism $k_{z_0,j}$ from Theorem~\ref{6A} can be extended to give 
an isomorphism of manifolds with corners
\[ 
k_j: p^{-1}(z_0) \times \coprod \FX_j / W_j \isoto p^{-1} \bigl( \coprod \FX_j/W_j \bigr) 
    \stackrel{\ref{4Lb}}{=} \coprod_{Q \in \adm^{-1}(Q_j)} e(Q) \; ,
\]
satisfying the following. \\[0.1cm]
(a)~The isomorphism $k_j$ is $\Cent_{C_j(\BQ)}(\coprod \FX_j/W_j)$-equivariant. \\[0.1cm]
(b)~The isomorphism $k_j$ restricts to give
\[
\bigl( p^{-1}(z_0) \cap e(Q) \bigr) \times \coprod \FX_j / W_j \isoto e(Q) 
\]
and
\[
\bigl( p^{-1}(z_0) \cap \overline{ e(Q) } \bigr) \times \coprod \FX_j / W_j
                    \isoto p^{-1} \bigl( \coprod \FX_j / W_j \bigr) \cap \overline{e(Q)} \; ,
\]
for any parabolic $Q \in \adm^{-1}(Q_j)$. \\[0.1cm]
(c)~The isomorphism $k_j$ identifies the projection to the second
factor (on $p^{-1}(z_0) \times \coprod \FX_j / W_j$) with $p$ 
(on $p^{-1} ( \coprod \FX_j/W_j )$). 
\end{Thm}

\begin{Proof}
First, for \emph{any} choice of points $z_1$ in all connected components of 
$\coprod \FX_j / W_j$
not containing $z_0$, we may apply Theorem~\ref{6A} 
and Corollary~\ref{6C} for $z_0$ and all $z_1$.
As a result, we get $\Cent_{C_j(\BQ)}(\coprod \FX_j/W_j)$-equivariant trivializations
over each connected component. 

In order to get a global trivialization,
use Corollary~\ref{5Sa} for the choice of the $z_1$. 
\end{Proof}

\forget{
\begin{Cor} \label{6B}
Let $Q$ be a parabolic subgroup of $P$ contained in $Q_j \,$.
The base changes of both
\[
p_{\tei e(Q)}: e(Q) \longinto \FX^{BS} \stackrel{p}{\longto} \FX^*
\]
and  
\[
p_{\tei \overline{e(Q)}}: \overline{e(Q)} \longinto \FX^{BS} \stackrel{p}{\longto} \FX^*
\] 
to any connected component of $\coprod \FX_j/W_j \subset \FX^*$ are trivial fibrations.
The fibres of these base changes are contractible if $Q \in \adm^{-1}(Q_j)$
(\emph{i.e.}, if $P_j \subset Q$). Else, they are empty.
\end{Cor}

\begin{Proof}
Apply Corollaries~\ref{6C} and \ref{5T}~(b), and Complement~\ref{4Lb}.
\end{Proof}
}
The reader will find no difficulty, proceeding as in the proof
of Corollary~\ref{5U}~(a), to show the analogues of Theorems~\ref{6A}
and \ref{6Ca} for the map $p^r: \FX^{rBS} \to \FX^*$.
We therefore content ourselves with their statement.

\begin{Cor} \label{6D}
Let $z_0 \in \coprod \FX_j / W_j \,$. \\[0.1cm]
(a)~The isomorphism $k_j$ from Theorem~\ref{6Ca}
induces a $\Cent_{C_j(\BQ)} (\coprod \FX_j/W_j)$-equivariant identification of 
$(p^r)^{-1}(z_0) \times \coprod \FX_j / W_j$
with the pre-image under $p^r: \FX^{rBS} \to \FX^*$ of 
$\coprod \FX_j / W_j \subset \FX^*$. \\[0.1cm]
(b)~Let $Q \in \adm^{-1}(Q_j)$. 
Then the isomorphism
\[
(p^r)^{-1}(z_0) \times \coprod \FX_j / W_j \isoto 
                                      (p^r)^{-1} \bigl( \coprod \FX_j / W_j \bigr) 
\]
from (a) restricts to give
\[
\bigl( (p^r)^{-1}(z_0) \cap e^r(Q) \bigr) \times \coprod \FX_j / W_j \isoto e^r(Q) 
\]
and 
\[
\bigl( (p^r)^{-1}(z_0) \cap \overline{ e^r(Q) } \bigr) \times \coprod \FX_j / W_j
                    \isoto (p^r)^{-1} \bigl( \coprod \FX_j / W_j \bigr) \cap \overline{e^r(Q)} 
\]
($\overline{e^r(Q)} :=$ the closure of $e^r(Q)$ in $\FX^{rBS}$). \\[0.1cm]
(c)~The isomorphism 
\[
(p^r)^{-1}(z_0) \times \coprod \FX_j / W_j \isoto 
                                      (p^r)^{-1} \bigl( \coprod \FX_j / W_j \bigr) 
\]
from (a) identifies the projection to the second
factor with $p^r$. 
\forget{
\\[0.1cm]
(c)~Let $Q \in \adm^{-1}(Q_j)$. 
Then the base changes of both
\[
p^r_{\tei e^r(Q)}: e^r(Q) \longinto \FX^{rBS} \stackrel{p^r}{\longto} \FX^*
\]
and  
\[
p^r_{\tei \overline{e^r(Q)}}: \overline{e^r(Q)} \longinto \FX^{rBS} \stackrel{p^r}{\longto} \FX^*
\] 
to $(\FX_j/W_j)^0 \subset \FX^*$ are trivial fibrations, whose fibres are connected. 
}
\end{Cor}

Let us finish the section by spelling out the consequences of the above results
for Shimura varieties and their compactifications. As before,
consider the disjoint union $\coprod \FX_j$ of the finitely many spaces 
underlying rational boundary components associated to $Q_j \,$,
and its quotient $\coprod \FX_j / W_j \,$.

\begin{Prop} \label{6E}
Let $K$ be an open compact subgroup of $G (\BA_f)$, and
$g \in G(\BA_f)$. Consider the image
$pr_{gK}(\coprod \FX_j/W_j)$ of $\coprod \FX_j/W_j \times \{gK\}$ under the projection 
\[
\FX^* \times G (\BA_f) / K 
\longonto G (\BQ) \backslash \bigl( \FX^* \times G (\BA_f) / K \bigr) = M^K (G,\FX)^* (\BC) \; .
\]
Define
\[
H_Q:= H_Q (gK) := Q_j(\BQ) \cap gKg^{-1} \subset Q_j(\BQ) \; .
\] 
(a)~As a subset of $M^K (G,\FX)^* (\BC)$, 
\[
pr_{gK} \bigl( \coprod \FX_j/W_j \bigr) 
= H_Q \backslash \Bigl(  \coprod \FX_j/W_j \times \{gK\} \Bigr)
\subset G (\BQ) \backslash \bigl( \FX^* \times G (\BA_f) / K \bigr) \; .
\]
Consequently, the pre-image of $pr_{gK}(\coprod \FX_j/W_j)$ under 
\[
p^K: M^K (G,\FX) (\BC)^{BS} \longto M^K (G,\FX)^* (\BC)
\]
equals
\[
H_Q \backslash \Bigl( p^{-1} \bigl( \coprod \FX_j/W_j \bigr) \times \{gK\} \Bigr)
\subset G (\BQ) \backslash \bigl( \FX^{BS} \times G (\BA_f) / K \bigr) \; ,
\]
and the pre-image of $pr_{gK}(\coprod \FX_j/W_j)$ under 
\[
p^{r,K}: M^K (G,\FX) (\BC)^{rBS} \to M^K (G,\FX)^* (\BC)
\]
equals
\[
H_Q \backslash \Bigl( (p^r)^{-1} \bigl( \coprod \FX_j/W_j \bigr) \times \{gK\} \Bigr)
\subset G (\BQ) \backslash \bigl( \FX^{rBS} \times G (\BA_f) / K \bigr) \; ,
\]
(b)~If $K$ is neat, then the action of 
\[
H_C:= H_C(gK) := C_j(\BQ) \cap gKg^{-1} \subset C_j(\BQ) 
\]
on $\coprod \FX_j/W_j$ is trivial. 
The induced action on
$\coprod \FX_j/W_j$ of the quotient $H_Q/H_C$ is free. \\[0.1cm]
(c)~If $K$ is neat, then the diagrams
\[
\vcenter{\xymatrix@R-10pt{
H_C \backslash \Bigl( p^{-1} \bigl( \coprod \FX_j/W_j \bigr) \times \{gK\} \Bigr)
                                                              \ar@{->>}[r] \ar[d]_-{p} & 
(p^K)^{-1} \bigl( pr_{gK} \bigl( \coprod \FX_j/W_j \bigr) \bigr) \ar[d]^-{p^K} \\
\coprod \FX_j/W_j \times \{gK\} \ar@{->>}[r] & pr_{gK} \bigl( \coprod \FX_j/W_j \bigr) 
\\}}
\]
and
\[
\vcenter{\xymatrix@R-10pt{
H_C \backslash \Bigl( (p^r)^{-1} \bigl( \coprod \FX_j/W_j \bigr) \times \{gK\} \Bigr)
                                                              \ar@{->>}[r] \ar[d]_-{p^r} & 
(p^{r,K})^{-1} \bigl( pr_{gK} \bigl( \coprod \FX_j/W_j \bigr) \bigr) \ar[d]^-{p^{r,K}} \\
\coprod \FX_j/W_j \times \{gK\} \ar@{->>}[r] & pr_{gK} \bigl( \coprod \FX_j/W_j \bigr) 
\\}}
\]
are Cartesian. 
\end{Prop}

\begin{Proof}
We leave the proof of part~(a) to the reader.

Part~(b) is Proposition~\ref{5X}~(e), applied to $\Gamma = H_Q$.

Under the identifications from (a), the first diagram in (c) transforms into
\[
\vcenter{\xymatrix@R-10pt{
H_C \backslash p^{-1} \bigl( \coprod \FX_j/W_j \bigr) \ar@{->>}[r] \ar[d]_-{p} & 
H_Q \backslash p^{-1} \bigl( \coprod \FX_j/W_j \bigr) \ar[d]^-{p} \\
\coprod \FX_j/W_j \ar@{->>}[r] & H_Q \backslash \bigl( \coprod \FX_j/W_j \bigr) 
\\}}
\]
Denote by $F$ the Cartesian product of 
$\coprod \FX_j/W_j$ and $H_Q \backslash p^{-1} (\coprod \FX_j/W_j )$
over $H_Q \backslash (\coprod \FX_j/W_j)$.
It is easy to see that the map 
\[
\alpha : H_C \backslash p^{-1} \bigl(\coprod \FX_j/W_j \bigr) \longto F
\]
is surjective. Let $x$ and $y$ be elements of $p^{-1} (\coprod \FX_j/W_j)$,
and assume that their images in $\coprod \FX_j/W_j$ and in 
$H_Q \backslash p^{-1} (\coprod \FX_j/W_j)$
are the same. In other words, there is an element $q \in H_Q$ such that $y = qx$, and
such that $q$ stabilizes $p(x)$. The action of $H_Q/H_C$ on $\coprod \FX_j/W_j$ being free
(according to (b)),
the latter property of $q$ implies 
$q \in H_C$. Therefore, the points $x$ and $y$ yield the same class
in $H_C \backslash p^{-1} (\coprod \FX_j/W_j)$, and $\alpha$ is injective.

The proof for the second diagram in (c) is formally identical. We leave it to the reader.
\end{Proof}

Combining Theorem~\ref{6Ca}, Corollary~\ref{6D} and Proposition~\ref{6E}, we get the following.

\begin{Cor} \label{6Ea}
Let $K$ be a neat open compact subgroup of $G (\BA_f)$, and $g \in G(\BA_f)$. 
Consider the map
\[
pr_{gK}: \coprod \FX_j/W_j  \longonto pr_{gK} \bigl( \coprod \FX_j/W_j \bigr) \into M^K (Q_j,\FX) \into M^K (G,\FX)^* (\BC) \; .
\]
The base changes \emph{via} $pr_{gK}$ of both
\[
p^K: M^K (G,\FX) (\BC)^{BS} \to M^K (G,\FX)^* (\BC)
\]
and
\[
p^{r,K}: M^K (G,\FX) (\BC)^{rBS} \to M^K (G,\FX)^* (\BC)
\] 
are trivial stratified fibrations. More precisely: \\[0.1cm]
(a)~The base changes \emph{via} $pr_{gK}$ of $p^K$ and of $p^{r,K}$ are given by the maps
\[
p: H_C \backslash p^{-1} \bigl( \coprod \FX_j/W_j \bigr) \longto \coprod \FX_j/W_j
\]
and 
\[
p^r: H_C \backslash (p^r)^{-1} \bigl( \coprod \FX_j/W_j \bigr) \longto \coprod \FX_j/W_j \; ,
\] 
respectively. \\[0.1cm]
(b)~The isomorphisms from Theorem~\ref{6Ca} and Corollary~\ref{6D}, associated to a choice of base point $z_0 \in \coprod \FX_j/W_j \,$, induce isomorphisms
\[ 
(p^K)^{-1} \bigl( [(z_0,gK)] \bigr) \times \coprod \FX_j / W_j \isoto H_C \backslash p^{-1} \bigl( \coprod \FX_j/W_j \bigr) 
\]
and
\[ 
(p^{r,K})^{-1} \bigl( [(z_0,gK)] \bigr) \times \coprod \FX_j / W_j \isoto H_C \backslash (p^r)^{-1} \bigl( \coprod \FX_j/W_j \bigr) \; ,
\]
respectively. These isomorphisms identify the projections to the second factor with the maps $p$ and $p^r$ from (b), respectively. \\[0.1cm]
(c)~The isomorphisms from (b) respect the stratifications induced by the canonical stratifications on $M^K (G,\FX) (\BC)^{BS}$
and on $M^K (G,\FX) (\BC)^{rBS}$, respectively (\emph{via} the inclusions of the fibres 
\[
(p^K)^{-1} \bigl( [(z_0,gK)] \bigr) \subset M^K (G,\FX) (\BC)^{BS} \; , 
\]
\[
(p^{r,K})^{-1} \bigl( [(z_0,gK)] \bigr) \subset M^K (G,\FX) (\BC)^{rBS}
\]
on the sources and the base change maps
\[
H_C \backslash p^{-1} \bigl( \coprod \FX_j/W_j \bigr) \longto M^K (G,\FX) (\BC)^{BS} \; , 
\]
\[
H_C \backslash (p^r)^{-1} \bigl( \coprod \FX_j/W_j \bigr) \longto M^K (G,\FX) (\BC)^{rBS}
\] 
on the targets).
\end{Cor}

\begin{Proof}
(a): since $K$ is assumed neat, Proposition~\ref{6E}~(c) applies.

\noindent (b): use parts~(a) and (c) of Theorem~\ref{6Ca} and Corollary~\ref{6D} respectively, and the canonical isomorphisms
\[
H_C \backslash p^{-1}(z_0) \isoto (p^K)^{-1} \bigl( [(z_0,gK)] \bigr)
\]
and
\[
H_C \backslash (p^r)^{-1}(z_0) \isoto (p^{r,K})^{-1} \bigl( [(z_0,gK)] \bigr) 
\]
(Theorem~\ref{5Z}~(a), Corollary~\ref{5Y}~(a)).

\noindent (c): we give the proof only for the Borel--Serre compactification,
and leave it to the reader to perform the necessary modifications for $M^K (G,\FX) (\BC)^{rBS}$. First, according to Theorem~\ref{4aG}~(a), 
\[
(p^K)^{-1} \bigl( M^K (Q_j,\FX)(\BC) \bigr) = e^K \bigl( Q_j, G(\BA_f) \bigr)' \; ,
\]
the latter subset of $M^K  (G,\FX) (\BC)^{BS}$ being (by Definition~\ref{4aC}) equal to  
\[
\bigcup_{Q \in \adm^{-1}(Q_j)} e^K \bigl( Q,G(\BA_f) \bigr) \; .
\]
Recall then (Definition~\ref{4aA}), that for each parabolic $Q \in \adm^{-1}(Q_j)$, we have
\[
e^K \bigl( Q,G(\BA_f) \bigr) = \bigcup_{g' \in G (\BA_f)} e^K(Q,g') \; ,
\]
where each $e^K(Q,g')$ is the image of
\[
e(Q) \times \{ g'K \} \subset \FX^{BS} \times G (\BA_f) / K
\]
under the projection from $\FX^{BS} \times G (\BA_f) / K$ to 
\[
M^K (G,\FX) (\BC)^{BS} = G (\BQ) \backslash \bigl( \FX^{BS} \times G (\BA_f) / K \bigr) \; .
\] 
The canonical stratification on $M^K (G,\FX) (\BC)^{BS}$
being thus coarser than the one by the $e^K(Q,g')$,
our claim follows from Theorem~\ref{6Ca}~(b).
\end{Proof}

Here is the second main result of this section.

\begin{Thm} \label{6F}
Let $K$ be a neat open compact subgroup of $G (\BA_f)$. 
Fix an admissible parabolic subgroup $Q_j$ of $G$, and
consider the canonical stratum $M^K (Q_j,\FX)$ of $M^K (G,\FX)^*$
(according to Definition~\ref{4aF}~(a), its space of complex points thus equals the image of 
$\coprod \FX_j / W_j \times G(\BA_f)$ under the projection
\[
\FX^* \times G (\BA_f) / K 
\longonto M^K (G,\FX)^* (\BC)) \; .
\]
(a)~The base change of the map $p^K: M^K (G,\FX) (\BC)^{BS} \to M^K (G,\FX)^* (\BC)$ to 
$M^K (Q_j,\FX)(\BC)$ 
is a locally trivial stratified fibration. \\[0.1cm]
(b)~The base change of the map $p^{r,K}: M^K (G,\FX) (\BC)^{rBS} \to M^K (G,\FX)^* (\BC)$ to 
$M^K (Q_j,\FX)(\BC)$ is a locally trivial stratified fibration.
\end{Thm}

\begin{Proof}
We may replace $M^K (Q_j,\FX)(\BC)$ by $pr_{gK}(\coprod \FX_j/W_j)$,
for $g$ running through the elements of $G(\BA_f)$. According to Proposition~\ref{6E}~(a), (b), it is sufficient to prove
the claims after base change \emph{via} $pr_{gK}$.
Now apply Corollary~\ref{6Ea}.
\forget{
This in turn will be achieved once we know that
\[
p: p^{-1} \bigl( \coprod \FX_j/W_j \bigr) 
\longonto \coprod \FX_j/W_j 
\]
is $H_C$-equivariantly locally trivial. Let $(\FX_j/W_j)^0$ be a connected component
of $\coprod \FX_j/W_j$. Choose a point $z_0$ in $(\FX_j/W_j)^0$. Then according to
Theorem~\ref{6A}~(d)--(f), the morphism of manifolds with corners
 \[
k_{z_0,j}: p^{-1}(z_0) \times (\FX_j / W_j)^0 \longto \FX^{BS}
\]
is a $\Cent_{C_j(\BQ)} ((\FX_j/W_j)^0)$-, hence $H_C$-equivariant trivialization of 
$p_{\tei (\FX_j / W_j)^0}$.
}
\end{Proof}

\begin{Rem}
Possibly up to the statement on the stratification,
Theorem~\ref{6F}~(b), together with the description of the 
fibres of $p^{r,K}$ from Corollary~\ref{5Y}~(c), follows from  
\cite[Thm.~5.8]{G}, at least if $\FX$ is connected. 
\end{Rem}
\forget{
\begin{Rem}
Let us reassure ourselves about the action of $H_C$ on 
$\FZ = A_{Q_j} \! \backslash \tilde{\pi}_j^{-1}(z_0)$. By construction (Proposition~\ref{5L}),
$\FZ$ carries the action of $C_j(\BR)$ obtained by restricting the extended action of
$G(\BR)$ on $\FX$ from Corollary~\ref{1G}. According to Remark~\ref{5Ma}~(b), the restriction
of this action to $\Cent_{C_j(\BR)} ( \FX_j/W_j )$ coincides with the restriction of the action 
underlying the Shimura data.  
But according to Proposition~\ref{6E}~(b),
\[
H_C \subset \Cent_{C_j(\BR)} ( \FX_j/W_j ) \; .
\]
Therefore, the action of $H_C$ on $\FZ$ is the one underlying the Shimura data.
\end{Rem}

Using Theorem~\ref{6Ca}~(b) as well, one obtains
the stratified version of Theorem~\ref{6F}~(a). 

\begin{Comp} \label{6G}
Let $K$ be a neat open compact subgroup of $G (\BA_f)$.
Fix an admissible parabolic subgroup $Q_j$ of $G$, and consider
the canonical stratum $M^K (Q_j,\FX)$ of $M^K (G,\FX)^*$.
Then the base change of the map $p^K: M^K (G,\FX) (\BC)^{BS} \to M^K (G,\FX)^* (\BC)$ to 
$M^K (Q_j,\FX)(\BC)$ 
is a locally trivial stratified fibration: local trivializations
of 
\[
p^K_{\tei e^K(Q_j,G(\BA_f))'}: e^K \bigl( Q_j,G(\BA_f) \bigr)' \longto M^K (Q_j,\FX)(\BC)
\]
(cmp.~Theorem~\ref{4aG}~(a)) can be chosen, inducing local trivializations of
\[
p^K_{\tei e^K(Q,G(\BA_f))}: e^K \bigl( Q,G(\BA_f) \bigr) \longto M^K (Q_j,\FX)(\BC)
\]
for every $Q \in \adm^{-1}(Q_j)$.
\end{Comp}

\begin{Rem}
The fibres of the maps $p^K_{\tei e^K(Q,G(\BA_f))}$ 
and $p^K_{\tei e^K(Q,G(\BA_f))'} \,$, and the stratification
of the fibres of $p^K_{\tei e^K(Q,G(\BA_f))'}$ induced by the canonical stratification 
of $M^K (G,\FX) (\BC)^{BS}$ (Definition~\ref{4aA}~(b)) are described in Theorem~\ref{5Z}. 
\end{Rem}

We leave it to the reader to formulate and prove the 
stratified version of Theorem~\ref{6F}~(b).
}


\bigskip

\bigskip

%
%

\end{document}